\documentclass[12pt]{amsart}
\usepackage{amssymb,amscd}
\usepackage[russian]{babel}
\usepackage{graphicx}
\usepackage{pdfpages}
\usepackage{placeins}
\newtheorem{theorem}{Теорема}[section]

\newtheorem{lemma}{Лемма}[section]

\theoremstyle{definition}
\newtheorem{definition}{Определение}[section]

\theoremstyle{remark}

\def \x#1{\begin{tabular}{c}#1\end{tabular}}

\ifx\pdfpagewidth\undefined \else
\fi
 \evensidemargin \oddsidemargin

\advance\hoffset by -1in
\advance\voffset by -1in

\usepackage[cp1251]{inputenc}
\author{И.~А.~Иванов-Погодаев, А.~Я.~Канель-Белов}

\title{Детерминированная раскраска семейства комплексов}

\thanks{Moscow Institute of Physics and Technology, Bar-Ilan University, Мехмат МГУ}

\address{Московский Физико-Технический Институт}
\email{ivanov.pogodaev@mail.ru}
\address{Bar-Ilan University, Israel, Мехмат МГУ}
\email{kanel@mccme.ru}

\begin{document}

\let \mathbf=\texttt
\tabcolsep 2pt

\begin{abstract}

Это вторая часть работы посвященной решению проблемы Л.~Н.~Шеврина и М.~В.~Сапира. Строится конечно определенная бесконечная нильполугруппа, удовлетворяющая тождеству $x^9=0$.

В первой части работы был построено семейство геометрических комплексов, удовлетворяющих набору свойств. В частности, каждый комплекс семейства является {\it равномерно-эллиптическим}:  любые две точки $A$ и $B$ на расстоянии $D$ соединяются системой геодезических,  образующих диск ширины $\lambda\cdot D$ для некоторой глобальной константы $\lambda>0$.

Во второй части доказательства вводится конечная система цветов, обладающих детерминированностью: для каждого минимального квадрата, из которых состоит комплекс, цвета трех углов определяет цвет четвертого угла. 

Настоящая работа посвящена второй части доказательства.

Работа была проведена с помощью Российского Научного Фонда, Грант 22-11-00177. Первый автор является победителем конкурса ``Молодая математика России''.

\end{abstract}

\maketitle

\medskip
%\tableofcontents

\medskip

\section{Введение}

Проблемы бернсайдовского типа внесли огромный вклад в развитие современной алгебры. Эта проблематика охватила большой круг вопросов, как в теории групп, так и в смежных областях, стимулировала алгебраические исследования.
Вопрос о локальной конечности групп с тождеством $x^n = 1$ был решен отрицательно в знаменитых работах П. С. Новикова и С. И. Адян. После этого оценка на экспоненту улучшалась. 
Работы П. С. Новикова и С. И. Адяна оказали огромное влияние на творчество И. А. Рипса, который в дальнейшем разработал метод канонической формы и построил примеры бесконечных периодических групп, обладающих дополнительными свойствами.

\medskip

Все имеющиеся примеры бесконечных периодических групп бесконечно определены.
В частности, вопрос о существовании конечно определенной бесконечной периодической группы открыт как для ограниченной так и для неограниченной экспоненты. 
В кольцах и полугруппах вопросы бернсайдовского типа затрагивают объекты с нильсвойством: то есть такие, в которых каждый элемент в некоторой степени равен нулю.
Также открытым остается поставленный Латышевым вопрос о существовании конечно определенного бесконечного нилькольца. В целом, в настоящее время известно мало примеров конечно определенных алгебраических объектов, обладающих бернсайдовскими свойствами. 

В этой связи интересен метод, позволяющий получать такие конструкции.
Доклад посвящен построению конечно определенной бесконечной нильполугруппы, удовлетворяющей тождеству $x^9 = 0$. Эта конструкция отвечает на проблему поставленную Л. Н. Шевриным и М. В. Сапиром.
 
\medskip

Построение состоит из трех этапов. 
На первом этапе определяется последовательность геометрических комплексов, представляющих собой несколько склееных между собой $4$-циклов (квадратов). 
На втором этапе вершины и ребра комплексов кодируются буквами конечного алфавита. Словам полугруппы будут соответствовать кодировки путей на комплексах. Определяющим соотношениям отвечают локальные эквивалентности двух путей длины $2$, соединяющих противоположные стороны каждого квадрата комплекса. Кроме того, вводятся мономиальные соотношения, приводящие к нулю кодировки невозможных путей, а также путей туда и обратно по некоторому ребру комплекса.
На третьем этапе строится алгоритм, приводящий произвольное слово к канонической форме. При этом помогает геометрический смысл введенных соотношений. Каждое слово локально представляет собой кодировку пути на комплексе. При этом применение определяющих соотношений приводит к локальному преобразованию этого пути, не меняющего его концов и его длину. В процессе локальных преобразований в слове может появиться запрещенный участок, в этом случае оно приводится к нулю. Иначе в итоге мы получаем вложение данного пути в некоторый комплекс.
 
 \medskip
 
От комплекса нужны следующие свойства:

1. Равномерная эллиптичность. Любой путь можно достаточно сильно шевелить локальными преобразованиями;

2. Локальная конечность. Все семейство комплексов допускает раскраску вершин и ребер в конечное множество цветов (букв)

3. Детерминированность. Если известны цвета ребер и вершин вдоль пути, длины соединяющего противоположные вершины некоторого квадрата, то однозначно вычисляются и цвета ребер и вершин вдоль парного к нему пути (по другим двум сторонам).
Это свойство позволяет корректно ввести определяющие соотношения.

4.Апериодичность. Кодировки путей на комплексе не содержат периодических слов периода $9$.
 
В итоге ненулевым словам в полугруппе соответствуют кратчайшие пути на комплексах. Кодировки, соответствующие некратчайшим и невозможным путям, приводятся к нулю, при этом остается бесконечное множество неэквивалентных ненулевых путей.

\section{Содержание первой работы}

В первой части цикла работ <<Комплексы со свойством равномерной эллиптичности>> была изложена геометрическая структура семейства комплексов, служащих базой для введения определяющих соотношений. 

Ниже кратко излагаются основные положения первой части работы. Точные определения и доказательства можно найти в \cite{complex}.

\begin{figure}[hbtp]
\centering
\includegraphics[width=0.5\textwidth]{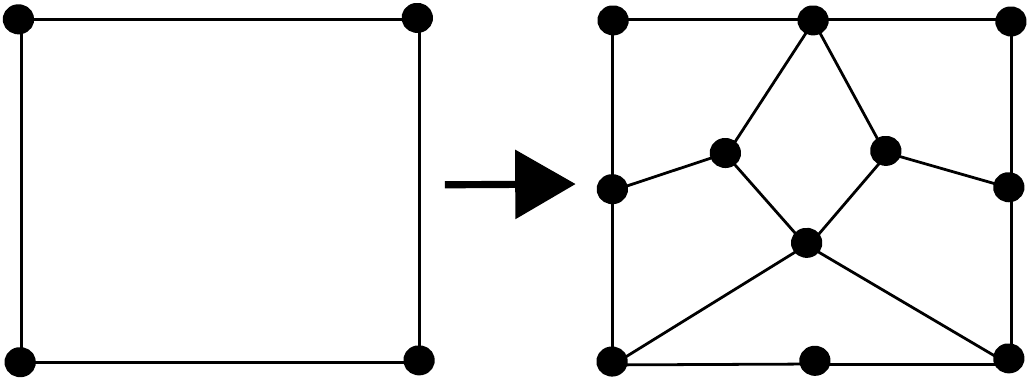}
\caption{Разбиение.}
\label{subst}
\end{figure}

Основной конструкцией работы является семейство геометрических комплексов, каждый из которых состоит из конечного числа квадратов, приклееных друг к другу сторона к стороне специальным образом. Нас будут интересовать свойства комплексов как графов. Комплекс $1$ уровня это просто квадрат из $4$ вершин (плитка). Комплекс $2$ и $3$ уровня изображены на рисунке. Комплекс $3$ уровня получается разбиением каждого квадрата на $6$ квадратов, при этом в середине каждого ребра появляется дополнительная вершина (рисунок~\ref{subst}. 
{\it Макроплиткой $n$ уровня} называется плоская часть комплекса -- результат применения 
$n-1$ уровней разбиения к квадрату из $4$ вершин. 

 Можно заметить, что при этом образуются различные виды вершин --   угловые, краевые, боковые, и внутренние. Угловые и краевые лежат на краю комплекса (в углах, или на сторонах), боковые лежат на ребре, проведенном при некотором разбиении, внутренние -- образуются при некотором разбиении.

\begin{figure}[hbtp]
\centering
\includegraphics[width=0.5\textwidth]{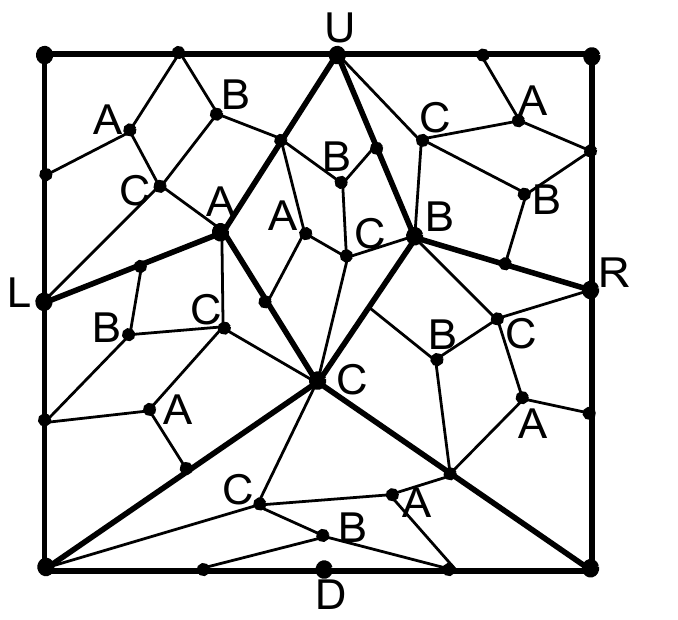}
\caption{Макроплитка третьего уровня. Отмечены типы внутренних вершин ($A$, $B$, $C$)}
\label{fig:level2}
\end{figure}

\begin{figure}[hbtp]
\centering
\includegraphics[width=0.5\textwidth]{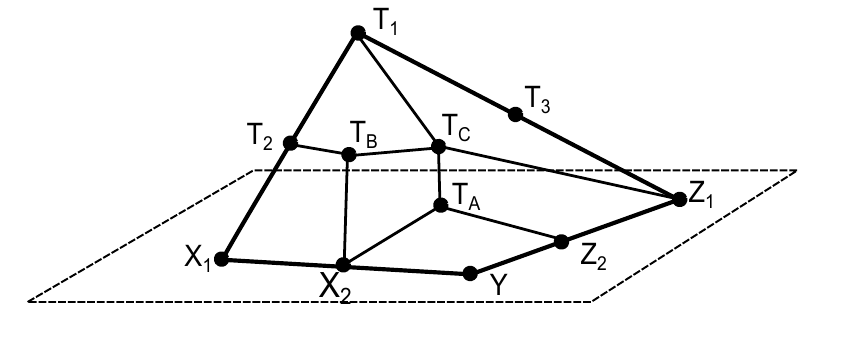}
\caption{Подклейка.}
\label{fig:pasting}
\end{figure}

Далее комплекс уровня $n$ определяется индуктивно. Сначала к комплексу $n$-уровня применяется разбиение -- каждый квадрат делится на $6$ квадратов по правилу, как на рисунке~\ref{fig:level2} потом
проводятся подклейки -- к некоторым путям длины $4$ на комплексе подклеивается $6$ квадратов. (рисунок ~\ref{fig:pasting}.

В результате образуется комплекс следующего уровня. Структура разбиений нужна чтобы обеспечить локальную преобразуемость -- любой путь с концами на границах комплекса может быть преобразован в путь только по границам локальными заменами (когда два соседних ребра некоторого квадрата меняются на два соседних ребра).

Подклейки нужны чтобы обеспечить равномерную эллиптичность комплексов: любой путь на комплексе может быть преобразован локальными заменами достаточно сильно в зависимости от его длины. Например, для естественной метрики, определенной на путях с общими концами, и некоторой константы $\alpha$, путь длины $n$ локально преобразуется в путь, отличающийся на 
$n\alpha$. 

В полугруппе, образованной кодировками путей на семействе комплексов, равномерная эллиптичность означает возможность достаточно сильно преобразововывать произвольные слова. Локальная преобразуемость же позволяет приводить слова к канонической форме.

Также нужны еще несколько свойств, позволяющих ввести раскраску всех вершин и ребер комплексов в конечное число цветов -- они будут играть роль букв в полугруппе. Также часть 
свойств посвящены геометрической структуре комплексов, помогающих приводить кодировки несуществующих на комплексе путей к нулю. 

 {\it Нулевой формой} мы называем путь по одному ребру туда-обратно.
В полугруппе слова, соответствующие {\it коротким} кодировкам перечисленных ниже типов, приравниваются к нулю:

 1) не встречающиеся на комплексах;
 
 2) пути по одному ребру туда-обратно (нулевая форма);
 
 3) принадлежащие к числу <<мертвых паттернов>> -- это пути, которые не могут быть частью достаточно длинного кратчайшего пути на комплексе.
 
 Если в процессе локальных преобразований в слове возникает короткое запрещенное подслово, мы приходим к нулю. Если путь лежит на комплексе, но является некратчайшим путем, соединяющим его концы, то можно локально преобразовать его к нулевой форме и обнулить.

Ниже приводится список основных предложений первой части работы:

\begin{lemma}[Об ограниченности роста степени вершины] \label{growth_bound}

Для каждой вершины $Z$ существует такое натуральное $N$, что
начиная с уровня макроплитки $N$, степень (число входящих ребер) вершины $Z$ не меняется, то есть она одинакова для макроплиток уровня $N$ и $N+k$ для любого натурального $k$.
\end{lemma}

{\bf Следствие.} 1. Каждая вершина заданной глубины $x$ выступает в качестве ядра подклейки для ограниченного количества вершин.

2. В каждую вершину входит ограниченное количество ребер различных уровней, включая ребра из подклеек.

Предложение нужно чтобы обосновать конечность числа возможных типов вершин. Это нужно для введения конечного алфавита полугруппы.

\begin{lemma}[О боковой вершине] \label{side_node}

1) Каждая боковая вершина лежит на середине стороны в какой-либо макроплитке или в двух макроплитках одного уровня, лежащих в одной плоскости;

2) Если боковая вершина не находится на границе исходного комплекса первого уровня, то она лежит на одном из восьми внутренних ребер в некоторой макроплитке, либо лежит на границе подклееной макроплитки.

\end{lemma}

Предложение нужно для описания структуры кодировки типов вершин комплекса.

\begin{lemma}[О выносе пути на границу]  \label{to_border}

Пусть начало и конец пути $P$, проходящего по макроплитке $T$, лежат на границе $T$. Тогда можно реализовать одну из двух возможностей:

1) $P$ может быть локально преобразован в нулевую форму.

2) $P$ может быть локально преобразован в форму $P'$ так, что $P'$ полностью лежит на границе $T$.

Кроме того, любой кратчайший путь, соединяющий противоположные углы или середины противоположных сторон макроплитки имеет длину $2^n$, и может быть локально преобразован в любой из двух путей по границе (полупериметр), с теми же концами.

\end{lemma}

Это удобное предложение для преобразований путей.

\begin{lemma}[О выделении локального участка] \label{longpath}

Пусть $n$ -- уровень макроплитки и путь $P$ лежит внутри нее.

1. Пусть оба края $P$ лежат на границе $T$. Тогда если путь $P$ имеет длину не менее $5\times 2^{n-2}$, он может быть локально преобразован в нулевую форму $P'$.

2. Пусть один край $P$ лежит в углу или на середине стороны $T$, а второй край - внутри $T$ (часть пути может проходить по границе $T$). Тогда если путь $P$ имеет длину $2^{n}$, он может быть локально преобразован в нулевую форму $P'$.

\end{lemma}

Это предложение используется для приведения к каноническому виду. Суть в том, что достаточно длинные пути на небольшом куске комплекса гарантированно являются некратчайшими.

\begin{lemma}[О мертвых паттернах] \label{DeadPaterns}

Рассмотрим некоторую макроплитку $T$ и обозначим в ней внутренние вершины $A$, $B$, $C$ и боковые $U$, $R$, $D$, $L$ (аналогично обозначениям при разбиениях).
Тогда паттерны $AUB$, $ACB$, $CXD$ (где $X$ -- нижняя левая, либо нижняя правая вершина) являются мертвыми.
\end{lemma}

\begin{lemma}[О мертвых путях в нижней подплитке] \label{DeadPaths}

Рассмотрим некоторую макроплитку $T$ уровня $n$. Пусть путь $XYZ$ лежит в $T$, причем $Y$ -- левый нижний угол $T$, $XY$ лежит на внутреннем ребре, идущем из левого нижнего угла $T$ во внутреннюю вершину $C$, а $YZ$ лежит на нижней стороне $T$. Тогда для любых плоских путей $W_1$, $W_2$, длины которых более $2^{n+1}$, путь $W_1XYZW_2$ можно локально преобразовать к нулевой форме.

\end{lemma}

\begin{lemma}[О некорректных участках] \label{uncorrect_sectors}

Пусть есть некорректный участок $XYZ$ в макроплитке $T$ уровня $n$, причем $T$ -- минимальная макроплитка, содержащая $XYZ$ в качестве некорректного участка. Тогда для любых плоских путей $W_1, W_2$, длины более $2^{n+2}$, путь $W_1XYZW_2$, может быть локально преобразован к нулевой форме.
\end{lemma}

Три предложения выше разбирают суть понятия мертвого паттерна и показывают, что некоторые геометрические пути на комплексе могут быть обнулены, как не являющиеся частями длинных кратчайших путей.

\begin{lemma}[О корректности путей] \label{correct_paths}

Пусть путь $P$ представляет собой проход по двум соседним сторонам некоторой макроплитки $T$. Тогда любые локальные преобразования не могут привести $P$ к нулевой форме, а также к форме, содержащей некорректный участок.

\end{lemma}

Это предложение нужно для объяснения бесконечности полугруппы, а именно что некоторые длинные пути на комплексе не могут быть преобразованы к нулю.

\begin{lemma}[О расстоянии от края до выхода в подклейку] \label{pasting_distance}

1. Пусть вершина $X$ лежит на краю некоторой макроплитки $T$ принадлежащей комплексу $K$, вершина $Y$ принадлежит $T$, но не находится на ее границе, а из $Y$ существует выход в подклееную макроплитку уровня $n \geq 2$, не содержащую вершин на границе $T$. Тогда расстояние от $X$ до $Y$ в комплексе $K$ (длина кратчайшего пути по ребрам) не менее $2^{n-1}$.

\end{lemma}

\begin{lemma}[Об ограниченности пути, уходящего в подклееную часть] \label{pastingpath}

1. Пусть путь $P$ имеет вид $V_1AV_2BV_3CV_4$, где ребра из вершин $A$, $B$, $C$ ведут в подклееные области, а участки $V_1$, $V_2$, $V_3$, $V_4$ -- плоские. Тогда путь $P$ не может быть плоским (лежать в одной макроплитке).

2.Пусть путь $P$ начинается в вершине $X$, выход из $X$ идет в подклееную макроплитку уровня $n$. Кроме того, пусть $P$ не содержит ребер-выходов из подклееных плиток.
Тогда, если длина пути не менее  $2^{n+1}$, то он может быть приведен к нулевой форме.

%Пусть путь $P$ имеет форму $XP_1YP_2$, где $P_1Y$ -- плоский участок (лежащий в какой-то макроплитке) и выходящие из вершин $X$ и $Y$ ребра идет в подклееные плитки, и больше выходов в подклееную область или входов в вершину из подклееной области нет. Пусть также $P$ не может быть приведен к нулевой форме.
%Тогда длина $YP_2$ не более, чем удвоенная длина $XP_1Y$.
\end{lemma}

Эти предложения нужны для доказательства, что периодические пути с периодом $9$ гарантированно приводятся к нулю.

\medskip

\section{Строение полугруппы путей} \label{pathsemigroup}

В предыдущей части работы мы описали геометрическую структуру семейства комплексов.  В следующей главе~\ref{coding_section} будет показано, как связать с ним специальную кодировку. Мы покажем, что с помощью конечного числа букв можно закодировать узел. Кроме того, конечным числом букв можно закодировать входящие и выходящие ребра в каждый узел, как плоские, так и выходящие в подклейку.

В этом разделе мы покажем, для чего нужна такая кодировка и как необходимая нам полугруппа строится с помощью кодировки путей на комплексе.

Пусть $X_i$ -- буквы, кодирующие входящие ребра, $Y_i$ -- буквы, кодирующие узлы (их типы, окружения и информации), $Z_i$ -- буквы, кодирующие выходящие ребра. (Детали кодирования будут показаны в последующих разделах.)

Будем говорить, что слово $W$ имеет  $\mathbf{CODE}$-{\it форму}, если в нем справа от любой (не последней в слове) буквы семейства $X$ обязательно стоит буква семейства $Y$, справа от любой (не последней в слове)  буквы семейства $Y$ стоит буква семейства $Z$, а справа от любой (не последней в слове)  буквы семейства $Z$ стоит буква семейства $X$.

Часть таких $\mathbf{CODE}$-слов соответствует путям на комплексе.
Рассмотрим полугруппу с нулем $S$ с порождающими $\{X_i,Y_i,Z_i \}$.

\medskip

\begin{lemma}[О стандартной форме слова] \label{standardform}

В полугруппе с нулем $S$ с порождающими $\{X_i,Y_i,Z_i \}$ можно ввести конечное число определяющих мономиальных соотношений так, чтобы все не $\mathbf{CODE}$-слова, можно было привести к нулевому слову.

\end{lemma}

\begin{proof} Действительно, введем соотношения $X_iZ_j=0$, $X_iX_j=0$, $Y_iX_j=0$, $Z_iY_j=0$, $Y_iY_j=0$, $Z_iZ_j=0$, для всевозможных пар $(i,j)$. В получившейся полугруппе, любое ненулевое слово будет иметь $\mathbf{CODE}$-форму: за буквой $X_i$ может следовать только $Y_j$, за $Y_j$ -- только $Z_k$, за $Z_k$ -- только $X_m$.

\medskip

Заметим также, что можно выписать все $\mathbf{CODE}$-слова длины не более $4$, которые кодируют какой-нибудь путь на комплексе. Пусть $F_1, \dots ,F_N$ -- все оставшиеся $\mathbf{CODE}$-слова длины не более $4$. Их мы будем считать {\it запрещенными} и введем в полугруппе $S$ соотношения $F_i=0$, для всех $i=1,\dots ,N$.

Таким образом, любое слово $W$ в полугруппе $S$ либо приводится к нулю, либо имеет $\mathbf{CODE}$-форму, где каждое подслово $W$ длины не более $4$ отвечает некоторому пути на комплексе.

\end{proof}

В последующих разделах будет показано, что для двух эквивалентных путей, не содержащих мертвых паттернов и нулевых форм, по коду одного из них можно восстановить код другого. Это дает возможность ввести конечное множество соотношений вида $W_i=W_j$, где $W_i,W_j$ -- коды двух эквивалентных путей.

Таким образом, появляется возможность преобразовывать слова, меняя подпути на им эквивалентные.
При этом если до преобразования слово являлось кодировкой некоторого пути на комплексе, то и после преобразования слово оно будет кодировкой эквивалентного ему пути. Допустим, после нескольких таких замен образуется запрещенное подслово. Это значит, что соответствующий этой кодировке путь не может быть реализован на комплексе. Но тогда и эквивалентный ему путь не мог быть реализован.
Таким образом, если, преобразовывая слова-кодировки по установленным правилам, мы получаем запрещенное слово, значит изначальное слово не соответствует пути на комплексе.

\medskip

{\bf Нулевая форма.} Если слово $W$ представляет собой кодировку пути, имеющего {\it нулевую форму}, мы вводим в полугруппе соотношение $W=0$. Рассмотрим слово-кодировку пути, концы которого лежат на периметре некоторой макроплитки $T$, и хотя бы один конец попадает в угол $T$. Заметим, что если путь не кратчайший (то есть существует путь короче по длине с теми же концами), то согласно лемме 3.6 о выделении локального участка из первой части работы, путь может быть преобразован в нулевую форму. Это означает, что соответствующее слово-кодировка приводится к нулю.

Фактически, обнуление нулевых форм позволяет приводить к нулю все слова, соответствующие некратчайшим путям.

\medskip

{\bf Мертвые паттерны.} Мертвые паттерны представляют собой некоторые типы слов-кодировок, которые никогда не могут встретиться в достаточно больших ненулевых словах. Нас будут интересовать только мертвые паттерны $AUB$, $ACB$, $CXD$, $BUA$, $BCA$, $DXC$, которые обсуждаются в лемме~3.8 о мертвых паттернах из первой части. Выпишем все слова, длины не более $4$, содержащие хотя бы один из перечисленных мертвых паттернов. Для каждого такого слова $W$, введем соотношение $W=0$.

Заметим, что все достаточно большие ненулевые слова не могут содержать мертвый паттерн (и не могут быть приведены к форме, его содержащей).

\medskip

Будем называть слово {\it регулярным}, если оно является кодировкой некоторого пути на построенном геометрическом комплексе и при этом не может быть преобразовано в нулевую форму или форму содержащую мертвый паттерн.

\medskip

Заметим, что замена некоторого подслова регулярного слова на эквивалентное ему, не нарушает регулярность. То есть, если слово можно преобразовать так, что оно будет содержать запрещенное подслово, нулевую форму или мертвый паттерн, то значит оно приводится к нулю, и не может быть регулярным.

Нашей основной целью будет показать, что ненулевыми словами в полугруппе являются только регулярные слова, то есть, что любое нерегулярное слово может быть приведено к нулю. Фактически, это будет означать, что ненулевыми элементами полугруппы являются только слова, кодирующие кратчайшие пути, не содержащие мертвые паттерны.

\medskip

{\bf Примечание о использовании мертвых паттернов для конструкции}. Ценой некоторого усложнения кодировки можно добиться того, чтобы пути, содержащие мертвые паттерны, допускали такие же локальные преобразования, как и обычные слова. В этом случае, можно не вводить мономиальные соотношения для слов, их содержащих. И тогда ненулевыми элементами полугруппы будут просто кратчайшие пути на построенном комплексе. Но в целях упрощения конструкции, представляется разумным ввести такие соотношения.

%Фактически, множество мертвых паттернов представляет собой радикал в полугруппе путей.

\medskip

\section{Кодировка вершин и путей на комплексе} \label{coding_section}

Рассмотрим путь на комплексе. Для каждой вершины из этого пути, кроме первой и последней, есть ребро входа и ребро выхода. Ребра входа и выхода могут вести в подклееные области, а могут принадлежать той же базовой плоскости, что и сама вершина.

Каждая вершина принадлежит некоторой макроплитке. Кроме того из нее могут выходить ребра ведущие в подклееные плитки/макроплитки. То есть, каждая вершина лежит на своей базовой плоскости и еще участвует в нескольких подклееных плоскостях. Фактически, вершина играет свою роль (занимает определенное положение) в своей макроплитке на базовой плоскости, и кроме этого выполняет другие какие-то роли в своих подклееных макроплитках. При этом в подклееных макроплитках эта вершина лежит всегда на краю.

Мы хотим определить систему кодирования путей на комплексе. Каждый путь это  конечная последовательность вершин, причем любые две соседние соединены ребром плитки минимального уровня (не макроплитки). Для кодирования путей нужно сначала закодировать все вершины, которые могут встретиться на этом пути, конечным числом букв.

Каждая вершина может принимать участие во множестве подклееных макроплиток, то есть входящие в нее ребра либо относятся к плоскости вершины либо классифицируются по принадлежности к различным подклееным макроплиткам. В общем случае, путь приходит в вершину из одной подклееной макроплитки и выходит в другую. Общая кодировка вершины будет состоять из трех {\it плоских} кодов: первый представляет собой плоский код вершины (без учета подклееных плиток) второй -- кодирует положение этой же вершины в плоскости подклееной макроплитки, откуда приходит входящее ребро, и третий кодирует положение этой же вершины в плоскости подклееной макроплитки, куда уходит выходящее ребро.

\medskip

В конечном счете, наша цель - задать раскраску вершин и ребер всего семейства комплексов в конечное множество цветов. Раскраска должна удовлетворять следующему свойству детерминированности:

\begin{theorem}[О введении раскраски со свойством детерминированности] \label{determ}

Введенная система цветов для ребер и вершин семейства комплексов обладает свойством детерминированности: если $A$, $B$, $ C$ - вершины некоторого минимального 4-цикла $ABCD$ на комплексе, и при этом путь $ABC$ не образует мертвый патерн, то по кодировке пути $ABC$ можно однозначно определить кодировку пути $ADC$.

\end{theorem}

Данное свойство детерминированности позволяет корректно ввести определяющие соотношения, приравнивающие друг к другу кодировки путей $ABC$ и $ADB$.

Вся данная работа посвящена введению такой раскраске и проверке этого свойства для различных случаев расположения путей.

\medskip

Приступим к введению раскраски вершин и ребер.

Если рассматривать вершину только в рамках одной макроплитки, в которую она входит, имеет смысл говорить о {\it плоском} коде. Например, если весь путь лежит в одной плоскости и не выходит в подклееные области, то код такого пути будет плоским. Плоский код одной вершины будет состоять из следующих частей:

1) тип вершины;

2) уровень вершины;

3) окружение;

4) информация.

Мы последовательно определим, что значит каждая из этих частей.

\medskip

\subsection{Параметр ``тип'' для вершин}

Тип вершины соответствует этому определению из первой части работы, вершины разделены на следующие категории:

1) {\it Угловые}. (Лежащие в углах подклееных макроплиток или всего комплекса). {\it Тип угловой вершины } определим как один из четырех вариантов углов, в котором она может находиться: $\mathbb{CUL}$, $\mathbb{CUR}$, $\mathbb{CDR}$, $\mathbb{CDL}$. (Corner Up-Left и так далее.) Угловые вершины {\it принадлежат} макроплитке, где они являются углами.

2) {\it Краевые}. (Лежащие на стороне подклееной макроплитки или всего комплекса). Каждая такая вершина лежит в середине стороны некоторой макроплитки, прилегающей к краю. Тип краевой вершины  определим в соответствии с тем, серединой какой стороны в этой макроплитке она является: $\mathbb{L}$, $\mathbb{R}$, $\mathbb{D}$, $\mathbb{U}$.
Краевые вершины {\it принадлежат} макроплитке, на краю которой они лежат.

3)  {\it Внутренние}. В этой категории определим три типа внутренних вершин: $\mathbb{A}$, $\mathbb{B}$, $\mathbb{C}$, отвечающих внутренним узлам макроплиток, эти вершины создаются {\it внутри} разбиваемой макроплитки. Будем считать, что эти вершины {\it принадлежат} данной макроплитке.

4)  {\it Боковые}. (Лежащие на границе между двумя макроплитками, на внутреннем ребре).
  Типы боковых вершин будут соответствовать всем упорядоченным парам из множества $\{ \mathbb{U},\mathbb{R},\mathbb{D},\mathbb{L}\}$. А именно: $\mathbb{DR}$, $\mathbb{RD}$, $\mathbb{DL}$, и так далее. Будем считать, что в упорядоченной паре первой называется буква, соответствующая $A$-стороне внутренего ребра, на котором лежит вершина.  Все боковые вершины создаются в середине стороны разбиваемой макроплитки. 
Тип боковой вершины определяет, серединой каких именно сторон она является в двух макроплитках, где она является серединой сторон. Это как раз те макроплитки, которые разбивались при создании данной вершины.
Будем считать, что боковые вершины {\it принадлежат} макроплитке, которой принадлежит данное внутреннее ребро.

\begin{figure}[hbtp]
\centering
\includegraphics[width=0.9\textwidth]{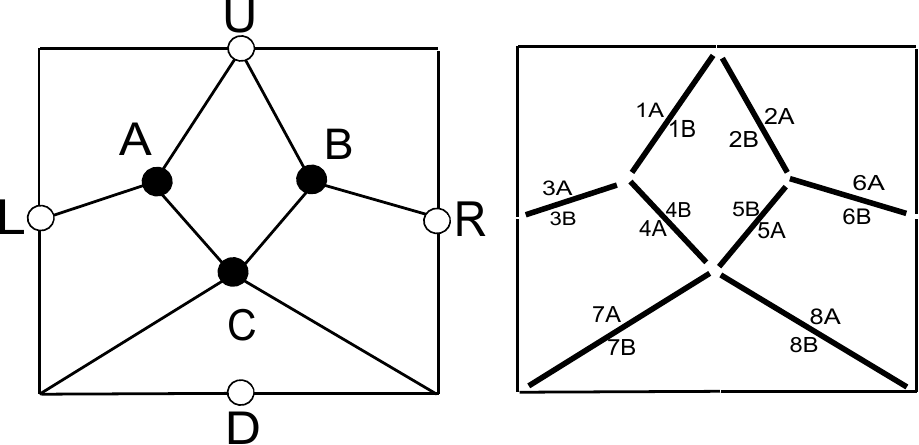}
\caption{Черные (внутренние) и белые (боковые) вершины; типы ребер}
\label{fig:blacknodes}
\end{figure}

В дальнейшем, вершину типа $\mathbb{A}$ будем, для простоты, называть $\mathbb{A}$-вершиной или $\mathbb{A}$-узлом. Аналогично для других типов. Вообще, вершины на графе мы также будем называть узлами.
В левой части рисунка~\ref{fig:blacknodes} представлены внутренние и боковые вершины. Следующая лемма показывает, что возможны только указанные сочетания сторон.

\begin{lemma}[О комбинациях сторон на границе макроплиток] \label{side_combos}

Пусть вершина $X$ лежит в середине сторон двух макроплиток, лежащих в одной плоскости. Тогда возможны только следующие сочетания ориентаций сторон: правая и нижняя, верхняя и правая,  левая и нижняя, либо верхняя и левая.
\end{lemma}

\begin{proof} По второй части леммы о боковой вершине, $X$ лежит на некотором ребре $YZ$, которое относится к одному из восьми типов. (Правая часть рисунка~\ref{fig:blacknodes}). Рассмотрим вершину $F$, лежащую в середине этого ребра.  $F$ тоже боковая вершина и сочетания ориентаций сторон у нее могут быть только право и низ, верх и право, лево и низ, верх и лево, то есть как в условии леммы. Если $F$ совпадает с $X$, то все доказано. Иначе рассмотрим разбиение наших соприкасающихся макроплиток.  Легко проверить, что, согласно правилам разбиения, середины ребер $YF$ и $FZ$ тоже могут иметь сочетание ориентаций только право и низ, верх и право, лево и низ, верх и лево ($\mathbb{RD}$, $\mathbb{UR}$, $\mathbb{LD}$, $\mathbb{UL}$ и транспозиции). Например, если макроплитки соприкасаются нижней и правой сторонами, то $F$ имеет тип $\mathbb{DR}$ или $\mathbb{RD}$, и тогда $YF$ и $FZ$ будут иметь типы $\mathbb{RU}$ и $\mathbb{UL}$ (или наоборот). Повторяя операцию разбиения, можно получить, что все боковые вершины на ребре $YZ$ будут относиться к одному из указанных четырех типов, то есть являются серединами правой и нижней, либо верхней и правой, либо левой и нижней, либо верхней и левой сторон.

\end{proof}

\subsection{Типы ребер и параметр ``уровень'' для вершин}

Будем считать, что вершина {\it принадлежит} макроплитке, если она образуется при первом разбиении этой макроплитки (когда она из плитки становится макроплиткой второго уровня). То есть, белые вершины принадлежат тем макроплиткам, в середине сторон которых они лежат, а черные -- макроплиткам минимального уровня, внутри которых они находятся.

\medskip

{\bf Типы ребер.} Всего существует $16$ внутренних ребер (с учетом стороны, рисунок~\ref{fig:blacknodes}), это ребра $\mathbf{1A}$, $\mathbf{1B}$, \dots ,$\mathbf{8B}$ и $4$ краевых ребра: $\mathbf{left}$, $\mathbf{top}$, $\mathbf{right}$, $\mathbf{bottom}$.

\medskip

{\bf Уровни вершин.}
Для черных, внутренних вершин уровень определять не будем (то есть, можно считать, что у всех один уровень).

%Для черных, внутренних вершин будет три возможных уровня. Первый будет у вершин, принадлежащих макроплиткам второго уровня (самым маленьким макроплиткам), второй у вершин, принадлежащих макроплиткам третьего уровня, и третий -- у всех остальных черных вершин.

%У угловых вершин уровень определим аналогично: первый у вершин, принадлежащих макроплиткам второго уровня (только что подклееная макроплитка), второй у вершин, принадлежащих макроплиткам третьего уровня, и третий у всех остальных угловых вершин (принадлежащих макроплиткам четвертого уровня и выше)

Для белых, боковых и краевых вершин уровня будет три: первый у вершин, принадлежащих макроплиткам второго уровня, второй -- у вершин, принадлежащих макроплиткам третьего уровня, третий -- у всех остальных белых вершин.

\medskip

Рассмотрим черную вершину $X$. Если $X$ имеет тип $\mathbb{A}$ и $\mathbb{B}$, то она находится на стыке трех макроплиток (разбивающих макроплитку, которой $X$ принадлежит), для $\mathbb{C}$ таких макроплиток четыре. В соответствии с определением 3.8 из первой части, ребра, исходящие из $X$ и при этом лежащие на границе каких-то двух из этих макроплиток, являются {\it главными ребрами} для $X$.

Если  $X$ -- белая вершина, лежащая на середине сторон двух макроплиток $T_1$ и $T_2$, то два ребра, лежащих на границе между $T_1$ и $T_2$, тоже называются {\it главными}. Если боковая вершина лежит на краю только одной макроплитки, то главные ребра -- это лежащие на ее границе.

 В зависимости от типа $X$, из нее выходят еще несколько неглавных ребер. 

%Если $X$ -- принадлежит макроплитке второго уровня, то все выходящие из нее ребра главные, их три для $\mathbb{A}$ и $\mathbb{B}$ и четыре для $\mathbb{C}$. Если $X$ -- принадлежит макроплитке третьего уровня, для $\mathbb{A}$, $\mathbb{B}$, $\mathbb{C}$ к главным ребрам добавляются $2$, $3$, $3$ неглавных соответственно.
%Для вершины $\mathbb{С}$ принадлежащей макроплитке четвертого уровня и выше, добавляются еще три неглавных ребра. Рассмотрим два соседних исходящих из $X$ ребра. Они ограничивают некоторую плитку или макроплитку, причем $X$ является одним из четырех ее углов. Легко проверяется, что этот угол при $X$ является простым, то есть не между нашими соседними ребрами не добавится еще ребро при дальнейших разбиениях. Таким образом, при дальнейших разбиениях исходящих ребер для $X$ не добавится, и поэтому все черные вершины можно классифицировать по описанным трем уровням.

\medskip

 %Мы можем разобрать все белые вершины третьего уровня и убедиться, что все прилежащие к $X$ углы будут простые и при дальнейших разбиениях новых исходящих ребер не добавляется.

\medskip

{\bf Ребра входа и выхода. Обозначения.}
{\it Первым главным ребром} для вершин типа $\mathbb{A}$ и $\mathbb{B}$ будем называть то главное ребро, которое уходит к середине верхней стороны той макроплитки, которой принадлежит данная вершина.
Для $\mathbb{C}$ первым главным ребром будем считать ребро, идущее в $\mathbb{A}$-узел.

\medskip

Для боковых вершин оба главных ребра являются двумя частями некоторого внутреннего ребра. Первым из них будем считать то ребро, относительно которого А-сторона (из определения типа внутреннего ребра) остается по правую сторону. Второе и третье главные ребра определяются по часовой стрелке после первого. Все главные ребра будем обозначать цифрами (от $1$ до $4$).

\medskip

Для краевых вершин первым главным ребром будем считать то ребро, которое соответствует обходу макроплитки по часовой стрелке. То есть, макроплитка остается по правую сторону от направления этого ребра. Ребро, идущее против часовой стрелке, будет вторым.

\medskip

У вершин типа $\mathbb{A}$ может быть два неглавных ребра, одно уходит в левую верхнюю подплитку, другое в левую нижнюю. Обозначим их, соответственно, $\mathbf{lu}$  и $\mathbf{ld}$.
У вершин типа $\mathbb{B}$ может быть три неглавных ребра, одно уходит в правую верхнюю подплитку, второе в правую нижнюю, третье в среднюю. Обозначим их, соответственно, $\mathbf{ru}$, $\mathbf{rd}$, $\mathbf{mid}$.
У вершин типа $\mathbb{C}$ может быть шесть неглавных ребра, два уходят в левую нижнюю подплитку, два в среднюю, и два в нижнюю. В каждой паре одно из ребер иерархически старше (появляется при разбиении раньше). Обозначим ребра, соответственно, $\mathbf{ld}_1$, $\mathbf{ld}_2$, $\mathbf{mid}_1$, $\mathbf{mid}_2$, $\mathbf{d}_1$, $\mathbf{d}_2$. Более старшему ребру даем первый номер, другому -- второй.

\begin{figure}[hbtp]
\centering
\includegraphics[width=0.5\textwidth]{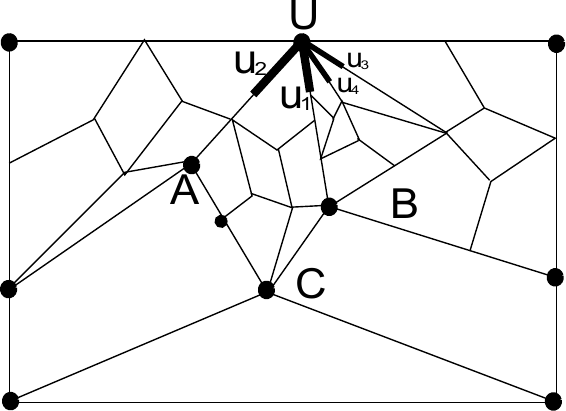}
\caption{Неглавные ребра выходящие из бокового узла}
\label{fig:upedges}
\end{figure}

Для боковых и краевых вершин будем давать неглавным ребрам имена в соответствии с типом вершины и стороны, откуда они уходят:

$\mathbb{U}$-узел или $U$-часть $\mathbb{UL}$-узла или $\mathbb{UR}$-узла -- это $\mathbf{u}_1$, $\mathbf{u}_2$, $\mathbf{u}_3$, $\mathbf{u}_4$;

$\mathbb{R}$-узел или $R$-часть $\mathbb{UR}$-узла или $\mathbb{DR}$-узла -- $\mathbf{r}$, $\mathbf{r}_2$, $\mathbf{r}_3$;

$\mathbb{L}$-узел или $L$-часть $\mathbb{UL}$-узла или $\mathbb{LD}$-узла -- $\mathbf{l}$, $\mathbf{l}_2$, $\mathbf{l}_3$.

В $\mathbb{D}$-узел неглавные ребра не приходят.

Ребра нумеруются по иерархическому старшинству, при равном старшинстве -- по часовой стрелке. Например, есть четыре неглавных ребра из узла типа $\mathbb{LU}$, входящих внутрь плитки $T$, примыкающей верхней стороной (рисунок~\ref{fig:upedges}). $\mathbf{u}_1$ будет ребром в сторону $\mathbb{B}$-узла $T$, $\mathbf{u}_2$ -- ребром в сторону $\mathbb{A}$-узла $T$. $\mathbf{u}_3$ и $\mathbf{u}_4$,  будут входить в верхнюю правую подплитку $T$.

\medskip

{\bf Типы ребер, ведущих в подклееные области.}
Рассмотрим вершину $X$. В подклееной макроплитке она может иметь один из следующих типов:
$\mathbb{L}$, $\mathbb{U}$, $\mathbb{CUR}$, $\mathbb{CUL}$, $\mathbb{CDL}$.

\begin{figure}[hbtp]
\centering
\includegraphics[width=0.9\textwidth]{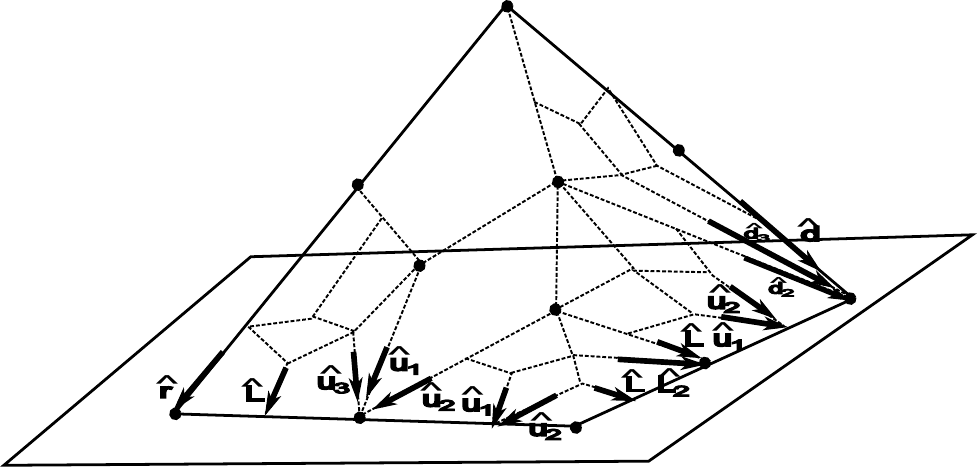}
\caption{Ребра в подклееные области}
\label{fig:pastingedges}
\end{figure}

Из вершины $\mathbb{L}$ могут выходить ребра $\mathbf{l}$, $ \mathbf{l}_2 $, $\mathbf{l}_3$ (как и в плоском случае). Чтобы показать, что имеются в виду ребра в подклееную область, будем записывать их как $\widehat{\mathbf{l}_1}$, $\widehat{\mathbf{l}_2}$, $\widehat{\mathbf{l}_3}$ (рисунок~\ref{fig:pastingedges}).

Аналогично, из вершины $\mathbb{U}$ могут выходить ребра $\widehat{\mathbf{u}_1}$, $\widehat{\mathbf{u}_2}$, $\widehat{\mathbf{u}_3}$, $\widehat{\mathbf{u}_4}$.

Из вершины $\mathbb{CUR}$ может выходить только одно ребро в подклееную область -- это ребро по правой стороне подклееной макроплитки. Будем обозначать его как $\widehat{\mathbf{r}}$.

Из вершины $\mathbb{CDL}$ может выходить два ребра в подклееную область -- по нижней стороне подклееной макроплитки, по ребру $7$ подклееной макроплитки, а также по $8$ ребру макроплитки, образующейся при разбиении нижней подплитки подклееной макроплитки. Будем обозначать их как $\widehat{\mathbf{d}}$, $\widehat{\mathbf{d}_2}$, $\widehat{\mathbf{d}_3}$ соответственно.

Из вершины $\mathbb{CUL}$ (она же -- ядро подклееной макроплитки) ребер в подклееные области не выходит.

\medskip

{\bf Цепи.}
Выберем некоторую вершину $X$. Обозначим как $T(X)$ множество всех макроплиток одного размера, для которых $X$ является левым верхним углом. ($T(X)$ может быть выбрано несколькими способами, для разных размеров макроплиток.)

\begin{definition}

Совокупность ребер $1$ и $3$ типов в макроплитках из $T(X)$ и {\bf боковых} узлов в их концах будем называть {\it цепью} для $X$. Будем называть вершину $X$ {\it центром} цепи.

\end{definition}

Из определения следует, что каждый узел может входить только в одну цепь, причем центр для этой цепи определяется однозначно.
Под {\it макроплитками из цепи} будем понимать макроплитки из соответствующего множества $T(X)$.
Заметим, что между ребрами, выходящими из $X$ и вершина цепи, а  то есть, количество узлов в любой цепи ограничено сверху (Напомним, что в цепь могут входить только боковые вершины).

\smallskip

Для вершин в одной цепи можно задать упорядоченность: первой будем считать тот узел, который лежит на первом главном ребре, выходящем из $X$, а остальные перечисляются в соответствии с выходящими из $X$ ребрами по часовой стрелке. Номер вершины в ее цепи будем называть {\it указателем}.

\medskip

Для выбранной вершины $X$ можно построить несколько цепей с центром $X$ (в каждой будут макроплитки одного размера). Каждой цепи можно присвоить уровень: для самой крупной цепи -- нулевой ($0$-цепь), далее первый и так далее. Если узел $X$ принадлежит некоторой макроплитке $T$, то ребра цепи уровня $k$ принадлежат макроплитке, получающейся после $k+1$ разбиений $T$.

\medskip

\begin{figure}[hbtp]
\centering
\includegraphics[width=1\textwidth]{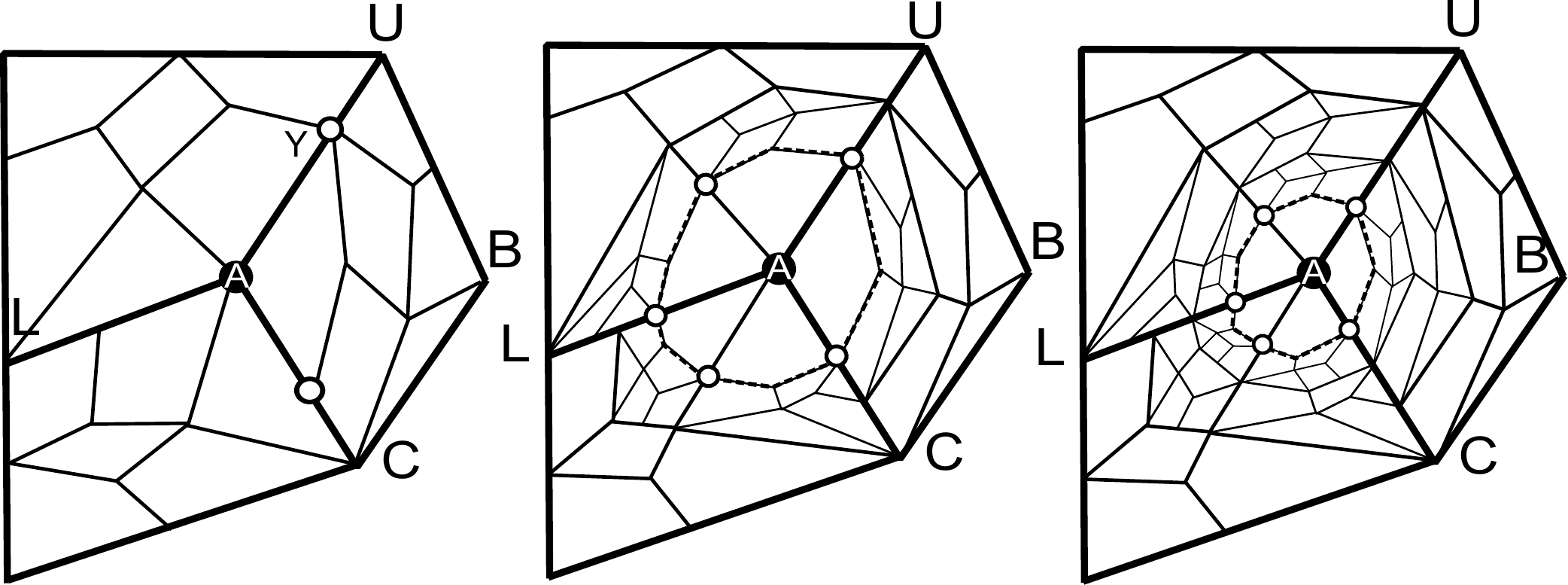}
\caption{Цепи узла $\mathbb{A}$.}
\label{fig:chainAdef}
\end{figure}

На рисунке показаны цепи нулевого, первого и второго уровней c центром типа $\mathbb{A}$.

\medskip

{\bf Параметр ``окружение'' для узлов.}

\begin{definition}

Для каждой макроплитки можно записать упорядоченную четверку типов ребер, на котором лежат ее стороны: левая, верхняя, правая, нижняя (строго в данном порядке). Такую четверку будем называть {\it окружением} данной макроплитки.

\end{definition}

\medskip

Рассмотрим некоторую вершину $X$. Она может входить в одну или несколько подклееных областей. Рассмотрим одну из таких подклееных областей $P$.

\begin{definition}

 {\it Меткой вершины относительно $P$} будем называть:

{\bf 1.} Для регулярных внутренних вершин (типы $\mathbb{A}$, $\mathbb{B}$, $\mathbb{C}$) -- окружение макроплитки, которой принадлежит данная вершина.

\smallskip
 {\bf Пример:}  {\it У узла типа $\mathbb{A}$ внутри средней подплитки метка $[\mathbf{1B},\mathbf{2B},\mathbf{4B},\mathbf{5B}]$}.

\smallskip

{\bf 2.} Для регулярных боковых вершин (типы $\mathbb{UL/LU}$, $\mathbb{UR/RU}$, $\mathbb{LD/DL}$, $\mathbb{RD/DR}$) -- упорядоченную пару окружений двух макроплиток, в середине сторон которых лежит данная вершина.

\smallskip

    {\bf Пример: } {\it У узла $Y$ на левой части рисунка~\ref{fig:chainAdef} метка $[x,y,\mathbf{1A},\mathbf{3A}]-[\mathbf{1B},\mathbf{2B},\mathbf{4B},\mathbf{5B}]$}, где $x$, $y$ -- {\it типы} левой и верхней сторон всей макроплитки рисунка.

\smallskip

{\bf 3.} Для угловых вершин -- метка будет ($\mathbf{left}$,$\mathbf{top}$,$\mathbf{right}$,$\mathbf{bottom}$), то есть как окружение подклееной макроплитки.

{\bf 4.} Для краевых вершин -- метка соответствует окружению макроплитки, в середине стороны которой лежит вершина.

\end{definition}

\medskip

{\bf Замечание}. У вершины может быть несколько меток -- по числу подклееных областей, в которые она входит. Ясно, что только в своей базовой плоскости вершина может быть регулярной. В подклееных областях она будет обязательно угловой или краевой.

\medskip

\begin{definition}

{\it Окружением цепи} будем называть упорядоченный набор меток всех вершин цепи.

\end{definition}

\smallskip

  {\bf Пример.} {\it У цепи в средней части рисунка~\ref{fig:chainAdef} окружение

   $$ [\mathbf{1A},\mathbf{8A},\mathbf{5A},\mathbf{6B}]:[\mathbf{2B},\mathbf{4B},\mathbf{1A},\mathbf{3A}]\times [\mathbf{2B},\mathbf{4B},\mathbf{1A},\mathbf{3A}]:[\mathbf{8B},\mathbf{4A},\mathbf{4A},\mathbf{7B}]\times $$

   $$ \times [\mathbf{8B},\mathbf{4A},\mathbf{4A},\mathbf{7B}]:[\mathbf{3B},\mathbf{8A},\mathbf{5A},\mathbf{6B}]\times [\mathbf{3B},\mathbf{8A},\mathbf{5A},\mathbf{6B}]:[\mathbf{8B},\mathbf{3A},\mathbf{3A},\mathbf{7B}] \times  $$

   $$\times [\mathbf{8B},\mathbf{3A},\mathbf{3A},\mathbf{7B}]:[\mathbf{1A},\mathbf{8A},\mathbf{5A},\mathbf{6B}]$$}

\smallskip

Ниже будет показано, что тип узла в центре цепи уровня выше нулевого полностью определяет окружение этой цепи.

\medskip

\begin{definition}

{\it Окружением вершины} относительно плоскости $P$ будем называть:

{\bf 1.} Для регулярных боковых вершин типов $\mathbb{UL/LU}$, $\mathbb{UR/RU}$, $\mathbb{LD/DL}$, а также для краевых вершин типов $\mathbb{U}$ и $\mathbb{L}$  -- упорядоченную пару $(x,y)$, где $x$ --  окружение цепи, в которую входит данная вершина, $y$ -- указатель, номер данной вершины в ее цепи.

\smallskip

{\bf 2.} Для регулярных боковых вершин типов $\mathbb{RD/DR}$ -- тип ребра, на которой находится данная вершина (это ребро $2$, $3$, $5$ или $6$).

{\bf 3.} Для всех остальных вершин (угловые, краевые, регулярные типов $\mathbb{A}$, $\mathbb{B}$, $\mathbb{C}$) -- метку этой вершины относительно плоскости $P$.

\end{definition}

Таким образом, если вершина входит в цепь, то зная окружение вершины, мы знаем окружения всех вершин в этой цепи.
Окружение вершины $X$ относительно плоскости $P$ будем обозначать как $\mathbf{SurrP}(X)$. Если из контекста ясно, о какой плоскости идет речь, будем использовать обозначение $\mathbf{Surr}(X)$.

\medskip

Пусть $X$ -- узел, входящий в некоторую цепь. Если $Y$ -- центр этой цепи, то $X$ соответствует одному из выходящих из $Y$ ребер. Номер этого ребра, (начиная с первого выходящего из $X$) мы называем {\it указателем} для вершины $X$. Очевидно, что по указателю и окружению цепи можно установить окружение самого узла.

\medskip

\begin{definition}

И наконец, определим {\it расширенное окружение} относительно области~$P$. Это упорядоченная пара окружения относительно $P$ и окружения относительно базовой плоскости вершины.
\end{definition}

Таким образом, расширенное окружение несет информацию не только об окрестности вершины в подклейке, но и об окрестности в базовой плоскости.

\medskip

\subsection{Параметр ``информация'' для узлов}

{\bf Начальники.}

Определим понятие {\it начальников} для регулярных вершин.
Пусть вершина $X$ -- регулярная ( то есть не угловая, не боковая). Тогда можно выбрать такую минимальную макроплитку $T$, что $X$  находится внутри (не на границе) $T$. В этом случае, $X$ либо является внутренней вершиной ($\mathbb{A}$, $\mathbb{B}$ или $\mathbb{C}$), либо лежит на одном из восьми внутренних ребер $T$.

\begin{definition}

{\it Первым начальником} вершины $X$ будем называть узел в середине верхней стороны макроплитки $T$.

Для боковых вершин, лежащих на ребрах $2$, $5$, $6$ макроплитки $T$, {\it вторым начальником} будем называть узел в правом нижнем углу макроплитки $T$.

Для боковых вершин, лежащих на ребрах $7$ и $8$, а также для для вершин типа $\mathbb{C}$ определим {\it второго} и {\it третьего начальников} как узлы в левом нижнем и правом нижнем углах $T$ соответственно.

\end{definition}

У краевых вершин начальников не будет. Из угловых вершин определим начальника только для $\mathbb{CDR}$-узла, и им будет $\mathbb{CDL}$-узел в той же макроплитке.

\medskip

Теперь дадим определение информации. Оно будет зависеть от расположения вершины внутри макроплитки.

\begin{definition}

Для регулярных боковых вершин, лежащих на внутренних ребрах $1$, $3$, $4$, а также для $\mathbb{A}$-узлов {\it информацией вершины} будем называть упорядоченную тройку состоящую из типа, уровня и расширенного окружения ее единственного начальника.

Для регулярных боковых вершин, лежащих на внутренних ребрах $2$, $5$, $6$, а также для $\mathbb{B}$-узлов {\it информацией вершины} будем называть упорядоченную четверку состоящую из типа, уровня и расширенного окружения ее первого начальника, а также типа ее второго начальника.

Для регулярных боковых вершин, лежащих на внутренних ребрах $7$, $8$, а также для $\mathbb{C}$-узлов {\it информацией вершины} будем называть упорядоченную девятку состоящую из типов, уровней и расширенных окружения всех трех ее начальников.

У краевых вершин информации не будет.

Из угловых вершин информацию определим только для $\mathbb{CDR}$-узла, и она будет равна упорядоченной тройке состоящей из типа, уровня и расширенного окружения ее единственного начальника, $\mathbb{CDL}$-узла в той же макроплитке.

\end{definition}

\medskip

{\bf Обозначения.} Информацию о первом начальнике узла $X$ будем обозначать как $\mathbf{FBoss}(X)$. Иногда это обозначение будет употребляться для типа или окружения начальника, тогда это будет понятно по контексту. Аналогичны обозначения для второго и третьего начальников:
$\mathbf{SBoss}(X)$, $\mathbf{TBoss}(X)$.

\medskip

Таким образом, информация -- это данные о начальниках, а для боковых вершин на ребрах $2$, $5$, $6$ -- это еще знание, какой тип у узла в правом нижнем углу.

\medskip

\subsection{Ядро подклейки и флаг макроплитки}

Рассмотрим некоторую макроплитку $T$, подклееную к ребрам $e_1$ и $e_2$ некоторой вершины $X$.
Вершина $X$ является {\it ядром подклейки}, в соответствии с определением из первой части. Она лежит в левом верхнем углу $T$. Упорядоченную пятерку параметров (тип, базовое окружение $X$, информация $X$, тип ребра $e_1$, тип ребра $e_2$) будем называть {\it флагом подклейки} для макроплитки $T$.

Для каждой вершины $Y$, лежащей внутри подклееной макроплитки $T$ (не на ее краю) определим параметр {\it флаг подклейки}, со значением, равным значению этого параметра для макроплитки. Таким образом, для всех таких вершин $Y$, лежащих внутри подклееной макроплитки значение этого параметра будет одинаковым.

\subsection{Кодировка путей}

\smallskip

Рассмотрим путь $X_1 e_1 X_2 e_2 X_3 e_3 \dots e_{n-1} X_n$, где $X_i$ -- вершины, $e_i$ -- ребра между соседними вершинами в пути.
Обозначим $\mathbf{In}(e_i)$ -- тип $e_i$ как выходящего ребра из вершины $X_i$. Также пусть  $\mathbf{Out}(e_i)$ -- тип $e_i$ как входящего ребра в вершину $X_{i+1}$.

\medskip

Рассмотрим некоторую вершину $X_i$. Если и входящее в $X_i$ ребро $e_{i-1}$, и выходящее ребро $e_i$ являются плоскими, то кодом $X_i$ будем считать упорядоченную четверку, состоящую из типа $X_i$, ее уровня, ее окружения и ее информации.

\smallskip

Если одно из ребер, например входящее ребро $e_{i-1}$ приходит из подклееной области $P_{i-1}$, а выходящее ребро $e_i$ -- плоское, то кодом $X_i$ будем считать упорядоченный набор из четырех элементов:
\begin{enumerate}

\item {\it первый элемент} -- упорядоченная пара из типа $X_i$ в области $P_{i-1}$ и типа $X_i$ в базовой плоскости;

\item {\it второй элемент} -- упорядоченная пара из уровня $X_i$ в области $P_{i-1}$ и уровня $X_i$ в базовой плоскости;

\item {\it третий элемент} -- упорядоченная пара из окружения $X_i$ в области $P_{i-1}$ и окружения $X_i$ в базовой плоскости;

\item {\it четвертый элемент} -- информация у $X_i$.

\end{enumerate}

\smallskip

Таким образом кодировка учитывает и положение вершины в макроплитке, откуда пришло ребро и в макроплитке базовой плоскости. Для случая, когда входящее ребро -- плоское, а выходящее ведет в подклейку, определение аналогично, только в упорядоченной паре (типов, уровней, окружений) сначала будет идти базовая плоскость, а потом подклееная.

\medskip

Если оба ребра, и входящее $e_{i-1}$, и выходящее $e_i$, ведут в подклееные области $P_{i-1}$ и $P_i$, то кодом $X_i$ будем считать упорядоченный набор из четырех элементов:

\begin{enumerate}

\item {\it первый элемент} -- упорядоченная тройка -- тип $X_i$ в области $P_{i-1}$, тип $X_i$ в базовой плоскости, тип $X_i$ в области $P_i$;

\item {\it второй элемент} -- упорядоченная тройка -- уровень $X_i$ в области $P_{i-1}$, уровень $X_i$ в базовой плоскости, уровень $X_i$ в области $P_i$;

\item {\it третий элемент } -- упорядоченная тройка -- окружение $X_i$ в области $P_{i-1}$, окружение $X_i$ в базовой плоскости, окружение $X_i$ в области $P_i$;

\item {\it четвертый элемент} -- информация у $X_i$.

\end{enumerate}

\medskip

{\bf Замечание}.  Как видно, информация не подчиняется общей логике. Это происходит по причине того, что любая вершина в своих подклееных макроплитках всегда занимает место на верхней или левой сторонах, либо в лежит в любом углу, кроме правого нижнего. Во всех этих случаях информация вершины в рамках такой подклееной макроплитки, по определению, пустая.

\smallskip

Код вершины $X$ будем обозначать как $\mathbf{Code}(X)$. Теперь определим код всего пути.

\begin{definition}
{\it Кодом пути} $X_1 e_1 X_2 e_2 X_3 e_3 \dots e_{n-1} X_n$ будем называть упорядоченный набор

\smallskip

$\mathbf{Code}(X_1)$, $\mathbf{Out}(e_1)$, $\mathbf{In}(e_1)$,  $\mathbf{Code}(X_2)$, $\mathbf{Out}(e_2)$, \dots , $\mathbf{In}(e_{n-1})$, $\mathbf{Code}(X_n)$.

\smallskip

Здесь $\mathbf{Code}(X_i)$ -- код вершины $X_i$, $\mathbf{Out}(e_i)$ -- тип $e_i$ как выходящего ребра из вершины $X_i$, $\mathbf{In}(e_i)$ -- тип $e_i$ как входящего ребра в вершину $X_{i+1}$.

\end{definition}

\medskip

Нашей основной целью будет определить локальные преобразования кодов путей так, чтобы переход от одного кода к другому соответствовал переходу к эквивалентному пути на комплексе.

При этом элементарным преобразованиям путей будут соответствовать замены одних слов-кодов путей на другие. То есть элементарные преобразования путей являются аналогами определяющих соотношений в полугруппе слов-кодов.

Для того, чтобы задать такие определяющие соотношения, нужно рассмотреть все локальные преобразования путей и для каждого из них указать, какой паре равных кодов соответствует данное локальное преобразование.

Часто знание типа, окружения или информации некоторой вершины помогает нам понять устройство ближайшей окрестности. Для удобства, мы зададим набор функций на узлах, аргументом которых являются типы, окружения или информации некоторых специально расположенных узлы, а значениями -- окружения или типы узлов в окрестности узла-аргумента. В следующем параграфе мы определим набор таких функций.

\medskip

\section{Функции на узлах и структура цепей} \label{functions}

Ниже приведены функции, которые облегчают нам работу по восстановлению параметров вершин.

\medskip

\begin{enumerate}

\item $\mathbf{TopFromCorner}$ Значение: окружение узла в середине верхней стороны $T$.

\item $\mathbf{RightCorner}$  Значение: окружение узла в правом нижнем углу $T$.

\item $\mathbf{TopRightType}$  Значение: тип узла в правом верхнем углу $T$, а также тип ребра для боковых узлов.

\item $\mathbf{BottomLeftType}$  Значение: тип узла в левом нижнем углу $T$,
  а также тип ребра для боковых узлов.

\item $\mathbf{LevelPlus}$  Значение: окружение цепи с уровнем на один больше, чем у $X$, с тем же центром и указателем.

\item $\mathbf{BottomRightTypeFromRight}$ Значение: тип узла в правом нижнем углу $T$, а также тип ребра для боковых узлов.

\item $\mathbf{TopFromRight}$ Значение: окружение узла в середине верхней стороны $T$.

\end{enumerate}

Функции 1 и 2 содержат два аргумента: первый -- это окружение узла $X$, являющимся левым нижним углом в макроплитке $T$, второй -- тип ребра выхода из $X$, соответствующего нижней стороне $T$.

Функции 3, 4, 5 содержат один аргумент: окружение узла $X$, являющегося серединой верхней стороны в макроплитке $T$

Функции 6 и 7 содержат два аргумента: окружение {\bf и информацию} узла $X$, являющимся серединой правой стороны в макроплитке~$T$.

\smallskip

Еще раз отметим, что под окружением узла в цепи мы понимаем окружение всей этой цепи вместе с указателем -- типом ребра входа-выхода, на котором лежит узел из цепи.

\medskip

\subsection{Свойства цепей}

В этом параграфе мы опишем свойства цепей для различных типов узлов.
В частности, мы установим, что зная тип центра, можно узнать окружение цепи, и наоборот. Помимо этого, мы используем свойства цепей для доказательства того, что по аргументу функции можно узнать ее значение.

\medskip

Мы разберем цепи с различными центрами и в каждом случае будем проверять, что если мы знаем аргумент функции, то можем узнать ее значение.

\medskip

\begin{figure}[hbtp]
\centering
\includegraphics[width=1\textwidth]{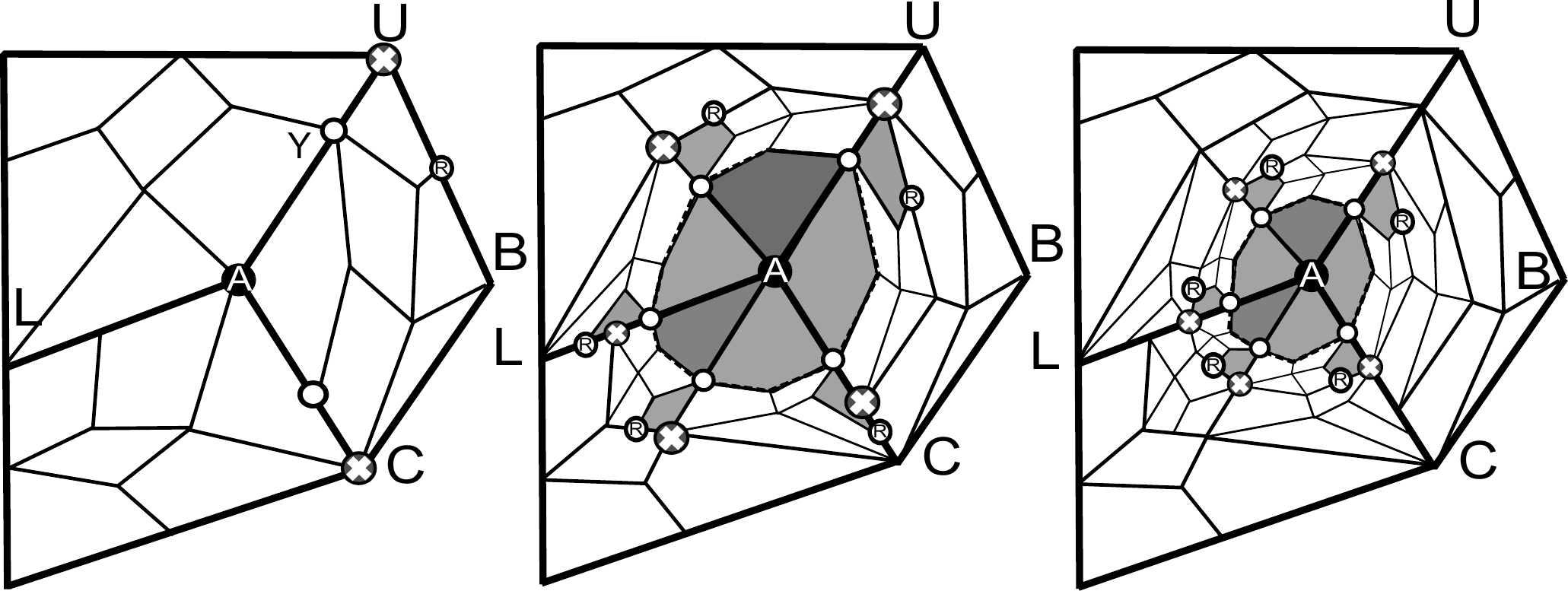}
\caption{Цепи узла $\mathbb{A}$.}
\label{fig:chainA}
\end{figure}

{\bf Цепи с центром в узле типа $\mathbb{A}$.}
На рисунке~\ref{fig:chainA} белыми точками отмечены цепи узла $\mathbb{A}$ нулевого, первого и второго уровней. Заметим, что на втором уровне все пять макроплиток, левый верхний угол которых попадает в узел $\mathbb{A}$, занимают левое верхнее положение в своих родительских макроплитках. Значит, если взять следующий уровень подразбиения, типы ребер-границ макроплиток не изменятся. То есть окружения вершин в цепях будут те же, то есть, окружение цепи третьего и последующих уровней совпадает с окружением цепи второго уровня.

Таким образом, существует только три возможных окружения цепи с центром в узле типа $\mathbb{A}$. Заметим, что по окружению цепи, мы можем установить, является ли центр цепи узлом типа $\mathbb{A}$, и какого уровня цепь. Действительно, одна из макроплиток в цепи с центром в $\mathbb{A}$ обязательно имеет окружение $(\mathbf{4B},\mathbf{1B},\mathbf{2B},\mathbf{5B})$ или $(\mathbf{4B},\mathbf{1B},\mathbf{1A},\mathbf{3A})$, причем, макроплитки с такими окружениями не могут входить в цепи с другими центрами. Цепь нулевого уровня содержит макроплитку с окружением $(\mathbf{4B},\mathbf{1B},\mathbf{2B},\mathbf{5B})$ (остальные -- нет). Для цепи второго уровня все четверки типов ребер имеют вид $(x,y,\mathbf{1A},\mathbf{3A})$, а для первого уровня это не так.

Забегая вперед, можно сказать, что возможность установить тип центра по окружению цепи будет и для центров других типов (не только $\mathbb{A}$).

Заметим также, что мы можем выписать полностью все окружение $1$-цепи и $2$-цепи с центром в узле типа $\mathbb{A}$, а зная окружение узла $\mathbb{A}$, можем выписать и окружение $0$-цепи.

\medskip

\medskip

Из сказанного выше следует, что для $\mathbb{A}$-узлов значение функции $\mathbf{LevelPlus}$ может быть установлено по ее аргументу.

На рисунке~\ref{fig:chainA} знаком ``$\otimes$'' отмечены узлы, являющиеся верхними правыми или левыми нижними углами в макроплитках, где середина верхней стороны попадает в узел цепи $X$ . Можно убедиться в том, что для цепей первого и второго уровней все типы узлов с крестами, а также типы ребер, на которых они лежат, мы можем выписать, для каждого заданного $X$. Для $0$-цепи мы можем выписать типы узлов со знаком ``$\otimes$'' (а также типы ребер), если знаем окружение центра (узла $\mathbb{A}$).
Таким образом, мы можем вычислить функции $\mathbf{TopRightType}$ и $\mathbf{BottomLeftType}$.

%Черными точками на рисунке~\ref{fig:chainA} помечены узлы являющиеся правыми нижними углами в макроплитках, где середина верхней стороны попадает в узел цепи $X$. Заметим, что зная окружение и информацию $X$, мы можем установить тип этого правого нижнего угла. Действительно, для всех узлов кроме одного (на $3$ выходящем из $\mathbb{A}$ ребре в $1$-цепи), мы сразу можем установить их тип. В оставшемся случае, этот угол может быть записан как $\mathbf{Next.FBoss}(X)$. Таким образом, мы можем вычислить функцию $\mathbf{BottomRightTypeFromTop}$.

\medskip

Серым цветом на рисунке выделены плитки, где левый нижний угол (при дальнейшем разбиении) будет попадать в вершину цепи. Можно заметить, что зная эту вершину (ее место в цепи) мы можем установить и окружение правого нижнего угла в соответствующей макроплитке. Кроме того, мы также можем установить тип левого верхнего угла, а значит и всю цепь, содержащую середину верхней стороны в соответствующей макроплитке. Все это значит, что мы можем вычислить функции
$\mathbf{TopFromCorner}$, $\mathbf{RightCorner}$.

 %Рассмотрим все макроплитки, где середина верхней стороны $X$ попадает в цепь. При этом середина
 %правой стороны попадет в один из узлов, помеченных $<<R>>$. Заметим, что зная окружение и информацию
 %$X$, мы можем вычислить окружение этого $R$-узла. (Во всех случаях, кроме одного, хватает просто
 %окружения $X$, для $X$ лежащей на $3$ ребре выхода из $\mathbb{A}$, окружение $R$-вершины будет
 %$\mathbf{Next.FBoss}(X)$). Таким образом, функция $\mathbf{RightFromTop}$ также вычисляется.

Функции $\mathbf{TopFromRight}$ и $\mathbf{BottomRightTypeFromRight}$ могут быть применимы только в одном случае: для $0$-цепи, если аргументом является узел в середине ребра типа $1$ (обозначим этот узел как $Y$, а макроплитку, где $Y$ является серединой правой стороны, как $T$). $\mathbf{BottomRightTypeFromRight}$ в этом случае это просто наш узел $\mathbb{A}$.

Для вычисления $\mathbf{TopFromRight}$ нам надо найти окружение узла в середине верхней стороны макроплитки $T$. Заметим, что этот узел входит в цепь с центром в левом верхнем углу $T$, причем уровень этой цепи на один выше, чем уровень цепи, куда входит правый верхний угол $T$. Окружение этой цепи мы можем установить по информации $Y$, так как правый верхний угол $T$ является серединой верхней стороны в макроплитке $T'$, родительской к $T$ и является первым начальником узла $Y$. По окружению этой цепи можно установить тип левого верхнего угла $T$ и окружение цепи на один уровень выше, что и требуется.

\medskip

Таким образом, мы вычислили значения всех функций для аргументов, входящих в $\mathbb{A}$-цепи.

\medskip

\begin{figure}[hbtp]
\centering
%\leftskip=-1cm
\includegraphics[width=1\textwidth]{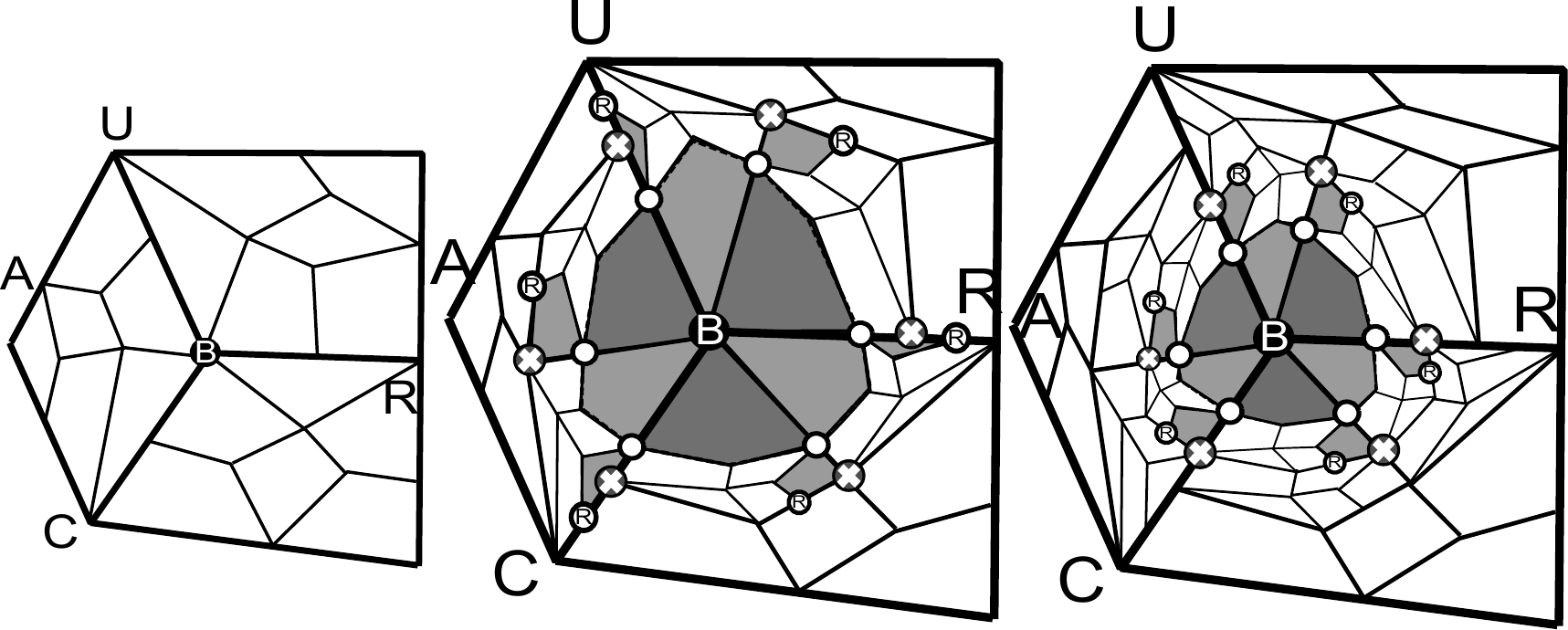}
\caption{Цепи узла $\mathbb{B}$.}
\label{fig:chainB}
\end{figure}

{\bf Цепи с центром в узле типа $\mathbb{B}$.}
На рисунке~\ref{fig:chainB} белыми точками отмечены цепи узла $\mathbb{B}$. $0$-цепи с центром в $\mathbb{B}$ не существует (все боковые узлы в соответствующей окрестности имеют тип $\mathbb{DR}$). Заметим, что в $2$-цепи все шесть макроплиток, левый верхний угол которых попадает в узел $\mathbb{B}$, занимают левое верхнее положение в своих родительских макроплитках. Значит, если взять следующий уровень подразбиения, типы ребер-границ макроплиток не изменятся. То есть окружения вершин в цепях будут те же, то есть, окружение цепи третьего и последующих уровней совпадает с окружением цепи второго уровня.

Таким образом, существует только два возможных окружения цепи с центром в узле типа $\mathbb{B}$. Заметим, что по окружению цепи, мы можем установить, является ли центр цепи узлом типа $\mathbb{B}$, и какого уровня цепь. Действительно, $1$-цепь с центром в $\mathbb{B}$ содержит узел с окружением $(\mathbf{8B},\mathbf{5B},\mathbf{5B},\mathbf{7B})$, который не может содержаться ни в какой другой цепи. $2$-цепь с центром в $\mathbb{B}$ содержит узел с окружением $(\mathbf{8B},\mathbf{5B},\mathbf{1A},\mathbf{3A})$, который также ни в какой другой цепи не встречается.

Заметим также, что мы можем выписать полностью все окружение $1$-цепи и $2$-цепи с центром в узле типа $\mathbb{B}$.

\medskip

\medskip

Из сказанного выше следует, что для $\mathbb{B}$-узлов значение функции $\mathbf{LevelPlus}$ может быть установлено по ее аргументу.
На рисунке~\ref{fig:chainB} знаком ``$\otimes$'' отмечены узлы, являющиеся верхними правыми или левыми нижними углами в макроплитках, где середина верхней стороны попадает в узел цепи $X$. Можно убедиться в том, что для цепей первого и второго уровней все типы узлов с крестами, а также типы ребер, на которых они лежат,  мы можем выписать, для каждого заданного $X$.
Таким образом, мы можем вычислить функции $\mathbf{TopRightType}$ и $\mathbf{BottomLeftType}$.

\medskip

Серым цветом на рисунке выделены плитки, где левый нижний угол (при дальнейшем разбиении) будет попадать в вершину цепи. Можно заметить, что зная эту вершину (ее место в цепи мы знаем, благодаря второму аргументу -- нижнему ребру) мы можем установить и окружение правого нижнего угла в соответствующей макроплитке. Кроме того, мы также можем установить тип левого верхнего угла, а значит и всю цепь, содержащую середину верхней стороны в соответствующей макроплитке. Все это значит, что мы можем вычислить функции
$\mathbf{TopFromCorner}$, $\mathbf{RightCorner}$.

Функции $\mathbf{TopFromRight}$ и $\mathbf{BottomRightTypeFromRight}$ не могут быть применимы, то есть аргументом этих функций никогда не может быть узел из B-цепи так как ни один узел из цепи с центром в $\mathbb{B}$ не является серединой правой стороны ни в какой макроплитке.

%Рассмотрим все макроплитки, где середина верхней стороны $X$ попадает в цепь. При этом середина правой стороны попадет в один из узлов, помеченных $<<R>>$. Заметим, что для $2$-цепи, зная окружение и информацию $X$, мы можем вычислить окружение этого R-узла.

%Для $1$-цепи возможно шесть случаев расположения $X$, в соответствии с выходящими из $\mathbb{B}$ ребрами. Для $\mathbf{ru}$, $\mathbf{rd}$, $\mathbf{mid}$ -- ребер $R$-вершина лежит в середине $5$ ребра, окружение легко вычисляется. Для $1$-ребра окружение $R$-вершины это $1$-цепь вокруг $\mathbf{FBoss}(X)$, с указателем $\mathbf{u}_2$. Для $3$-ребра окружение $R$-вершины это $1$-цепь вокруг $\mathbb{C}$-узла, с указателем $2$. Для $2$-ребра $R$-вершины входит в $1$-цепь вокруг узла в середине правой стороны макроплитки $\mathbb{B}$-узла с указателем $\mathbf{r}$. Тип этого узла может быть установлен по окружению макроплитки, то есть $U$-части окружения узла в середине верхней стороны: $\mathbf{U.FBoss}(X)$.  Для левого верхнего положения тип будет $\mathbb{UR}$, для левого нижнего $\mathbb{DR}$, для среднего $\mathbb{DR}$, правого верхнего $\mathbb{RD}$, для правого нижнего $\mathbb{RD}$. Для нижнего ребра тип будет $\mathbb{RU}$ или $\mathbb{UR}$ в соответствии с ориентацией нижнего ребра ($А$ или $B$). Таким образом, функция $\mathbf{RightFromTop}$ также вычисляется.

\medskip

Таким образом, мы вычислили значения всех функций для аргументов, входящих в $\mathbb{B}$-цепи.

\medskip

\begin{figure}[hbtp]
\centering
%\leftskip=-1cm
\includegraphics[width=1\textwidth]{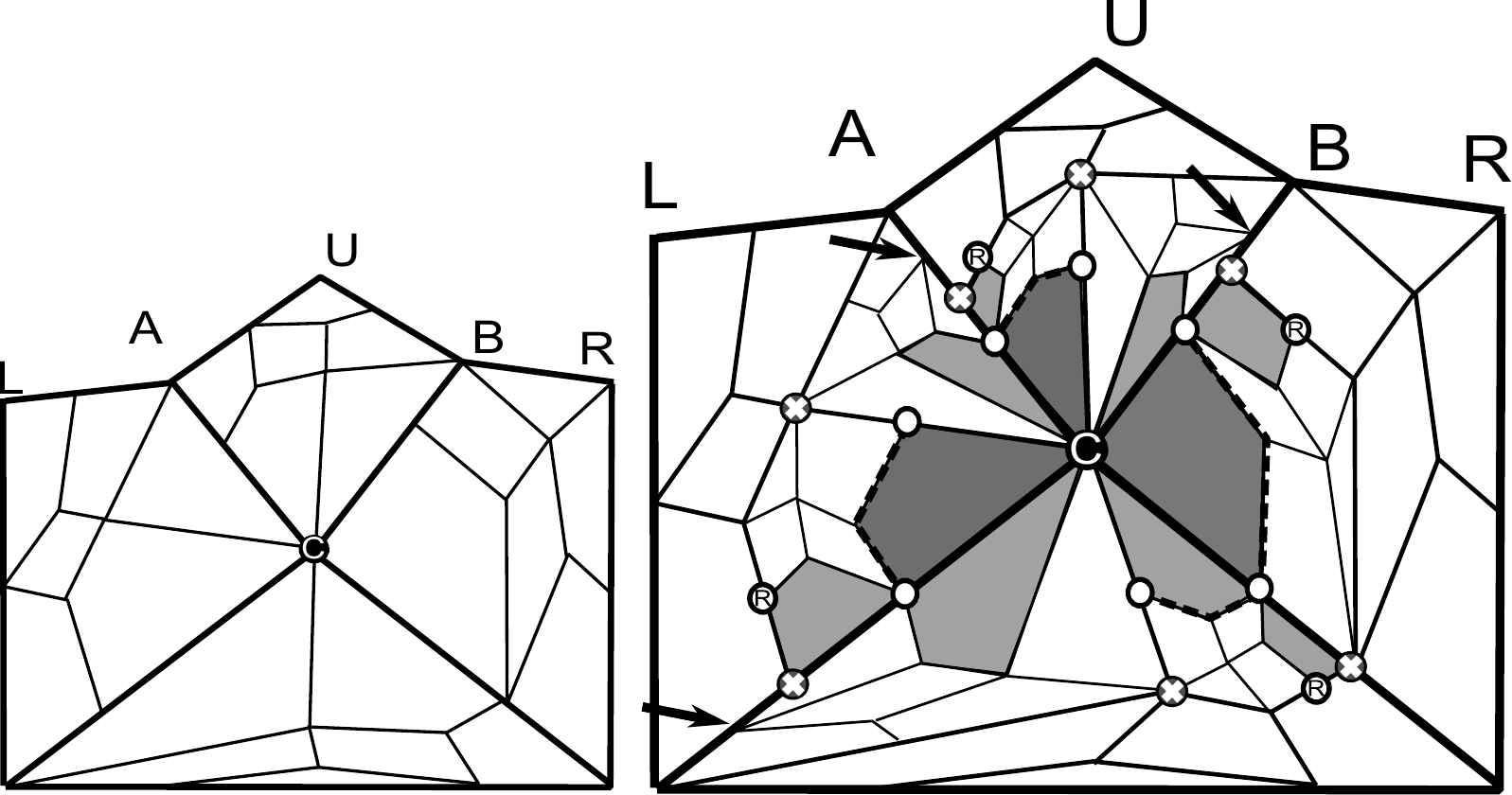}
\caption{$1$-цепь узла $\mathbb{C}$.}
\label{fig:chainC1}
\end{figure}

\begin{figure}[hbtp]
\centering
\includegraphics[width=0.8\textwidth]{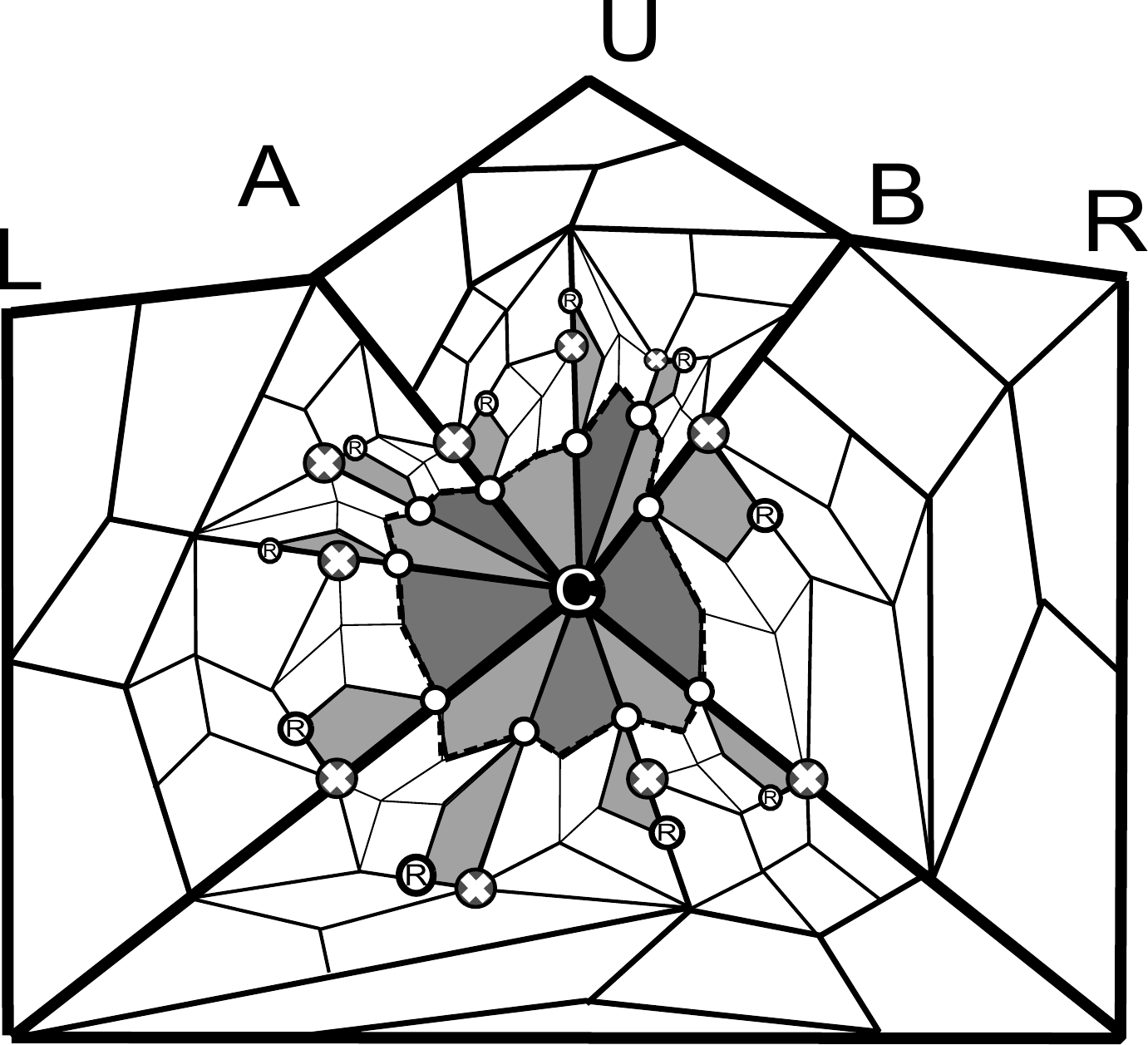}
\caption{$2$-цепь узла $\mathbb{C}$.}
\label{fig:chainC2}
\end{figure}

{\bf Цепи с центром в узле типа $\mathbb{C}$.}
На рисунках~\ref{fig:chainC1} и~\ref{fig:chainC2} белыми точками отмечены $1$-цепи и $2$-цепи узла $\mathbb{C}$. $0$-цепи с центром в $\mathbb{C}$ не существует.

Кроме указанных на рисунках цепей первого и второго уровней существует также $3$-цепь, получаемая применением операции разбиения к макроплиткам на рисунке~\ref{fig:chainC2}.

На получающемся третьем уровне все  макроплитки, левый верхний угол которых попадает в узел $\mathbb{C}$, будут занимать левое верхнее положение в своих родительских макроплитках. Значит, на следующих уровнях подразбиения, типы ребер-границ макроплиток не изменятся. То есть окружения вершин в цепях будут те же, то есть, окружение цепи четвертого и последующих уровней совпадает с окружением цепи третьего уровня.

Таким образом, существует только три возможных окружения цепи с центром в узле типа $\mathbb{C}$ (для цепей первого, второго и третьего уровней). Заметим, что по окружению цепи, мы можем установить, является ли центр цепи узлом типа $\mathbb{C}$, и какого уровня цепь. Действительно, $1$-цепь с центром в $\mathbb{C}$ содержит узел с окружением $(\mathbf{8A},\mathbf{5A},\mathbf{1A},\mathbf{3A})-(\mathbf{8B},\mathbf{5B},\mathbf{5B},\mathbf{7B})$ (на ребре $5$), который не может содержаться ни в какой другой цепи. $2$-цепь с центром в $\mathbb{C}$ содержит узел с окружением $(\mathbf{8A},\mathbf{5A},\mathbf{1A},\mathbf{3A})-(\mathbf{5B},\mathbf{8A},\mathbf{5A},\mathbf{6B})$, который также ни в какой другой цепи не встречается.
Для $3$-цепи таким узлом будет $(\mathbf{8A},\mathbf{5A},\mathbf{1A},\mathbf{3A})-(\mathbf{5B},\mathbf{8A},\mathbf{1A},\mathbf{3A})$.

Заметим также, что мы можем выписать полностью все окружение $1$-цепей, $2$-цепей, $3$-цепей с центром в узле типа $\mathbb{C}$.

\medskip

\leftskip=0cm
\medskip

Из сказанного выше следует, что для $\mathbb{C}$-узлов значение функции $\mathbf{LevelPlus}$ может быть установлено по ее аргументу.
На рисунках~\ref{fig:chainC1} и~\ref{fig:chainC2} знаком ``$\otimes$'' отмечены узлы, являющиеся верхними правыми или левыми нижними углами в макроплитках, где середина верхней стороны попадает в узел цепи $X$ . Можно убедиться в том, что для цепей первого, второго и третьего уровней все типы узлов с крестами, а также типы ребер, на которых они лежат,  мы можем выписать, для каждого заданного $X$.
Таким образом, мы можем вычислить функции $\mathbf{TopRightType}$ и $\mathbf{BottomLeftType}$.

Серым цветом на рисунке выделены плитки, где левый нижний угол (при дальнейшем разбиении) будет попадать в вершину цепи. Можно заметить, что зная эту вершину (ее место в цепи) мы можем установить и окружение правого нижнего угла в соответствующей макроплитке. Кроме того, мы также можем установить тип левого верхнего угла, а значит и всю цепь, содержащую середину верхней стороны в соответствующей макроплитке. (Все это также очевидно проверяется для $3$-цепи.) Все это значит, что мы можем вычислить функции $\mathbf{TopFromCorner}$, $\mathbf{RightCorner}$.

\medskip

\leftskip=0cm

%Рассмотрим все макроплитки, где середина верхней стороны $X$ попадает в цепь. При этом середина правой стороны попадет в один из узлов, помеченных <<$R$>>. Заметим, что зная окружение и информацию $X$, мы можем вычислить окружение этого $R$-узла.
%Таким образом, функция $\mathbf{RightFromTop}$ также вычисляется.

Функции $\mathbf{TopFromRight}$ и $\mathbf{BottomRightTypeFromRight}$ могут быть применимы только к $1$-цепи (в остальные входят только узлы типов $\mathbb{UL}$ и $\mathbb{LU}$). В случае $1$-цепи аргументом могут быть узлы, лежащие на ребрах $1$, $2$ и $4$ (относительно $\mathbb{C}$). Пусть $X$ такой узел.  Во всех случаях
$\mathbf{BottomRightTypeFromRight}$ будет $\mathbb{C}$. Узлы, на которые указывает $\mathbf{TopFromRight}$ отмечены стрелками на рисунке~\ref{fig:chainC1}. Таким образом, $\mathbf{TopFromRight}$ принимает значения: $1$-цепи вокруг $\mathbb{A}$ для $1$ ребра (с указателем $2$), $1$-цепи вокруг $\mathbb{B}$ для ребра $2$ (с указателем $3$),
$2$-цепь вокруг $\mathbf{SBoss}(X)$ для ребра $4$ (с указателем, соответствующим входу ребра $7$ в левый нижний угол (см параграф~\ref{pointers} ``Указатели'').

\medskip

Таким образом, мы вычислили значения всех функций для аргументов, входящих в $\mathbb{C}$-цепи.

\medskip

\begin{figure}[hbtp]
\centering
\includegraphics[width=1\textwidth]{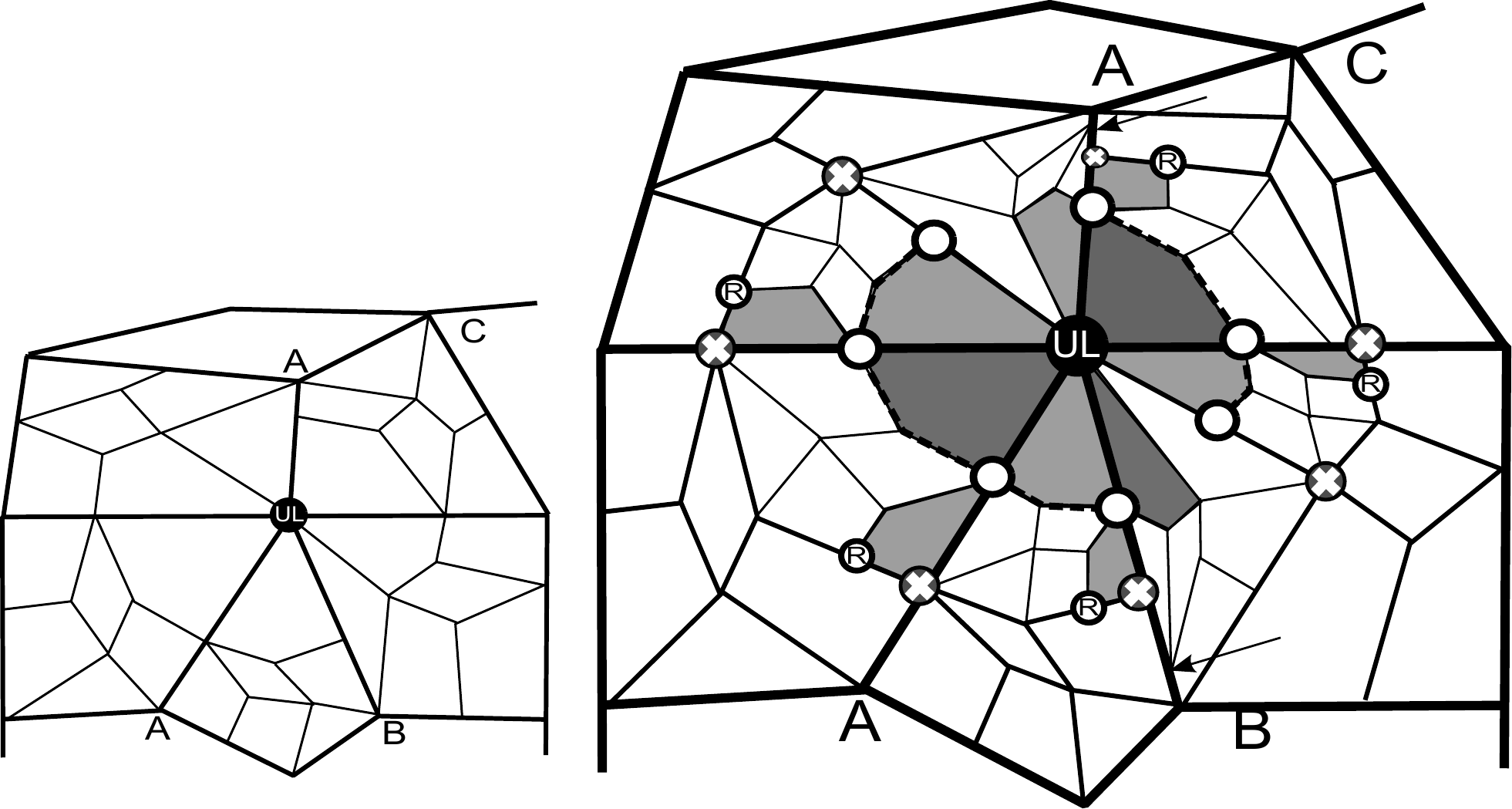}
\caption{$1$-цепь узла $\mathbb{UL}$.}
\label{fig:chainUL1}
\end{figure}

\begin{figure}[hbtp]
\centering
\includegraphics[width=0.8\textwidth]{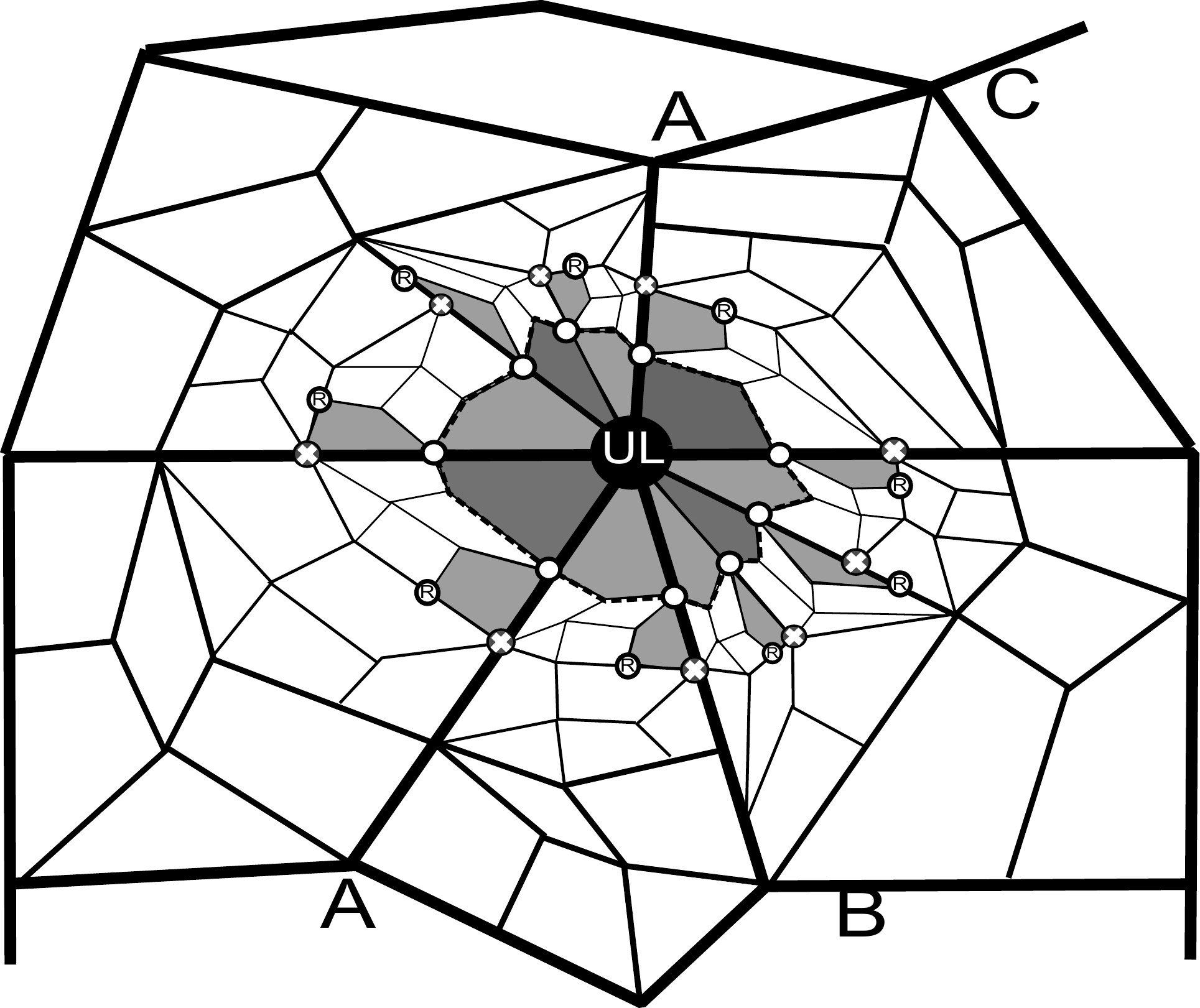}
\caption{$2$-цепь узла $\mathbb{UL}$.}
\label{fig:chainUL2}
\end{figure}

{\bf Цепи с центром в узлах типа $\mathbb{UL}$ и $\mathbb{LU}$.}
Случаи $\mathbb{UL}$ и $\mathbb{LU}$ узлов симметричны, для определенности будем далее разбирать случай $\mathbb{UL}$ узла. Будем считать, что главное ребро имеет тип $t$, то есть с $U$-стороны это $tA$, а с $L$-стороны $tB$.
На рисунках~\ref{fig:chainUL1} и~\ref{fig:chainUL2} белыми точками отмечены $1$-цепи и $2$-цепи узла $\mathbb{UL}$. $0$-цепи с центром в $\mathbb{UL}$ не существует.
Кроме указанных на рисунках цепей первого и второго уровней существует также $3$-цепь, получаемая применением операции разбиения к макроплиткам на рисунке~\ref{fig:chainUL2}.

На получающемся третьем уровне все  макроплитки, левый верхний угол которых попадает в узел $\mathbb{UL}$, будут занимать левое верхнее положение в своих родительских макроплитках. Значит, на следующих уровнях подразбиения, типы ребер-границ макроплиток не изменятся. То есть окружения вершин в цепях будут те же, то есть, окружение цепи четвертого и последующих уровней совпадает с окружением цепи третьего уровня.

Таким образом, существует только три возможных конфигурации окружения цепи с центром в узле типа $\mathbb{UL}$ (для цепей первого, второго и третьего уровней, при этом еще может быть выбрано разное главное ребро $t$). Заметим, что по окружению цепи, мы можем установить, является ли центр цепи узлом типа $\mathbb{UL}$, какого уровня цепь, а также сам параметр $t$. Действительно, $1$-цепь с центром в $\mathbb{UL}$ содержит узлы с окружением $(tA,\mathbf{1A},\mathbf{6A},\mathbf{2A})-(\mathbf{1B},\mathbf{2B},\mathbf{6A},\mathbf{2A})$ и $(\mathbf{8B},\mathbf{3A},\mathbf{3A},\mathbf{7B})-(tB,\mathbf{3B},\mathbf{6A},\mathbf{2A})$, которые вместе ни в какой другой цепи не встречаются. $2$-цепь с центром в $\mathbb{UL}$ содержит одновременно узлы с окружением $(tA,\mathbf{1A},\mathbf{1A},\mathbf{3A})-(\mathbf{1B},\mathbf{2B},\mathbf{1A},\mathbf{3A})$ и $(\mathbf{3A},\mathbf{8A},\mathbf{5A},\mathbf{6B})-(tB,\mathbf{3B},\mathbf{1A},\mathbf{3A})$, которые вместе также ни в какой другой цепи не встречаются.
Для $3$-цепи такой парой будет $(tA,\mathbf{1A},\mathbf{1A},\mathbf{3A})-(\mathbf{1B},\mathbf{2B},\mathbf{1A},\mathbf{3A})$ и $(\mathbf{3A},\mathbf{8A},\mathbf{1A},\mathbf{3A})-(tB,\mathbf{3B},\mathbf{1A},\mathbf{3A})$

Заметим также, что зная тип ребра $t$, мы можем выписать полностью все окружение $1$-цепей, $2$-цепей, $3$-цепей с центром в узле типа $\mathbb{UL}$.

\medskip

\leftskip=0cm

Из сказанного выше следует, что для $\mathbb{UL}$-узлов значение функции $\mathbf{LevelPlus}$ может быть установлено по ее аргументу.

На рисунках~\ref{fig:chainUL1} и~\ref{fig:chainUL2} знаком ``$\otimes$'' отмечены узлы, являющиеся верхними правыми или левыми нижними углами в макроплитках, где середина верхней стороны попадает в узел цепи $X$ . Можно убедиться в том, что для цепей первого, второго и третьего уровней все типы узлов с крестами, а также типы ребер на которых они лежат, мы можем выписать, для каждого заданного $X$.
Таким образом, мы можем вычислить функции $\mathbf{TopRightType}$ и $\mathbf{BottomLeftType}$.

Серым цветом на рисунке выделены плитки, где левый нижний угол (при дальнейшем разбиении) будет попадать в вершину цепи. Можно заметить, что зная эту вершину (ее место в цепи) мы можем установить и окружение правого нижнего угла в соответствующей макроплитке. Кроме того, мы также можем установить тип левого верхнего угла, а значит и всю цепь, содержащую середину верхней стороны в соответствующей макроплитке. (Все это также очевидно проверяется для $3$-цепи.) Все это значит, что мы можем вычислить функции $\mathbf{TopFromCorner}$, $\mathbf{RightCorner}$.

\medskip

\leftskip=0cm

%Рассмотрим все макроплитки, где середина верхней стороны $X$ попадает в цепь. При этом середина правой стороны попадет в один из узлов, помеченных <<$R$>>. Заметим, что зная окружение и информацию $X$, мы можем вычислить окружение этого $R$-узла.
%Таким образом, функция $\mathbf{RightFromTop}$ также вычисляется.

Функции $\mathbf{TopFromRight}$ и $\mathbf{BottomRightTypeFromRight}$ могут быть применимы только к $1$-цепи (в остальные входят только узлы типов $\mathbb{UL}$ и $\mathbb{LU}$). В случае $1$-цепи аргументом могут быть узлы, лежащие на ребрах $\mathbf{u}_1$ и $\mathbf{l}$ (относительно центрального $\mathbb{UL}$-узла). Пусть $X$ такой узел.  Во всех случаях
$\mathbf{BottomRightTypeFromRight}$ будет наш же узел $\mathbb{UL}$. Узлы, на которые указывает $\mathbf{TopFromRight}$ отмечены стрелками на рисунке~\ref{fig:chainUL1}. Таким образом, $\mathbf{TopFromRight}$ принимает значения: $1$-цепи вокруг $\mathbb{B}$ (указатель $1$) для ребра $\mathbf{u}_1$, $1$-цепи вокруг $\mathbb{A}$ (указатель $3$) для ребра $\mathbf{l}$.

\medskip

Таким образом, мы вычислили значения всех функций для аргументов, входящих в $\mathbb{UL}$-цепи.

\medskip

{\bf Цепи с центром в узлах типа $\mathbb{UR}$ и $\mathbb{RU}$.}
Случаи узлов $\mathbb{UR}$ и $\mathbb{RU}$ симметричны, для определенности будем далее разбирать случай $\mathbb{UR}$ узла. Будем считать, что главное ребро имеет тип $t$, то есть с $U$-стороны это $tA$, а с $R$-стороны $tB$.

\begin{figure}[hbtp]
\centering
\includegraphics[width=0.9\textwidth]{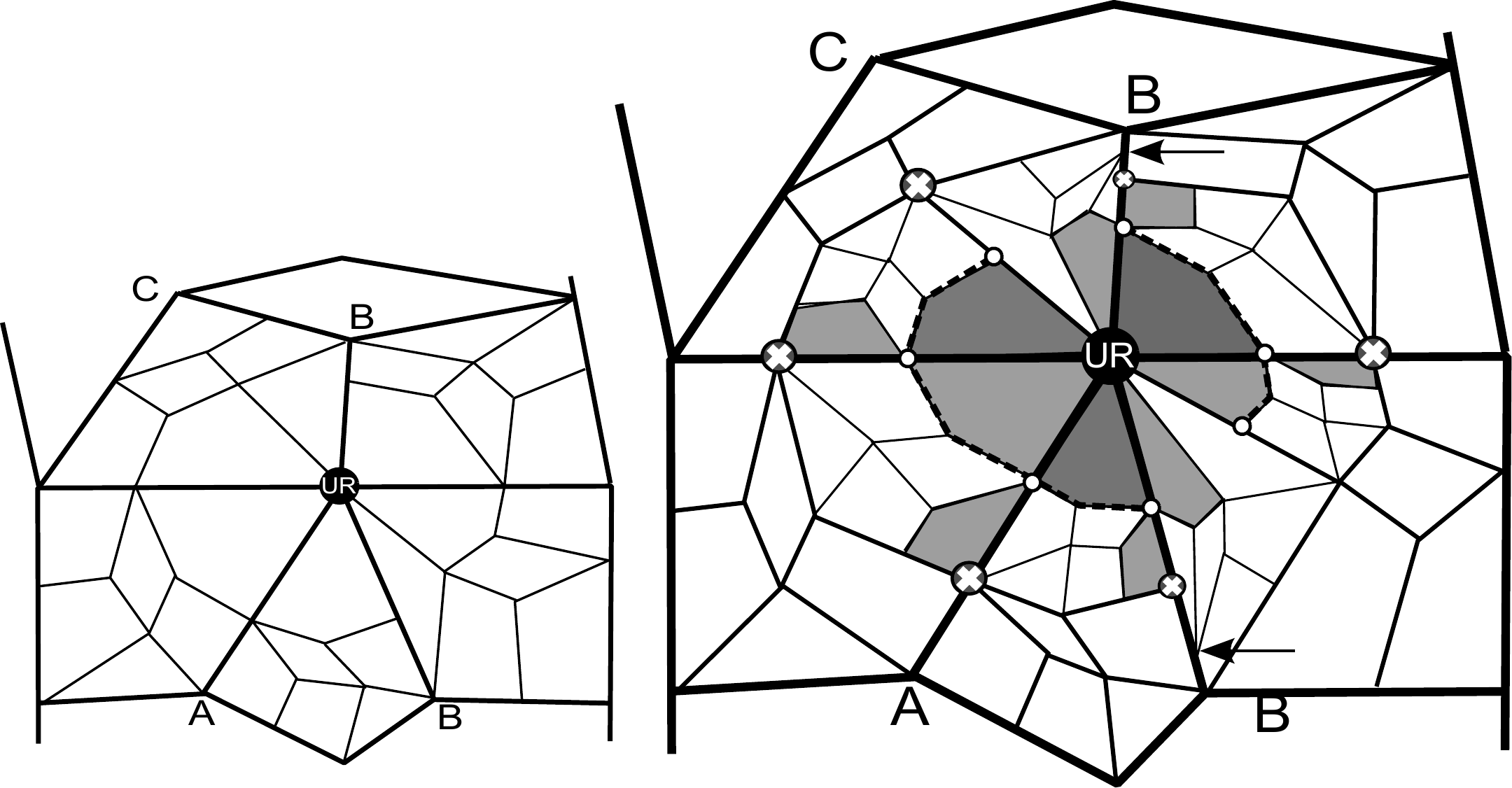}
\caption{$1$-цепь узла $\mathbb{UR}$.}
\label{fig:chainUR1}
\end{figure}

\begin{figure}[hbtp]
\centering
\includegraphics[width=0.8\textwidth]{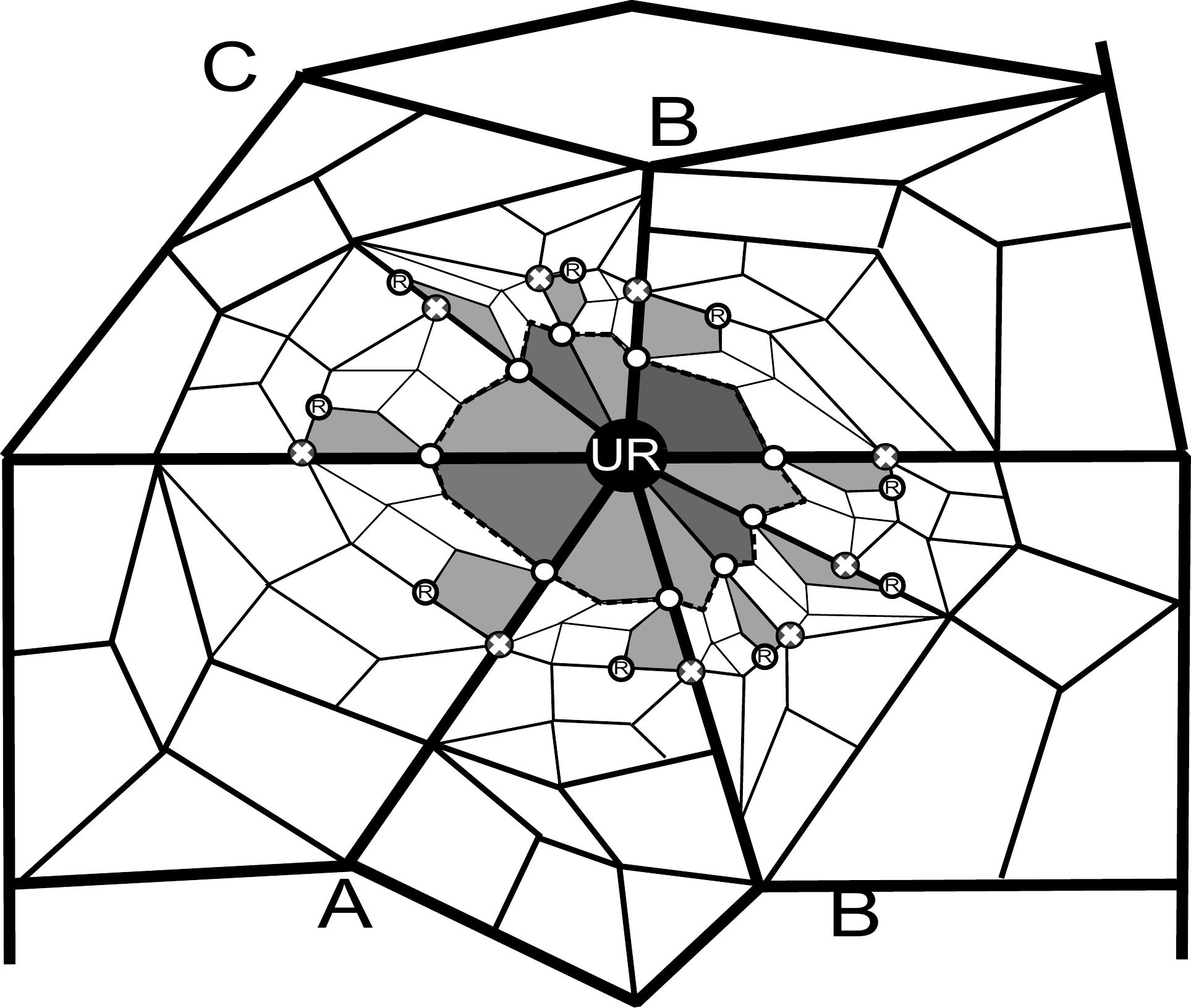}
\caption{$2$-цепь узла $\mathbb{UR}$.}
\label{fig:chainUR2}
\end{figure}

На рисунках~\ref{fig:chainUR1} и~\ref{fig:chainUR2} белыми точками отмечены $1$-цепи и $2$-цепи узла $\mathbb{UR}$. $0$-цепи с центром в $\mathbb{UR}$ не существует.
Кроме указанных на рисунках цепей первого и второго уровней существует также $3$-цепь, получаемая применением операции разбиения к макроплиткам на рисунке~\ref{fig:chainUR2}.

На получающемся третьем уровне все  макроплитки, левый верхний угол которых попадает в узел $\mathbb{UR}$, будут занимать левое верхнее положение в своих родительских макроплитках. Значит, на следующих уровнях подразбиения, типы ребер-границ макроплиток не изменятся. То есть окружения вершин в цепях будут те же, то есть, окружение цепи четвертого и последующих уровней совпадает с окружением цепи третьего уровня.

Таким образом, существует только три возможных конфигурации окружения цепи с центром в узле типа $\mathbb{UR}$ (для цепей первого, второго и третьего уровней, при этом еще может быть выбрано разное главное ребро $t$). Заметим, что по окружению цепи, мы можем установить, является ли центр цепи узлом типа $\mathbb{UR}$, какого уровня цепь, а также сам параметр $t$. Действительно, $1$-цепь с центром в $\mathbb{UR}$ содержит одновременно узлы с окружением $(tA,\mathbf{1A},\mathbf{6A},\mathbf{2A})-(\mathbf{1B},\mathbf{2B},\mathbf{6A},\mathbf{2A})$ и $(tB,\mathbf{6A},\mathbf{6A},\mathbf{2A})-(\mathbf{8B},\mathbf{6B},\mathbf{6B},\mathbf{7B})$, которые вместе ни в какой другой цепи не встречаются. $2$-цепь с центром в $\mathbb{UR}$ содержит одновременно узлы с окружением $(tA,\mathbf{1A},\mathbf{1A},\mathbf{3A})-(\mathbf{1B},\mathbf{2B},\mathbf{1A},\mathbf{3A})$ и $(tB,\mathbf{6A},\mathbf{1A},\mathbf{3A})-(\mathbf{6B},\mathbf{8A},\mathbf{5A},\mathbf{6B})$, которые вместе также ни в какой другой цепи не встречаются.
Для $3$-цепи такой парой будет $(tA,\mathbf{1A},\mathbf{1A},\mathbf{3A})-(\mathbf{1B},\mathbf{2B},\mathbf{1A},\mathbf{3A})$ и $(tB,\mathbf{6A},\mathbf{1A},\mathbf{3A})-(\mathbf{6B},\mathbf{8A},\mathbf{1A},\mathbf{3A})$

Заметим также, что зная тип ребра $t$, мы можем выписать полностью все окружение $1$-цепей, $2$-цепей, $3$-цепей с центром в узле типа $\mathbb{UR}$.

\medskip

\leftskip=0cm

\medskip
Из сказанного выше следует, что для $\mathbb{UR}$-узлов значение функции $\mathbf{LevelPlus}$ может быть установлено по ее аргументу.

\medskip

\leftskip=0cm

На рисунках~\ref{fig:chainUR1} и~\ref{fig:chainUR2} знаком ``$\otimes$'' отмечены узлы, являющиеся верхними правыми или левыми нижними углами в макроплитках, где середина верхней стороны попадает в узел цепи $X$. Можно убедиться в том, что для цепей первого, второго и третьего уровней все типы узлов с крестами, а также типы ребер на которых они лежат, мы можем выписать для каждого заданного $X$.
Таким образом, мы можем вычислить функции $\mathbf{TopRightType}$ и $\mathbf{BottomLeftType}$.

Серым цветом на рисунке выделены плитки, где левый нижний угол (при дальнейшем разбиении) будет попадать в вершину цепи. Можно заметить, что зная эту вершину (ее место в цепи) мы можем установить и окружение правого нижнего угла в соответствующей макроплитке. Кроме того, мы также можем установить тип левого верхнего угла, а значит и всю цепь, содержащую середину верхней стороны в соответствующей макроплитке. (Все это также очевидно проверяется для $3$-цепи.) Все это значит, что мы можем вычислить функции $\mathbf{TopFromCorner}$, $\mathbf{RightCorner}$.

\medskip

Функции $\mathbf{TopFromRight}$ и $\mathbf{BottomRightTypeFromRight}$ могут быть применимы только к $1$-цепи (в остальные входят только узлы типов $\mathbb{UL}$ и $\mathbb{LU}$). В случае $1$-цепи аргументом могут быть узлы, лежащие на ребрах $\mathbf{u}_1$ и $\mathbf{r}$ (относительно $\mathbb{UR}$). Пусть $X$ такой узел.  Во всех случаях
$\mathbf{BottomRightType}$ будет наш же узел $\mathbb{UR}$. Узлы, на которые указывает $\mathbf{TopFromRight}$ отмечены стрелками на рисунке~\ref{fig:chainUR1}. Таким образом, $\mathbf{TopFromRight}$ принимает значения: $1$-цепи вокруг $\mathbb{B}$ (указатель $1$) для $\mathbf{u}_1$ ребра, $1$-цепи вокруг $\mathbb{B}$ (указатель $2$) для $\mathbf{l}$ ребра.

\medskip

Таким образом, мы вычислили значения всех функций для аргументов, входящих в $\mathbb{UR}$-цепи.

\medskip

\medskip

\begin{figure}[hbtp]
\centering
\includegraphics[width=0.9\textwidth]{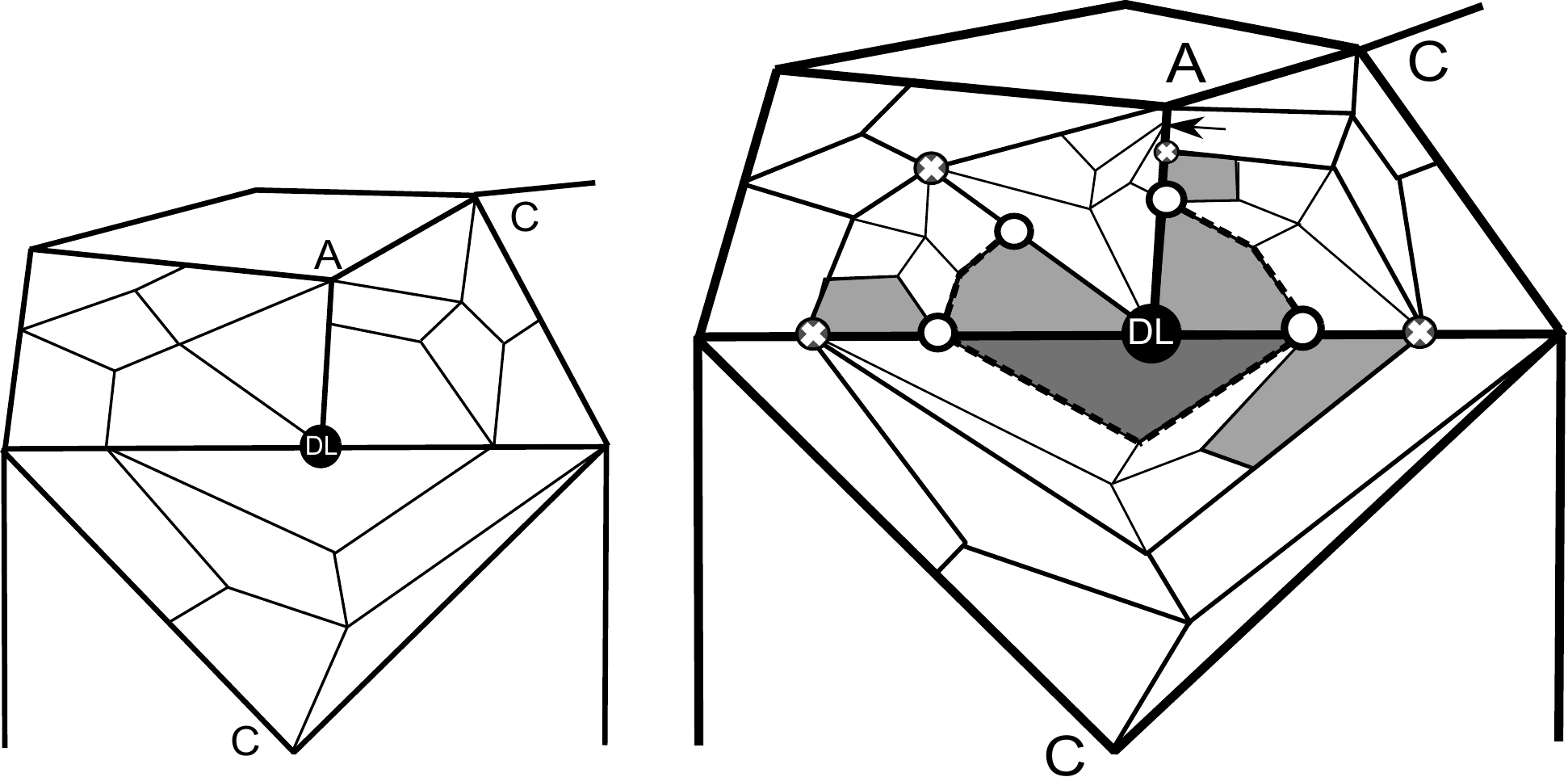}
\caption{$1$-цепь узла $\mathbb{DL}$.}
\label{fig:chainDL1}
\end{figure}

\begin{figure}[hbtp]
\centering
\includegraphics[width=0.8\textwidth]{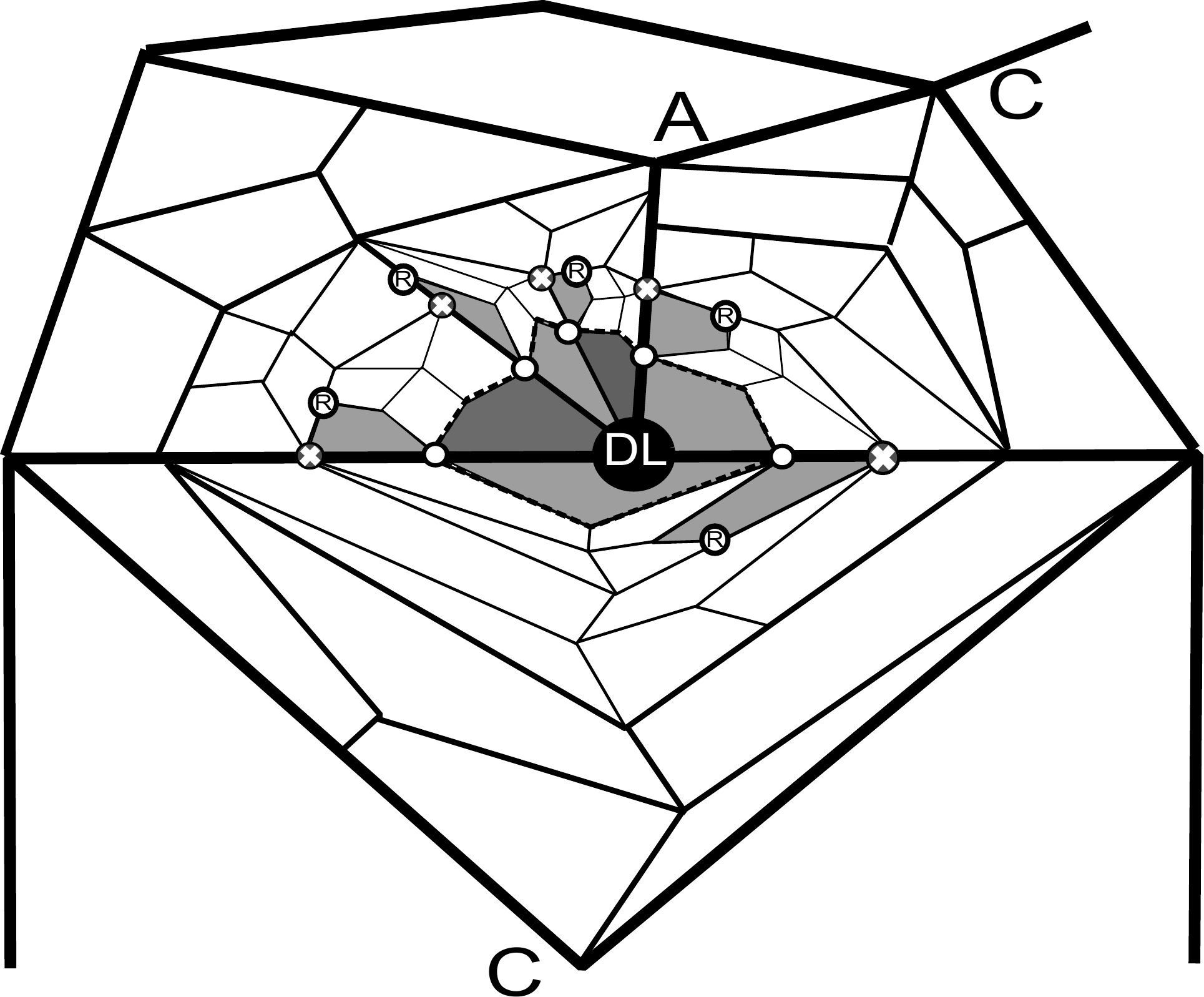}
\caption{$2$-цепь узла $\mathbb{DL}$.}
\label{fig:chainDL2}
\end{figure}

{\bf Цепи с центром в узлах типа $\mathbb{DL}$ и $\mathbb{LD}$.}
Случаи узлов $\mathbb{DL}$ и $\mathbb{LD}$ симметричны, для определенности будем далее разбирать случай $\mathbb{DL}$ узла. Будем считать, что главное ребро имеет тип $t$, то есть с $D$-стороны это $tA$, а с $L$-стороны $tB$.
На рисунках~\ref{fig:chainDL1} и~\ref{fig:chainDL2} белыми точками отмечены $1$-цепи и $2$-цепи узла $\mathbb{DL}$. $0$-цепи с центром в $\mathbb{DL}$ не существует.
Кроме указанных на рисунках цепей первого и второго уровней существует также $3$-цепь, получаемая применением операции разбиения к макроплиткам на рисунке~\ref{fig:chainDL2}.

\medskip

\leftskip=0cm

\medskip

На получающемся третьем уровне все  макроплитки, левый верхний угол которых попадает в узел $\mathbb{DL}$, будут занимать левое верхнее положение в своих родительских макроплитках. Значит, на следующих уровнях подразбиения, типы ребер-границ макроплиток не изменятся. То есть окружения вершин в цепях будут те же, то есть, окружение цепи четвертого и последующих уровней совпадает с окружением цепи третьего уровня.

Таким образом, существует только три возможных конфигурации окружения цепи с центром в узле типа $\mathbb{DL}$ (для цепей первого, второго и третьего уровней, при этом еще может быть выбрано разное главное ребро $t$). Заметим, что по окружению цепи, мы можем установить, является ли центр цепи узлом типа $\mathbb{DL}$, какого уровня цепь, а также сам параметр $t$. Действительно, $1$-цепь с центром в $\mathbb{DL}$ содержит одновременно узлы с окружением $(tA,tA,\mathbf{6A},\mathbf{2A})-(tB,\mathbf{3B},\mathbf{6A},\mathbf{2A})$ и $(tA,tA,\mathbf{6A},\mathbf{2A})-(\mathbf{7A},tB,\mathbf{3B},\mathbf{4A})$, которые вместе ни в какой другой цепи не встречаются. $2$-цепь с центром в $\mathbb{DL}$ содержит одновременно узлы с окружением $(tA,tA,\mathbf{1A},\mathbf{3A})-(tB,\mathbf{3B},\mathbf{1A},\mathbf{3A})$ и $(\mathbf{3A},\mathbf{8A},\mathbf{5A},\mathbf{6B})-(\mathbf{8B},\mathbf{7B},\mathbf{7B},\mathbf{7B})$, которые вместе также ни в какой другой цепи не встречаются.
Для $3$-цепи такой парой будет $(tA,tA,\mathbf{1A},\mathbf{3A})-(tB,\mathbf{3B},\mathbf{1A},\mathbf{3A})$ и $(\mathbf{3A},\mathbf{8A},\mathbf{1A},\mathbf{3A})-(\mathbf{8B},\mathbf{7B},\mathbf{1A},\mathbf{3A})$

Заметим также, что зная тип ребра $t$, мы можем выписать полностью все окружение $1$-цепей, $2$-цепи, $3$-цепей с центром в узле типа $\mathbb{DL}$.

\medskip

\leftskip=0cm

Из сказанного выше следует, что для $\mathbb{DL}$-узлов значение функции $\mathbf{LevelPlus}$ может быть установлено по ее аргументу.
На рисунках~\ref{fig:chainDL1} и~\ref{fig:chainDL2} знаком ``$\otimes$'' отмечены узлы, являющиеся верхними правыми или левыми нижними углами в макроплитках, где середина верхней стороны попадает в узел цепи $X$ . Можно убедиться в том, что для цепей первого, второго и третьего уровней все типы узлов с крестами, а также типы ребер на которых они лежат, мы можем выписать для каждого заданного $X$. Таким образом, мы можем вычислить функции $\mathbf{TopRightType}$ и $\mathbf{BottomLeftType}$.

Серым цветом на рисунке выделены плитки, где левый нижний угол (при дальнейшем разбиении) будет попадать в вершину цепи. Можно заметить, что зная эту вершину (ее место в цепи) мы можем установить и окружение правого нижнего угла в соответствующей макроплитке. Кроме того, мы также можем установить тип левого верхнего угла, а значит и всю цепь, содержащую середину верхней стороны в соответствующей макроплитке. (Все это также очевидно проверяется для $3$-цепи.) Все это значит, что мы можем вычислить функции $\mathbf{TopFromCorner}$, $\mathbf{RightCorner}$.

Функции $\mathbf{TopFromRight}$ и $\mathbf{BottomRightTypeFromRight}$ могут быть применимы только к $1$-цепи (в остальные входят только узлы типов $\mathbb{UL}$ и $\mathbb{LU}$). В случае $1$-цепи аргументом могут быть только узел, лежащий на $\mathbf{l}$ ребре (относительно $\mathbb{DL}$ узла). Пусть $X$ этот узел.  Тогда
$\mathbf{BottomRightTypeFromRight}$ будет наш же узел $\mathbb{DL}$. Узел, на которые указывает $\mathbf{TopFromRight}$ отмечен стрелкой на рисунке~\ref{fig:chainDL1}. Таким образом, $\mathbf{TopFromRight}$ это $1$-цепь вокруг $\mathbb{A}$ (указатель $3$).

\medskip

\begin{figure}[hbtp]
\centering
\includegraphics[width=0.9\textwidth]{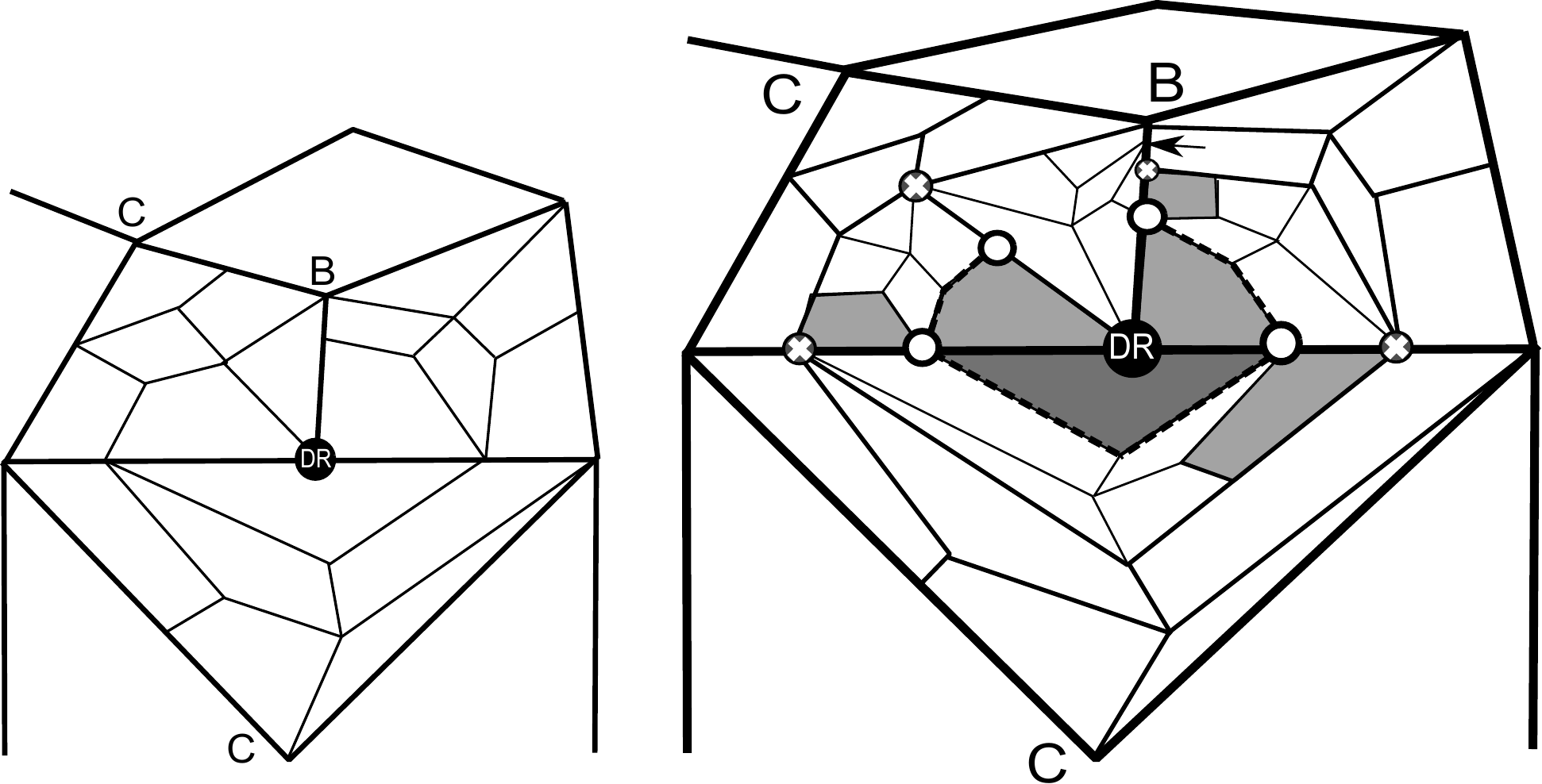}
\caption{$1$-цепь узла $\mathbb{DR}$.}
\label{fig:chainDR1}
\end{figure}

\begin{figure}[hbtp]
\centering
\includegraphics[width=0.8\textwidth]{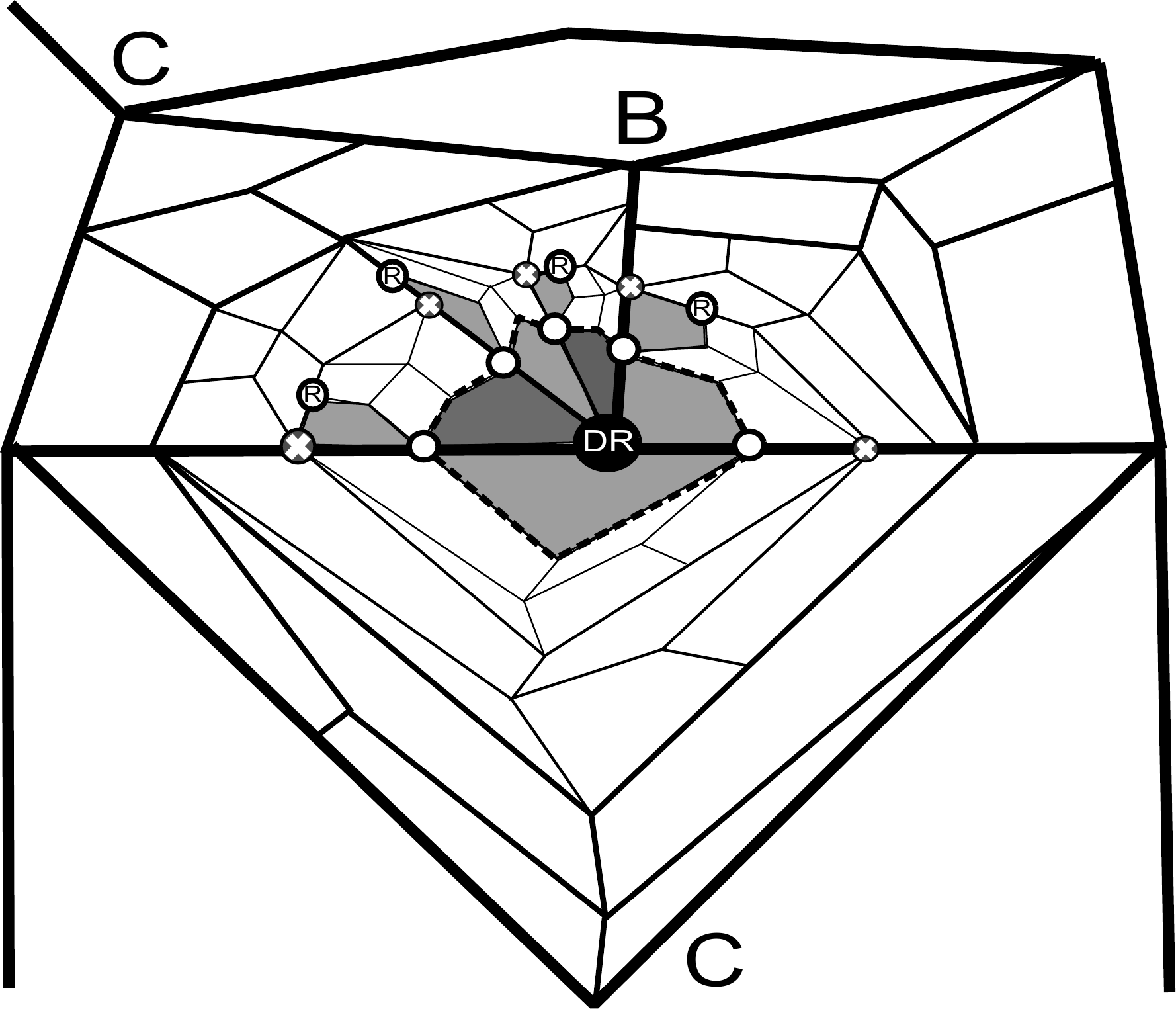}
\caption{$2$-цепь узла $\mathbb{DR}$.}
\label{fig:chainDR2}
\end{figure}

{\bf Цепи с центром в узлах типа $\mathbb{DR}$ и $\mathbb{RD}$.}
На рисунках~\ref{fig:chainDR1} и~\ref{fig:chainDR2} белыми точками отмечены $1$-цепи и $2$-цепи узла $\mathbb{DR}$. $0$-цепи с центром в $\mathbb{DR}$ не существует.

Кроме указанных на рисунках цепей первого и второго уровней существует также $3$-цепь, получаемая применением операции разбиения к макроплиткам на рисунке~\ref{fig:chainDL2}.
Полностью аналогично цепи $\mathbb{DL}$, можно показать, что существует только три возможных конфигурации окружения цепи с центром в узле типа $\mathbb{DR}$.
Кроме того, зная тип ребра $t$, мы можем выписать полностью все окружение $1$-цепи, $2$-цепи, $3$-цепи с центром в узле типа $\mathbb{DR}$.

\medskip

\leftskip=0cm

\medskip

Таким образом, что для $\mathbb{DR}$-узлов значение функции $\mathbf{LevelPlus}$ может быть установлено по ее аргументу.

На рисунках~\ref{fig:chainDR1} и~\ref{fig:chainDR2} знаком ``$\otimes$'' отмечены узлы, являющиеся верхними правыми или левыми нижними углами в макроплитках, где середина верхней стороны попадает в узел цепи $X$.
Вычисление функций $\mathbf{TopRightType}$, $\mathbf{BottomLeftType}$, $\mathbf{TopFromCorner}$, $\mathbf{RightCorner}$ полностью аналогично случаю $\mathbb{DL}$-цепи.

\medskip

\leftskip=0cm

\medskip

Функции $\mathbf{TopFromRight}$ и $\mathbf{BottomRightTypeFromRight}$ могут быть применимы только к $1$-цепи (в остальные входят только узлы типов $\mathbb{UL}$ и $\mathbb{LU}$). В случае $1$-цепи аргументом могут быть только узел, лежащий на $\mathbf{r}$ ребре (относительно $\mathbb{DR}$ узла). Пусть $X$ этот узел.  Тогда
$\mathbf{BottomRightType}$ будет наш же узел $\mathbb{DR}$. Узел, на которые указывает $\mathbf{TopFromRight}$ отмечен стрелкой на рисунке~\ref{fig:chainDR1}. Таким образом, $\mathbf{TopFromRight}$ это $1$-цепь вокруг $\mathbb{B}$ (указатель $2$).

\medskip

Нам остается разобрать цепи с центрами в краевых и угловых вершинах. Важное замечание: при разборе цепи с центром в краевой вершине часть узлов в попадает на край макроплитки. Под их окружением, в данном случае, мы понимаем базовое окружение в рамках нашей макроплитки (эта макроплитка может быть подклееной и тогда вершина на ее краю одновременно имеет какое-то другое окружение в основной своей области, вот это окружение мы не рассматриваем).

С учетом этого замечания, для рассмотрения цепей с центрами в вершинах типов $\mathbb{U}$, $\mathbb{L}$, $\mathbb{R}$, $\mathbb{D}$ достаточно применить те же рассуждения как и для цепей $\mathbb{UL}$, $\mathbb{LD}$, $\mathbb{RD}$, $\mathbb{DL}$, но рассматривать только половину картинки для каждой цепи. В силу полной аналогичности этого рассмотрения, мы не будем тут приводить этот разбор.

Цепи с центрами в угловых вершинах рассмотрим отдельно.

\medskip

\begin{figure}[hbtp]
\centering
\includegraphics[width=0.9\textwidth]{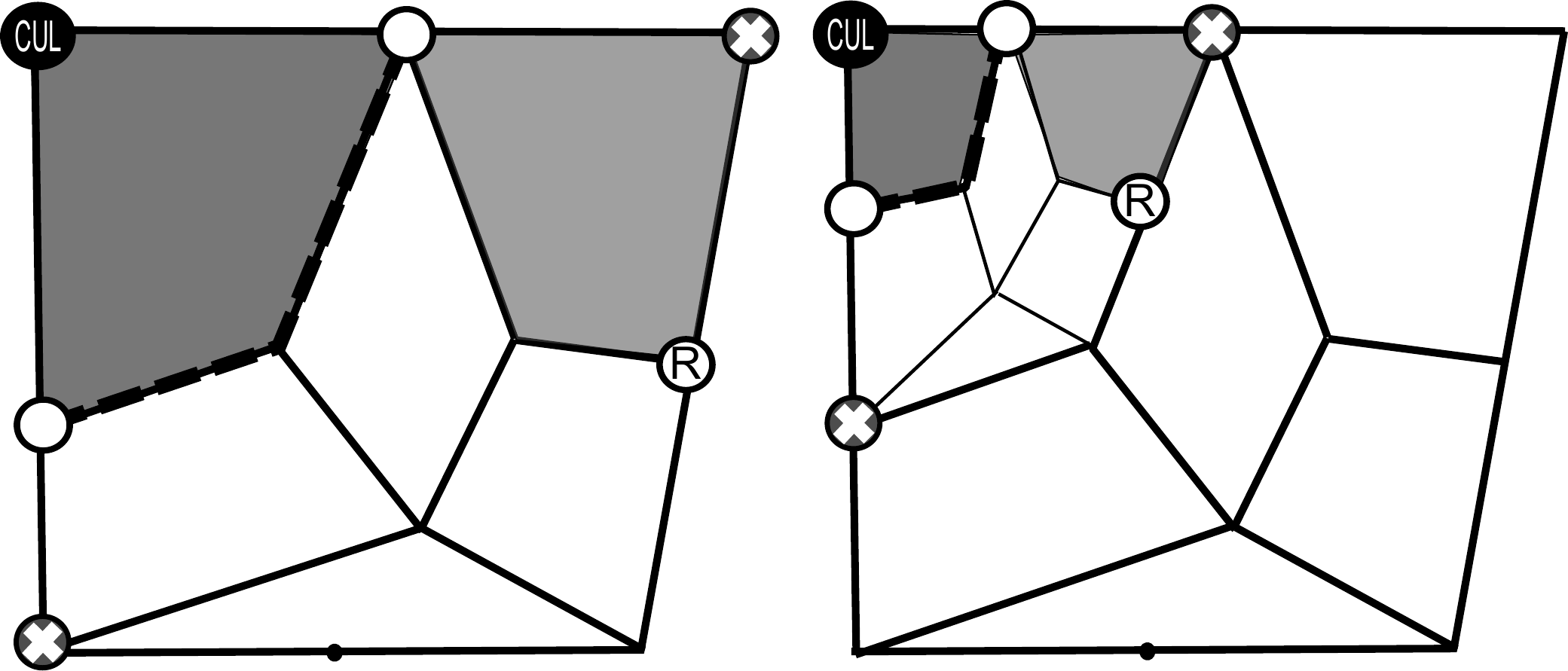}
\caption{$0$-цепь и $1$-цепь узла $\mathbb{CUL}$.}
\label{chainCUL01}
\end{figure}

\begin{figure}[hbtp]
\centering
\includegraphics[width=0.9\textwidth]{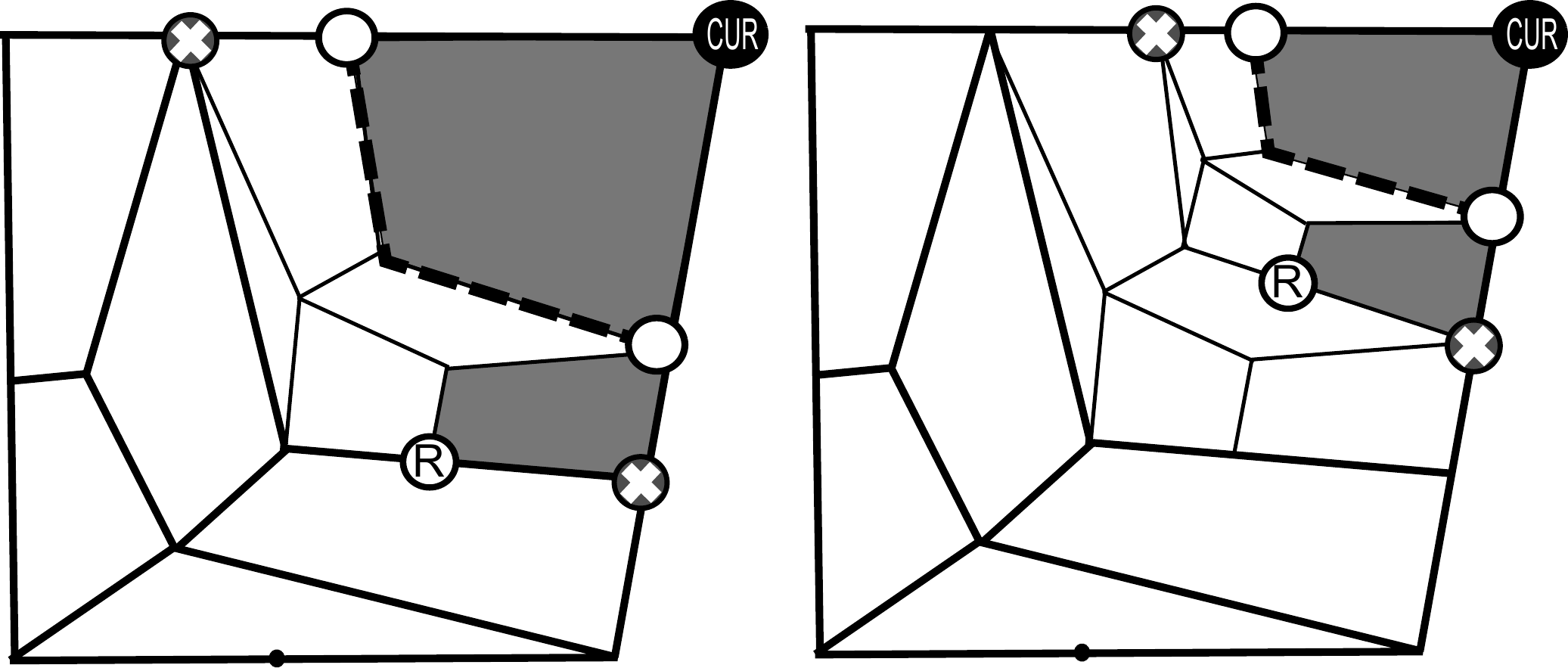}
\caption{$1$-цепь и $2$-цепь узла $\mathbb{CUR}$.}
\label{chainCUR12}
\end{figure}

\begin{figure}[hbtp]
\centering
\includegraphics[width=0.9\textwidth]{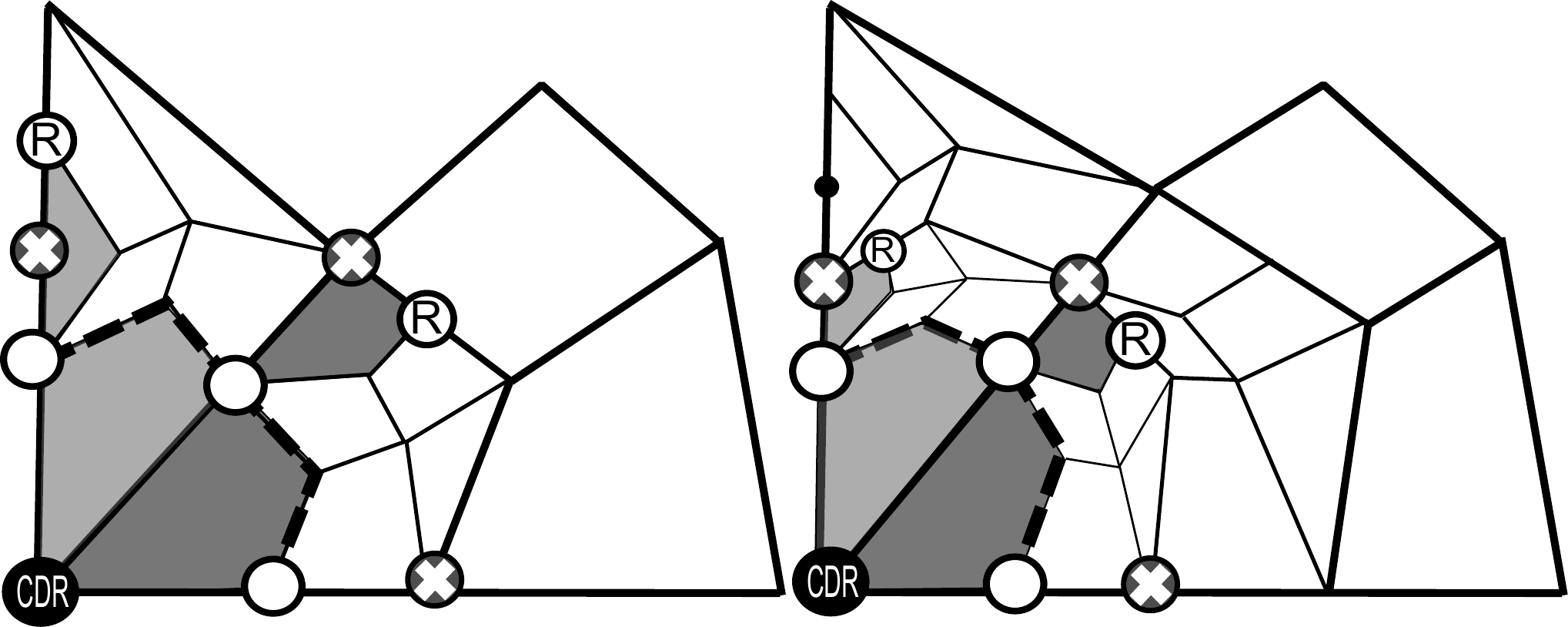}
\caption{$1$-цепь и $2$-цепь узла $\mathbb{CDR}$.}
\label{chainCDR12}
\end{figure}

\begin{figure}[hbtp]
\centering
\includegraphics[width=0.9\textwidth]{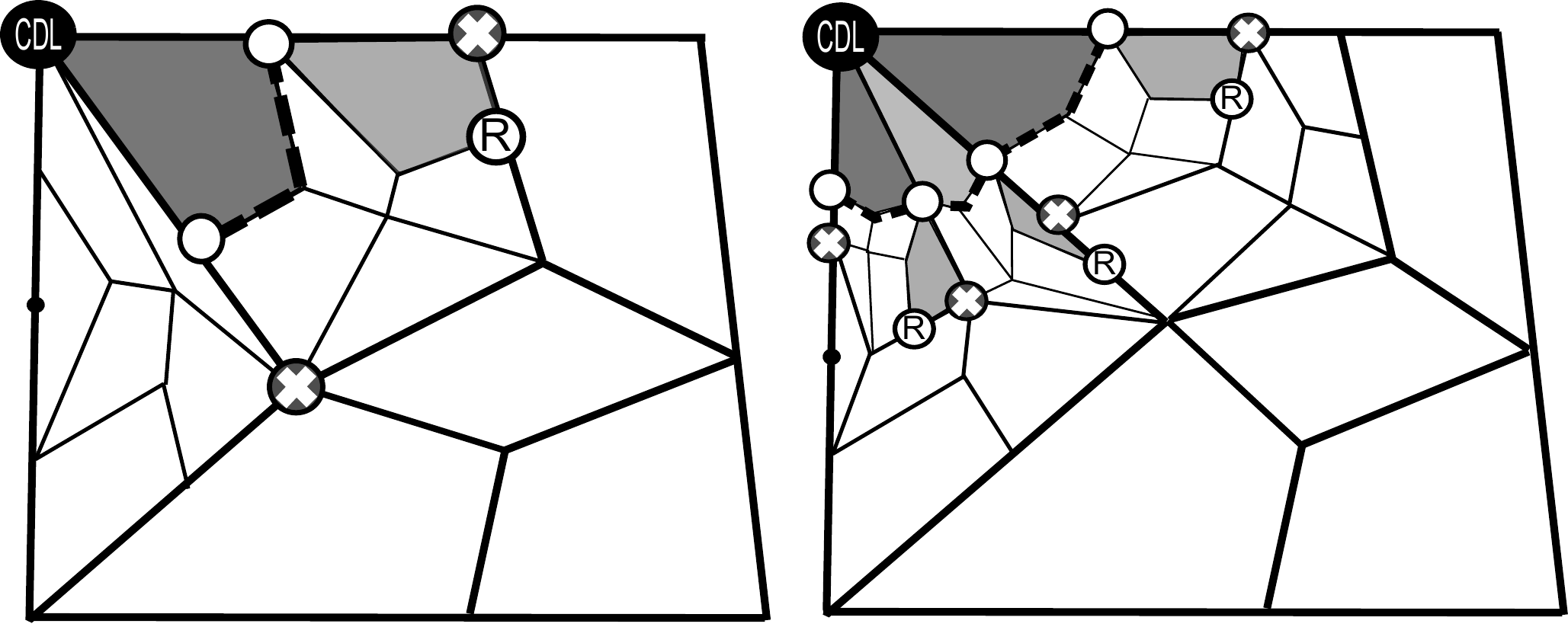}
\caption{$1$-цепь и $2$-цепь узла $\mathbb{CDL}$.}
\label{chainCDL12}
\end{figure}

{\bf Цепи с центром в узлах типа $\mathbb{CUL}$, $\mathbb{CUR}$, $\mathbb{CDR}$, $\mathbb{CDL}$.}
На рисунках~\ref{chainCUL01}, ~\ref{chainCUR12},~\ref{chainCDR12},~\ref{chainCDL12}    белыми точками отмечены цепи c вершинами в узлах типов  $\mathbb{CUL}$, $\mathbb{CUR}$, $\mathbb{CDR}$, $\mathbb{CDL}$.

Аналогично предыдущим цепям, можно показать, что существует только две возможных конфигурации окружения цепи в каждом из этих четырех случаев.

\medskip

\leftskip=0cm

\medskip

\leftskip=0cm

\medskip

\leftskip=0cm

\medskip

\leftskip=0cm

\medskip

Мы используем те же обозначения, что и в предыдущих случаях, крестами отмечены вершины, являющиеся правыми верхними или левыми нижними.

Вычисление функций $\mathbf{LevelPlus}$, $\mathbf{TopRightType}$, $\mathbf{BottomLeftType}$, $\mathbf{TopFromCorner}$, $\mathbf{RightCorner}$ очевидно, а функции
 $\mathbf{TopFromRight}$ и $\mathbf{BottomRightTypeFromRight}$ не могут быть применимы к вершинам данным цепям.

\medskip

{\bf Узлы типа $\mathbb{DR}$ и $\mathbb{RD}$.}
Для некоторых функций нужно также рассмотреть ситуации, когда аргументом является узел $\mathbb{DR}$ (или симметричный случай $\mathbb{RD}$).

\begin{figure}[hbtp]
\centering
\includegraphics[width=0.5\textwidth]{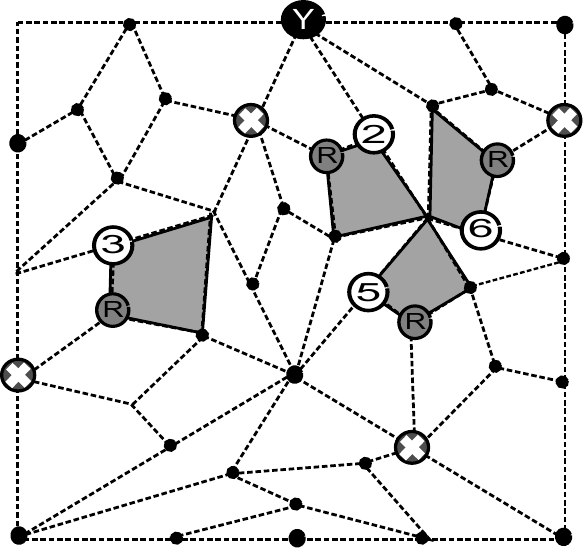}
\caption{Случаи расположения $\mathbb{DR}$-узла.}
\label{RDnode}
\end{figure}

Окружение $\mathbb{DR}$-узла может быть четырех видов, по числу внутренних ребер, на которых он может располагаться. На рисунке~\ref{RDnode} отмечены эти четыре возможные ситуации. Во всех случаях начальником нашего $\mathbb{DR}$-узла является вершина $Y$ в середине верхней стороны. На рисунке отмечены серым цветом макроплитки, левый нижний угол которых попадает в $\mathbb{DR}$-узел. Таким образом, значение функции $\mathbf{RightCorner}$ попадает в узлы отмеченные как ``R''. Это $\mathbb{B}$-узел, и его окружение во всех случаях можно вычислить, учитывая, что мы знаем окружение узла в середине верхней стороны.

\medskip

\leftskip=0cm

\medskip

Значение $\mathbf{TopFromCorner}$  указывает на $1$ цепь вокруг $\mathbb{A}$-узла с указателем $\mathbf{ld}$ для одного из случаев расположения, и $1$-цепь вокруг $\mathbb{B}$-узла для трех остальных случаев, указатели ясны из рисунка. Значение $\mathbf{BottomRightTypeFromRight}$ также очевидно во всех четырех случаях.

Теперь рассмотрим функцию $\mathbf{TopFromRight}$. Вершины, являющиеся ее значениями, отмечены знаком ``$\otimes$'' на рисунке. Поскольку мы знаем окружение вершины $Y$ в середине верхней стороны, мы можем вычислить и все требуемые окружения.

\medskip

\subsection{Указатели} \label{pointers}
Пусть узел $X$ принадлежит некоторой цепи, а узел $Y$ является ее центром, причем известно окружение $X$ и тип внутреннего ребра, на котором лежит $X$ (один из восьми типов внутренних ребер), а также известно окружение макроплитки $T$, которой принадлежит это ребро. Покажем, как восстановить {\it указатель} $X$, то есть тип ребра входа-выхода для $Y$.

Почти для всех ребер {\it указатель}, являющийся ребром выхода из $Y$ устанавливается сразу по типу ребра:
Для  ребра типа $1$, {\it указатель} будет очевидно $\mathbf{u}_2$, для ребра типа $2$ будет $\mathbf{u}_1$, для ребра типа $3$ будет $\mathbf{l}$, для ребра типа $4$ -- в зависимости от типа $Y$: для $\mathbb{A}$ указатель $2$, для $\mathbb{C}$ будет $1$. Для $5$ ребра опять в зависимости от типа $Y$ -- для $\mathbb{B}$ -- указатель $3$, для $\mathbb{C}$ будет $2$. Для $6$ ребра указатель $\mathbf{r}$.

Пусть ребро имеет тип $7$. Зная окружение макроплитки $T$, мы можем установить ее положение в родительской макроплитке. Далее вычисляем значение указателя: для левого-верхнего положения $T$ это $\mathbf{l}_2$, для левого нижнего $\mathbf{ld}$, для среднего $\mathbf{mid}$, для правого верхнего $\mathbf{u}_3$, для правого нижнего -- $\mathbf{rd}$, для нижнего -- $\mathbf{d}$.

Пусть ребро имеет тип $8$. Зная окружение макроплитки $T$, мы можем установить ее положение в родительской макроплитке. Далее вычисляем значение указателя: для левого-верхнего $\mathbf{lu}$ , для левого нижнего $\mathbf{ld}$, для среднего $\mathbf{mid}$, для правого верхнего $\mathbf{ru}$, для правого нижнего -- $\mathbf{rd}$. В случае нижнего положения $T$,  по типу верхней (и правой) стороны $T$ можно узнать, какое положение уже родительская макроплитка для $T$  занимает в своей родительской макроплитке  и в зависимости от этого установить значение указателя:  для левого-верхнего положения  -- $\mathbf{l}_3$, для левого нижнего $\mathbf{ld}_2$, для среднего $\mathbf{mid}_2$, для правого верхнего $\mathbf{u}_4$, для правого нижнего -- $\mathbf{r}_3$, для нижнего -- $\mathbf{d}_2$.

\medskip

\subsection{Дополнительные функции} \label{addfunc}

\begin{figure}[hbtp]
\centering
\includegraphics[width=0.4\textwidth]{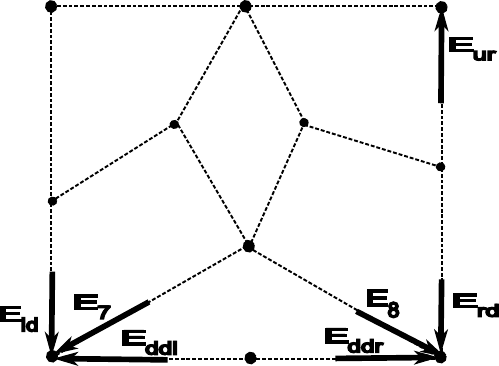}
\caption{Функции $E_7$, $E_8$, $E_{\mathbf{ld}}$, $E_{\mathbf{rd}}$, $E_{\mathbf{ur}}$, $E_{\mathbf{ddl}}$, $E_{\mathbf{ddr}}$ обозначают входящие или выходящие ребра в углы.}
\label{fig:e7e8}
\end{figure}

Символами $E_7$ и $E_8$ будем обозначать код соответствующего входящего и выходящего ребра. Очевидно, что зная по окружению узла в левом нижнем углу можно установить $E_7$, а по окружению узла в правом нижнем углу можно установить $E_8$ -- (рисунок~\ref{fig:e7e8}).

\medskip

Аналогично, обозначим символами  $E_{\mathbf{ld}}$,  $E_{\mathbf{rd}}$, $E_{\mathbf{ur}}$ коды входящих или выходящих ребер в нижние углы по левой или правой сторонам и в верхний угол по правой стороне.  Также символами $E_{\mathbf{ddl}}$,  $E_{\mathbf{ddr}}$ обозначим коды ребер в нижние углы по нижней стороне, слева и справа. Очевидно, что зная окружение соответствующего узла, можно установить букву, обозначающую нужное ребро (рисунок~\ref{fig:e7e8}).

\medskip

Пусть узел $X$ принадлежит некоторой цепи. Поскольку нам известны окружения всех узлов в цепи, а также типы указателей и их порядок следования, мы можем рассмотреть следующий по часовой стрелке узел в цепи, после $X$. Будем обозначать его как $\mathbf{Next}(X)$. Аналогично можно определить следующий после $X$ узел против часовой стрелки: $\mathbf{Prev}(X)$. Например, если в некоторой макроплитке известно окружение $X$ узла в середине верхней стороны, то окружение узла в середине левой стороны будет $\mathbf{Next}(X)$.

\medskip

Также отдельно определим функцию $\mathbf{BottomRightType}$, аргументами которой являются окружение и информация узла $X$, являющимся серединой верхней стороны в некоторой макроплитке $T$,
а значением --  тип узла в правом нижнем углу $T$ и ребро, на котором он расположен, если это боковой узел.

Значение этой функции вычисляется следующим образом. Мы можем установить положение $T$ в родительской макроплитке $T'$. Далее, для левого верхнего и левого нижнего положений тип правого нижнего угла будет $\mathbb{A}$, для среднего, правого верхнего и правого нижнего положений тип правого нижнего угла будет $\mathbb{B}$.

В случае нижнего положения $T$, правый нижний угол попадает в левый нижний угол $T'$. По типу правого и верхнего ребер $T$ можно установить, какое положение занимает уже $T'$ в своей родительской макроплитке $T''$.

Для левого верхнего положения правый нижний угол $T$ попадет в середину левой стороны $T''$ и его тип будет вычисляться функцией $\mathbf{Next.FBoss}(X)$.

Для левого нижнего, среднего и нижнего положений правый нижний угол $T$ попадет в $\mathbb{C}$-узел $T''$.

Для правого верхнего положения правый нижний угол $T$ попадет в середину верхней стороны $T''$ и его тип будет вычисляться как $\mathbf{FBoss}(X)$.

Для правого нижнего положения правый нижний угол $T$ попадет в $\mathbb{B}$ узел $T''$.

\medskip

\begin{figure}[hbtp]
\centering
\includegraphics[width=0.7\textwidth]{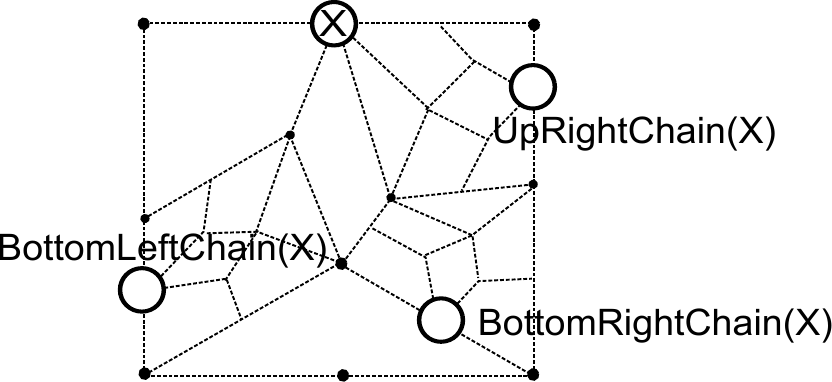}
\caption{Узлы, на которые указывают функции $\mathbf{BottomLeftChain}$, $\mathbf{BottomRightChain}$, $\mathbf{UpRightChain}$.}
\label{fig:addfunctions}
\end{figure}

Если узел $X$ в середине верхней стороны макроплитки известен, то функции $\mathbf{BottomRightType}$, $\mathbf{BottomLeftType}$, $\mathbf{TopRightType}$ дают нам типы узлов в трех соответствующих углах. Учитывая, что мы можем установить ребра входа в углы (символы $E_{\mathbf{ld}}$, $E_{\mathbf{rd}}$, $E_{\mathbf{ur}}$), мы можем установить окружения отмеченных на рисунке~\ref{fig:addfunctions} трех узлов. Обозначим соответствующие функции как $\mathbf{BottomLeftChain}$, $\mathbf{BottomRightChain}$, $\mathbf{UpRightChain}$.

\medskip

\medskip

Также определим функцию $\mathbf{RightFromB}(X)$, которая вычисляет окружение узла в середине правой стороны в макроплитке, по известному $\mathbb{B}$-узлу $X$. Значение этой функции вычисляется следующим образом.

Мы можем установить положение $T$ в родительской макроплитке $T'$.
Для левого верхнего положения, значением будет $\mathbb{UR}$-узел c окружением, соответствующим A0-цепи, с указателем $1$. Типы левой и верхней стороны в $T'$ те же что и в $T$, то есть мы их знаем по окружению $X$.
Для среднего, левого нижнего, правого верхнего, правого нижнего положений значением будет $\mathbb{RD}$-узел с окружением соответственно $3$, $2$, $6$, $5$ внутреннем ребрам.
Для нижнего положения значением будет $1$-цепь вокруг второго начальника $X$ (мы знаем его тип), с указателем $E_{\mathbf{ddl}}$ (вход по нижнему ребру в левый нижний угол).

\medskip

\medskip

\section{Оценка количества букв} \label{count}

Алфавит состоит из вершинных букв и реберных букв. Реберные буквы кодируют всевозможные ребра входа и выхода, включая возможные выходы в подклееные области. Плоских ребер, выходящих из вершины -- не более $10$ (больше всего в $\mathbb{C}$-вершине), всех подклееных ребер  не более $11$ ($\widehat{\mathbf{r}}$, $\widehat{\mathbf{u}_i}$; где $i=1,\dots, 4$, $\widehat{\mathbf{l}}$, $\widehat{\mathbf{l}_2}$, $\widehat{\mathbf{l}_3}$, $\widehat{\mathbf{d}}$, $\widehat{\mathbf{d}_2}$, $\widehat{\mathbf{d}_3}$). То есть реберных букв всего $21$.
Каждая вершинная буква кодирует одну комбинацию значений параметров вершины, то есть тип, уровень, расширенное окружение, информация и флаг подклейки.

\medskip

Посчитаем количество вершинных букв, используемых при кодировании. Четверок типов ребер всего существует $210$, подробно этот подсчет произведен в приложении~\ref{Appendix}. То есть, существует $210$ вариантов базовых окружений для типов $\mathbb{A}$, $\mathbb{B}$ и $\mathbb{C}$.
В цепях с центром в $\mathbb{UL}$-вершине существует $25$ возможных вершин ($7$ в $1$-цепи, и по $9$ в $2$-цепи и $3$-цепи). Каждая может быть одного из трех уровней, то есть всего $75$ возможных сочетаний базового окружения и уровня.

В таблице~\ref{chains} ниже представлено, сколько различных базовых окружений может быть у вершины заданного типа, входящей в цепь указанного уровня.

  \begin{table}[hbtp]
  \caption{Количество вершин с различными окружениями, сгруппированных по цепям}
\centering
 \begin{tabular}{|c|c|c|c|c|c|}   \hline
  тип центра цепи & $0$-цепь & $1$-цепь & $2$-цепь & $3$-цепь & всего \cr \hline
$\mathbb{UL}$/$\mathbb{LU}$ & -- & $7$ & $9$ & $9$ & $50$ \cr \hline
$\mathbb{UR}$/$\mathbb{RU}$ & -- & $7$ & $9$ & $9$ & $50$ \cr \hline
$\mathbb{DL}$/$\mathbb{LD}$ & -- & $4$ & $5$ & $5$ & $28$ \cr \hline
$\mathbb{DR}$/$\mathbb{RD}$ & -- & $4$ & $5$ & $5$ & $28$ \cr \hline
$\mathbb{A}$ & $420$ & $5$ & $5$ & -- & $430$ \cr \hline
$\mathbb{B}$ & -- & $6$ & $6$ & -- & $12$ \cr \hline
$\mathbb{C}$ & -- & $7$ & $10$ & $10$ & $27$ \cr \hline
$\mathbb{D}$ & -- & $2$ & $2$ & -- & $4$ \cr \hline
$\mathbb{L}$ & -- & $4$ & $5$ & $5$ & $14$ \cr \hline
$\mathbb{R}$ & -- & $4$ & $5$ & $5$ & $14$ \cr \hline
$\mathbb{U}$ & -- & $5$ & $6$ & $6$ & $17$ \cr \hline
$\mathbb{CUL}$ & -- & $2$ & $2$ & -- & $4$ \cr \hline
$\mathbb{CUR}$ & -- & $2$ & $2$ & -- & $4$ \cr \hline
$\mathbb{CDL}$ & -- & $2$ & $4$ & $4$ & $10$ \cr \hline
$\mathbb{CDR}$ & -- & $3$ & $3$ & -- & $6$ \cr \hline
  \end{tabular}
\label{chains}
\end{table}

Таким образом, всего возможно не более $698$ окружений для вершин, входящих в некоторую цепь, причем для каждой из этих вершин также известен ее тип. Каждая из них может быть одного из трех уровней, итого $2094$ сочетаний {\it тип-уровень-окружение}.

\medskip

Теперь посчитаем количество вершин, не входящих в цепи. В таблице~\ref{nonchains} ниже подсчитаны все возможные окружения вершин, не входящих в цепи.

  \begin{table}[hbtp]
  \caption{Количество вершин с разными окружениями, сгруппированных по типам}
\centering
 \begin{tabular}{|c|c|c|c|}   \hline
  тип вершины & число окружений & число уровней & всего  \cr \hline
$\mathbb{DR}$/$\mathbb{RD}$ & $4$ & $3$ & $12$ \cr \hline
$\mathbb{A}$ & $210$ & -- & $210$  \cr \hline
$\mathbb{B}$ & $210$ & -- & $210$  \cr \hline
$\mathbb{C}$ & $210$ & -- & $210$   \cr \hline
$\mathbb{D}$ & $1$ & $3$ & $3$  \cr \hline
$\mathbb{R}$ & $2$ & $3$ & $6$ \cr \hline
$\mathbb{CUL}$ & $1$ & -- & $1$  \cr \hline
$\mathbb{CUR}$ & $1$ & -- & $1$ \cr \hline
$\mathbb{CDL}$ & $1$ & -- & $1$  \cr \hline
$\mathbb{CDR}$ & $1$ & -- & $1$  \cr \hline
  \end{tabular}
\label{nonchains}
\end{table}

Заметим, что у $\mathbb{D}$-вершины может быть только окружение $(\mathbf{left},\mathbf{top},\mathbf{right},\mathbf{bottom})$, так как нижней стороной макроплитка примыкает к краю только являясь при этом подклееной плиткой.
Аналогично, у $\mathbb{R}$-вершины может быть два окружения $(\mathbf{left},\mathbf{top},\mathbf{right},\mathbf{bottom})$ и $(\mathbf{8B},\mathbf{bottom},\mathbf{bottom},\mathbf{7B})$, так как она макроплитка может примыкать правой стороной только к правому или нижнему краям и только в двух этих случаях.

То есть всего существует не более $655$ различных окружений для вершин не из цепей. Итого всех сочетаний тип-уровень-базовое окружение- не более $698\times 3 + 655 = 2749$.

\medskip

Рассмотрим теперь тип в подклееной области. Он может быть один из следующих: $\mathbb{CUL}$, $\mathbb{CUR}$, $\mathbb{CDL}$, $\mathbb{U}$, $\mathbb{L}$.  У первых трех по одному возможному подклееному окружению. У $\mathbb{U}$-вершины всего $6$ вариантов (по одному узлу в $\mathbb{U}1$, $\mathbb{U}2$, $\mathbb{U}3$, $\mathbb{L}1$, $\mathbb{L}2$, $\mathbb{L}3$, цепях). Аналогично у $\mathbb{L}$-вершины тоже $6$ вариантов подклееных окружений. С учетом трех возможных уровней, получается $3+6\times 3+ 6\times 3=39$ сочетаний.

Из определения операции подклейки следует, что к вершине типов $\mathbb{A}$, $\mathbb{B}$, $\mathbb{C}$ и к угловым вершинм подклееные макроплитки не могут примыкать левой или верхней стороной (только левым-верхним углом). Значит, вершина с нетривиальным подклееным окружением, не являющаяся ядром, в базовой плоскости либо имеет тип $\mathbb{DR}$/$\mathbb{RD}$, $\mathbb{D}$, $\mathbb{R}$, либо входит в некоторую цепь.
То есть для таких вершин базовых окружений всего $698\times 3 + 12 +3 +6 =2115$. То есть расширенных окружений (сочетаний базового и подклееного) имеется не более $2115\times 39 =82485$.

\medskip

{\bf Число значений параметра ``информация''}.
Для простоты будем говорить не ``число значений параметра информация'', а просто ``число информаций''.
Оценим это число. Пусть сначала все начальники не лежат на левом или верхнем краю подклееных плиток (то есть у них только базовое окружение).

Тогда первый начальник обязательно является вершиной, входящей в некоторую цепь, то есть $2094$ варианта. Второй начальник (левый нижний угол) либо попадает в вершину из цепи, либо $\mathbb{DR}$/$\mathbb{RD}$-вершину, либо $\mathbb{C}$-вершину либо $\mathbb{CDL}$-вершину. Всего $2094+12+210+1=2317$ вариантов. Окружение третьего начальника (правый нижний угол) может быть вычислено по окружению второго (функция $\mathbf{RightCorner}$).

Итак, если начальник один, есть $2094$ вариантов информации. Можно проверить, что тип правого нижнего угла можно восстановить по окружению узла в середине верхней стороны, для всех случаев, кроме двух -- цепи $\mathbb{A}0$ с указателем $3$ и цепи  $\mathbb{B}1$ c указателем $2$. То есть, всего возможных сочетаний окружения первого начальника и типа второго существует не более $2092+19+19=2130$.

Если начальника три, то окружение первого и третьего восстанавливаются по окружению второго (функции $\mathbf{TopFromCorner}$ и $\mathbf{RightCorner}$), то есть всего вариантов информации не более $2317$. Итого вариантов информации для плоского случая не более $2094+2317+2130=6541$.

\medskip

Теперь рассмотрим случай, когда хотя бы один начальник лежит на левом или верхнем краю подклееной макроплитки. Тогда в случае одного начальника вариантов информации будет не более $2115\times 39 =82485$ (любое сочетание базового и подклееного окружения). Для вершин с двумя начальниками заметим, что когда первый начальник имеет краевой тип, мы можем вычислить подклееный тип правого нижнего угла по подклееному окружению первого начальника. То есть тут столько же вариантов информаций, $82485$.

Для вершин с тремя начальниками, третий начальник не может попасть на левый или верхней край подклееной макроплитки, то есть у него будет только базовое окружение, которое можно вычислить по функции $\mathbf{RightCorner}$.
Возможны три случая:

{\bf 1.} И первый и второй начальники имеют и подклееное и базовое окружение (то есть они оба попали на левую и верхнюю стороны подклееной макроплитки). Тогда
расширенное окружение второго начальника можно выбрать не более чем $82485$ способами. Подклееное окружение первого начальника можно вычислить по функции $\mathbf{TopFromCorner}$. Если считать, что базовое окружение первого начальника -- любое,  то всего
вариантов информации $82485\cdot 2115$.

{\bf 2.} Только первый начальник попадает на край. Тогда число вариантов информации не превосходит $2749 \cdot 2115$. (Произвольное базовое окружение первого начальника и произвольное базовое окружение второго).

{\bf 3.} Только второй начальник попадает на край.  Окружение первого начальника можно вычислить по функции $\mathbf{TopFromCorner}$. Тогда число вариантов информации не превосходит $82485$.
Таким образом, общее число всех вариантов не более $82485\cdot 2115+ 2115\cdot 2749+82485=180352395$.

Любых возможных информаций не более $181\cdot 10^6$. Всего сочетаний (тип, уровень, расширенное окружение, информация) не более чем $181\cdot 10^6 \cdot (82485+2749)$ (то есть либо произвольное расширенное окружение, либо базовое окружение, без подклееного). Это число не превосходит $16 \cdot 10^{12}$.

\medskip

{\bf Число флагов макроплиток}.
Ядро, то есть вершина в левом верхнем углу подклееной макроплитки, может иметь произвольное базовое окружение и произвольную информацию. Флагом макроплитки называется сочетание типа, уровня, базового окружения и информации ядра этой макроплитки и упорядоченной пары двух выходящих ребер, которые определяют стороны макроплитки. Одно из ребер может вести в подклееные области.
Плоских выходящих ребер у любой вершины не более $10$. Различных подклееных ребер существует не более $11$. То есть возможных флагов существует не более $2749\cdot 11\cdot 10 \cdot 16 \cdot 10^{12}<5\cdot 10^{18}$.

\medskip

Таким образом, число всех возможных букв в алфавите не превосходит $82485 \cdot 16\cdot 10^{12} \cdot 5\cdot 10^{18}< 7\cdot 10^{36}$.

\section{Разбор случаев расположения путей} \label{flip_section}

\begin{figure}[hbtp]
\centering
\includegraphics[width=0.6\textwidth]{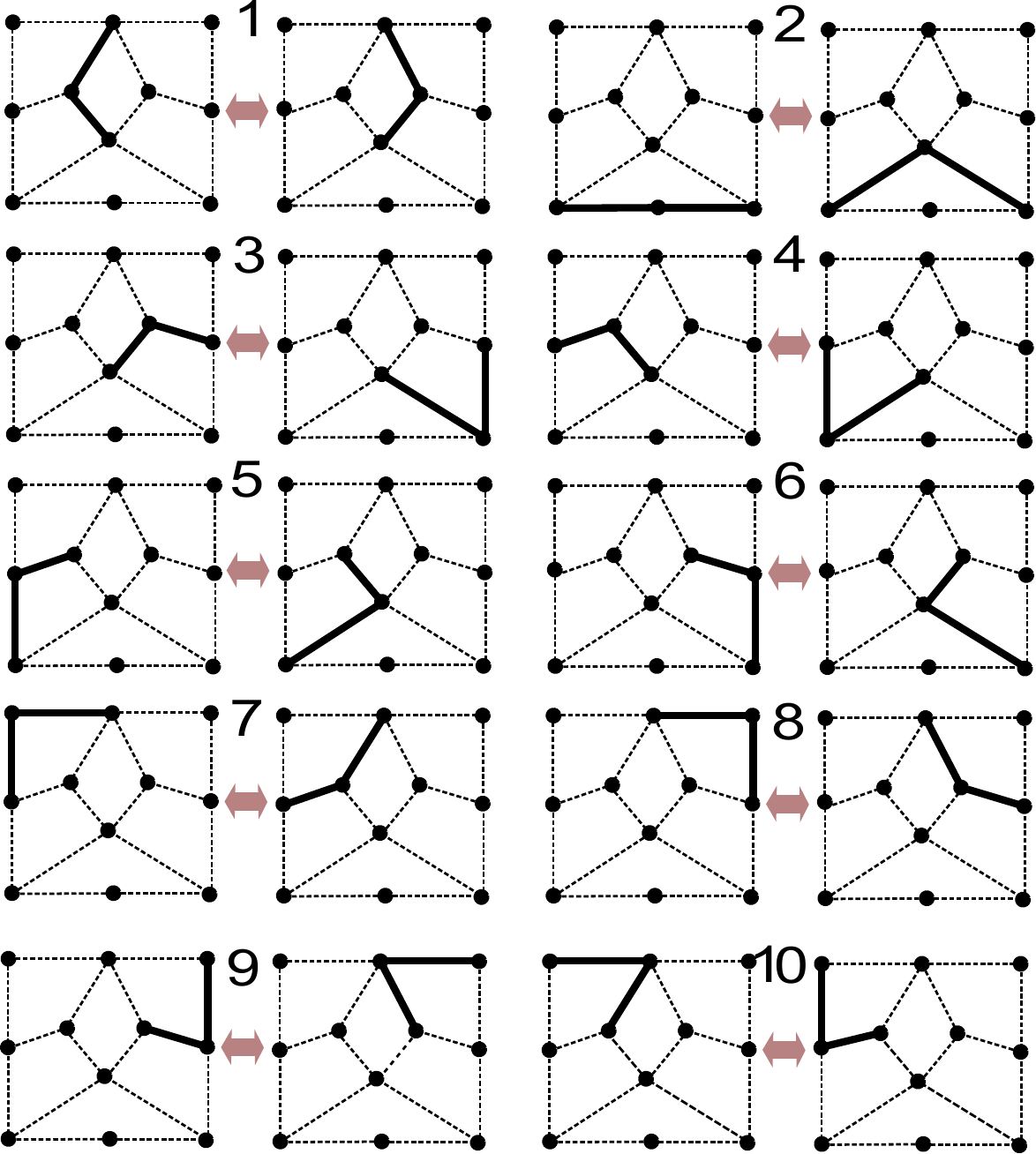}
\caption{По известному коду одного пути из пары можно восстановить код другого.}
\label{fig:flippaths}
\end{figure}

На рисунке~\ref{fig:flippaths} изображены пары путей. Далее мы покажем, как, зная код одного пути из пары, можно определить код другого пути либо установить, что код является мертвым.
В целом для этого требуется уметь вычислять код любой вершины в любой макроплитки по известным кодам трех остальных вершин. Рассмотрим произвольный путь $X_1e_1e_2X_2e_3e_4X_3$, где $X_1$, $X_2$, $X_3$ -- буквы, отвечающие кодам вершин, а $e_1$, $e_2$, $e_3$, $e_4$ -- буквы, отвечающие ребрам входа и выхода. Мы должны установить, к какой из десяти конфигураций относится наш путь, а также показать, как провести локальное преобразование, то есть получить код другого пути из соответствующей пары.

\medskip

В этой главе мы рассмотрим случай, когда все буквы $e_1$, $e_2$, $e_3$, $e_4$ отвечают плоским ребрам, то есть наш путь не выходит в подклееные области.
Случай выхода в подклееные области, когда хотя бы одна буква отвечает ребру, выходящему в подклейку, мы разберем в следующей главе.
Фактически, в этой и следующей главах мы вводим определяющие соотношения в полугруппе.

\medskip

\subsection{Обзор перебора случаев: сочетания типов сторон макроплитки}

Чтобы провести восстановление вершин для случаев $(1)-(10)$, требуется знать что расположено вокруг плитки, в которой идет наш путь.

\medskip

Итак, мы рассматриваем путь из двух последовательных ребер, проходящий внутри либо по границе некоторой макроплитки $T$. На рисунке~\ref{fig:flippaths} изображены макроплитки $T$ c проходящими внутри них нашими путями. Основная задача состоит в том, чтобы в каждом случае {\it восстановить путь}, то есть по известным значениям параметров трех вершин и входящим ребрам вычислить параметры четвертой вершины и код парного пути.

Для плоского пути, когда среди ребер вдоль пути не встречаются выходы в подклейки, параметр {\it флаг подклейки} восстанавливается тривиально -- для четвертой вершины он будет такой же, как и для остальных вершин (среди которых можно выбрать гарантированно лежащую не на границе подклееной макроплитки).
То есть нужно вычислить остальные параметры --  тип, уровень, окружение и информацию.

\medskip

Локальные преобразования с $1$ по $6$ мы рассмотрим каждое отдельно. При этом мы исследуем различные варианты расположения макроплитки, в которой проходит путь, и покажем, как вводятся определяющие соотношения в различных случаях расположения.

Локальные преобразования с $7$ по $10$ мы будем рассматривать все вместе для каждого случая расположения макроплитки.

\medskip

{\bf Обозначения.}
Буквами $e_1$, $e_2$, $e_3$, $e_4$, $e_{\mathbf{u}_1}$, $e_{\mathbf{u}_2}$, $e_{\mathbf{u}_3}$, $e_{\mathbf{ld}}$, $e_{\mathbf{r}}$, $e_{\mathbf{r}_2}$, $e_{\mathbf{r}_3}$ (и все подобные) будем обозначать входящие и выходящие ребра, соответствующие названию.  Символами $E_{\mathbf{l}}$, $E_{\mathbf{ld}}$, $E_{\mathbf{ddl}}$, (и подобными, указанными в разделе~\ref{addfunc} ``дополнительные функции'') будем обозначать ребра, получаемые в результате применения указанных функций. Все эти буквы мы будем использовать эти буквы при введении определяющих соотношений.

\smallskip

Иногда нам нужно будет зафиксировать часть параметров некоторой вершины $X$, чтобы потом использовать их для вычисления параметров других вершин. В этом случае мы будем использовать слово {\it назначение} $X$.

\smallskip

Пусть часть параметров вершины $X$ зафиксированы и известны. Пусть также оперируя известными параметрами $X$, мы вычисляем некоторые параметры вершин $X_1$, $X_2$, $X_3$. {\it Разрешенной} комбинацией будем называть такую упорядоченную четверку букв ($x$, $x_1$, $x_2$, $x_3$), что $x$ кодирует $X$ c зафиксированными параметрами, $x_i$ кодирует $X_i$, причем параметры совпадают с вычисленными. Заметим, что возможно несколько комбинаций разрешенных четверок, например, при разных параметрах флага подклейки. Часто кодирующие буквы мы будем обозначать так же, как сами вершины, например, $Z$, $J$, $Y$, $F$.

\medskip

{\bf Замечание.}
Рассматриваемые нами пути могут быть кусками более длинных путей. В том числе возможна ситуация, когда после прохождения нашего подпути, далее путь уходит по подклееному ребру, например, после вершины $J$. В этой ситуации, согласно определению, окружением $J$ является расширенное окружение, состоящее из базового окружения (в плоскости нашего пути) и подклееного окружения (в плоскости той подклееной плитки, куда уходит наш путь). Это отражается в том, что даже если тип, базовое окружение, флаг подклейки и информация у $J$ вычислены нами из условия расположения пути и остальных вершин, подклееное окружение $J$ часто может быть произвольным. В этом случае в качестве {\it разрешенных} букв, кодирующих $J$, можно выбрать буквы из множества с фиксированными параметрами типа, базового окружения, флаг подклейки и информации, но при этом с разными подклееными окружениями.

В этом случае мы будем говорить, что $J$ выбрана с точностью до подклееного окружения. Указанная ситуация может быть применима и к другим параметрам, например, информации.

\medskip

{\bf Комментарий о вкладе авторов.}

Первому автору (Иванов-Погодаев) принадлежит конструкция подстановочной системы семейства комплексов, система кодировки вершин, перебор для проверки детерминированности и  финальная схема приведения периодического слова к нулю.

\medskip

\subsection{Локальное преобразование 1}

\begin{figure}[hbtp]
\centering
\includegraphics[width=0.7\textwidth]{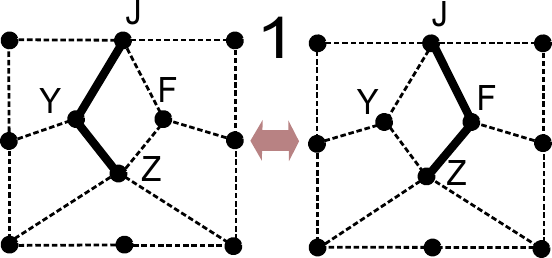}
\caption{Локальное преобразование $1$}
\label{flip1}
\end{figure}

Рассмотрим пару путей на рисунке~\ref{flip1}.

\medskip

\medskip

\medskip

{\bf Определяющие соотношения}.

\medskip

{\bf Назначение $Z$}. Пусть $Z$ -- вершина типа $\mathbb{C}$. Зафиксируем ее флаг подклейки, базовое окружение, подклееное окружение и информацию.

Теперь можно вычислить базовые окружения вершин $J$, $Y$, $F$, так как $J$ -- это просто первый начальник $Z$, а базовые окружения $Y$ и $F$ совпадают с окружением макроплитки, которое содержится в базовом окружении $J$. Зафиксируем некоторую информацию $J$. Информацию $Y$ и $F$ известна, так как известны окружения $J$ и вершины в правом нижнем углу (третий начальник $Z$).

Таким образом, выбрав произвольно базовое и подклееное окружение $Z$, подклееное окружение и информацию $J$, а также флаг макроплитки, мы можем вычислить значение остальных параметров. То есть, мы можем выписать множество разрешенных четверок букв ($Z$, $J$, $Y$, $F$), кодирующих описанное положение вершин.

Заметим, что если буква $Z$ зафиксирована, то буквы $Y$, $F$  -- точно определены, а буква $J$ -- определена с точностью до информации, то есть $J$ могут кодировать разные буквы, отличающиеся информацией.

\medskip

Для каждой разрешенной комбинации четырех букв $J$, $Z$, $Y$, $F$, введем следующие определяющие соотношения:

$Ze_1e_2Ye_1e_{\mathbf{u}_2}J=Ze_2e_3Fe_1e_{\mathbf{u}_1}J$

$Je_{\mathbf{u}_2}e_1Ye_2e_{1}J=Ze_{\mathbf{u}_1}e_1Fe_3e_{2}J$

\medskip

Два соотношения вводятся, потом что путь может быть пройден как в прямом порядке, так и в обратном.

Итак, $Z$ -- это произвольная буква типа $\mathbb{C}$ (произвольные расширенное окружение, флаг подклейки и информация), а у $J$ может быть произвольная информация и произвольное расширенное окружение. Остальные буквы в этом соотношении заданы однозначно.

\medskip

\subsection{Локальное преобразование 2}

\begin{figure}[hbtp]
\centering
\includegraphics[width=0.7\textwidth]{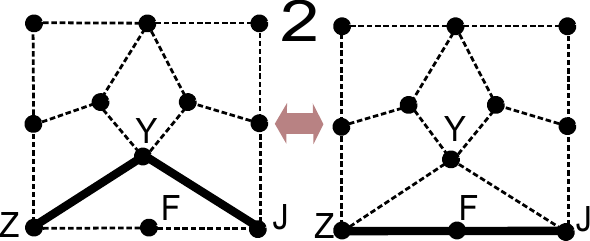}
\caption{Локальное преобразование $2$}
\label{flip2}
\end{figure}

Рассмотрим пару путей на рисунке~\ref{flip2}.

\medskip

\medskip

{\bf Определяющие соотношения}.
Мы рассмотрим по очереди шесть случаев расположения. Для каждого случая мы определим разрешенные множества букв, которыми могут быть вершины $Z$, $J$, $F$, $Y$. После чего для каждой комбинации букв, введем два соотношения, одно отвечает прямому обходу данного пути (по стрелкам), а другое -- обратному (против стрелок).

\medskip

\begin{figure}[hbtp]
\centering
\includegraphics[width=0.9\textwidth]{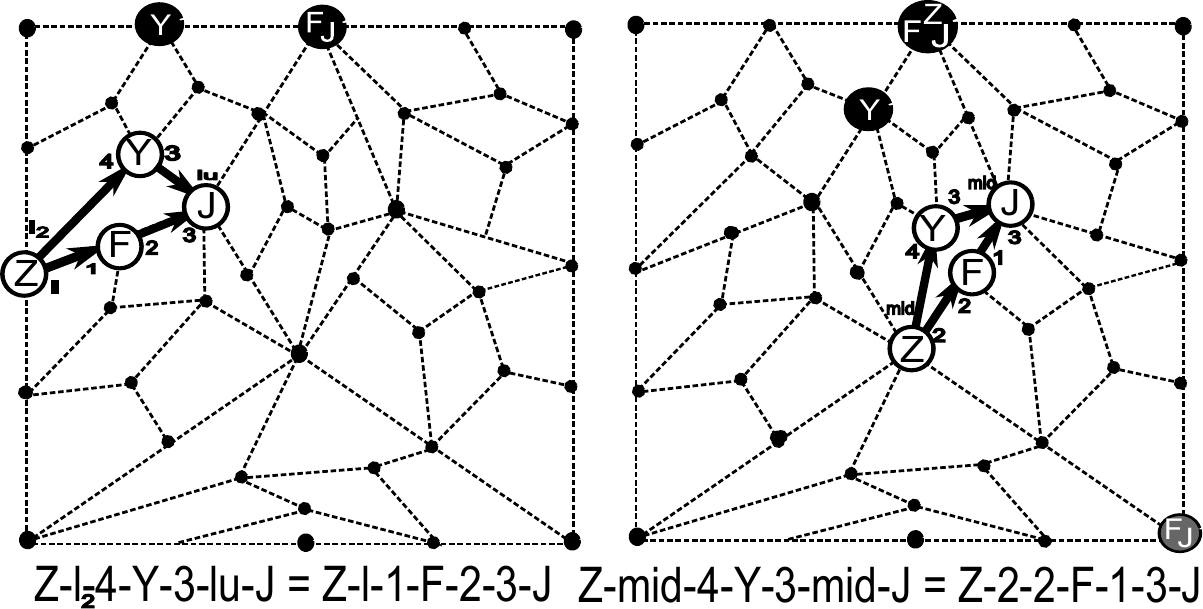}
\caption{Левое-верхнее и среднее положения. При входе и выходе в вершины обозначены коды входящих и выходящих ребер}
\label{pathplace2-1}
\end{figure}

{\bf 1. Левое-верхнее расположение}  (левая часть рисунка~\ref{pathplace2-1}).

\medskip

{\bf Назначение $Z$}. Пусть $Z$ -- вершина типа $\mathbb{LU/UL}$, $\mathbb{LD/DL}$ или $\mathbb{L}$  c произвольным расширенным окружением и информацией.

Начальники вершин $J$, $Y$, $F$ отмечены на рисунке~\ref{pathplace2-1} (левая часть) черными кругами.
Заметим, что мы можем вычислить их окружения. Одна из них входит в ту же цепь что и $Z$, а другая вычисляется с помощью функции $\mathbf{LevelPlus}$.  Значит, мы можем определить базовые окружения и информации вершин $J$, $Y$, $F$. Флаг подклейки может быть назначен произвольно (одинаковое значение для всех четырех вершин).

Таким образом, можно сформировать множество разрешенных четверок букв $Z$, $J$, $Y$, $F$.

\medskip

Для каждой разрешенной комбинации заданных букв, введем следующее определяющие соотношения:

\smallskip

$Ze_{\mathbf{l}_2}e_4Ye_3e_{\mathbf{lu}}J=Ze_{\mathbf{l}}e_1Fe_2e_3J$

$Je_{\mathbf{lu}}e_3Ye_4e_{\mathbf{l}_2}Z=Je_3e_2Fe_1e_{\mathbf{l}}Z$

\smallskip

Два соотношения отвечают двум направлениям обхода заданного участка, по стрелкам или против стрелок.

\medskip

{\bf Восстановление кода.} Зная коды $J$ и $Z$ можно вычислить код как $Y$, так и $F$, так как их начальники и базовые окружения легко вычисляются по известному окружению и начальникам $Z$, в частности, у $F$ общий набор начальников с $J$.

\medskip

\medskip

\medskip

{\bf 2. Среднее расположение} (правая часть рисунка~\ref{pathplace2-1}).

\medskip

{\bf Назначение $Z$}. Пусть $Z$ -- вершина с типом $\mathbb{C}$ и произвольным расширенным окружением, флагом подклейки и информацией.

Начальники вершин $J$, $Y$, $F$ отмечены черными кругами. Зная код $Z$, можно вычислить их окружения. То есть можно вычислить коды $J$, $Y$, $F$.

\medskip

Заметим, что если буква $Z$ зафиксирована, то $Y$, $F$ -- точно определены, а $J$ -- с точностью до произвольного подклееного окружения, при заданном базовом. То есть можно выделить множество разрешенных четверок букв, кодирующих наши вершины.

\medskip

Для каждой разрешенной комбинации заданных букв, введем следующие определяющие соотношения:

\smallskip

$Ze_{\mathbf{mid}}e_4Ye_3e_{\mathbf{mid}}J=Ze_2e_2Fe_1e_3J$

$Je_{\mathbf{mid}}e_3Ye_4e_{\mathbf{mid}}Z=Je_3e_1Fe_2e_2Z$

\smallskip

Два соотношения отвечают двум направлениям обхода заданного участка, по стрелкам или против стрелок.

\medskip

{\bf Восстановление кода.} Аналогично левому-верхнему расположению.

\medskip

\begin{figure}[hbtp]
\centering
\includegraphics[width=0.9\textwidth]{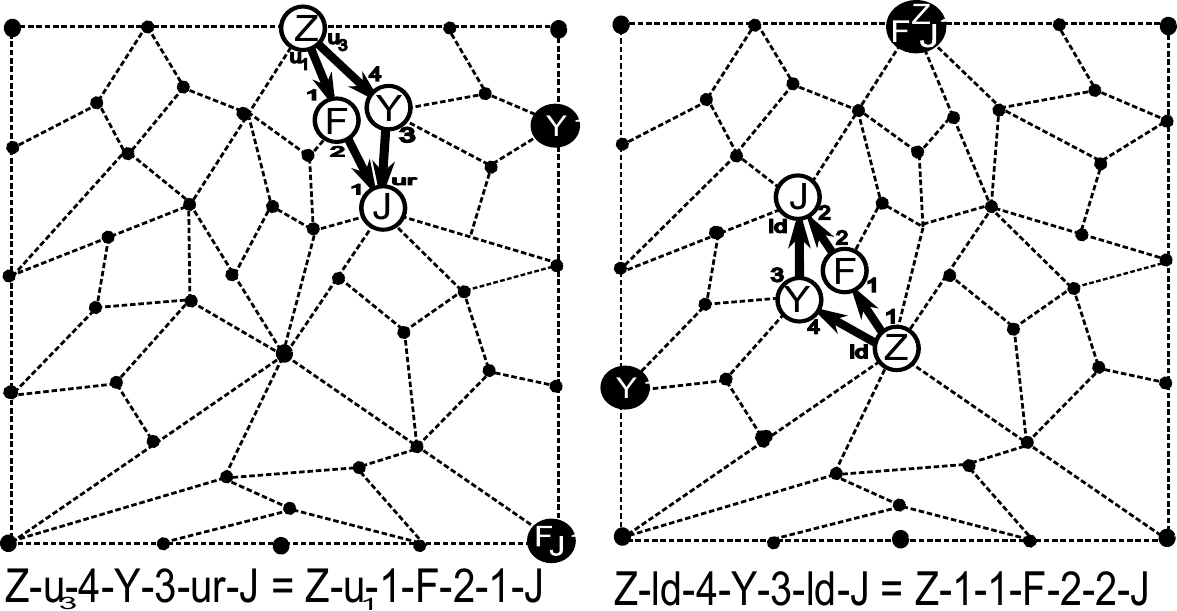}
\caption{Правое-верхнее и левое нижнее положения. При входе и выходе в вершины обозначены коды входящих и выходящих ребер}
\label{pathplace2-2}
\end{figure}

{\bf 3. Правое-верхнее расположение} (рисунок~\ref{pathplace2-2}, левая часть).

\medskip

{\bf Назначение $Z$}. Пусть $Z$ --  вершина типа $\mathbb{UL/LU}$ или $\mathbb{UR/RU}$  c произвольным расширенным окружением, флагом подклейки и информацией.

Окружение первого начальника $Y$ мы полуаем, используя функцию $\mathbf{UpRightChain}(Z)$.
Тип узла в правом нижнем углу -- $\mathbf{BottomRightType}(Z)$. То есть, коды всех остальных вершин можно вычислить.

Заметим, что если буква $Z$ зафиксирована, то буквы $F$, $Y$  -- точно определены, а буква $J$ -- с точностью до произвольного подклееного окружения, при заданном базовом.

\medskip

Для каждой разрешенной комбинации заданных букв, введем следующие определяющие соотношения:

\smallskip

$Ze_{\mathbf{u}_3}e_4Ye_3e_{\mathbf{ur}}J=Ze_{u_1}e_1Fe_2e_1J$

$Je_{\mathbf{ur}}e_3Ye_4e_{\mathbf{u}_3}Z=Je_1e_2Fe_1e_{\mathbf{u}_1}Z$

\smallskip

{\bf Восстановление кода.} Аналогично левому-верхнему расположению.

\medskip

\medskip

{\bf 4. Левое-нижнее расположение}  (рисунок~\ref{pathplace2-2}, правая часть).

\medskip

{\bf Назначение $Z$}. Пусть $Z$ вершина типа $\mathbb{C}$ c произвольным расширенным окружением, флагом подклейки и информацией.

Цепь вокруг левого нижнего угла $T$ можно получить как $\mathbf{BottomLeftChain}(Z)$. То есть, окружения вершин, отмеченных кругами можно вычислить, и коды всех вершин $J$, $Y$, $F$ можно вычислить, зная код $Z$.

Заметим, что если буква $Z$ зафиксирована, то буквы $F$, $Y$  -- точно определены, а буква $J$ -- с точностью до произвольного подклееного окружения, при заданном базовом.

\medskip

Для каждой разрешенной комбинации заданных букв, введем следующие определяющие соотношения:

\smallskip
$Ze_{\mathbf{ld}}e_4Ye_3e_{\mathbf{ld}}J=Ze_1e_1Fe_2e_2J$

$Je_{\mathbf{ld}}e_3Ye_4e_{\mathbf{ld}}Z=Je_2e_2Fe_1e_1Z$
\smallskip

{\bf Восстановление кода.} Аналогично левому-верхнему расположению.

\medskip

\begin{figure}[hbtp]
\centering
\includegraphics[width=0.9\textwidth]{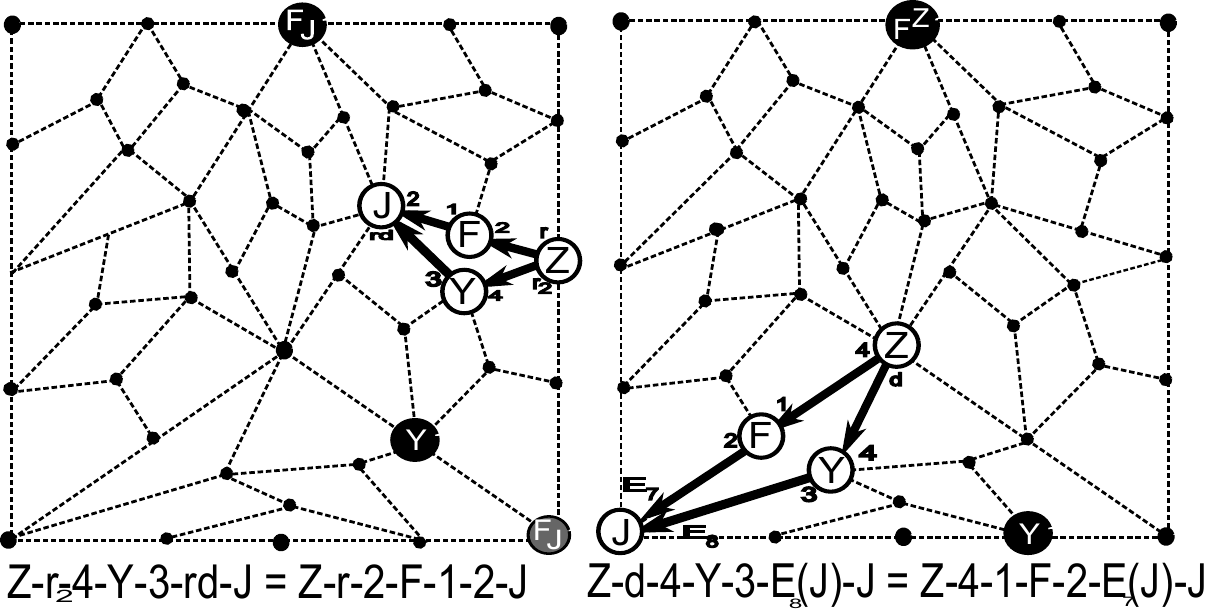}
\caption{Правое-нижнее и нижнее положения. При входе и выходе в вершины обозначены коды входящих и выходящих ребер}
\label{pathplace2-3}
\end{figure}

{\bf 5. Правое-нижнее расположение}  (рисунок~\ref{pathplace2-3}, левая часть).

\medskip

{\bf Назначение $Z$}. Пусть $Z$ -- вершина типа $\mathbb{RD/DR}$, $\mathbb{RU/UR}$  или $\mathbb{R}$  c произвольным расширенным окружением, флагом подклейки и информацией.

Цепь вокруг правого нижнего угла $T$ можно получить как $\mathbf{BottomRightChain}(Z)$. То есть, можно вычислить окружения всех начальников, значит и все коды вершин $J$, $Y$, $F$.

Заметим, что если буква $Z$ зафиксирована, то $F$, $Y$  -- точно определены, а $J$ -- с точностью до одного начальника и произвольного подклееного окружения, при заданном базовом.

\medskip

Для каждой разрешенной комбинации заданных букв, введем следующие определяющие соотношения:

\smallskip

$Ze_{\mathbf{r}_2}e_4Ye_3e_{\mathbf{rd}}J=Ze_{\mathbf{r}}e_2Fe_1e_2J$

$Je_{\mathbf{rd}}e_3Ye_4e_{\mathbf{r}_2}Z=Je_2e_1Fe_2e_{\mathbf{r}}Z$

\smallskip

{\bf Восстановление кода.} Аналогично левому-верхнему расположению.

\medskip

\medskip

{\bf 6. Нижнее расположение} (рисунок~\ref{pathplace2-3}, правая часть).

\medskip

{\bf Назначение $Z$}. Пусть $Z$ -- вершина типа $\mathbb{C}$ c произвольным расширенным окружением, флагом подклейки и информацией.

Заметим, что $J$ -- второй начальник $Z$. Окружения и информации $F$ и $Y$  мы также можем вычислить.

\medskip

{\bf Назначение $J$}. Тип, окружение и флаг подклейки для $J$ уже определены. Назначим произвольную информацию.

Теперь буквы $F$, $Y$  -- точно определены, а буква $J$ -- с точностью до произвольной информации.

\medskip

$E_7(J)$ и  $E_8(J)$ --  реберные буквы, полученные с помощью применения функций $E_7$ и $E_8$ к узлу $J$, то есть это тоже буквы, кодирующие ребра входов и выходов, но в зависимости от узла $J$.

\medskip

Для каждой разрешенной комбинации заданных букв, введем следующие определяющие соотношения:

\smallskip
$Ze_{\mathbf{d}}e_4Ye_3E_8(J)J=Ze_{4}e_1Fe_2E_7(J)J$

$JE_8(J)e_3Ye_4e_{\mathbf{d}}Z=JE_7(J)e_2Fe_1e_4Z$
\smallskip

{\bf Восстановление кода.} Аналогично левому-верхнему расположению.

\medskip

\subsection{Локальное преобразование 3}

\begin{figure}[hbtp]
\centering
\includegraphics[width=0.7\textwidth]{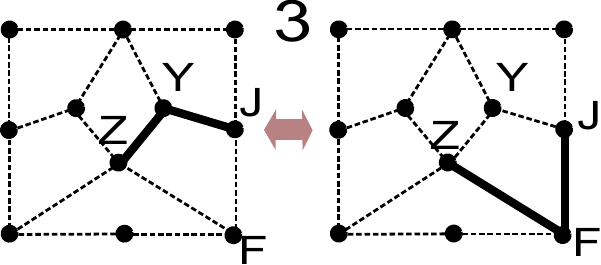}
\caption{Локальное преобразование $3$}
\label{flip3}
\end{figure}

Рассмотрим пару путей на рисунке~\ref{flip3}.

\medskip

\medskip

\begin{figure}[hbtp]
\centering
\includegraphics[width=0.9\textwidth]{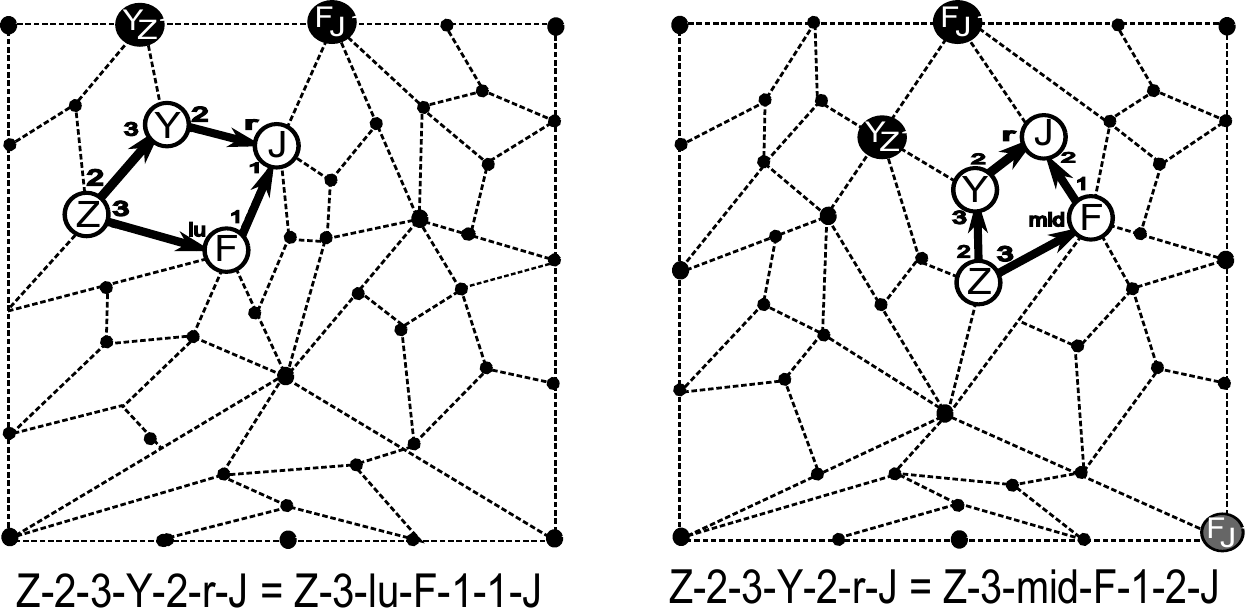}
\caption{Левое верхнее и среднее положения. При входе и выходе в вершины обозначены коды входящих и выходящих ребер}
\label{pathplace3-1}
\end{figure}

{\bf 1. Левое-верхнее расположение} (рисунок~\ref{pathplace3-1}, левая часть).

\medskip

{\bf Назначение $F$}. Пусть $F$ -- вершина типа $\mathbb{A}$ c произвольным базовым окружением, флагом подклейки и информацией.

Аналогично предыдущим случаям, этого достаточно, чтобы вычислить коды других вершин.

Заметим, что буквы $F$, $Y$ -- точно определены, а множества разрешенных букв $J$ и $Z$ -- различаются произвольным подклееным окружением, при заданном базовом.

\medskip

Для каждой разрешенной комбинации заданных букв, введем следующие определяющие соотношения:

\smallskip
$Ze_{2}e_3Ye_2e_{\mathbf{r}}J=Ze_{3}e_{\mathbf{lu}}Fe_1e_1J$

$Je_{\mathbf{r}}e_2Ye_3e_{2}Z=Je_1e_1Fe_{\mathbf{lu}}e_{3}Z$
\smallskip

\medskip

{\bf Восстановление кода.} Зная коды $J$ и $Z$ можно вычислить код как $Y$, так и $F$: у $Y$ базовое окружение и первый начальник как у $Z$, а $F$ является третьим начальником $Z$ и имеет общий набор начальников с $J$.

\medskip

\medskip

{\bf 2. Среднее расположение}  (рисунок~\ref{pathplace3-1}, правая часть).

\medskip

{\bf Назначение $F$}. Пусть $F$ --  вершина типа $\mathbb{B}$ c произвольным базовым окружением, флагом подклейки и информацией.

Аналогично предыдущим случаям, мы можем вычислить коды всех вершин ($J$ и $Z$ -- с точностью до произвольного подклееного окружения, при заданном базовом).

\medskip

Для каждой разрешенной комбинации заданных букв, введем следующие определяющие соотношения:

\smallskip
$Ze_{2}e_3Ye_2e_{\mathbf{r}}J=Ze_{3}e_{\mathbf{mid}}Fe_1e_2J$

$Je_{\mathbf{r}}e_2Ye_3e_{2}Z=Je_2e_1Fe_{\mathbf{mid}}e_{3}Z$
\smallskip

{\bf Восстановление кода.} Аналогично левому-верхнему расположению.

\medskip

\begin{figure}[hbtp]
\centering
\includegraphics[width=0.9\textwidth]{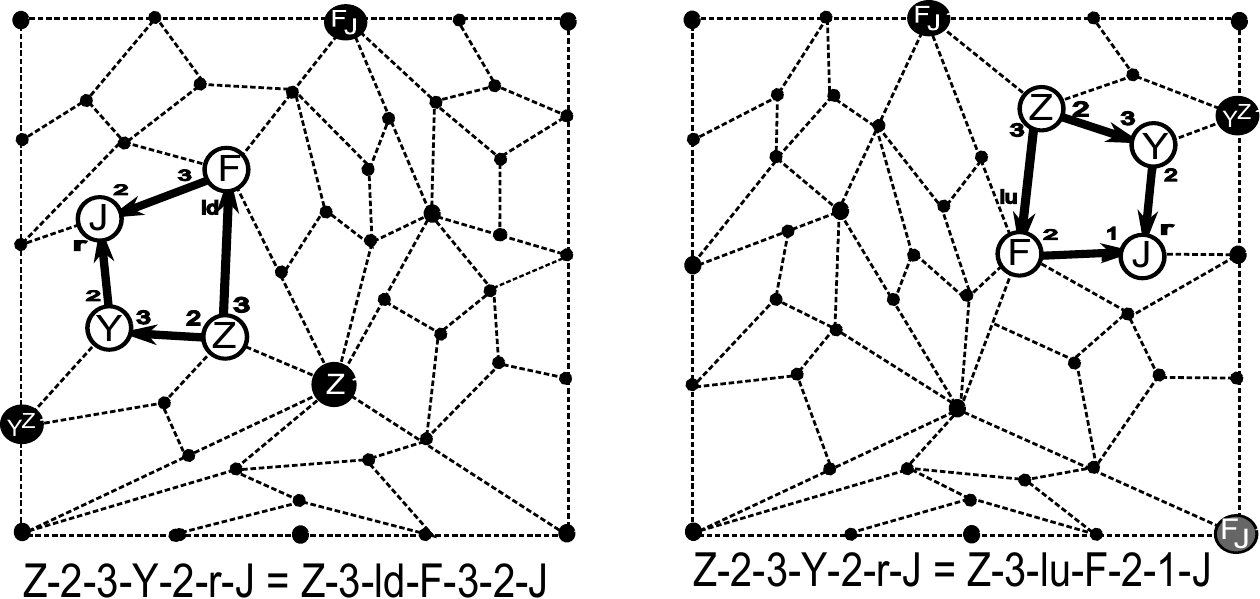}
\caption{Левое-нижнее и правое-верхнее положения. При входе и выходе в вершины обозначены коды входящих и выходящих ребер}
\label{pathplace3-2}
\end{figure}

{\bf 3. Левое нижнее расположение}  (рисунок~\ref{pathplace3-2}, левая часть).

\medskip

{\bf Назначение $F$}. Пусть $F$ -- вершина типа $\mathbb{A}$ c произвольным базовым окружением, флагом подклейки и информацией.

Цепь вокруг левого нижнего угла $T$ можно получить как $\mathbf{BottomLeftChain.FBoss}(F)$. Учитывая это, можно вычислить коды всех вершин ($J$ и $Z$ -- с точностью до произвольного подклееного окружения, при заданном базовом).

\medskip

Для каждой разрешенной комбинации заданных букв, введем следующие определяющие соотношения:

$Ze_{2}e_3Ye_2e_{\mathbf{r}}J=Ze_{3}e_{\mathbf{ld}}Fe_3e_2J$

$Je_{\mathbf{r}}e_2Ye_3e_{2}Z=Je_2e_3Fe_{\mathbf{ld}}e_{3}Z$

\medskip

{\bf Восстановление кода.} Аналогично левому-верхнему расположению.

\medskip

\medskip

{\bf 4. Правое-верхнее расположение} (рисунок~\ref{pathplace3-2}, правая часть).

\medskip

{\bf Назначение $F$}. Пусть $F$ -- вершина типа $\mathbb{B}$ c произвольным базовым окружением, флагом подклейки и информацией.

Цепь вокруг правого верхнего угла можно получить как  $\mathbf{UpRightChain.FBoss}(F)$.
Учитывая это, можно вычислить коды всех вершин ($J$ и $Z$ -- с точностью до произвольного подклееного окружения, при заданном базовом).

\medskip

Для каждой разрешенной комбинации заданных букв, введем следующие определяющие соотношения:

\smallskip
$Ze_{2}e_3Ye_2e_{\mathbf{r}}J=Ze_{3}e_{\mathbf{lu}}Fe_2e_1J$

$Je_{\mathbf{r}}e_2Ye_3e_{2}Z=Je_1e_2Fe_{\mathbf{lu}}e_{3}Z$

\medskip

{\bf Восстановление кода.} Аналогично левому-верхнему расположению.
\medskip

\begin{figure}[hbtp]
\centering
\includegraphics[width=0.9\textwidth]{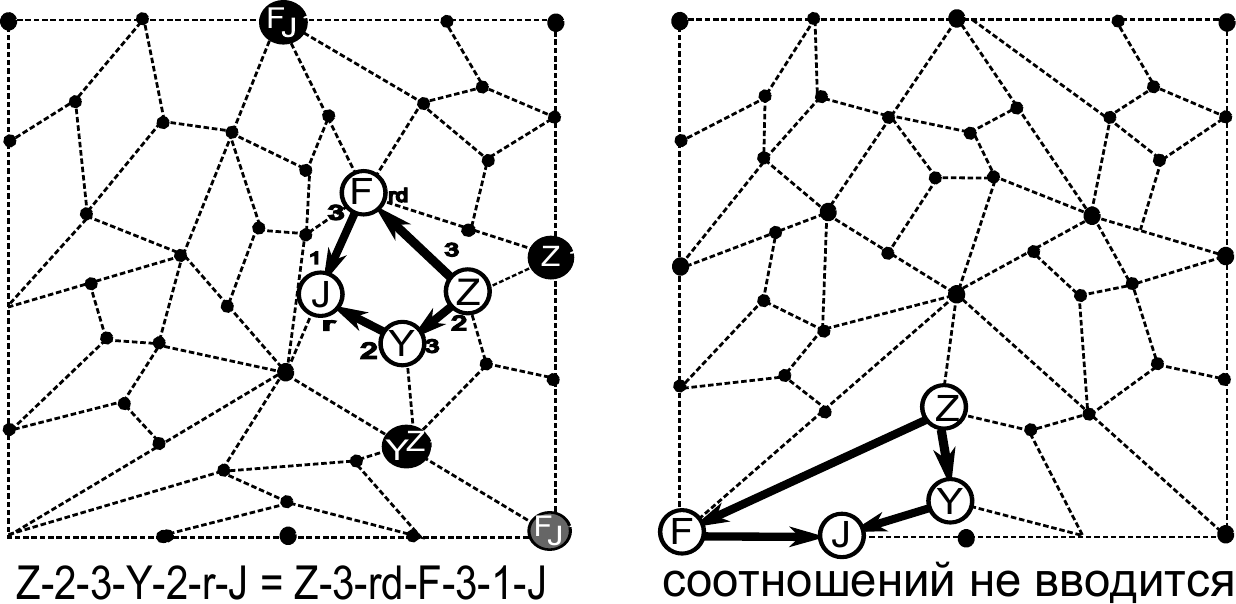}
\caption{Правое-нижнее и нижнее положения. При входе и выходе в вершины обозначены коды входящих и выходящих ребер}
\label{pathplace3-3}
\end{figure}

{\bf 5. Правое-нижнее расположение} (рисунок~\ref{pathplace3-3}, левая часть).

\medskip

{\bf Назначение $F$}. Пусть $F$ -- вершина типа $\mathbb{A}$ c произвольным базовым окружением, флагом подклейки и информацией.

Окружение вершины в середине правой стороны можно вычислить с помощью функции $\mathbf{RightFromB}$. Учитывая это, мы можем вычислить коды всех вершин ($J$ и $Z$ -- с точностью до произвольного подклееного окружения, при заданном базовом).

\medskip

Для каждой разрешенной комбинации заданных букв, введем следующие определяющие соотношения:

\smallskip
$Ze_{2}e_3Ye_2e_{\mathbf{r}}J=Ze_{3}e_{\mathbf{rd}}Fe_3e_1J$

$Je_{\mathbf{r}}e_2Ye_3e_{2}Z=Je_1e_3Fe_{\mathbf{rd}}e_{3}Z$

\medskip

{\bf Восстановление кода.} Аналогично левому-верхнему расположению.

\medskip

\medskip

{\bf 6. Нижнее расположение} (рисунок~\ref{pathplace3-3}, правая часть).

В этом случае соотношений мы не вводим, так как расположение пути удовлетворяет признакам мертвого паттерна. То есть участок пути с данным кодом не может являться подпутем достаточно длинного ненулевого пути.

\medskip

\subsection{Локальное преобразование 4}

\begin{figure}[hbtp]
\centering
\includegraphics[width=0.7\textwidth]{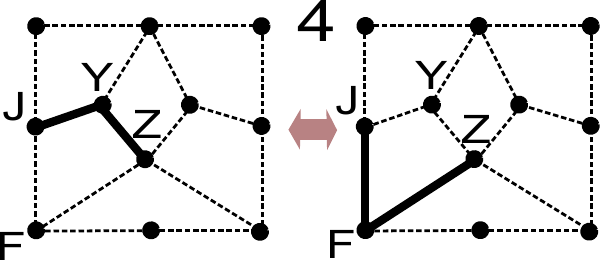}
\caption{Локальное преобразование $4$}
\label{flip4}
\end{figure}

Рассмотрим пару путей на рисунке~\ref{flip4}.

\medskip

\begin{figure}[hbtp]
\centering
\includegraphics[width=0.9\textwidth]{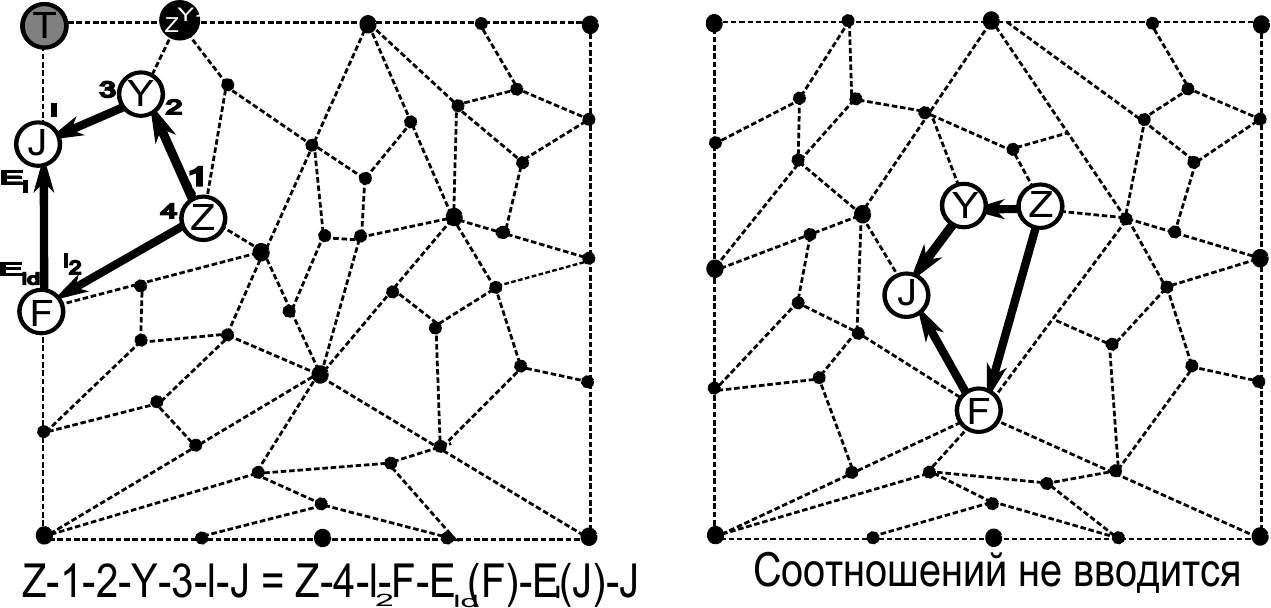}
\caption{Левое верхнее и среднее положения. При входе и выходе в вершины обозначены коды входящих и выходящих ребер}
\label{pathplace4-1}
\end{figure}

{\bf 1. Левое-верхнее расположение} (рисунок~\ref{pathplace4-1}, левая часть).

В этом случае сначала определим тип узла $T$ в левом верхнем углу. При этом $J$ и $F$ будут лежать на ребре выхода из $T$.

\medskip

{\bf Назначение $F$}. Пусть $F$ -- буква, кодирующая вершину типа $\mathbb{L}$, $\mathbb{LD}$ или $\mathbb{LU}$, лежащую в цепи вокруг $T$.

Заметим, учитывая расположение на рисунке, по известному коду $F$, мы можем вычислить окружения и начальников остальных трех вершин. То есть мы можем выписать все комбинации четырех букв $Z$, $Y$, $J$, $F$ такие, что их окружения и начальники соответствуют окружению и начальникам изображенных вершин. Будем называть такие комбинации {\it разрешенными}. При этом, буквы $Y$ и $F$ будут определены однозначно, а $Z$ и $J$ могут быть любыми буквами с заданным базовым окружением и начальниками, но с различными подклееными окружениями.

\medskip

Кроме того, буквы $E_\mathbf{l}$ и $E_{\mathbf{ld}}$ (кодирующие соответствующие на рисунке ребра входа и выхода) легко вычисляются по известным кодам $J$ и $F$.

\medskip

Для каждой разрешенной комбинации заданных букв, введем следующие определяющие соотношения:

\smallskip
$Ze_{1}e_2Ye_3e_{\mathbf{l}}J=Ze_{4}e_{\mathbf{l}_2}FE_{\mathbf{ld}}E_{\mathbf{l}}J$

$Je_{\mathbf{l}}e_{3}Ye_2e_{1}Z=JE_{\mathbf{l}}E_{\mathbf{ld}}Fe_{\mathbf{l}_2}e_{4}Z$

\medskip

{\bf Восстановление кода.} Зная коды $J$ и $Z$ можно вычислить код как $Y$, так и $F$: для $Y$ это очевидно, а $F$ является вторым начальником $Z$ и имеет общий набор начальников с $J$.

\medskip

\medskip

{\bf 2. Среднее расположение}, рисунок~\ref{pathplace4-1}, правая часть.

В этом случае соотношений мы не вводим, так как расположение пути удовлетворяет признакам мертвого паттерна. То есть участок пути с данным кодом не может являться подпутем достаточно длинного ненулевого пути.

\medskip

\begin{figure}[hbtp]
\centering
\includegraphics[width=0.9\textwidth]{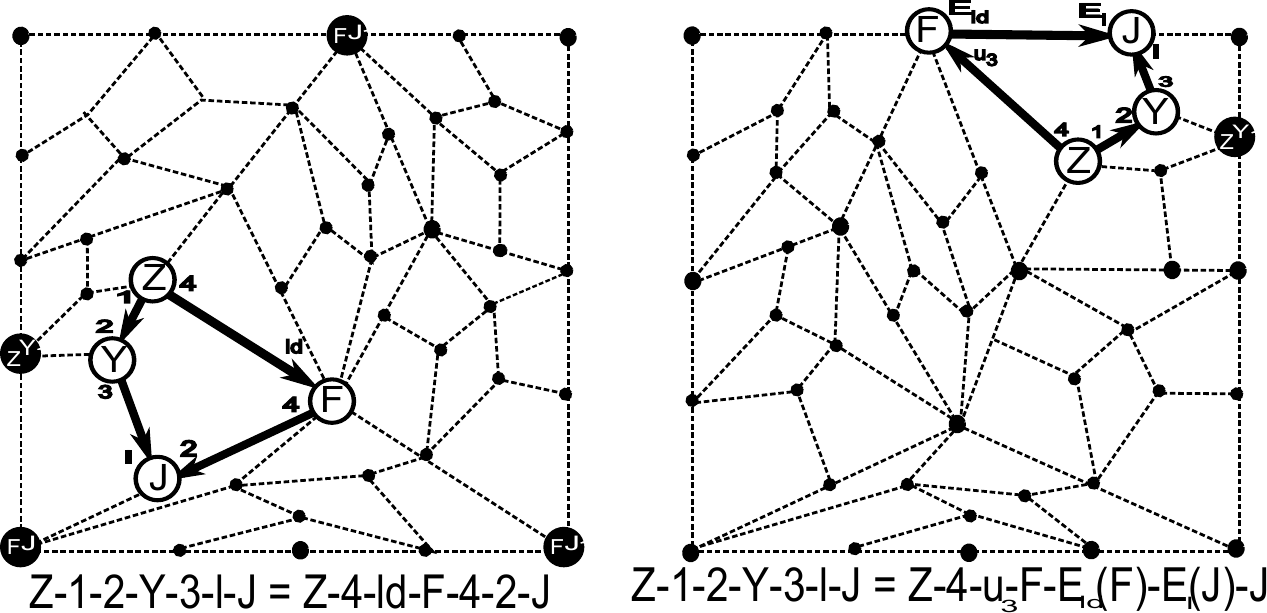}
\caption{Левое-нижнее и правое-верхнее положения. При входе и выходе в вершины обозначены коды входящих и выходящих ребер}
\label{pathplace4-2}
\end{figure}

{\bf 3. Левое нижнее расположение} (рисунок~\ref{pathplace4-2}, левая часть).

\medskip

{\bf Назначение $F$}. Пусть $F$ -- буква, кодирующая вершину типа $\mathbb{C}$, с произвольным окружением и начальниками.

Тогда базовые окружения вершин $Z$, $Y$, $J$ легко вычисляются. Начальники всех вершин отмечены черными кругами, и очевидно все они вычисляются, если мы знаем начальников $F$. Таким образом, мы можем аналогично определить все разрешенные четверки букв $Z$, $Y$, $J$, $F$, соответствующих вычисленным значениям окружений и начальников. У вершин $Z$ и $J$ могут быть разные подклееные окружения.

\medskip

Для каждой разрешенной комбинации заданных букв, введем следующие определяющие соотношения:

\smallskip
$Ze_{1}e_2Ye_3e_{\mathbf{l}}J=Ze_{4}e_{\mathbf{ld}}Fe_4e_2J$

$Je_{\mathbf{l}}e_3Ye_2e_{1}Z=Je_2e_4Fe_{\mathbf{ld}}e_{4}Z$

\medskip

{\bf Восстановление кода.} Аналогично левому-верхнему расположению.

\medskip

\medskip

{\bf 4. Правое верхнее расположение} (рисунок~\ref{pathplace4-2}, правая часть).

\medskip

{\bf Назначение $F$}. Пусть $F$ -- буква, кодирующая вершину типа $\mathbb{U}$, $\mathbb{UL}$ или $\mathbb{UR}$.

Тогда базовые окружения вершин $Z$, $Y$, $J$ легко вычисляются, учитывая, что мы можем применить функцию $\mathbf{UpRightChain}(F)$ и узнать тип вершины в правом верхнем углу.

Ясно также что начальники всех вершин также вычисляются, например у $J$ они такие же, как у $F$.  Таким образом, мы можем аналогично определить все разрешенные четверки букв $Z$, $Y$, $J$, $F$, соответствующих вычисленным значениям окружений и начальников. У вершин $Z$ и $J$ могут быть разные подклееные окружения.

\medskip

Кроме того, буквы $E_{\mathbf{l}}$ и $E_{\mathbf{ld}}$ (кодирующие соответствующие на рисунке ребра входа и выхода) легко вычисляются по известным кодам $J$ и $F$.

\medskip

Для каждой разрешенной комбинации заданных букв, введем следующие определяющие соотношения:

\smallskip
$Ze_{1}e_2Ye_3e_{\mathbf{l}}J=Ze_{4}e_{\mathbf{u}_3}FE_{\mathbf{ld}}E_{\mathbf{l}}J$

$Je_{\mathbf{l}}e_3Ye_2e_{1}Z=JE_1E_{\mathbf{ld}}Fe_{\mathbf{u}_3}e_{4}Z$

\medskip

{\bf Восстановление кода.} Аналогично левому-верхнему расположению.

\medskip

\begin{figure}[hbtp]
\centering
\includegraphics[width=0.9\textwidth]{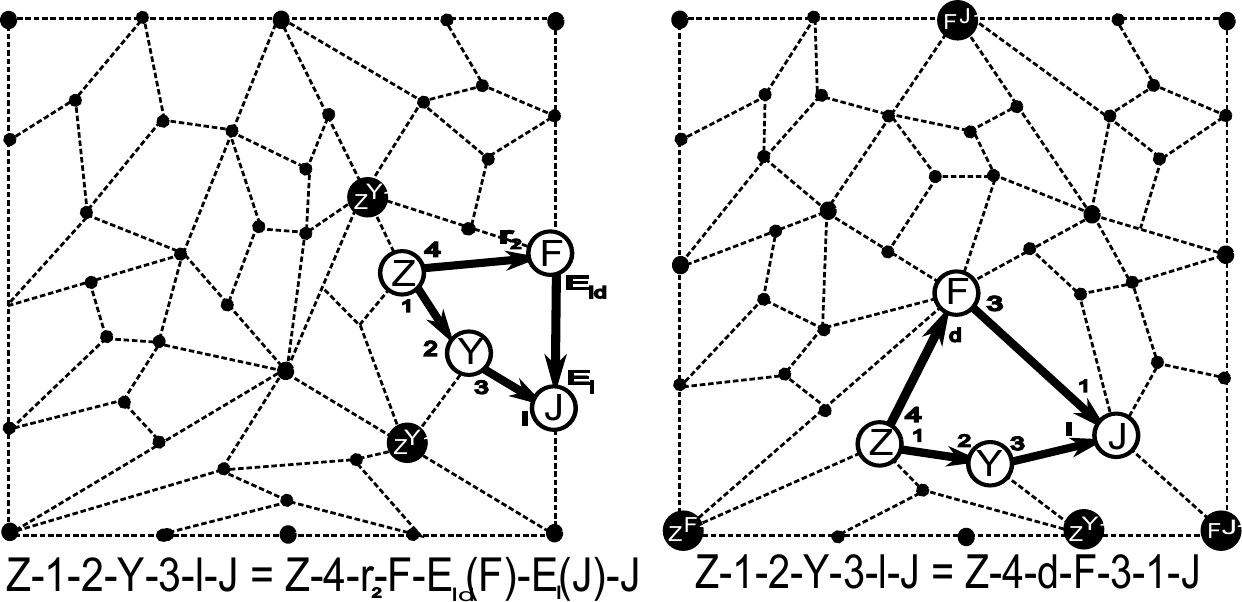}
\caption{Правое-нижнее и нижнее положения. При входе и выходе в вершины обозначены коды входящих и выходящих ребер}
\label{pathplace4-3}
\end{figure}

{\bf 5. Правое-нижнее расположение} (рисунок~\ref{pathplace4-3}, левая часть).

\medskip

{\bf Назначение $F$}. Пусть $F$ -- буква, кодирующая вершину типа $\mathbb{R}$, $\mathbb{RD}$ или $\mathbb{RU}$  c произвольным окружением и начальниками.

Тогда базовые окружения вершин $Z$, $Y$, $J$ легко вычисляются. Начальники всех вершин отмечены черными кругами, и очевидно все они вычисляются, если мы знаем начальников $F$. Таким образом, мы можем аналогично определить все разрешенные четверки букв $Z$, $Y$, $J$, $F$, соответствующих вычисленным значениям окружений и начальников. У вершин $Z$ и $J$ могут быть разные подклееные окружения.

\medskip
Кроме того, буквы $E_{\mathbf{l}}$ и $E_{\mathbf{ld}}$ (кодирующие соответствующие на рисунке ребра входа и выхода) легко вычисляются по известным кодам $J$ и $F$.

\medskip

Для каждой разрешенной комбинации заданных букв, введем следующие определяющие соотношения:

\smallskip

$Ze_{1}e_2Ye_3e_{\mathbf{l}}J=Ze_{4}e_{\mathbf{r}_2}FE_{\mathbf{ld}}E_{\mathbf{l}}J$

$Je_{\mathbf{l}}e_3Ye_2e_{1}Z=JE_{\mathbf{l}}E_{\mathbf{ld}}Fe_{\mathbf{r}_2}e_{4}Z$

\medskip

{\bf Восстановление кода.} Аналогично левому-верхнему расположению.

\medskip

\medskip

{\bf 6. Нижнее расположение} (рисунок~\ref{pathplace4-3}, правая часть).

\medskip

{\bf Назначение $F$}. Пусть $F$ -- буква, кодирующая вершину типа $\mathbb{C}$, c произвольным окружением и начальниками.

Тогда базовые окружения вершин $Z$, $Y$, $J$ легко вычисляются. Начальники всех вершин отмечены черными кругами, и очевидно все они вычисляются, если мы знаем начальников $F$. Таким образом, мы можем аналогично определить все разрешенные четверки букв $Z$, $Y$, $J$, $F$, соответствующих вычисленным значениям окружений и начальников. У вершин $Z$ и $J$ могут быть разные подклееные окружения.

\medskip

Для каждой разрешенной комбинации заданных букв, введем следующие определяющие соотношения:

\smallskip
$Ze_{1}e_2Ye_3e_{\mathbf{l}}J=Ze_{4}e_{\mathbf{d}}Fe_3e_1J$

$Je_{\mathbf{l}}e_3Ye_2e_{1}Z=Je_1e_3Fe_{\mathbf{d}}e_{4}Z$

\medskip

{\bf Восстановление кода.} Аналогично левому-верхнему расположению.

\medskip

\subsection{Локальное преобразование 5}

\begin{figure}[hbtp]
\centering
\includegraphics[width=0.7\textwidth]{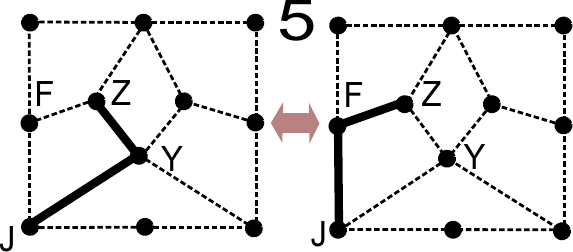}
\caption{Локальное преобразование $5$}
\label{flip5}
\end{figure}

Рассмотрим пару путей на рисунке~\ref{flip5}.

\medskip

\medskip

\begin{figure}[hbtp]
\centering
\includegraphics[width=0.9\textwidth]{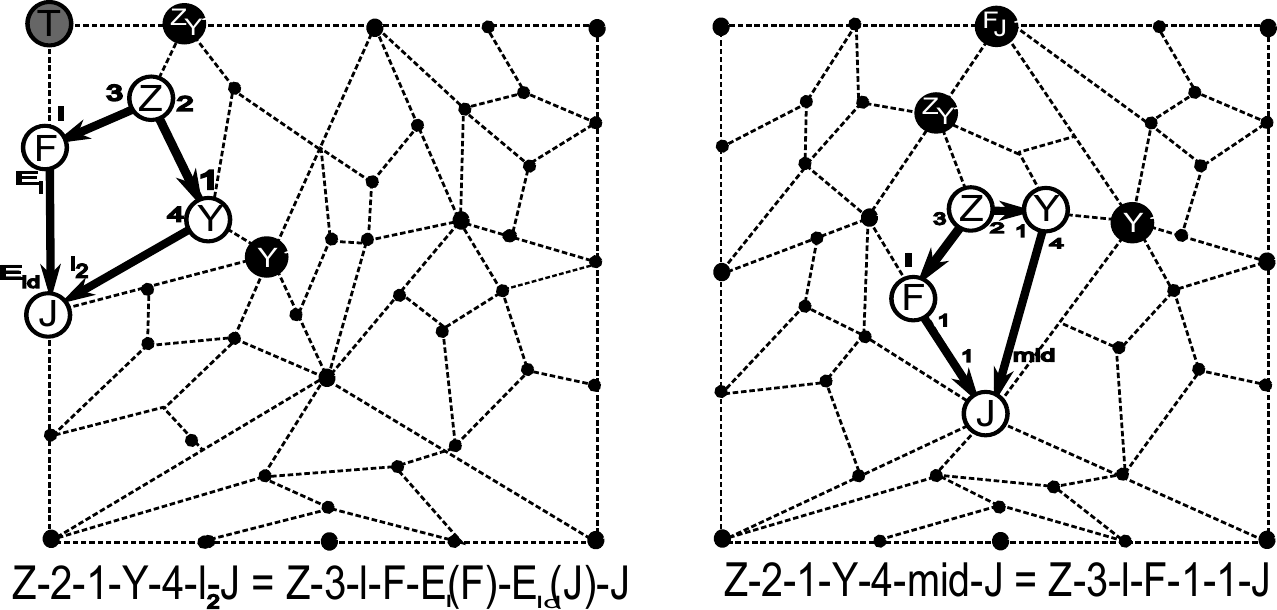}
\caption{Левое верхнее и среднее положения. При входе и выходе в вершины обозначены коды входящих и выходящих ребер}
\label{pathplace5-1}
\end{figure}

{\bf 1. Левое-верхнее расположение} (рисунок~\ref{pathplace5-1}, левая часть).

В этом случае сначала определим тип узла $T$ в левом верхнем углу. При этом $J$ и $F$ будут лежать на ребре выхода из $T$.

\medskip

{\bf Назначение $J$}. Пусть $J$ -- буква, кодирующая вершину типа $\mathbb{L}$, $\mathbb{LD}$ или $\mathbb{LU}$, лежащую в цепи вокруг $T$.

Аналогично предыдущему случаю, по известному коду $J$, мы можем вычислить окружения и начальников остальных трех вершин. То есть мы можем выписать все {\it разрешенные} комбинации четырех букв $Z$, $Y$, $J$, $F$ то есть такие, что их окружения и начальники соответствуют окружению и начальникам изображенных вершин. При этом, буквы $Y$ и $F$ будут определены однозначно, а $Z$ и $J$ могут быть любыми буквами с заданным базовым окружением и начальниками, но с различными подклееными окружениями.

\medskip

Кроме того, буквы $E_{\mathbf{l}}$ и $E_{\mathbf{ld}}$ (кодирующие соответствующие на рисунке ребра входа и выхода) легко вычисляются по известным кодам $J$ и $F$.

\medskip

Для каждой разрешенной комбинации заданных букв, введем следующие определяющие соотношения:

\smallskip
$Ze_{2}e_1Ye_4e_{\mathbf{l}_2}J=Ze_{3}e_{\mathbf{l}}FE_{\mathbf{l}}E_{\mathbf{ld}}J$

$Je_{\mathbf{l}_2}e_{4}Ye_1e_{2}Z=JE_{\mathbf{ld}}E_{\mathbf{l}}Fe_{\mathbf{l}}e_{3}Z$

\medskip

{\bf Восстановление кода.} Зная коды $J$ и $Z$ можно вычислить код $Y$ и $F$: для $Y$ это очевидно, а $F$ соответствует $\mathbf{Next.FBoss}(Z)$ и имеет общий набор начальников с $J$.

\medskip

\medskip

{\bf 2. Среднее расположение} (рисунок~\ref{pathplace5-1}, правая часть).

\medskip

{\bf Назначение $J$}. Пусть $J$ -- буква, кодирующая вершину типа $\mathbb{C}$, с произвольным окружением и начальниками.

Тогда базовые окружения вершин $Z$, $Y$, $F$ легко вычисляются. Начальники всех вершин отмечены черными кругами, и очевидно все они вычисляются, если мы знаем начальников $J$. Таким образом, мы можем аналогично определить все разрешенные четверки букв $Z$, $Y$, $J$, $F$, соответствующих вычисленным значениям окружений и начальников. У вершин $Z$ и $J$ могут быть разные подклееные окружения.

\medskip

Для каждой разрешенной комбинации заданных букв, введем следующие определяющие соотношения:

\smallskip
$Ze_{2}e_1Ye_4e_{\mathbf{mid}}J=Ze_{3}e_{\mathbf{l}}Fe_1e_1J$

$Je_{\mathbf{mid}}e_4Ye_1e_{2}Z=Je_1e_1Fe_{\mathbf{l}}e_{3}Z$

\medskip

{\bf Восстановление кода.} Аналогично левому-верхнему расположению.

\medskip

\begin{figure}[hbtp]
\centering
\includegraphics[width=0.9\textwidth]{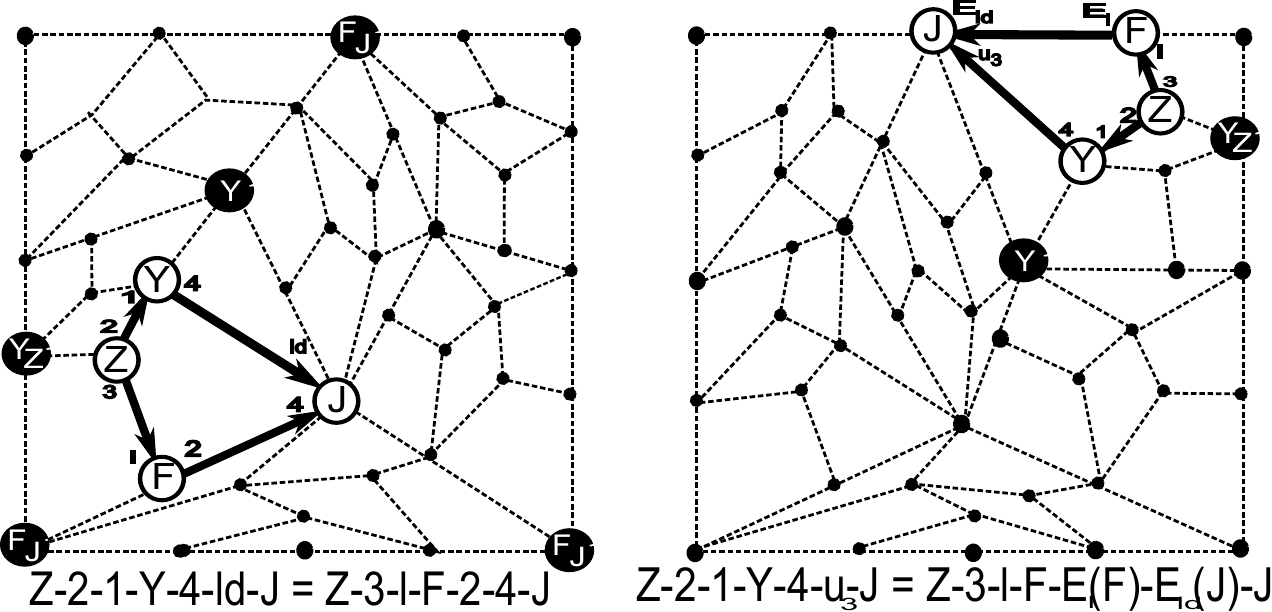}
\caption{Левое-нижнее и правое-верхнее положения. При входе и выходе в вершины обозначены коды входящих и выходящих ребер}
\label{pathplace5-2}
\end{figure}

{\bf 3. Левое нижнее расположение} (рисунок~\ref{pathplace5-2}, левая часть).

\medskip

{\bf Назначение $J$}. Пусть $J$ -- буква, кодирующая вершину типа $\mathbb{C}$, с произвольным окружением и начальниками.

Тогда базовые окружения вершин $Z$, $Y$, $F$ легко вычисляются. Начальники всех вершин отмечены черными кругами, и очевидно все они вычисляются, если мы знаем начальников $J$. Таким образом, мы можем аналогично определить все разрешенные четверки букв $Z$, $Y$, $J$, $F$, соответствующих вычисленным значениям окружений и начальников. У вершин $Z$ и $J$ могут быть разные подклееные окружения.

\medskip

Для каждой разрешенной комбинации заданных букв, введем следующие определяющие соотношения:

\smallskip
$Ze_{2}e_1Ye_4e_{\mathbf{ld}}J=Ze_{3}e_{\mathbf{l}}Fe_2e_4J$

$Je_{\mathbf{ld}}e_4Ye_1e_{2}Z=Je_4e_2Fe_{\mathbf{l}}e_{3}Z$

\medskip

{\bf Восстановление кода.} Аналогично левому-верхнему расположению.

\medskip

\medskip

{\bf 4. Правое верхнее расположение} (рисунок~\ref{pathplace5-2}, правая часть).

\medskip

{\bf Назначение $J$}. Пусть $J$ -- буква, кодирующая вершину типа $\mathbb{U}$, $\mathbb{UL}$ или $\mathbb{UR}$.

Тогда базовые окружения вершин $Z$, $Y$, $J$ легко вычисляются, учитывая, что мы можем применить функцию $\mathbf{UpRightChain}(F)$ и узнать тип вершины в правом верхнем углу.

Начальники всех вершин также вычисляются, например у $F$ они такие же, как у $J$.  Таким образом, мы можем аналогично определить все разрешенные четверки букв $Z$, $Y$, $J$, $F$, соответствующих вычисленным значениям окружений и начальников. У вершин $Z$ и $J$ могут быть разные подклееные окружения.

\medskip

Кроме того, буквы $E_{\mathbf{l}}$ и $E_{\mathbf{ld}}$ (кодирующие соответствующие на рисунке ребра входа и выхода) легко вычисляются по известным кодам $J$ и $F$.

\medskip

Для каждой разрешенной комбинации заданных букв, введем следующие определяющие соотношения:

\smallskip
$Ze_{2}e_1Ye_4e_{\mathbf{u}_3}J=Ze_{3}e_{\mathbf{l}}FE_{\mathbf{l}}E_{\mathbf{ld}}J$

$Je_{\mathbf{u}_3}e_4Ye_1e_{2}Z=JE_{\mathbf{ld}}E_{\mathbf{l}}Fe_{\mathbf{l}}e_{3}Z$

\medskip

{\bf Восстановление кода.} Аналогично левому-верхнему расположению.

\medskip

\begin{figure}[hbtp]
\centering
\includegraphics[width=0.9\textwidth]{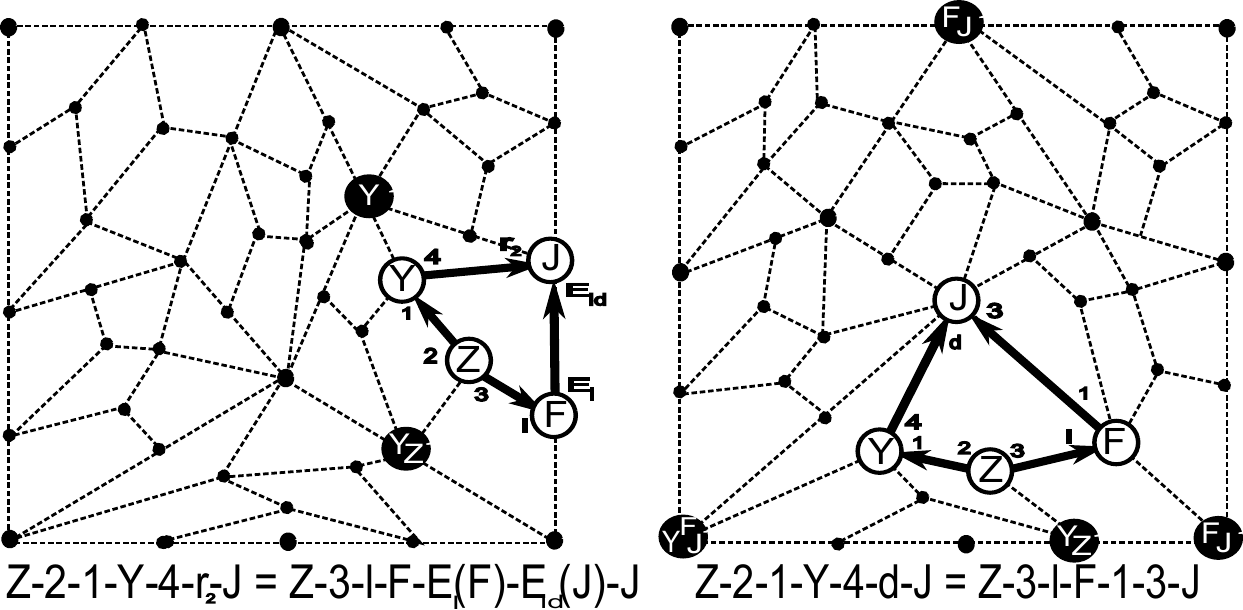}
\caption{Правое-нижнее и нижнее положения. При входе и выходе в вершины обозначены коды входящих и выходящих ребер}
\label{pathplace5-3}
\end{figure}

{\bf 5. Правое-нижнее расположение} (рисунок~\ref{pathplace5-3}, левая часть).

\medskip

{\bf Назначение $J$}. Пусть $J$ -- буква, кодирующая вершину типа $\mathbb{R}$, $\mathbb{RD}$ или $\mathbb{RU}$  c произвольным окружением и начальниками.

Тогда базовые окружения вершин $Z$, $Y$, $F$ (рисунок~\ref{pathplace5-3}, левая часть), очевидны. Начальники всех вершин отмечены черными кругами, и очевидно все они вычисляются, если мы знаем начальников $J$. Таким образом, мы можем аналогично определить все разрешенные четверки букв $Z$, $Y$, $J$, $F$, соответствующих вычисленным значениям окружений и начальников. У вершин $Z$ и $J$ могут быть разные подклееные окружения.

\medskip
Кроме того, буквы $E_{\mathbf{l}}$ и $E_{\mathbf{ld}}$ (кодирующие соответствующие на рисунке ребра входа и выхода) легко вычисляются по известным кодам $J$ и $F$.

\medskip

Для каждой разрешенной комбинации заданных букв, введем следующие определяющие соотношения:

\smallskip
$Ze_{2}e_1Ye_4e_{\mathbf{r}_2}J=Ze_{3}e_{\mathbf{l}}FE_{\mathbf{l}}E_{\mathbf{ld}}J$

$Je_{\mathbf{r}_2}e_4Ye_1e_{2}Z=JE_{\mathbf{ld}}E_{\mathbf{l}}Fe_{\mathbf{l}}e_{3}Z$

\medskip

{\bf Восстановление кода.} Аналогично левому-верхнему расположению.

\medskip

\medskip

{\bf 6. Нижнее расположение} (рисунок~\ref{pathplace5-3}, правая часть).

\medskip

{\bf Назначение $J$}. Пусть $J$ -- буква, кодирующая вершину типа $\mathbb{C}$, c произвольным окружением и начальниками.

Тогда базовые окружения вершин $Z$, $Y$, $F$ на вышеуказанном рисунке очевидны. Начальники всех вершин отмечены черными кругами, и очевидно все они вычисляются, если мы знаем начальников $J$. Таким образом, мы можем аналогично определить все разрешенные четверки букв $Z$, $Y$, $J$, $F$, соответствующих вычисленным значениям окружений и начальников. У вершин $Z$ и $J$ могут быть разные подклееные окружения.

\medskip

Для каждой разрешенной комбинации заданных букв, введем следующие определяющие соотношения:

\smallskip
$Ze_{2}e_1Ye_4e_{\mathbf{d}}J=Ze_{3}e_{\mathbf{l}}Fe_1e_3J$

$Je_{\mathbf{d}}e_4Ye_1e_{2}Z=Je_3e_1Fe_{\mathbf{l}}e_{3}Z$

\medskip

{\bf Восстановление кода.} Аналогично левому-верхнему расположению.

\medskip

\subsection{Локальное преобразование 6}

\begin{figure}[hbtp]
\centering
\includegraphics[width=0.7\textwidth]{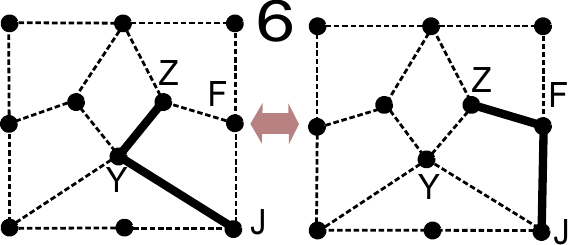}
\caption{Локальное преобразование $6$}
\label{flip6}
\end{figure}

Рассмотрим пару путей на рисунке~\ref{flip6}.

\medskip

\medskip

\begin{figure}[hbtp]
\centering
\includegraphics[width=0.9\textwidth]{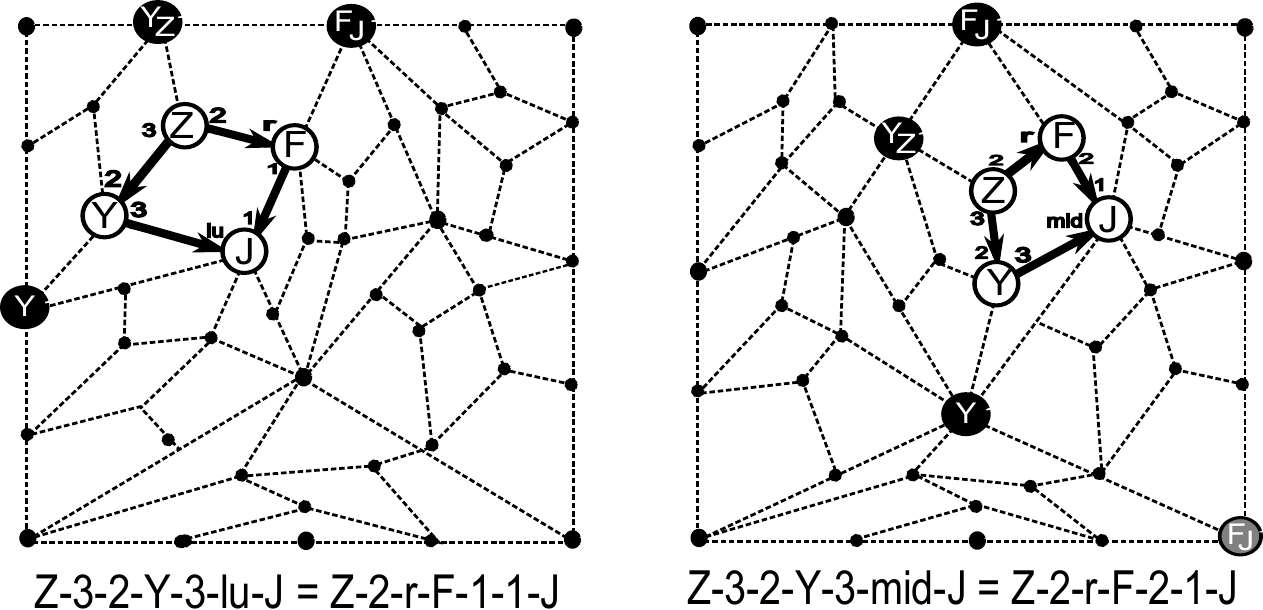}
\caption{Левое верхнее и среднее положения. При входе и выходе в вершины обозначены коды входящих и выходящих ребер}
\label{pathplace6-1}
\end{figure}

{\bf 1. Левое-верхнее расположение} (рисунок~\ref{pathplace6-1}, левая часть).

\medskip

{\bf Назначение $J$}. Пусть $J$ -- буква, кодирующая вершину типа $\mathbb{A}$, c произвольным окружением и произвольным начальником.

По известному коду $J$, мы можем выписать все {\it разрешенные} комбинации четырех букв $Z$, $Y$, $J$, $F$, то есть такие, что их окружения и начальники соответствуют окружению и начальникам изображенных на рисунке~\ref{pathplace6-1} вершин. При этом, буквы $Y$ и $F$ будут определены однозначно, а $Z$ и $J$ могут быть любыми буквами с заданным базовым окружением и начальниками, но с различными подклееными окружениями.

\medskip

Для каждой разрешенной комбинации заданных букв, введем следующие определяющие соотношения:

\smallskip
$Ze_{3}e_2Ye_3e_{\mathbf{lu}}J=Ze_{2}e_{\mathbf{r}}Fe_{1}e_{1}J$

$Je_{\mathbf{lu}}e_{3}Ye_2e_{3}Z=Je_{1}e_{1}Fe_{\mathbf{r}}e_{2}Z$

\medskip

{\bf Восстановление кода.} Зная коды $J$ и $Z$ можно вычислить код как $Y$, так и $F$. Окружение $Y$ и первый начальник такие же как у $Z$, второй начальник соответствует $\mathbf{Next.FBoss}(J)$, а третий -- вершина $J$. Окружение $F$ соответствует $0$-цепи вокруг $J$ с указателем $1$, начальник тот же что и у $J$.

\medskip

\medskip

{\bf 2. Среднее расположение} (рисунок~\ref{pathplace6-1}, правая часть).

\medskip

{\bf Назначение $J$}. Пусть $J$ -- буква, кодирующая вершину типа $\mathbb{B}$, с произвольным окружением и начальниками.

Тогда базовые окружения вершин $Z$, $Y$, $F$ очевидны. Начальники всех вершин отмечены черными кругами, и очевидно все они вычисляются, если мы знаем начальников $J$. Таким образом, мы можем аналогично определить все разрешенные четверки букв $Z$, $Y$, $J$, $F$, соответствующих вычисленным значениям окружений и начальников. У вершин $Z$ и $J$ могут быть разные подклееные окружения.

\medskip

Для каждой разрешенной комбинации заданных букв, введем следующие определяющие соотношения:

\smallskip

$Ze_{3}e_2Ye_3e_{\mathbf{mid}}J=Ze_{2}e_{\mathbf{r}}Fe_2e_1J$

$Je_{\mathbf{mid}}e_3Ye_2e_{3}Z=Je_1e_2Fe_{\mathbf{r}}e_{2}Z$

\medskip

{\bf Восстановление кода.} Зная коды $J$ и $Z$ можно вычислить код как $Y$, так и $F$. Окружение $Y$ и первый начальник такие же как у $Z$, второй начальник это вершина типа $\mathbb{C}$ с окружением как у $J$, а третий -- вершина $J$. Окружение $F$ вычисляется по известному окружению $J$, начальник тот же что и у $J$.

\medskip

\begin{figure}[hbtp]
\centering
\includegraphics[width=0.9\textwidth]{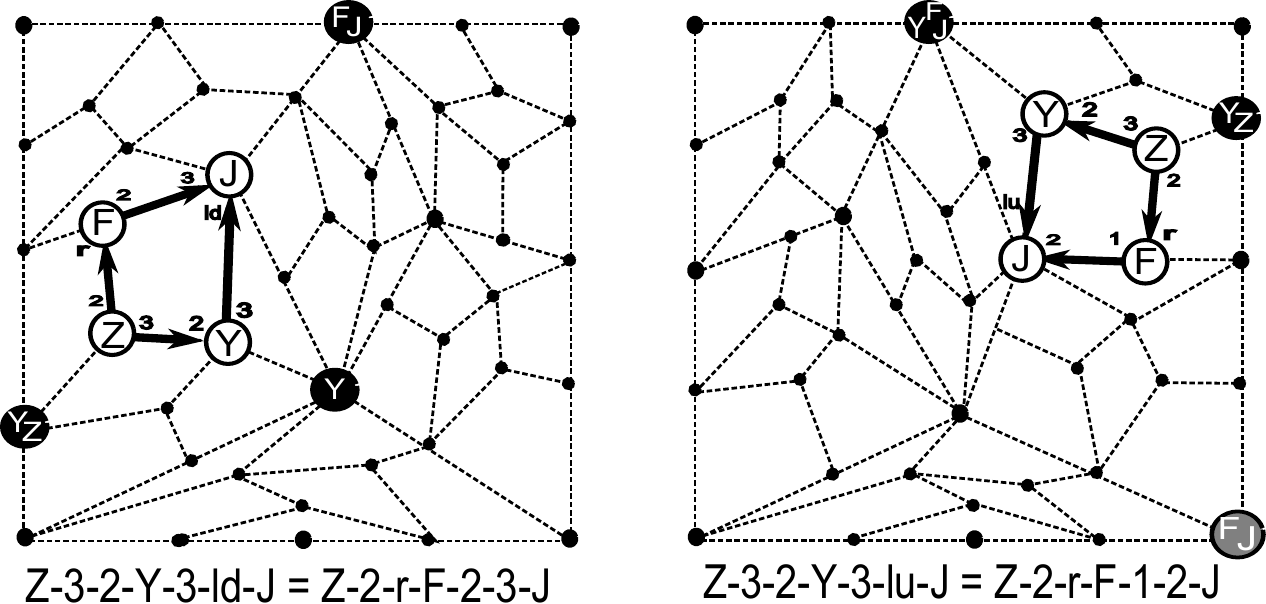}
\caption{Левое-нижнее и правое-верхнее положения. При входе и выходе в вершины обозначены коды входящих и выходящих ребер}
\label{pathplace6-2}
\end{figure}

{\bf 3. Левое нижнее расположение}  (рисунок~\ref{pathplace6-2}, левая часть).

\medskip

{\bf Назначение $J$}. Пусть $J$ -- буква, кодирующая вершину типа $\mathbb{A}$, с произвольным окружением и начальником.

Тогда базовые окружения вершин $Z$, $Y$, $F$ очевидны. Начальники всех вершин отмечены черными кругами, и очевидно все они вычисляются, если мы знаем код $J$. Таким образом, мы можем аналогично определить все разрешенные четверки букв $Z$, $Y$, $J$, $F$, соответствующих вычисленным значениям окружений и начальников. У вершин $Z$ и $J$ могут быть разные подклееные окружения.

\medskip

Для каждой разрешенной комбинации заданных букв, введем следующие определяющие соотношения:

\smallskip
$Ze_{3}e_2Ye_3e_{\mathbf{ld}}J=Ze_{2}e_{\mathbf{r}}Fe_2e_3J$

$Je_{\mathbf{ld}}e_3Ye_2e_{3}Z=Je_3e_2Fe_{\mathbf{r}}e_{2}Z$

\medskip

{\bf Восстановление кода.} Зная коды $J$ и $Z$ можно вычислить код как $Y$, так и $F$. Окружение $Y$ и первый начальник такие же как у $Z$, второй начальник соответствует узлу $\mathbb{C}$ c окружением как у $J$, а третий -- сама вершина $J$. Окружение $F$ вычисляется по известному окружению $J$, начальник тот же что и у $J$.

\medskip

\medskip

{\bf 4. Правое верхнее расположение} (рисунок~\ref{pathplace6-2}, правая часть).

\medskip

{\bf Назначение $J$}. Пусть $J$ -- буква, кодирующая вершину типа $\mathbb{B}$, c произвольным окружением и начальниками.

Тогда базовые окружения вершин $Z$, $Y$, $J$, очевидны.
Заметим, что мы можем применить функцию $\mathbf{UpRightChain}(F)$ и узнать тип вершины в правом верхнем углу. То есть, начальники всех вершин, также вычисляются. Таким образом, мы можем аналогично определить все разрешенные четверки букв $Z$, $Y$, $J$, $F$, соответствующих вычисленным значениям окружений и начальников. У вершин $Z$ и $J$ могут быть разные подклееные окружения.

\medskip

Для каждой разрешенной комбинации заданных букв, введем следующие определяющие соотношения:

\smallskip

$Ze_{3}e_2Ye_3e_{\mathbf{lu}}J=Ze_{2}e_{\mathbf{r}}Fe_1e_{2}J$

$Je_{\mathbf{lu}}e_3Ye_2e_{3}Z=Je_{2}e_{1}Fe_{\mathbf{r}}e_{2}Z$

\medskip

{\bf Восстановление кода.} Зная коды $J$ и $Z$ можно вычислить код как $Y$, так и $F$. Окружение $Y$ и первый начальник такие же как у $Z$, второй начальник соответствует $\mathbf{FBoss}(J)$, а третий -- вершина $J$. Окружение $F$ вычисляется по известному окружению $J$, начальники те же что и у $J$.

\medskip

\begin{figure}[hbtp]
\centering
\includegraphics[width=0.9\textwidth]{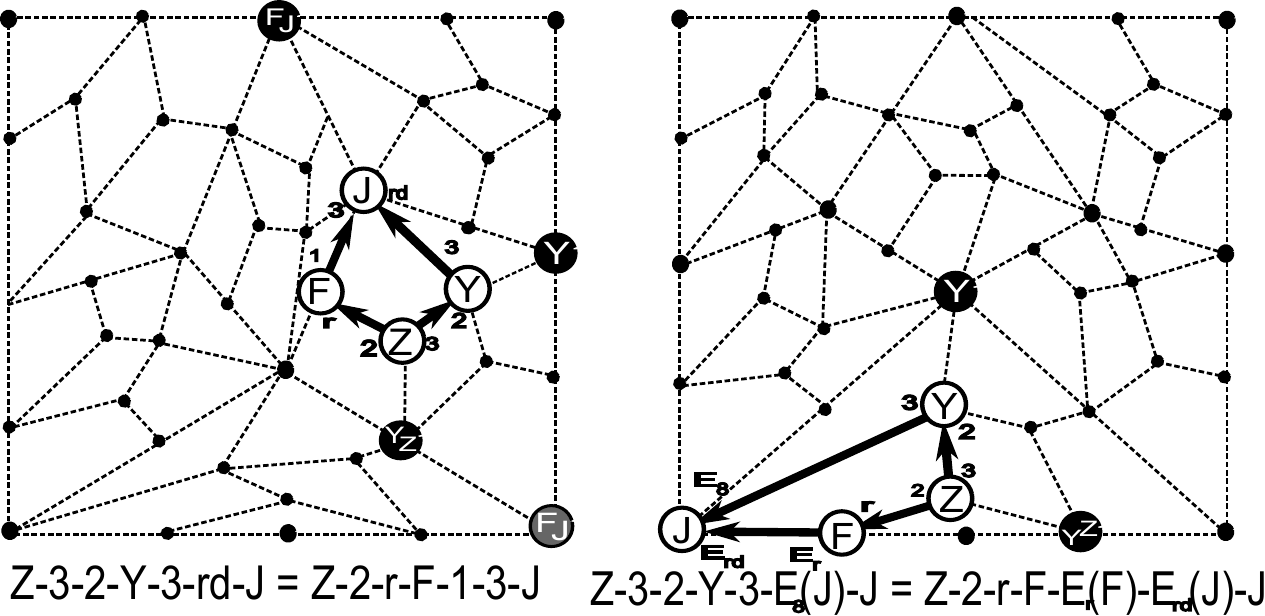}
\caption{Правое-нижнее и нижнее положения. При входе и выходе в вершины обозначены коды входящих и выходящих ребер}
\label{pathplace6-3}
\end{figure}

{\bf 5. Правое-нижнее расположение} (рисунок~\ref{pathplace6-3}, левая часть).

\medskip

{\bf Назначение $J$}. Пусть $J$ -- буква, кодирующая вершину типа $\mathbb{B}$, c произвольным окружением и начальниками.

Тогда базовые окружения вершин $Z$, $Y$, $F$ (рисунок~\ref{pathplace6-3}, левая часть), очевидны. Начальники всех вершин отмечены черными кругами, и очевидно все они вычисляются, если мы знаем начальников $J$. Таким образом, мы можем аналогично определить все разрешенные четверки букв $Z$, $Y$, $J$, $F$, соответствующих вычисленным значениям окружений и начальников. У вершин $Z$ и $J$ могут быть разные подклееные окружения.

\medskip

Для каждой разрешенной комбинации заданных букв, введем следующие определяющие соотношения:

\smallskip
$Ze_{3}e_2Ye_3e_{\mathbf{rd}}J=Ze_{2}e_{\mathbf{r}}Fe_{1}e_{3}J$

$Je_{\mathbf{rd}}e_3Ye_2e_{3}Z=Je_{3}e_{1}Fe_{\mathbf{r}}e_{2}Z$

\medskip

{\bf Восстановление кода.} Зная коды $J$ и $Z$ можно вычислить код как $Y$, так и $F$. Окружение $Y$ и первый начальник такие же как у $Z$, второй начальник вычисляется с помощью процедуры $\mathbf{RightFromB}(J)$, а третий -- вершина $J$. Окружение $F$ вычисляется по известному окружению $J$, начальники те же что и у $J$.

\medskip

\medskip

{\bf 6. Нижнее расположение} (рисунок~\ref{pathplace6-3}, правая часть).

\medskip

{\bf Назначение $J$}. Пусть $J$ -- буква, кодирующая вершину типа $\mathbb{C}$, c произвольным окружением и начальниками.

Тогда базовые окружения вершин $Z$, $Y$, $F$ на вышеуказанном рисунке очевидны. Начальники всех вершин отмечены черными кругами, и очевидно все они вычисляются, если мы знаем начальников $J$. Таким образом, мы можем аналогично определить все разрешенные четверки букв $Z$, $Y$, $J$, $F$, соответствующих вычисленным значениям окружений и начальников. У вершин $Z$ и $J$ могут быть разные подклееные окружения.

\medskip

Для каждой разрешенной комбинации заданных букв, введем следующие определяющие соотношения:

\smallskip
$Ze_{2}e_1Ye_4e_{\mathbf{d}}J=Ze_{3}e_{\mathbf{l}}Fe_1e_3J$

$Je_{\mathbf{d}}e_4Ye_1e_{2}Z=Je_3e_1Fe_{\mathbf{l}}e_{3}Z$

\medskip

{\bf Восстановление кода.} Зная коды $J$ и $Z$ можно вычислить код как $Y$, так и $F$. Окружение $Y$ и первый начальник такие же как у $Z$, второй начальник -- вершина $J$, третий соответствует $\mathbf{SBoss}(J)$. Окружение $F$ соответствует $1$-цепи вокруг третьего начальника $J$ с указателем, соответствующим входу в узел по нижней стороне, функция $E_{ddr}(\mathbf{TBoss}(J))$. Начальники те же что и у $J$.

\bigskip

\medskip

{\bf Локальные преобразования 7, 8, 9, 10.}

Рассмотрим узел, лежащий в середине верхней стороны макроплитки, где расположен наш путь. Этот узел должен входить в некоторую цепь. То есть может быть три варианта цепи вокруг центра, который может иметь следующий тип:

\smallskip

$\mathbb{C}$, $\mathbb{B}$, $\mathbb{A}$, $\mathbb{UL}/ \mathbb{UL}$,  $\mathbb{UR}/ \mathbb{RU}$, $\mathbb{DL}/ \mathbb{LD}$,  $\mathbb{DR}/ \mathbb{RD}$,  $\mathbb{U}$, $\mathbb{L}$, $\mathbb{D}$, $\mathbb{R}$, $\mathbb{CUL}$, $\mathbb{CUR}$, $\mathbb{CDL}$, $\mathbb{CDR}.$

\smallskip

 Мы разберем все эти случаи и введем соотношения в соответствии с устройством всех возможных цепей.

\medskip

{\bf Замечание.} Ниже мы не будем выписывать симметричное к введенному соотношение, отвечающее проходу пути в обратном порядке. Просто будем считать, что соотношений вводится в два раза больше.

\medskip

\subsection{Случай цепи $\mathbb{C}1$; преобразования 8 и 9}

\begin{figure}[hbtp]
\centering
\leftskip=-1.0cm
\includegraphics[width=1\textwidth]{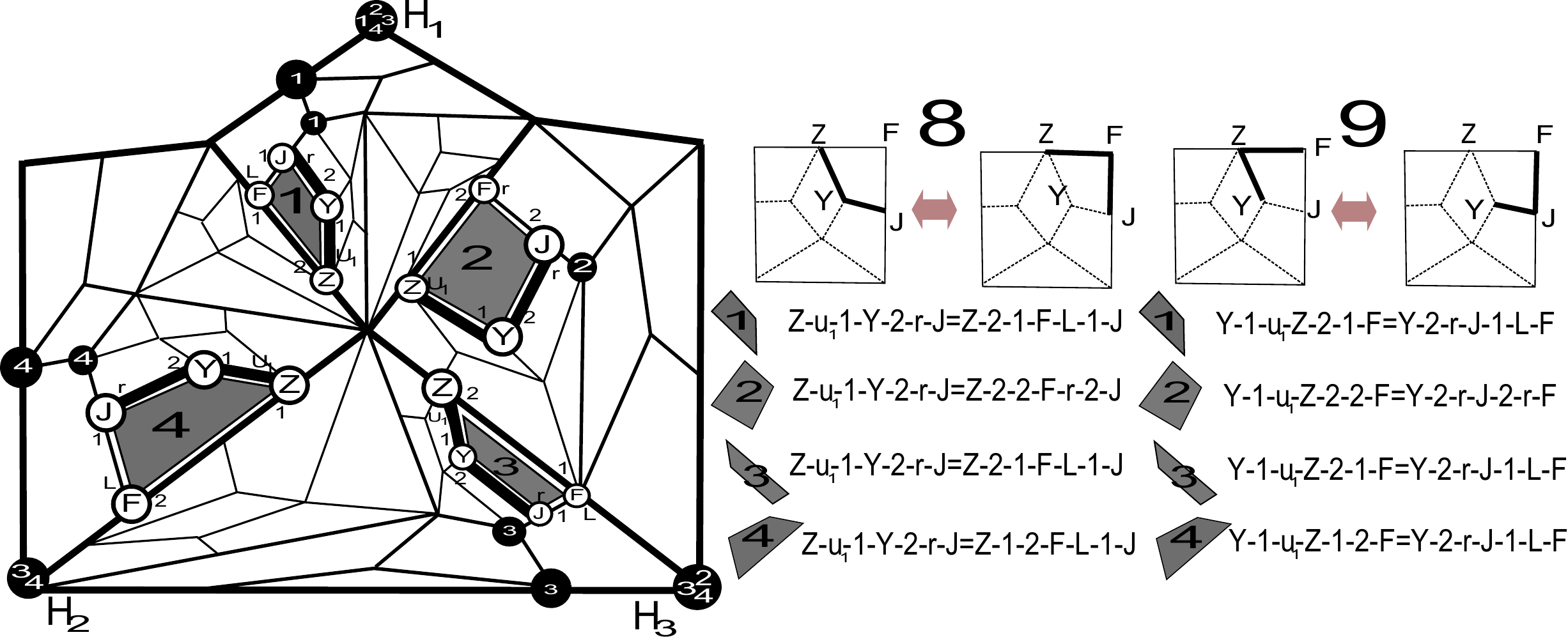}
\caption{Случаи расположения пути вокруг C1-цепи и соответствующие им определяющие соотношения для локальных преобразований 8 и 9}
\label{C1b}
\end{figure}

В правой части рисунка~\ref{C1b} изображены локальные преобразования $8$ и $9$.
Обозначим макроплитку, в которой проходят пути как $T$. Вершина $Z$ попадает в середину верхней стороны $T$.

\medskip

Итак, пусть середина верхней стороны $T$ входит в $\mathbb{C}1$-цепь. Различные случаи расположения интересующих нас путей показаны в левой части рисунка~\ref{C1b}.

\medskip

Сначала введем определяющие соотношения для данного случая.
Зафиксируем $\mathbb{C}$-вершину с некоторым базовым окружением и начальниками $\mathbf{FBoss}(Z)$, $\mathbf{SBoss}(Z)$, $\mathbf{TBoss}(Z)$. Исходя из этих данных, мы можем определить коды вершин $Z$, $Y$, $F$, $J$ во всех четырех случаях расположения и выписать определяющие соотношения для локальных преобразований $8$ и $9$.

\medskip

Покажем, как можно определить эти коды.
Вершины $H_1$, $H_2$, $H_3$  являются начальниками центральной $\mathbb{C}$-вершины. Зная типы, окружения и уровни вершин $H_1$, $H_2$, $H_3$, можно выписать типы, окружения и уровни $Z$, $Y$, $F$, $J$ во всех четырех случаях расположения, что видно на рисунке~\ref{C1b}. Например, окружение $F$ в случае $4$ соответствует $1$-цепи вокруг $H_2$ с указателем $E_7(H_2)$, и уровень $F$ равен $2$.

На рисунке~\ref{C1b} черными кругами отмечены вершины, являющиеся начальниками хотя бы одной из вершин $Z$, $Y$, $F$, $J$ в каждом из четырех случаев, число в круге обозначает номер случая. Окружение всех вершин, отмеченных черными кругами, также можно выписать, зная начальников центральной $\mathbb{C}$-вершины. То есть, мы можем вычислить начальников каждой из четырех вершин во всех случаях.

Уровни вершин $F$ во всех случаях равны $2$, а остальных боковых вершин ($J$ и $Z$) равны~$1$.

\medskip

Таким образом, мы можем определить множество букв $Z$, кодирующих вершины с заданным базовым окружением и начальниками, аналогично с остальными вершинами. Для локального преобразования $8$ буквы $F$ и $Y$ будут выбраны однозначно, а  $Z$ и $J$ -- с точностью до подклееного окружения. Для локального преобразования $9$ все наоборот.

\medskip

{\bf Соотношения для преобразования $8$.}\

\medskip

{\bf Назначение $Z$}. Для каждого из четырех случаев расположения, обозначим символом $Z$ буквы в алфавите, соответствующие кодам вершин с заданными типом, уровнем, базовым окружением и информацией. То есть эти буквы отличаются друг от друга только различными подклееными окружениями.

\medskip

{\bf Назначение $J$}. Аналогично определим буквы $J$, как буквы в алфавите, соответствующие кодам вершин $J$ с заданными типом, уровнем, базовым окружением и информацией, и произвольным подклееным окружением.

\medskip

{\bf Назначение $Y$ и $F$}. Буквы $Y$ и $F$ мы определим как конкретные буквы в алфавите с заданными типом, уровнем, базовым окружением и информацией, при пустом подклееном окружении. Мы это делаем, так как путь проходит через $F$ и $Y$ по плоским ребрам.

Теперь для каждого случая и для каждого разрешенного набора букв $Z$, $F$, $Y$, $J$ выпишем соотношение, представленное в правой части рисунка~\ref{C1b}.
Буквы $Z$, $F$, $Y$, $J$ мы только что определили, а остальные символы в соотношении обозначают входящие и выходящие ребра.

\medskip

{\bf Соотношения для преобразования $9$.}
Тут уже $Z$ и $J$ определены однозначно, а $Y$ и $F$ с точностью для подклееного окружения.
В остальном, все аналогично.

\medskip

{\bf Оценка числа соотношений.}
Так как в каждом соотношении две буквы определены однозначно, а другие две -- с точностью для подклееных окружений, то для каждого случая мы выписываем число соотношений, равное квадрату числа возможных подклееных окружений. Кроме того, все соотношения могут быть выписаны для любого значения параметра {\it флаг макроплитки}.

Мы вводим не более $8FP^2\mathbf{Num}(\mathbb{C})$ соотношений, где $F$ -- число различных флагов макроплиток, $P$ -- число различных подклееных окружений, $\mathbf{Num}(\mathbb{C})$ -- число вершин типа $\mathbb{C}$, то есть число сочетаний ``тип-уровень-окружение-информация''.

\medskip

{\bf Характеризация.}
Пусть есть слово $W$ представляющее собой код пути $X_1 e_1 e_2 X_2 e_3 e_4 X_3$. Покажем, как по нему установить, имеем ли мы дело с локальным преобразованием $8$ или $9$, а также как его провести.

На рисунке~\ref{C1b} отмечены входящие и выходящие ребра в каждом из четырех случаев. Исходя из этого легко выписать свойства пути, позволяющие нам установить, с каким случаем мы имеем дело. Все показано в таблицах~\ref{tableC1ba} и~\ref{tableC1bb}.

\medskip

\begin{table}[hbtp]
\caption{Характеристические условия на принадлежность пути к случаю $\mathbb{C}1$ цепи. }
\centering
\begin{tabular}{|c|c|c|c|}   \hline
Условие на буквы  & \x{Симметричное условие \cr (проход в обратном порядке)} & \x{ случай}  & \x{Локальное \cr преобр} \cr \hline
$\mathbf{Surr}(X_1)=\mathbb{C}11$, $e_1 e_2 e_3 e_4 = u_1 1 2 r$  & $\mathbf{Surr}(X_3)=\mathbb{C}11$, $e_1 e_2 e_3 e_4 = r 2 1 u_1$ & 1  & $8$ (левая) \cr \hline
$\mathbf{Surr}(X_1)=\mathbb{C}11$, $e_1 e_2 e_3 e_4 = 2 1 l 1$ & $\mathbf{Surr}(X_3)=\mathbb{C}11$, $e_1 e_2 e_3 e_4 = 1 l 1 2$  & 1 & $8$ (правая) \cr \hline
$\mathbf{Surr}(X_1)=\mathbb{C}12$, $e_1 e_2 e_3 e_4 = u_1 1 2 r$ & $\mathbf{Surr}(X_3)=\mathbb{C}12$, $e_1 e_2 e_3 e_4 = r 2 1 u_1$  & 2  & $8$ (левая) \cr \hline
$\mathbf{Surr}(X_1)=\mathbb{C}12$, $e_1 e_2 e_3 e_4 = 1 2 r 2$ & $\mathbf{Surr}(X_3)=\mathbb{C}12$, $e_1 e_2 e_3 e_4 = 2 r 2 1$  & 2 & $8$ (правая)\cr \hline
$\mathbf{Surr}(X_1)=\mathbb{C}13$, $e_1 e_2 e_3 e_4 = u_1 1 2 r$ & $\mathbf{Surr}(X_3)=\mathbb{C}13$, $e_1 e_2 e_3 e_4 = r 2 1 u_1$  & 3  & $8$ (левая) \cr \hline
$\mathbf{Surr}(X_1)=\mathbb{C}13$, $e_1 e_2 e_3 e_4 = 2 1 l 1$ & $\mathbf{Surr}(X_3)=\mathbb{C}13$, $e_1 e_2 e_3 e_4 = 1 l 1 2$  & 3 & $8$ (правая)\cr \hline
$\mathbf{Surr}(X_1)=\mathbb{C}14$, $e_1 e_2 e_3 e_4 = u_1 1 2 r$ & $\mathbf{Surr}(X_3)=\mathbb{C}14$, $e_1 e_2 e_3 e_4 = r 2 1 u_1$  & 4  & $8$ (левая) \cr \hline
$\mathbf{Surr}(X_1)=\mathbb{C}14$, $e_1 e_2 e_3 e_4 = 1 2 l 1$ & $\mathbf{Surr}(X_3)=\mathbb{C}14$, $e_1 e_2 e_3 e_4 = 1 l 2 1$  & 4 & $8$ (правая)\cr \hline
  \end{tabular}
\label{tableC1ba}
\end{table}

\medskip

\begin{table}[hbtp]
\caption{Характеристические условия на принадлежность пути к случаю $\mathbb{C}1$ цепи. }
\centering
\begin{tabular}{|c|c|c|c|}   \hline
Условие на буквы  & \x{Симметричное условие \cr (проход в обратном порядке)} & \x{ случай}  & \x{Локальное \cr преобр} \cr \hline
\x{ $\mathbf{Type}(X_1)=\mathbb{B}$, $\mathbf{FBoss}(X_1)=\mathbb{C}11$, \cr $e_1 e_2 e_3 e_4 = 1 u_1 2 1$}  & \x{$\mathbf{Type}(X_3)=\mathbb{B}$,  $\mathbf{FBoss}(X_3)=\mathbb{C}11$, \cr $e_1 e_2 e_3 e_4 = 1 2 u_1 1$} & 1  & $9$ (левая) \cr \hline
\x{ $\mathbf{Type}(X_1)=\mathbb{B}$, $\mathbf{FBoss}(X_1)=\mathbb{C}11$, \cr $e_1 e_2 e_3 e_4 = 2 r 1 l$} & \x{ $\mathbf{Type}(X_3)=\mathbb{B}$, $\mathbf{FBoss}(X_3)=\mathbb{C}11$, \cr $e_1 e_2 e_3 e_4 = l 1 r 2$}  & 1 & $9$ (правая) \cr \hline
\x{ $\mathbf{Type}(X_1)=\mathbb{B}$, $\mathbf{FBoss}(X_1)=\mathbb{C}12$, \cr $e_1 e_2 e_3 e_4 = 1 u_1 1 2$}  & \x{$\mathbf{Type}(X_3)=\mathbb{B}$, $\mathbf{FBoss}(X_3)=\mathbb{C}12$, \cr $e_1 e_2 e_3 e_4 = 2 1 u_1 1$} & 2  & $9$ (левая) \cr \hline
\x{ $\mathbf{Type}(X_1)=\mathbb{B}$, $\mathbf{FBoss}(X_1)=\mathbb{C}12$, \cr $e_1 e_2 e_3 e_4 = 2 r 2 r$} & \x{$\mathbf{Type}(X_3)=\mathbb{B}$, $\mathbf{FBoss}(X_3)=\mathbb{C}12$, \cr $e_1 e_2 e_3 e_4 = r 2 r 2$}  & 2 & $9$ (правая) \cr \hline
\x{ $\mathbf{Type}(X_1)=\mathbb{B}$, $\mathbf{FBoss}(X_1)=\mathbb{C}13$, \cr $e_1 e_2 e_3 e_4 = 1 u_1 2 1$}  & \x{$\mathbf{Type}(X_3)=\mathbb{B}$, $\mathbf{FBoss}(X_3)=\mathbb{C}13$, \cr $e_1 e_2 e_3 e_4 = 1 2 u_1 1$} & 3  & $9$ (левая) \cr \hline
\x{ $\mathbf{Type}(X_1)=\mathbb{B}$, $\mathbf{FBoss}(X_1)=\mathbb{C}13$, \cr $e_1 e_2 e_3 e_4 = 2 r 1 l$} & \x{$\mathbf{Type}(X_3)=\mathbb{B}$, $\mathbf{FBoss}(X_3)=\mathbb{C}13$, \cr $e_1 e_2 e_3 e_4 = l 1 r 2$}  & 3 & $9$ (правая) \cr \hline
\x{$\mathbf{Type}(X_1)=\mathbb{B}$, $\mathbf{FBoss}(X_1)=\mathbb{C}14$, \cr $e_1 e_2 e_3 e_4 = 1 u_1 1 2$}  & \x{$\mathbf{Type}(X_3)=\mathbb{B}$, $\mathbf{FBoss}(X_3)=\mathbb{C}14$, \cr $e_1 e_2 e_3 e_4 = 2 1 u_1 1$} & 4  & $9$ (левая) \cr \hline
\x{$\mathbf{Type}(X_1)=\mathbb{B}$, $\mathbf{FBoss}(X_1)=\mathbb{C}14$, \cr $e_1 e_2 e_3 e_4 = 2 r 1 l$} & \x{$\mathbf{Type}(X_3)=\mathbb{B}$, $\mathbf{FBoss}(X_3)=\mathbb{C}14$, \cr $e_1 e_2 e_3 e_4 = l 1 r 2$}  & 4 & $9$ (правая) \cr \hline
  \end{tabular}
\label{tableC1bb}
\end{table}

\medskip

Ясно, что указанные условия на буквы полностью определяют конфигурацию пути, а также его код, с точностью до подклееных окружений. Иначе говоря, не бывает никакого другого пути с заданными условиями на буквы, кроме пути, указанного нами.

То есть, для заданной конфигурации есть конечное множество слов $Q$, которые могли бы кодировать такой путь. Два слова в этом множестве отличаются только кодами подклееных окружений для крайних вершинных букв.

\medskip

{\bf Восстановление кода.}
Зафиксируем некоторый случай расположения $n$, один из четырех. Допустим, мы знаем коды трех из четырех вершин (из числа $Z$, $Y$, $J$, $F$). Тогда, учитывая характеризацию, мы знаем с каким из случаев расположения мы имеем дело. Далее, окружение вершины $H_1$ мы вычисляем, так как она первый начальник $Z$ и $F$. В первом и втором случаях расположения знания окружения $H_1$ достаточно, чтобы определить окружение оставшейся вершины. Для третьего и четвертого случаев, мы знаем также окружение $H_2$ и $H_3$ -- второй и третий начальник $F$ и $Z$. Этого достаточно, чтобы определить оставшееся окружение.

Теперь покажем, как определить начальников.
У $F$ и $Z$ общее множество начальников. Первым начальником $Y$ всегда является $Z$, а тип второго очевиден из расположения. В  случаях расположения $1$, $3$, $4$ начальники $J$ лежат в той же цепи, что и $F$.

В случае $2$, первый начальник $J$ соответствует $1$-цепи вокруг $H_3$. Тип $H_3$ мы знаем, так как это второй начальник $F$ и $Z$. Указатель можно установить по окружению макроплитки, $\mathbf{FBoss}(Z)$. Тип второго начальника $J$ -- это $\mathbb{B}$.

Таким образом, зная код пути $ZYJ$, можно вычислить код пути $ZFJ$, и наоборот. А также зная код пути $FZY$, можно вычислить код пути $FJY$, и наоборот.

\medskip

Пусть слово $W$ подходит под один из перечисленных случаев. Тогда мы
знаем с каким случаем мы имеем дело, а также какую конфигурацию имеет путь. Слово $W$ должно кодировать этот путь, то есть $W$ должно входить в $Q$.  Если слово $W$ не удовлетворяет этому условию, например, если код средней вершины не такой, какой должен быть исходя из конфигурации пути, то $W$ кодирует невозможный путь. В этом случае $W=0$, так как мы ввели обнуляющие соотношения для всех достаточно коротких невозможных путей.

\medskip
Если же $W$ входит в $Q$, то к $W$ можно применить соотношение введенное нами, из числа указанных справа на рисунке~\ref{C1b}, и получить код для другой части локального преобразования.
Таким образом, мы можем осуществить локальное преобразование пути через операцию с его кодом.

\medskip

\subsection{Случай цепи $\mathbb{C}1$; преобразования 7 и 10}

\begin{figure}[hbtp]
\centering
\leftskip=-1.0cm
\includegraphics[width=1\textwidth]{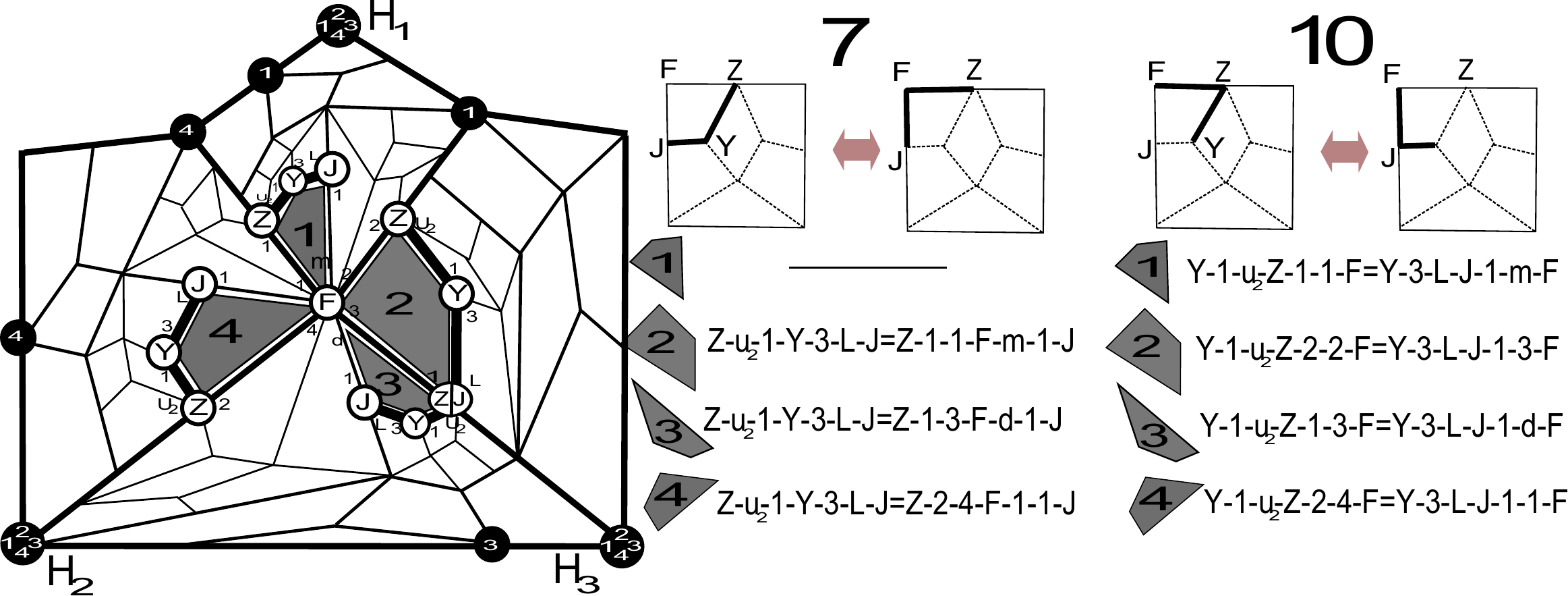}
\caption{Случаи расположения пути вокруг C1-цепи и соответствующие им определяющие соотношения для локальных преобразований 7 и 10}
\label{C1a}
\end{figure}

В правой части рисунка~\ref{C1a} изображены локальные преобразования $7$ и $10$.

\medskip

Итак, зафиксируем $\mathbb{C}$-узел с некоторым окружением и начальниками. Различные случаи расположения интересующих нас путей показаны в левой части рисунка~\ref{C1a}.

\medskip

Зная окружение и начальников центрального $\mathbb{C}$-узла, можно вычислить базовые окружения и начальников вершин $Z$, $Y$, $F$, $J$ во всех случаях расположения.
То есть мы можем определить коды вершин $Z$, $Y$, $F$, $J$ во всех четырех случаях расположения и выписать определяющие соотношения для локальных преобразований $7$ и $10$.

Уровни боковых вершин (то есть $J$ и $Z$) во всех случаях равны $1$.

\medskip

{\bf Соотношения для  преобразования 7.}

\medskip

{\bf Назначение $Z$}. Для каждого из четырех случаев расположения, обозначим символом $Z$ буквы в алфавите, соответствующие кодам вершин $Z$ с заданными типом, уровнем, базовым окружением и информацией. То есть эти буквы отличаются друг от друга только различными подклееными окружениями.

\medskip

{\bf Назначение $J$}. Аналогично определим буквы $J$, как буквы в алфавите, соответствующие кодам вершин $J$ с заданными типом, уровнем, базовым окружением и информацией, и произвольным подклееным окружением.

\medskip

{\bf Назначение $Y$ и $F$}. Буквы $Y$ и $F$ мы определим как конкретные буквы в алфавите с заданными типом, уровнем, базовым окружением и информацией, при пустом подклееном окружении. Мы это делаем, так как путь проходит через $F$ и $Y$ по плоским ребрам.

Теперь для каждого случая и для каждого разрешенного набора букв $Z$, $F$, $Y$, $J$ выпишем соотношение, представленное в правой части рисунка~\ref{C1a}.
Буквы $Z$, $F$, $Y$, $J$ мы только что определили, а остальные символы в соотношении обозначают входящие и выходящие ребра.

\medskip

{\bf Соотношения для преобразования 10.}
Тут уже $Z$ и $J$ определены однозначно, а $Y$ и $F$ с точностью для подклееного окружения.
В остальном, все аналогично.

\medskip

{\bf Оценка числа соотношений.}
Мы вводим не более $7FP^2\mathbf{Num}(\mathbb{C})$ соотношений, где $F$ -- число различных флагов макроплиток, $P$ -- число различных подклееных окружений, $\mathbf{Num}(\mathbb{C})$ -- число вершин типа $\mathbb{C}$.

В последующих случаях отношения считаются аналогично.

\medskip

{\bf Характеризация.}
Пусть есть слово $W$ представляющее собой код пути $X_1 e_1 e_2 X_2 e_3 e_4 X_3$. Покажем, как по нему установить, имеем ли мы дело с локальным преобразованием $7$ или $10$, а также как его провести.

На рисунке~\ref{C1a} отмечены входящие и выходящие ребра в каждом из четырех случаев. Исходя из этого легко выписать свойства пути, позволяющие нам установить, с каким случаем мы имеем дело. Все показано в таблицах~\ref{tableC1aa} и~\ref{tableC1ab}.

\begin{table}[hbtp]
\caption{Характеристические условия на принадлежность пути к случаю $\mathbb{C}1$ цепи. }
\centering
\begin{tabular}{|c|c|c|c|}   \hline
Условие на буквы  & \x{Симметричное условие \cr (проход в обратном порядке)} & \x{ случай}  & \x{Локальное \cr преобр} \cr \hline
$\mathbf{Surr}(X_1)=\mathbb{C}11$, $e_1 e_2 e_3 e_4 = u_1 1 3 l$  & $\mathbf{Surr}(X_3)=\mathbb{C}11$, $e_1 e_2 e_3 e_4 = l 3 1 u_1$ & 1  & $7$ (левая) \cr \hline
$\mathbf{Surr}(X_1)=\mathbb{C}11$, $e_1 e_2 e_3 e_4 = 1 1 m 1$ & $\mathbf{Surr}(X_3)=\mathbb{C}11$, $e_1 e_2 e_3 e_4 = 1 m 1 1$  & 1 & $7$ (правая) \cr \hline
$\mathbf{Surr}(X_1)=\mathbb{C}12$, $e_1 e_2 e_3 e_4 = u_2 1 3 l$ & $\mathbf{Surr}(X_3)=\mathbb{C}12$, $e_1 e_2 e_3 e_4 = l 3 1 u_2$  & 2  & $7$ (левая) \cr \hline
$\mathbf{Surr}(X_1)=\mathbb{C}12$, $e_1 e_2 e_3 e_4 = 2 2 3 1$ & $\mathbf{Surr}(X_3)=\mathbb{C}12$, $e_1 e_2 e_3 e_4 = 1 3 2 2$  & 2 & $7$ (правая)\cr \hline
$\mathbf{Surr}(X_1)=\mathbb{C}13$, $e_1 e_2 e_3 e_4 = u_2 1 3 l$ & $\mathbf{Surr}(X_3)=\mathbb{C}13$, $e_1 e_2 e_3 e_4 = l 3 1 u_2$  & 3  & $7$ (левая) \cr \hline
$\mathbf{Surr}(X_1)=\mathbb{C}13$, $e_1 e_2 e_3 e_4 = 1 3 d 1$ & $\mathbf{Surr}(X_3)=\mathbb{C}13$, $e_1 e_2 e_3 e_4 = 1 d 3 1$  & 3 & $7$ (правая)\cr \hline
$\mathbf{Surr}(X_1)=\mathbb{C}14$, $e_1 e_2 e_3 e_4 = u_2 1 3 l$ & $\mathbf{Surr}(X_3)=\mathbb{C}14$, $e_1 e_2 e_3 e_4 = l 3 1 u_2$  & 4  & $7$ (левая) \cr \hline
$\mathbf{Surr}(X_1)=\mathbb{C}14$, $e_1 e_2 e_3 e_4 = 2 4 1 1$ & $\mathbf{Surr}(X_3)=\mathbb{C}14$, $e_1 e_2 e_3 e_4 = 1 1 4 2$  & 4 & $7$ (правая)\cr \hline
  \end{tabular}
\label{tableC1aa}
\end{table}

\begin{table}[hbtp]
\caption{Характеристические условия на принадлежность пути к случаю $\mathbb{C}1$ цепи.}
\centering
\begin{tabular}{|c|c|c|c|}   \hline
Условие на буквы  & \x{Симметричное условие \cr (проход в обратном порядке)} & \x{ случай}  & \x{Локальное \cr преобр} \cr \hline
\x{ $\mathbf{Type}(X_1)=\mathbb{B}$, $\mathbf{FBoss}(X_1)=\mathbb{C}11$, \cr $e_1 e_2 e_3 e_4 = 1 u_1 1 1$}  & \x{$\mathbf{Type}(X_3)=\mathbb{A}$,  $\mathbf{FBoss}(X_3)=\mathbb{C}11$, \cr $e_1 e_2 e_3 e_4 = 1 1 u_1 1$} & 1  & $10$ (левая) \cr \hline
\x{ $\mathbf{Type}(X_1)=\mathbb{A}$, $\mathbf{FBoss}(X_1)=\mathbb{C}11$, \cr $e_1 e_2 e_3 e_4 = 3 l 1 m$} & \x{ $\mathbf{Type}(X_3)=\mathbb{A}$, $\mathbf{FBoss}(X_3)=\mathbb{C}11$, \cr $e_1 e_2 e_3 e_4 = m 1 l 3$}  & 1 & $10$ (правая) \cr \hline
\x{ $\mathbf{Type}(X_1)=\mathbb{A}$, $\mathbf{FBoss}(X_1)=\mathbb{C}12$, \cr $e_1 e_2 e_3 e_4 = 1 u_2 2 2$}  & \x{$\mathbf{Type}(X_3)=\mathbb{A}$, $\mathbf{FBoss}(X_3)=\mathbb{C}12$, \cr $e_1 e_2 e_3 e_4 = 2 2 u_2 1$} & 2  & $10$ (левая) \cr \hline
\x{ $\mathbf{Type}(X_1)=\mathbb{A}$, $\mathbf{FBoss}(X_1)=\mathbb{C}12$, \cr $e_1 e_2 e_3 e_4 = 3 l 1 3$} & \x{$\mathbf{Type}(X_3)=\mathbb{A}$, $\mathbf{FBoss}(X_3)=\mathbb{C}12$, \cr $e_1 e_2 e_3 e_4 = 3 1 l 3$}  & 2 & $10$ (правая) \cr \hline
\x{ $\mathbf{Type}(X_1)=\mathbb{A}$, $\mathbf{FBoss}(X_1)=\mathbb{C}13$, \cr $e_1 e_2 e_3 e_4 = 1 u_2 1 3$}  & \x{$\mathbf{Type}(X_3)=\mathbb{A}$, $\mathbf{FBoss}(X_3)=\mathbb{C}13$, \cr $e_1 e_2 e_3 e_4 = 3 1 u_2 1$} & 3  & $10$ (левая) \cr \hline
\x{ $\mathbf{Type}(X_1)=\mathbb{A}$, $\mathbf{FBoss}(X_1)=\mathbb{C}13$, \cr $e_1 e_2 e_3 e_4 = 3 l 1 d$} & \x{$\mathbf{Type}(X_3)=\mathbb{A}$, $\mathbf{FBoss}(X_3)=\mathbb{C}13$, \cr $e_1 e_2 e_3 e_4 = d 1 l 3$}  & 3 & $10$ (правая) \cr \hline
\x{$\mathbf{Type}(X_1)=\mathbb{A}$, $\mathbf{FBoss}(X_1)=\mathbb{C}14$, \cr $e_1 e_2 e_3 e_4 = 1 u_2 2 4$}  & \x{$\mathbf{Type}(X_3)=\mathbb{A}$, $\mathbf{FBoss}(X_3)=\mathbb{C}14$, \cr $e_1 e_2 e_3 e_4 = 4 2 u_2 1$} & 4  & $10$ (левая) \cr \hline
\x{$\mathbf{Type}(X_1)=\mathbb{B}$, $\mathbf{FBoss}(X_1)=\mathbb{C}14$, \cr $e_1 e_2 e_3 e_4 = 3 l 1 ld$} & \x{$\mathbf{Type}(X_3)=\mathbb{A}$, $\mathbf{FBoss}(X_3)=\mathbb{C}14$, \cr $e_1 e_2 e_3 e_4 = ld 1 l 3$}  & 4 & $10$ (правая) \cr \hline
  \end{tabular}
 \label{tableC1ab}
\end{table}

Ясно, что указанные условия на буквы полностью определяют конфигурацию пути, а также его код, с точностью до подклееных окружений. Иначе говоря, не бывает никакого другого пути с заданными условиями на буквы, кроме пути, указанного нами.

То есть, для заданной конфигурации есть конечное множество слов $Q$, которые могли бы кодировать такой путь. Два слова в этом множестве отличаются только кодами подклееных окружений для крайних вершинных букв.

\medskip

{\bf Восстановление кода.}
Зная код $F$, все остальные коды легко вычисляются. То есть остается разобраться, как определить код $F$ по известным кодам остальных вершин, в $7$ локальном преобразовании,  случаи расположения $2$, $3$, $4$.
Окружение $F$ совпадает с $U$-частью окружения $H_1$ первого начальника $Z$. Заметим, что $Z$ в случае $3$ и $4$, и $J$ в случае $2$ имеют тех же начальников, что и $F$.

Таким образом, зная код пути $ZYJ$, можно вычислить код пути $ZFJ$, и наоборот. А также зная код пути $FZY$, можно вычислить код пути $FJY$, и наоборот.

\medskip

Пусть слово $W$ подходит под один из перечисленных случаев.
Аналогично предыдущему параграфу,  мы можем определить с каким случаем мы имеем дело, а также какую конфигурацию имеет путь. Слово $W$ должно кодировать этот путь, то есть $W$ должно входить в $Q$. В противном случае $W$ приводится к нулю.

\medskip
Если же $W$ входит в $Q$, то к $W$ можно применить соотношение введенное нами, из числа указанных справа на рисунке~\ref{C1a}, и получить код для другой части локального преобразования.

Таким образом, мы можем осуществить локальное преобразование пути через операцию с его кодом.

\medskip

\subsection{Случай цепи $\mathbb{C}2$; преобразования 8 и 9}

\begin{figure}[hbtp]
\centering
\includegraphics[width=1\textwidth]{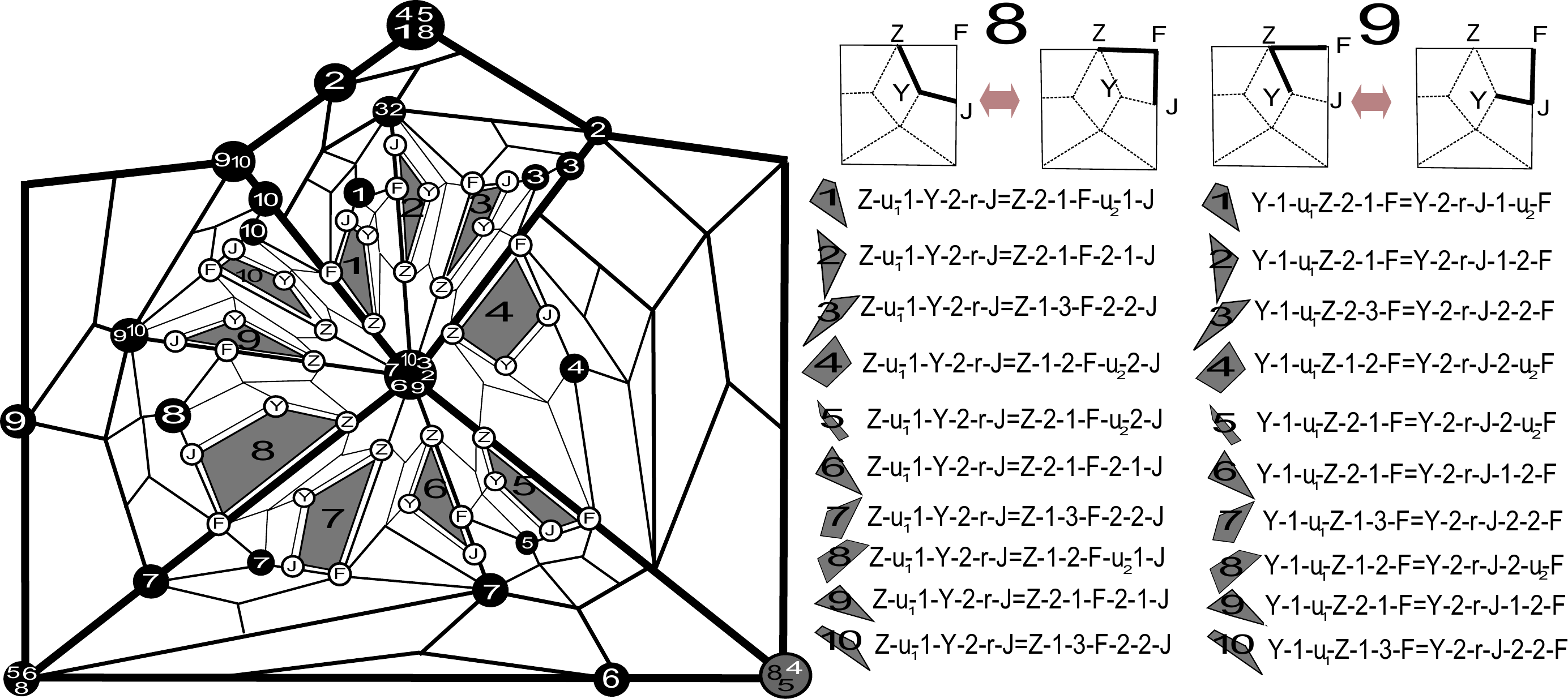}
\caption{Случаи расположения пути вокруг C2-цепи и соответствующие им определяющие соотношения для локальных преобразований 8 и 9}
\label{C2b}
\end{figure}

В правой части рисунка~\ref{C2b} изображены локальные преобразования $8$ и $9$.

\medskip

Аналогично предыдущим случаям, мы можем ввести определяющие соотношения. Фиксируем $\mathbb{C}$-узел с базовым окружением и тремя начальниками.
Ясно, что зная код центрального узла, типы, уровни и окружения вершин $Z$, $Y$, $F$, $J$ можно вычислить (с точностью до подклееными окружениями крайних вершин в пути).

Черными кругами отметим вершины, являющиеся начальниками вершин  $Z$, $Y$, $F$, $J$. Заметим, что зная код центрального $\mathbb{C}$-узла мы можем найти их типы, уровни и окружения. То есть коды вершин $Z$, $Y$, $F$, $J$ во всех десяти случаях мы можем назвать явно, с точностью до подклееных окружений.

Это позволяет ввести определяющие отношения, записанные в правой части рисунка~\ref{C2b}.

\medskip

{\bf Характеризация.}
Пусть есть слово $W$ представляющее собой код пути $X_1 e_1 e_2 X_2 e_3 e_4 X_3$.
Аналогично предыдущим случаям, мы можем выписать характеризующую последовательность кодов вершин и входящих-выходящих ребер. Это позволит для любого слова определить, является ли он кодом локального преобразования $8$ или $9$.

Мы не будем приводить здесь таблицу как для случая $\mathbb{C}1$-цепи, ее можно составить полностью аналогично, опираясь на рисунок~\ref{C2b}.

\medskip

{\bf Восстановление кода.} Учитывая характеризацию, мы знаем, с каким случаем расположения имеем дело. Кроме того, ясно, что окружение каждой из четырех вершин можно определить во всех десяти случаях. Заметим, что у $Z$ и $F$ одно и то же множество начальников во всех случаях. Первым начальником $Y$ всегда является $Z$, а тип второго очевиден из расположения. Остается определить начальников $J$. В  случаях $1$, $4$, $5$, $8$ это просто $F$, а в остальных случаях у $J$ либо те же начальники что и у $F$, либо один начальник, совпадающий с первым начальником $F$.

Таким образом, зная код пути $ZYJ$, можно вычислить код пути $ZFJ$, и наоборот. А также зная код пути $FZY$, можно вычислить код пути $FJY$, и наоборот.

\medskip

Таким образом, если слово $W$ подходит под один из перечисленных случаев то к нему можно применить соотношение введенное нами, из числа указанных справа на рисунке~\ref{C2b}, и получить код для другой части локального преобразования.

\medskip

\subsection{Случай цепи $\mathbb{C}2$; преобразования 7 и 10}

\begin{figure}[hbtp]
\centering
\includegraphics[width=1\textwidth]{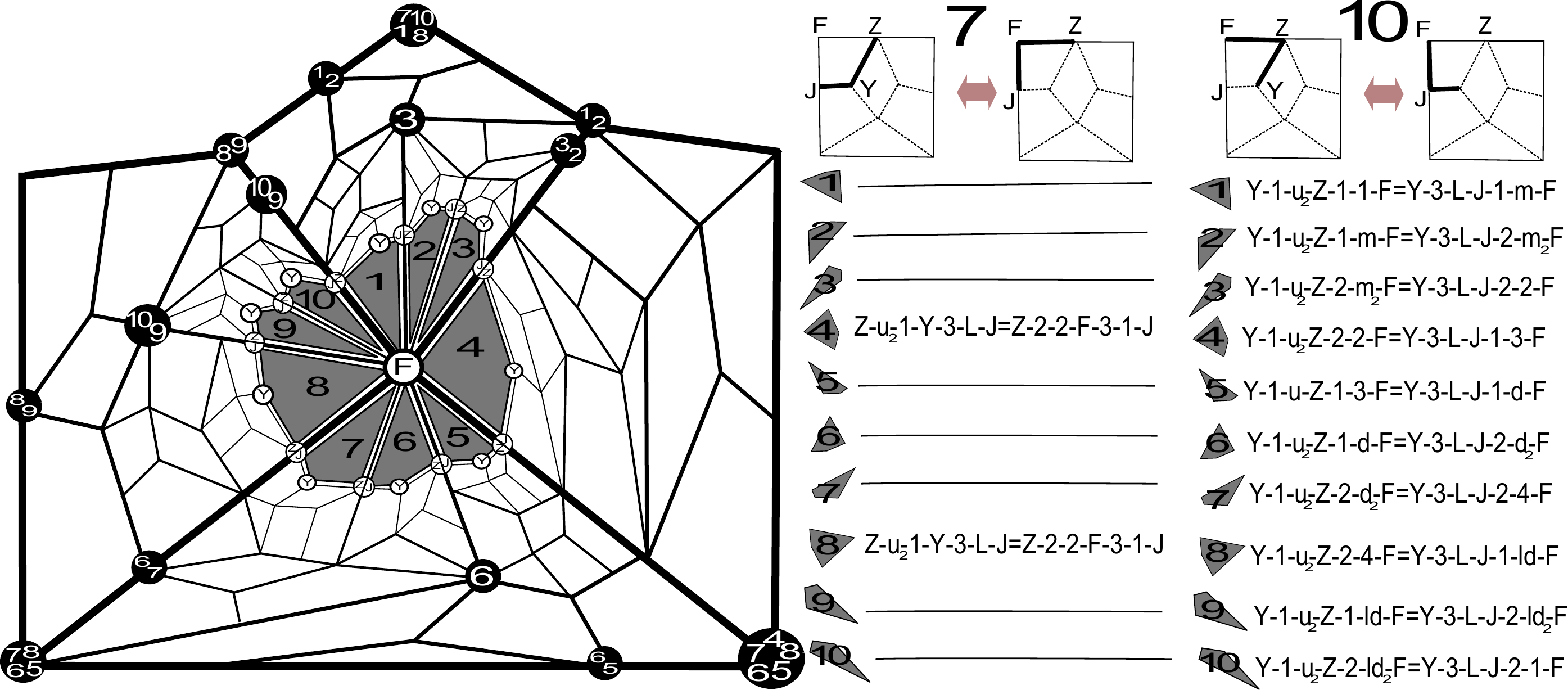}
\caption{Случаи расположения пути вокруг C2-цепи и соответствующие им определяющие соотношения для локальных преобразований 7 и 10}
\label{C2a}
\end{figure}

В правой части рисунка~\ref{C2a} изображены локальные преобразования $7$ и $10$.

\medskip

Черными кругами отметим вершины, являющиеся начальниками вершин  $Z$, $Y$, $F$, $J$. Заметим, что зная код центрального $\mathbb{C}$-узла мы можем найти их окружения. То есть коды вершин $Z$, $Y$, $F$, $J$ (включающие типы, уровни, окружения и информацию), во всех десяти случаях мы можем назвать явно, с точностью до подклееных окружений.
Это позволяет ввести определяющие отношения, записанные в правой части рисунка~\ref{C2a}.

\medskip

{\bf Характеризация.} Составление таблиц полностью аналогично случаю $\mathbb{C}1$-цепи.

\medskip

{\bf Восстановление кода.}
Рассмотрим случаи расположения $4$ и $8$ для преобразования $7$. Заметим, что $Z$ в  случае $8$ и $J$ в случае $4$ имеют тех же начальников, что и $F$. Окружение $F$ совпадает с $U$-частью его первого начальника. Значит, зная коды $Z$, $J$, $Y$, можно вычислить код $F$. Зная же код $F$, очевидно, все остальные вершины тоже вычисляются.
Для преобразования $10$, зная коды $F$, $J$, $Y$, можно вычислить $Z$, а зная $F$, $Z$, $Y$, можно вычислить $J$.

Таким образом, зная код пути $ZYJ$, можно вычислить код пути $ZFJ$, и наоборот. А также зная код пути $FZY$, можно вычислить код пути $FJY$, и наоборот.

\medskip

Таким образом, к $W$ можно применить соотношение введенное нами, из числа указанных справа на рисунке~\ref{C2a}, и получить код для другой части локального преобразования.

\medskip

\subsection{Случай цепи $\mathbb{C}3$; преобразования 7, 8, 9, 10}

Случай $\mathbb{C}3$-цепи полностью аналогичен случаю $\mathbb{C}2$-цепи, соотношения выглядят идентично, только кодировки вершин $J$, $F$, $Z$, $Y$ отвечают $\mathbb{C}3$-цепи. Все рассуждения о вычислении путей полностью аналогичны. Соотношений вводится столько же, сколько для $\mathbb{C}2$ случая, то есть $32FP^2\mathbf{Num}(\mathbb{C})$.

\medskip

\subsection{Случай цепи $\mathbb{B}1$; преобразования 8 и 9}

\begin{figure}[hbtp]
\centering
\includegraphics[width=1\textwidth]{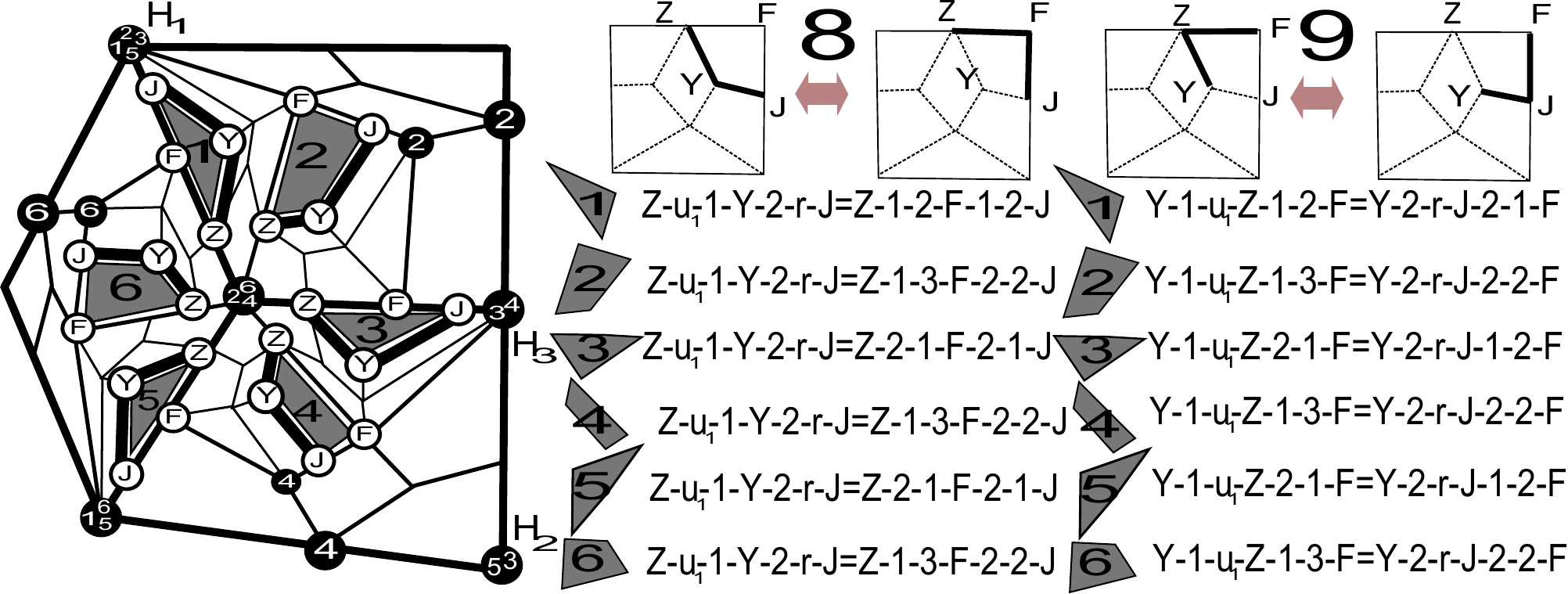}
\caption{Случаи расположения пути вокруг B1-цепи и соответствующие им определяющие соотношения для локальных преобразований 8 и 9}
\label{B1b}
\end{figure}

В правой части рисунка~\ref{B1b} изображены локальные преобразования $8$ и $9$.

\medskip

Аналогично предыдущим случаям, мы можем ввести определяющие соотношения. Фиксируем вершину типа $\mathbb{B}$ и ее начальников, вершины $H_1$ и $H_2$, причем от $H_2$ мы фиксируем только тип.
Черными кругами отметим вершины, являющиеся начальниками вершин  $Z$, $Y$, $F$, $J$.
Числа в них обозначают номер случая расположения. Заметим, что мы можем найти их окружения (кроме вершины $H_2$). Например, $H_3$ вычисляется с помощью $\mathbf{RightFromB}$. Окружения для вершин $Z$, $Y$, $F$, $J$ в каждом случае вычисляются с учетом известных нам окружения $H_1$ и типа $H_2$. То есть что коды вершин $Z$, $Y$, $F$, $J$ во всех случаях мы можем назвать явно, с точностью до подклееных окружений.
Это позволяет ввести определяющие отношения, записанные в правой части рисунка~\ref{B1b}.

\medskip

{\bf Характеризация.} Составление таблиц полностью аналогично случаю $\mathbb{C}1$-цепи.

\medskip

{\bf Восстановление кода.}
Пусть мы знаем коды трех вершин из $Z$, $J$, $Y$, $F$. Поскольку мы можем установить, с каким именно случаем расположения мы имеем дело, то и окружение оставшейся вершины мы можем вычислить.

$F$ и $Z$ имеют общее множество начальников во всех случаях. Первым начальником $Y$ всегда является $Z$, тип второго очевиден из расположения в каждом случае. В случаях~$1$, $3$, $5$ у $Y$ те же начальники что и у $F$, в случаях $2$, $4$, $6$ у $J$ первый начальник как у $F$, а второй -- как третий у $F$.

Таким образом, зная код пути $ZYJ$, можно вычислить код пути $ZFJ$, и наоборот. А также зная код пути $FZY$, можно вычислить код пути $FJY$, и наоборот.

\medskip

\subsection{Случай цепи $\mathbb{B}1$; преобразования 7 и 10}

\begin{figure}[hbtp]
\centering
\includegraphics[width=1\textwidth]{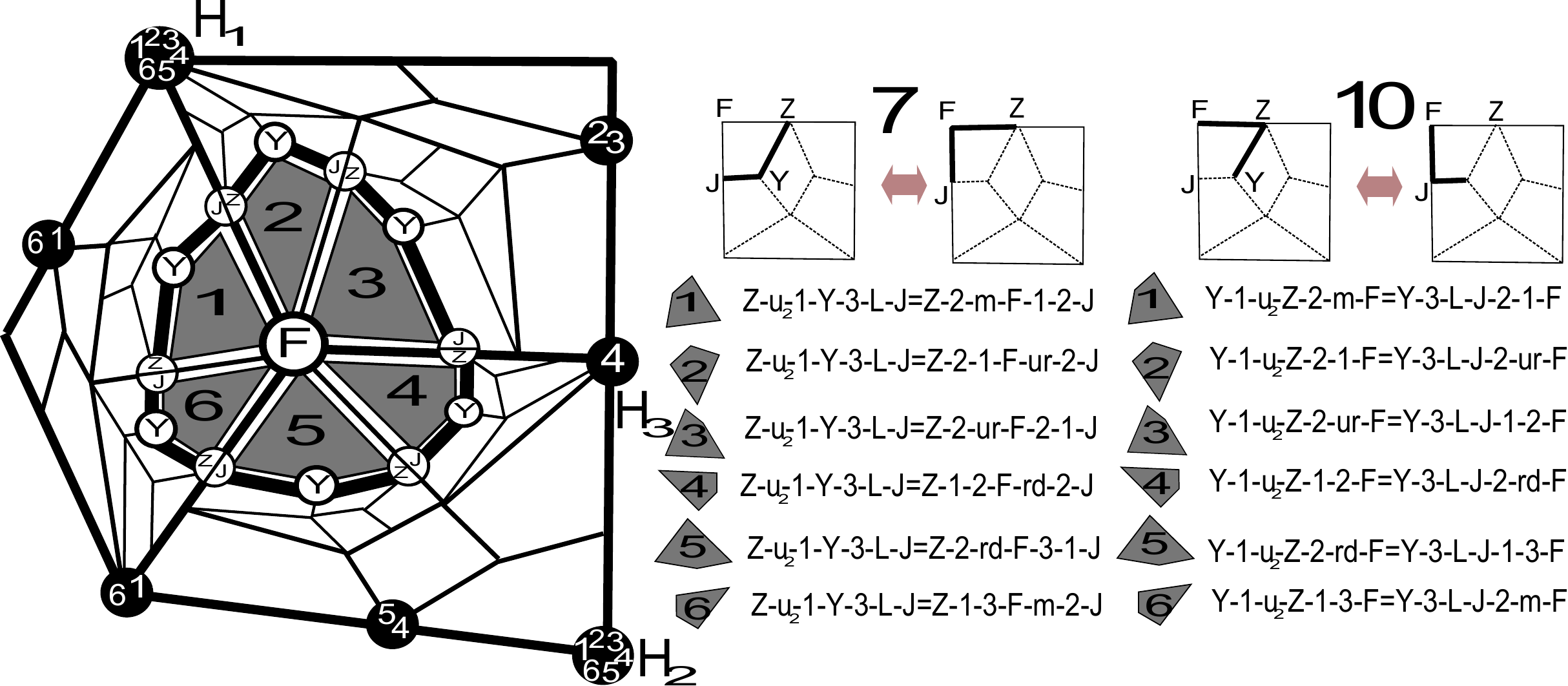}
\caption{Случаи расположения пути вокруг B1-цепи и соответствующие им определяющие соотношения для локальных преобразований 7 и 10}
\label{B1a}
\end{figure}

В правой части рисунка~\ref{B1a} изображены локальные преобразования $7$ и $10$.

\medskip

Аналогично предыдущим случаям, мы можем ввести определяющие соотношения. Фиксируем вершину типа $\mathbb{B}$ и ее начальников, вершины $H_1$ и $H_2$, причем от $H_2$ мы фиксируем только тип.
Черными кругами отметим вершины, являющиеся начальниками вершин  $Z$, $Y$, $F$, $J$.
Числа в них обозначают номер случая расположения. Заметим, что мы можем найти их типы, уровни и окружения (кроме вершины $H_2$, для нее только тип). Типы, уровни и окружения для вершин $Z$, $Y$, $F$, $J$, в каждом случае вычисляются с учетом известных нам окружения $H_1$ и типа $H_2$. То есть что коды вершин $Z$, $Y$, $F$, $J$ во всех случаях мы можем назвать явно, с точностью до подклееных окружений.
Это позволяет ввести определяющие отношения, записанные в правой части рисунка~\ref{B1a}.

\medskip

{\bf Характеризация.} Составление таблиц полностью аналогично случаю $\mathbb{C}1$-цепи.

\medskip

{\bf Восстановление кода.}
Для локального преобразования $7$, если нам нужно вычислить код $F$ по известным кодам $Z$, $J$, $Y$, можно заметить, что хотя бы у одной вершины из  $Z$, $J$, $Y$ начальником будет $F$, а у какой-либо другой начальники будут совпадать с $F$. Если же $F$ нам известно, все остальные вершины, очевидно, вычисляются.

Таким образом, зная код пути $ZYJ$, можно вычислить код пути $ZFJ$, и наоборот. А также зная код пути $FZY$, можно вычислить код пути $FJY$, и наоборот.

\medskip

\subsection{Случай цепи $\mathbb{B}2$; преобразования 8 и 9}

\begin{figure}[hbtp]
\centering
\includegraphics[width=1\textwidth]{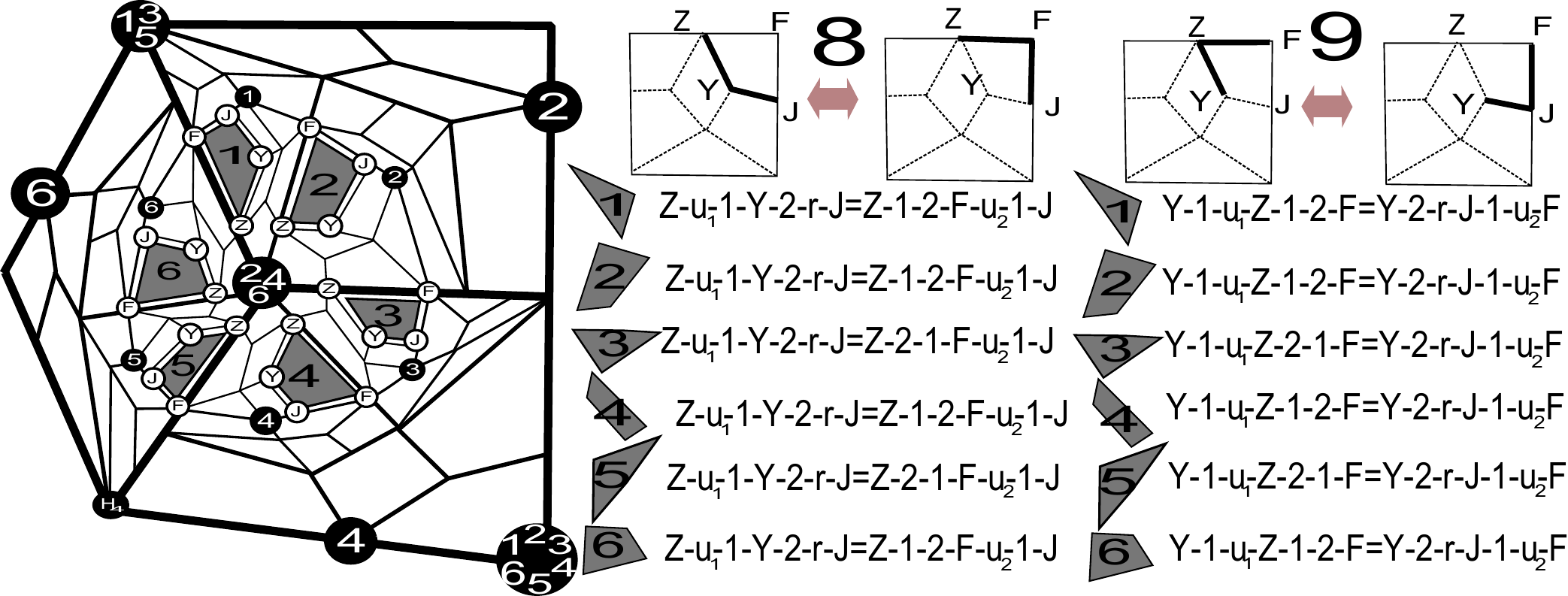}
\caption{Случаи расположения пути вокруг B2-цепи и соответствующие им определяющие соотношения для локальных преобразований 8 и 9}
\label{B2b}
\end{figure}

В правой части рисунка~\ref{B2b} изображены локальные преобразования $8$ и $9$.

\medskip

Введение определяющих соотношений полностью аналогично предыдущим случаям, можно легко проверить, что зная центральный $\mathbb{B}$-узел, коды всех четырех вершин в каждом из случаев расположения легко можно выписать. Это позволяет ввести определяющие отношения, записанные в правой части рисунка~\ref{B2b}.

\medskip

{\bf Характеризация.} Составление таблиц полностью аналогично случаю $\mathbb{C}1$-цепи.

\medskip

{\bf Восстановление кода.}
Поскольку мы можем установить, с каким именно случаем расположения мы имеем дело, то и окружение каждой из четырех вершин мы можем вычислить.
Кроме того, заметим, что каждая из четырех вершин  $Z$, $J$, $Y$, $F$ в каждом из случаев либо имеет начальников среди других трех вершин, либо у нее общие начальники с какой-то из трех других вершин. Значит, в каждом из случаев, зная коды трех вершин, можно вычислить код четвертой.

Таким образом, зная код пути $ZYJ$, можно вычислить код пути $ZFJ$, и наоборот. А также зная код пути $FZY$, можно вычислить код пути $FJY$, и наоборот.

\medskip

\subsection{Случай цепи $\mathbb{B}2$; преобразования 7 и 10}

\begin{figure}[hbtp]
\centering
\includegraphics[width=1\textwidth]{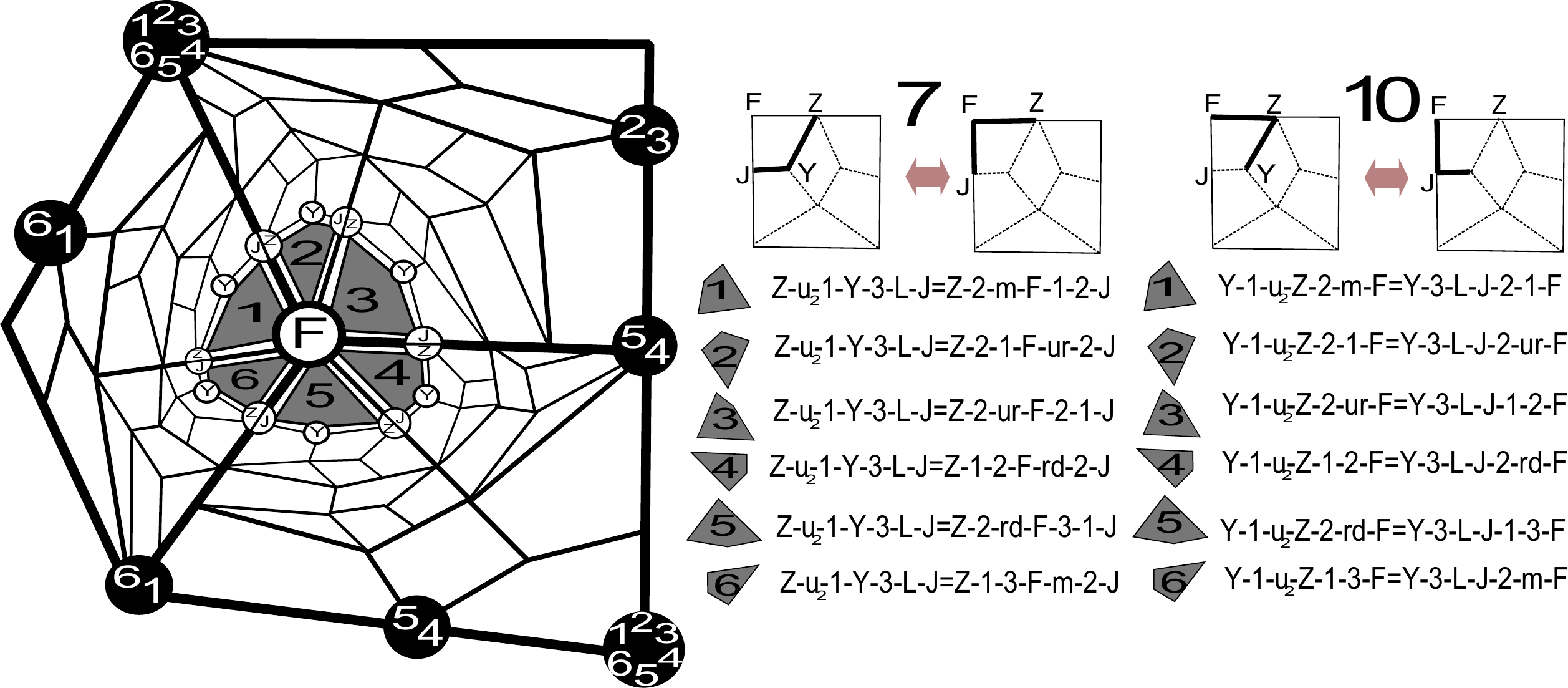}
\caption{Случаи расположения пути вокруг B2-цепи и соответствующие им определяющие соотношения для локальных преобразований 7 и 10}
\label{B2a}
\end{figure}

В правой части рисунка~\ref{B2a} изображены локальные преобразования $7$ и $10$.

\medskip

Зафиксируем вершину типа $\mathbb{B}$ и ее начальников. Аналогично предыдущим случаям, мы можем вычислить коды всех вершин и ввести определяющие соотношения.

\medskip

{\bf Характеризация.} Составление таблиц полностью аналогично случаю $\mathbb{C}1$-цепи.

\medskip

{\bf Восстановление кода.}
Для локального преобразования $7$, если нам нужно вычислить код $F$ по известным кодам $Z$, $J$, $Y$ можно заметить, что хотя бы у одной вершины из  $Z$, $J$, $Y$ начальники будут совпадать с $F$. Если же $F$ нам известно, все остальные вершины, очевидно, вычисляются.

Таким образом, зная код пути $ZYJ$, можно вычислить код пути $ZFJ$, и наоборот. А также зная код пути $FZY$, можно вычислить код пути $FJY$, и наоборот.

\medskip

\subsection{Случай цепи $\mathbb{A}0$; преобразования 8 и 9}

\begin{figure}[hbtp]
\centering
\includegraphics[width=0.8\textwidth]{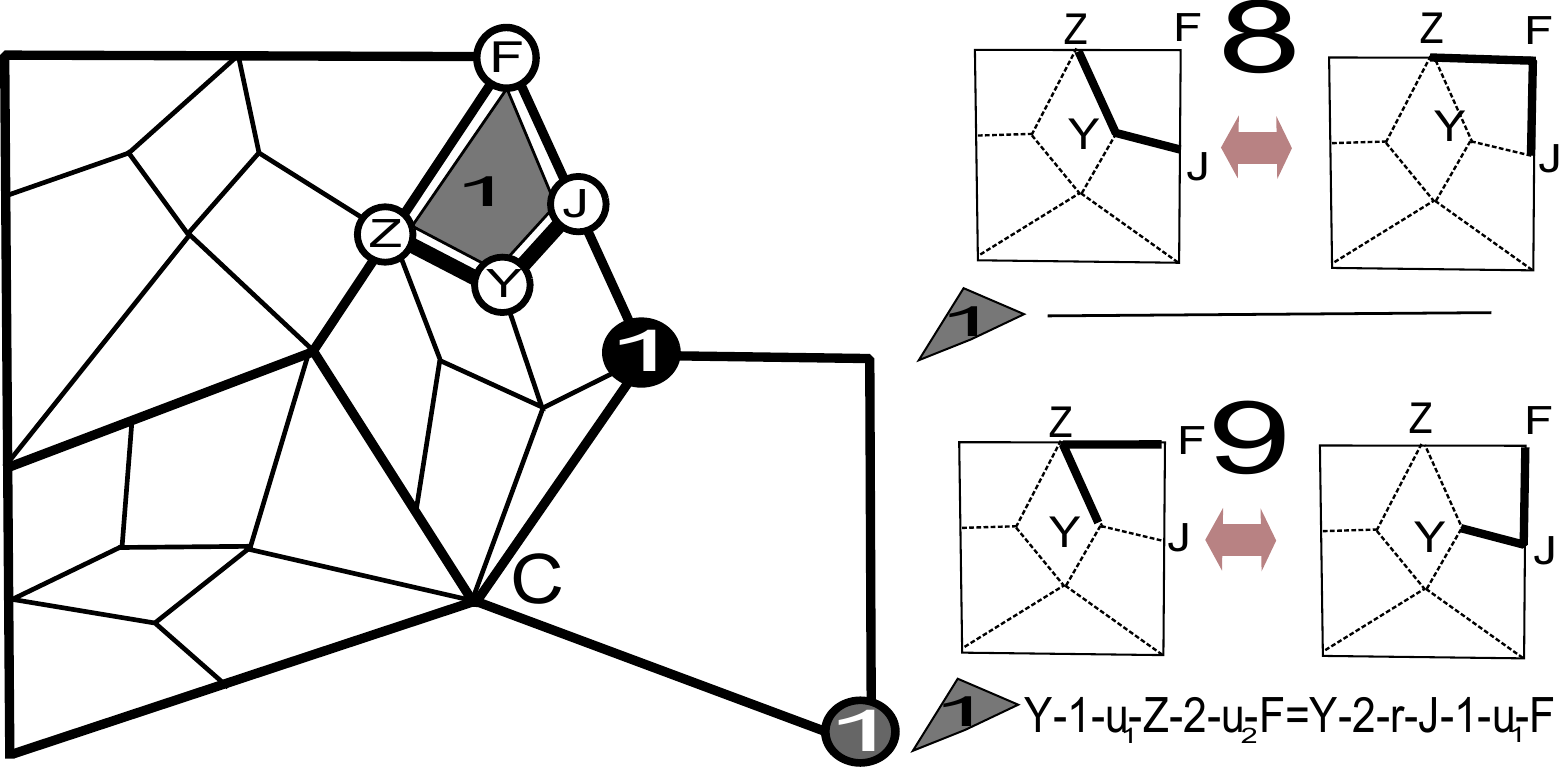}
\caption{Случаи расположения пути вокруг A0-цепи и соответствующие им определяющие соотношения для локальных преобразований 8 и 9}
\label{A1b}
\end{figure}

В правой части рисунка~\ref{A1b} изображены локальные преобразования $8$ и $9$.

\medskip

Зафиксируем вершину $X$ второго уровня с типом $\mathbb{UL}$, $\mathbb{UR}$ или $\mathbb{U}$. Мы будем вводить соотношения только для преобразования $9$. Путь в преобразовании $8$ подпадает под мертвый паттерн, для него определяющие соотношения вводить не требуется. По коду $X$ вычисляются коды остальных трех вершин, например тип второго начальника $J$ -- вершины в правом нижнем углу может быть найден с помощью функции $\mathbf{BottomRightType}(X)$.

\medskip

{\bf Характеризация.} Составление таблиц полностью аналогично случаю $\mathbb{C}1$-цепи.

\medskip

{\bf Восстановление кода.}
В преобразовании $9$, зная код $F$, мы можем вычислить и $J$ и $Z$.

\medskip

\subsection{Случай цепи $\mathbb{A}0$; преобразования 7 и 10}

\begin{figure}[hbtp]
\centering
\includegraphics[width=0.8\textwidth]{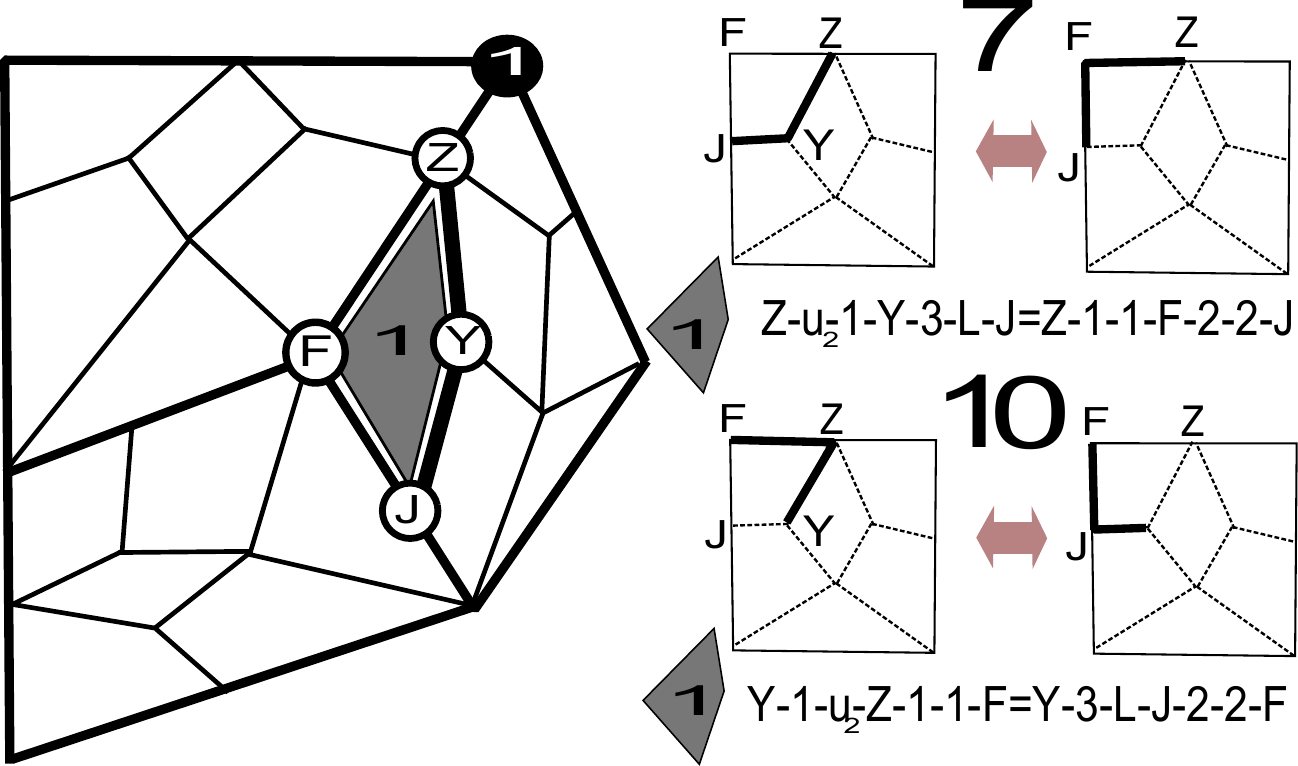}
\caption{Случаи расположения пути вокруг A0-цепи и соответствующие им определяющие соотношения для локальных преобразований 7 и 10}
\label{A1a}
\end{figure}

В правой части рисунка~\ref{A1a} изображены локальные преобразования $7$ и $10$.

\medskip

Зафиксируем вершину типа $\mathbb{A}$ и ее начальника (вершина, отмеченная черным кругом с ``$1$'', уровень которой равен $2$).
Заметим, что зная окружение этой вершины, мы можем вычислить коды вершин $Z$, $Y$, $F$, $J$ с точностью до подклееных окружений. Это позволяет ввести определяющие отношения, записанные в правой части рисунка~\ref{A1a}.

\medskip

{\bf Характеризация.} Составление таблиц полностью аналогично случаю $\mathbb{C}1$-цепи.

\medskip

{\bf Восстановление кода.}
Вершина, отмеченная черным кругом с $1$, является единственным начальником  $Z$, $F$, $J$.
Таким образом, зная код пути $ZYJ$, можно вычислить код пути $ZFJ$, и наоборот. А также зная код пути $FZY$, можно вычислить код пути $FJY$, и наоборот.

\medskip

\subsection{Случай цепи $\mathbb{A}1$; преобразования 8 и 9}

\begin{figure}[hbtp]
\centering
\includegraphics[width=1\textwidth]{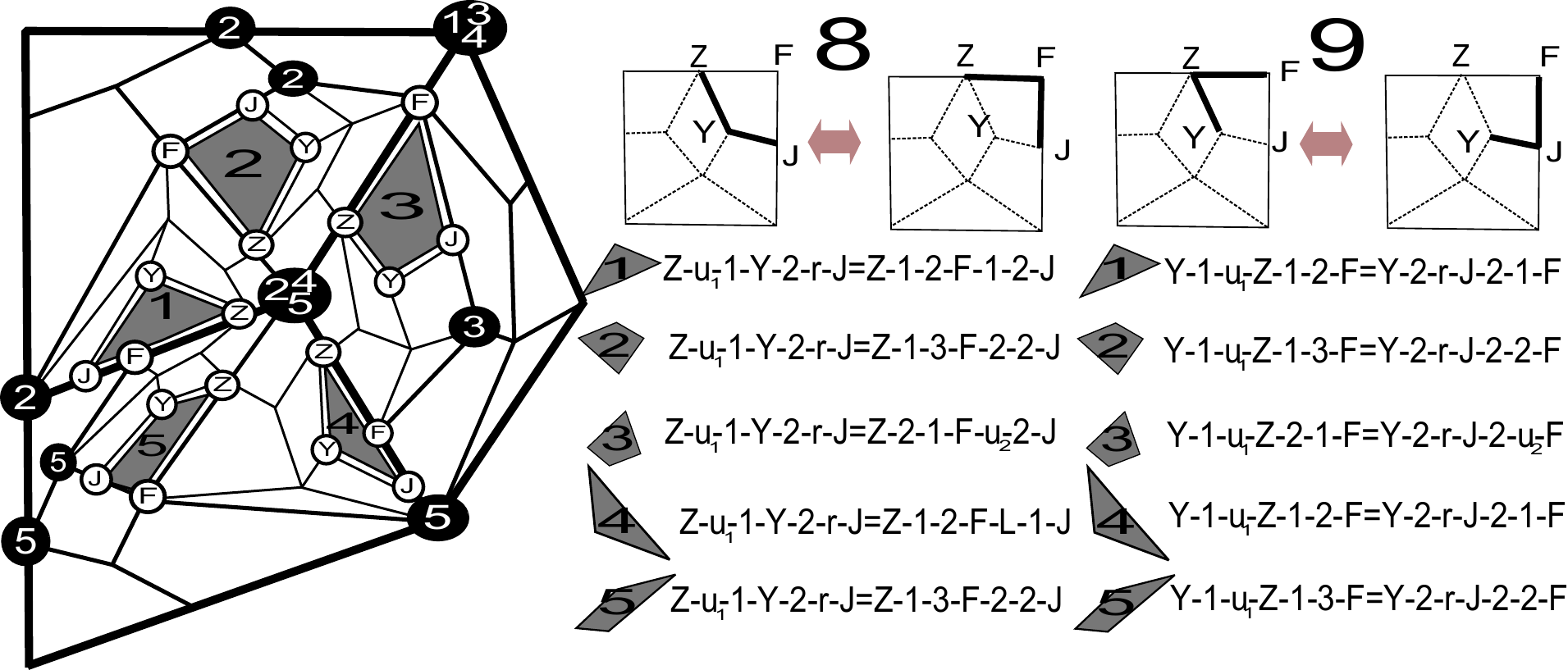}
\caption{Случаи расположения пути вокруг A1-цепи  и соответствующие им определяющие соотношения для локальных преобразований 8 и 9}
\label{A2b}
\end{figure}

В правой части рисунка~\ref{A2b} изображены локальные преобразования $8$ и $9$.

\medskip

Зафиксируем вершину типа $\mathbb{A}$ и ее начальника (третьего уровня).
Черными кругами отметим вершины, являющиеся начальниками вершин  $Z$, $Y$, $F$, $J$.
Заметим, что зная окружение центрального $\mathbb{A}$-узла, мы можем найти их типы, уровни и окружения. Так как начальники вершин  $Z$, $Y$, $F$, $J$ содержатся среди вершин, отмеченных черными кругами, то коды вершин $Z$, $Y$, $F$, $J$ во всех случаях мы можем назвать явно, с точностью до подклееных окружений. Это позволяет ввести определяющие отношения, записанные в правой части рисунка~\ref{A2b}.

\medskip

{\bf Характеризация.} Составление таблиц полностью аналогично случаю $\mathbb{C}1$-цепи.

\medskip

{\bf Восстановление кода.}
Поскольку мы можем установить, с каким именно случаем расположения мы имеем дело, то и окружение каждой из вершин мы можем вычислить, зная остальные три. Кроме того, заметим, что каждая из четырех вершин  $Z$, $J$, $Y$, $F$ в каждом из случаев либо имеет начальников среди других трех вершин, либо у нее общие начальники с какой-то из трех других вершин. Значит, в каждом из случаев, зная коды трех вершин, можно вычислить код четвертой.

Таким образом, зная код пути $ZYJ$, можно вычислить код пути $ZFJ$, и наоборот. А также зная код пути $FZY$, можно вычислить код пути $FJY$, и наоборот.

\medskip

\subsection{Случай цепи $\mathbb{A}1$; преобразования 7 и 10}

\begin{figure}[hbtp]
\centering
\includegraphics[width=1\textwidth]{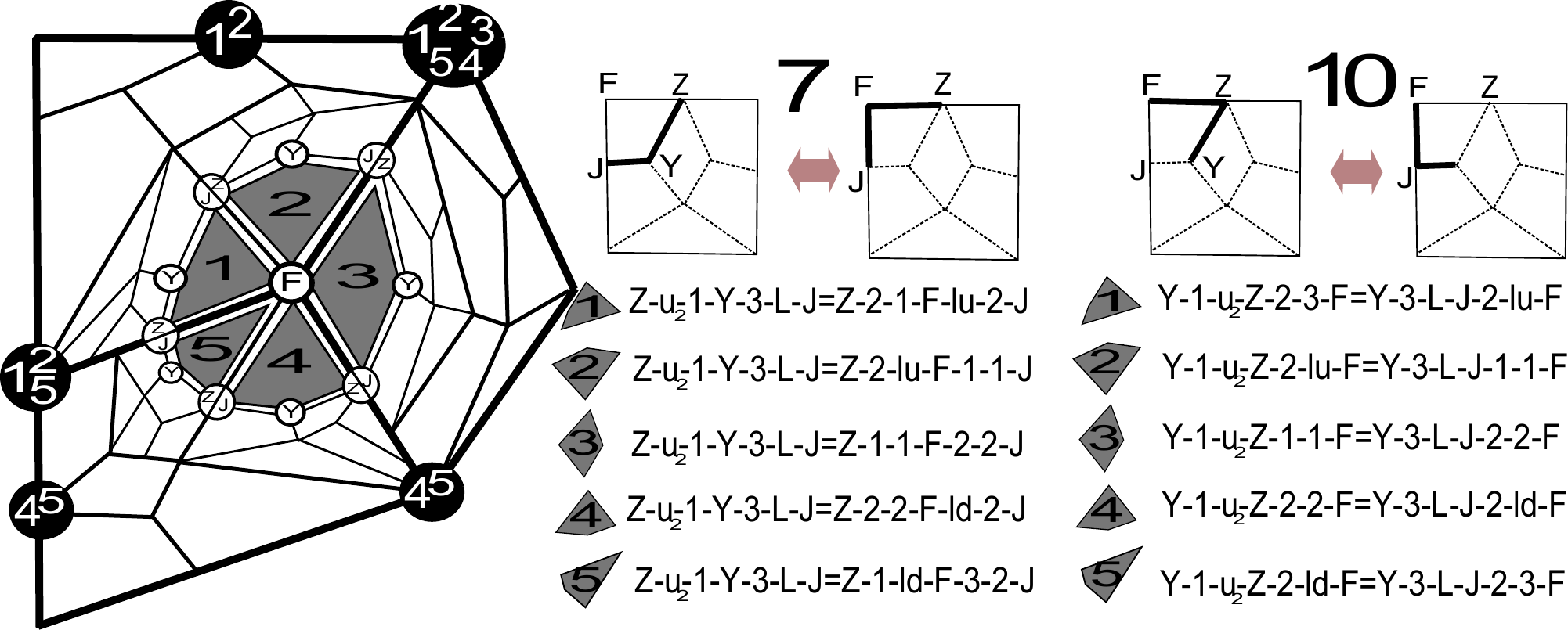}
\caption{Случаи расположения пути вокруг A1-цепи и соответствующие им определяющие соотношения для локальных преобразований 7 и 10}
\label{A2a}
\end{figure}

В правой части рисунка~\ref{A2a} изображены локальные преобразования $7$ и $10$.

\medskip

Фиксируем некоторую вершину типа $\mathbb{A}$ и ее начальника.
Черными  кругами отметим вершины, являющиеся начальниками вершин  $Z$, $Y$, $F$, $J$.
Заметим, что зная окружение этой центральной $\mathbb{A}$-вершины, мы можем вычислить коды вершин $Z$, $Y$, $F$, $J$ с точностью до подклееных окружений. Это позволяет ввести определяющие отношения, записанные в правой части рисунка~\ref{A2a}.

\medskip

{\bf Характеризация.} Составление таблиц полностью аналогично случаю $\mathbb{C}1$-цепи.

\medskip

{\bf Восстановление кода.}
Вершина $F$ имеет общего начальника с одной из вершин $J$ или $Z$. То есть, зная код пути $ZYJ$, можно вычислить код пути $ZFJ$. В обратную сторону, а также для локального преобразования 10: очевидно, что зная $F$, можно вычислить коды остальных вершин. Таким образом, мы можем осуществить локальное преобразование пути через операцию с его кодом.

\medskip

\subsection{Случай цепи $\mathbb{A}2$; преобразования 8 и 9}

\begin{figure}[hbtp]
\centering
\includegraphics[width=1\textwidth]{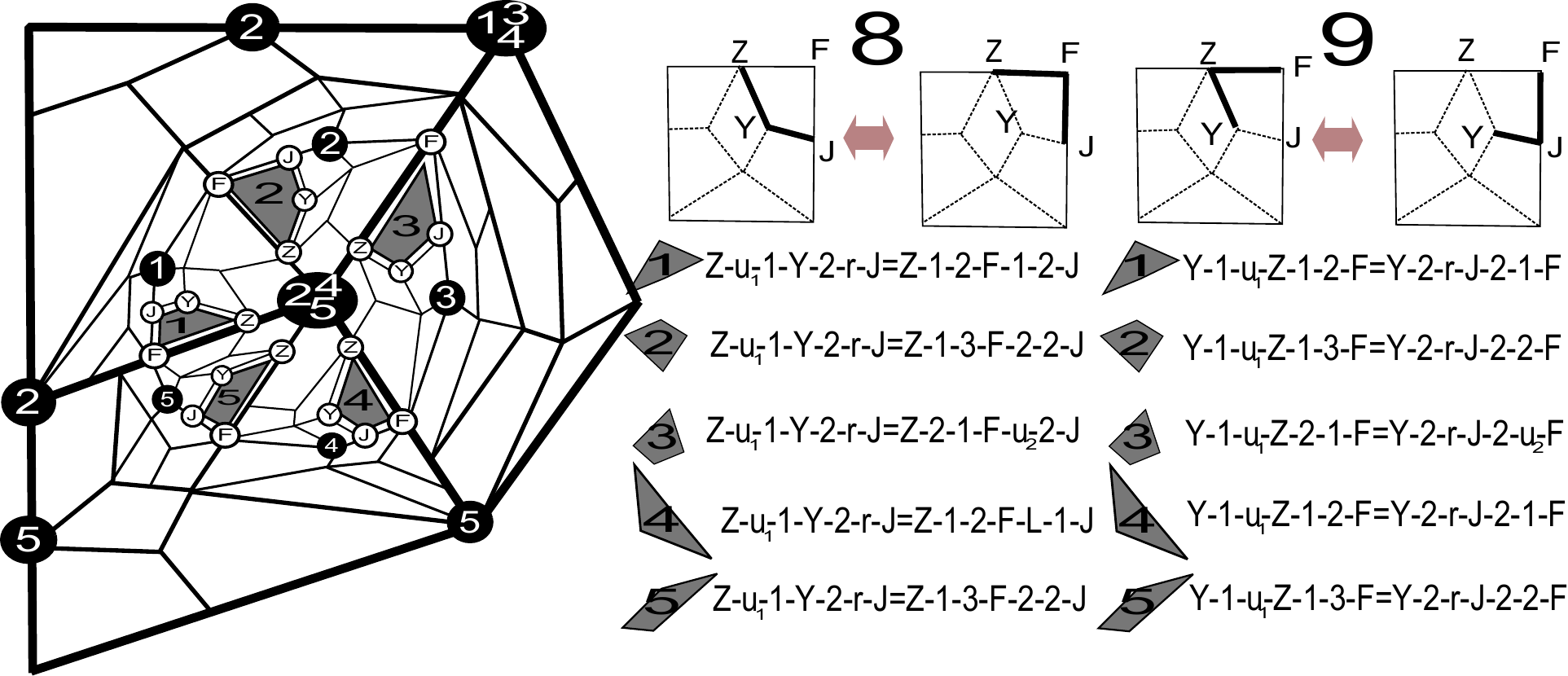}
\caption{Случаи расположения пути вокруг A2-цепи  и соответствующие им определяющие соотношения для локальных преобразований 8 и 9}
\label{A3b}
\end{figure}

В правой части рисунка~\ref{A3b} изображены локальные преобразования $8$ и $9$.

\medskip

Введение определяющих соотношений полностью аналогично случаю с $\mathbb{A}1$-цепью. Также аналогично можно установить, что зная код пути $ZYJ$, можно вычислить код пути $ZFJ$, и наоборот. А также зная код пути $FZY$, можно вычислить код пути $FJY$, и наоборот. Аналогично можно установить последовательности букв, характеризующие данную цепь и локальные преобразования.

\medskip

\subsection{Случай цепи $\mathbb{A}2$; преобразования 7 и 10}

\begin{figure}[hbtp]
\centering
\includegraphics[width=1\textwidth]{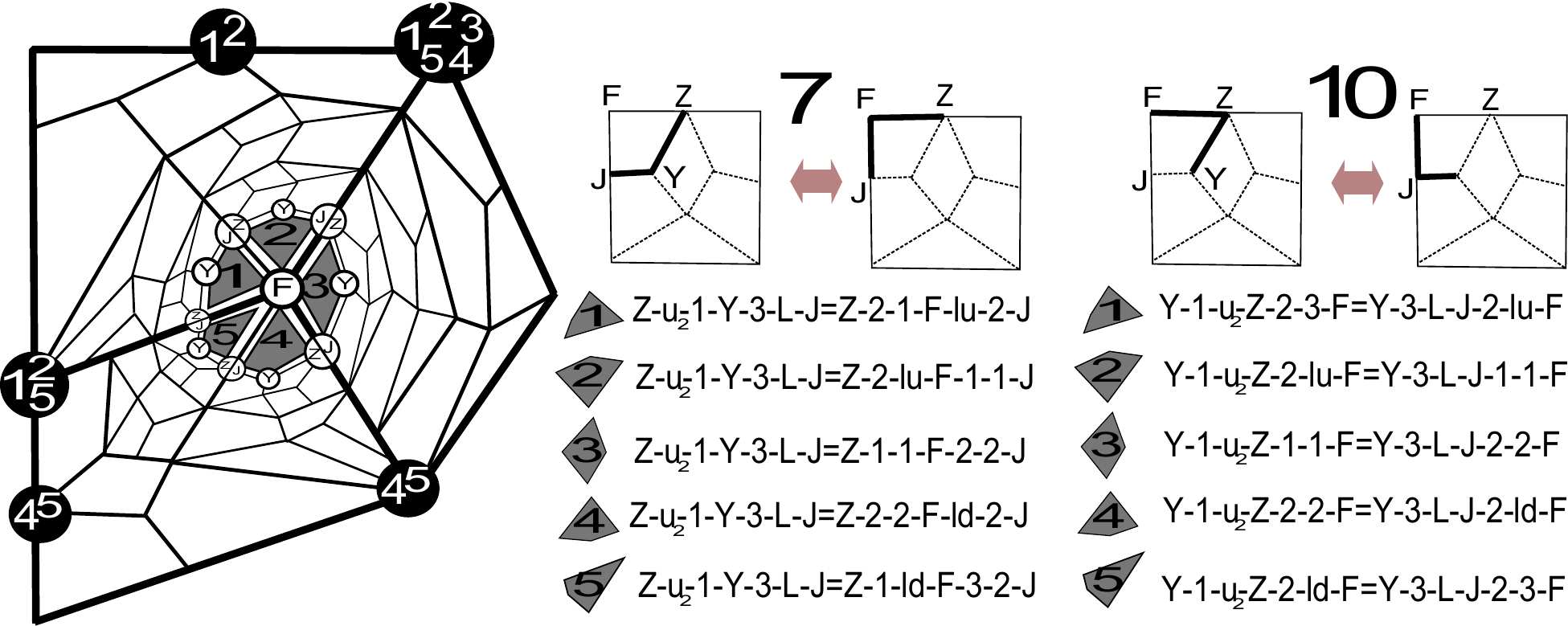}
\caption{Случаи расположения пути вокруг A2-цепи  и соответствующие им определяющие соотношения для локальных преобразований 7 и 10}
\label{A3a}
\end{figure}

В правой части рисунка~\ref{A3a} изображены локальные преобразования $7$ и $10$.

\medskip

Введение определяющих соотношений полностью аналогично случаю с $\mathbb{A}1$-цепью. Также аналогично можно установить, что зная код пути $ZYJ$, можно вычислить код пути $ZFJ$, и наоборот. А также зная код пути $FZY$, можно вычислить код пути $FJY$, и наоборот. Аналогично можно установить последовательности букв, характеризующие данную цепь и локальные преобразования.

\medskip

\subsection{Случай цепи $\mathbb{UL}1$; преобразования 8 и 9}

\begin{figure}[hbtp]
\centering
\includegraphics[width=1\textwidth]{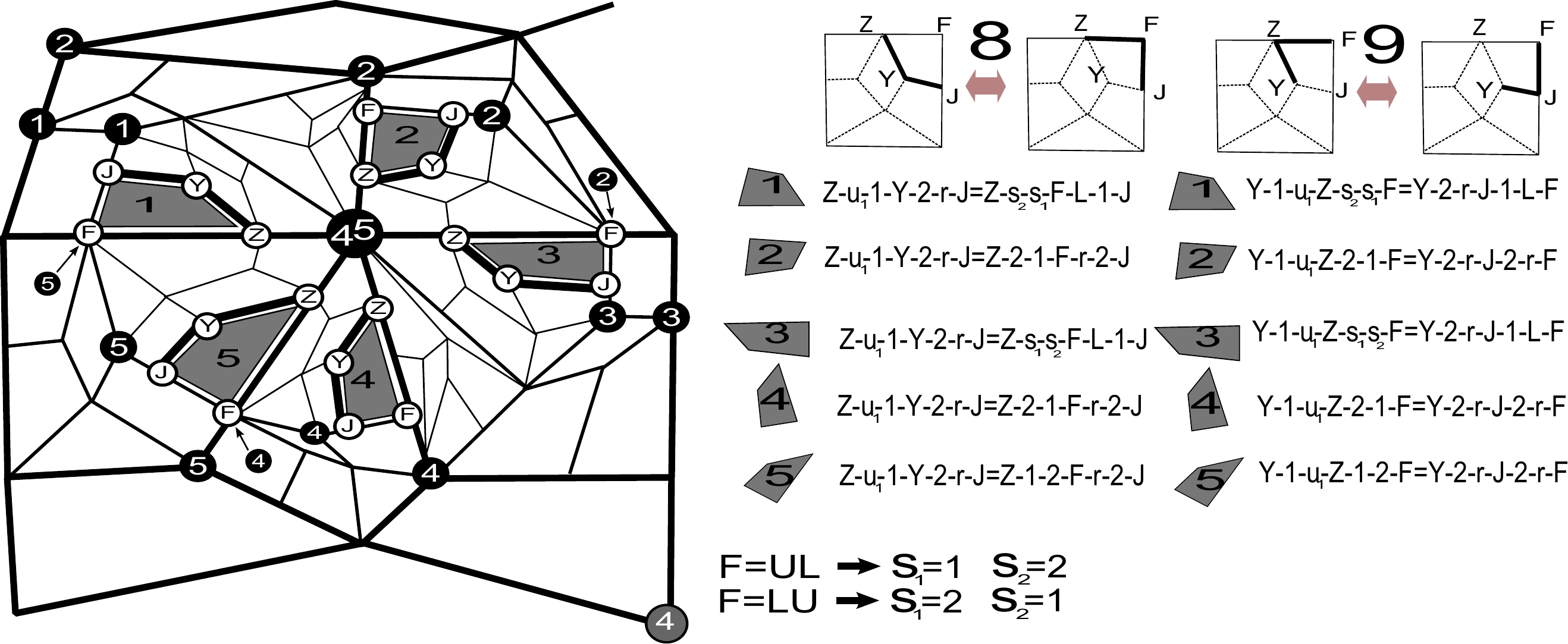}
\caption{Случаи расположения пути вокруг UL1-цепи и определяющие соотношения для  локальных преобразований 8 и 9}
\label{UL1b}
\end{figure}

В правой части рисунка~\ref{UL1b} изображены локальные преобразования $8$ и $9$.

\medskip

Зафиксируем вершину типа $\mathbb{UL}$ третьего уровня и ее начальников.
Черными кругами отметим вершины, являющиеся начальниками вершин  $Z$, $Y$, $F$, $J$. Три круга попадают в $F$-узлы других расположений, это отмечено стрелками.
Заметим, что зная окружение центрального $\mathbb{UL}$-узла, мы можем найти типы, уровни и окружения вершин с черными кругами (кроме вершины в правом нижнем углу, для нее только тип). Так как начальники вершин  $Z$, $Y$, $F$, $J$ содержатся среди вершин, отмеченных черными кругами, то коды вершин $Z$, $Y$, $F$, $J$ во всех случаях мы можем назвать явно, с точностью до подклееных окружений. Это позволяет ввести определяющие отношения, записанные в правой части рисунка~\ref{UL1b}.

\medskip

{\bf Характеризация.} Составление таблиц полностью аналогично случаю $\mathbb{C}1$-цепи.

\medskip

{\bf Восстановление кода.}
Поскольку мы можем установить, с каким именно случаем расположения мы имеем дело, то окружение оставшейся вершины мы можем вычислить, зная коды остальных трех. Кроме того, заметим, что у $Z$ и $F$ в каждом из случаев общие начальники, то есть код каждой из этих вершин может быть вычислен, исходя из кодов остальных трех вершин. У $Y$ во всех случаях первым начальником является $Z$, а тип второго в каждом случае ясен из расположения.

Начальники $J$ вычисляется следующим образом:  случай $1$ -- это $\mathbf{Prev}(F)$;  случай $2$ -- это $1$-цепь вокруг $\mathbf{TypeBottomLeft.FBoss}(Z)$, указатель $\mathbf{l}$ (первый начальник), и $\mathbb{A}$-тип второго начальника; случай $3$ -- это $\mathbf{Prev}(F)$;
 случай $4$ -- это $0$-цепь с указателем $1$ вокруг $\mathbb{A}$, окружения которого как $U$-часть $\mathbf{Fboss}(Z)$ (первый начальник), и тип второго начальника $\mathbb{B}$;
 случай $5$ -- это первый $\mathbf{Plus.FBoss}(Z)$, тип второго $\mathbb{A}$.

Таким образом, зная код пути $ZYJ$, можно вычислить код пути $ZFJ$, и наоборот. А также зная код пути $FZY$, можно вычислить код пути $FJY$, и наоборот. Таким образом, мы можем осуществить локальное преобразование пути через операцию с его кодом.

\medskip

\subsection{Случай цепи $\mathbb{UL}1$; преобразования 7 и 10}

\begin{figure}[hbtp]
\centering
\includegraphics[width=1\textwidth]{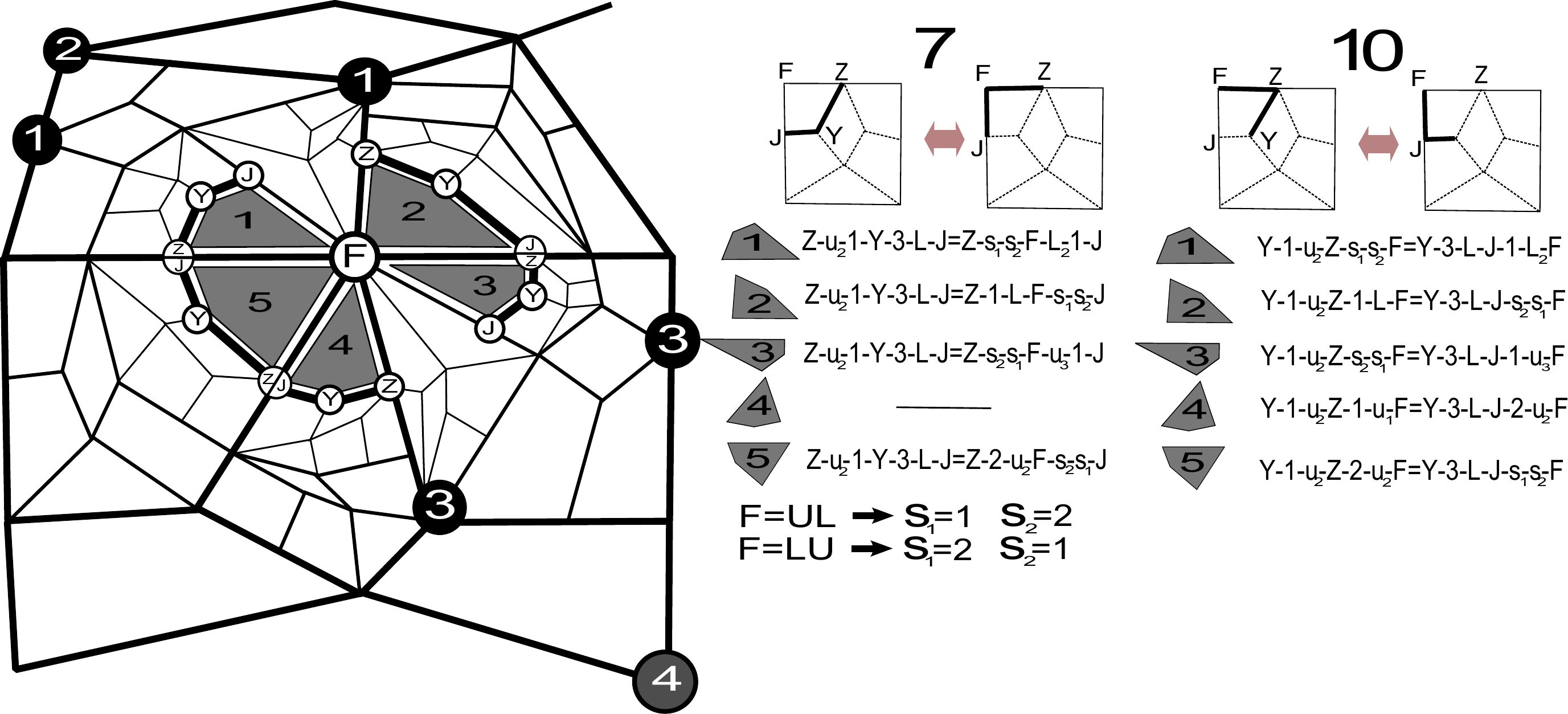}
\caption{Случаи расположения пути вокруг UL1-цепи и определяющие соотношения для локальных преобразований 7 и 10}
\label{UL1a}
\end{figure}

В правой части рисунка~\ref{UL1a} изображены локальные преобразования $7$ и $10$.

\medskip

Фиксируем некоторую вершину типа $\mathbb{UL}$ третьего уровня и ее начальников.
Черными кругами отметим вершины, являющиеся начальниками вершин  $Z$, $Y$, $F$, $J$.
Заметим, что зная окружение центральной вершины (которая совпадает с $F$), мы можем вычислить коды вершин $Z$, $Y$, $F$, $J$ с точностью до подклееных окружений во всех случаях. Это позволяет ввести определяющие отношения, записанные в правой части рисунка~\ref{UL1a}.

\medskip

{\bf Характеризация.} Составление таблиц полностью аналогично случаю $\mathbb{C}1$-цепи.

\medskip

{\bf Восстановление кода.}
Вершина $F$ имеет общего начальника с одной из вершин $J$ или $Z$, в каждом из случаев, кроме случая $4$ (а в  случае $4$ нам ее не надо вычислять).  То есть, зная код пути $ZYJ$, можно вычислить код пути $ZFJ$. В обратную сторону, а также для локального преобразования $10$: очевидно, что зная $F$, можно вычислить коды остальных вершин. Таким образом, мы можем осуществить локальное преобразование пути через операцию с его кодом.

\medskip

\subsection{Случай цепи $\mathbb{UL}2$; преобразования 8 и 9}

\begin{figure}[hbtp]
\centering
\includegraphics[width=0.9\textwidth]{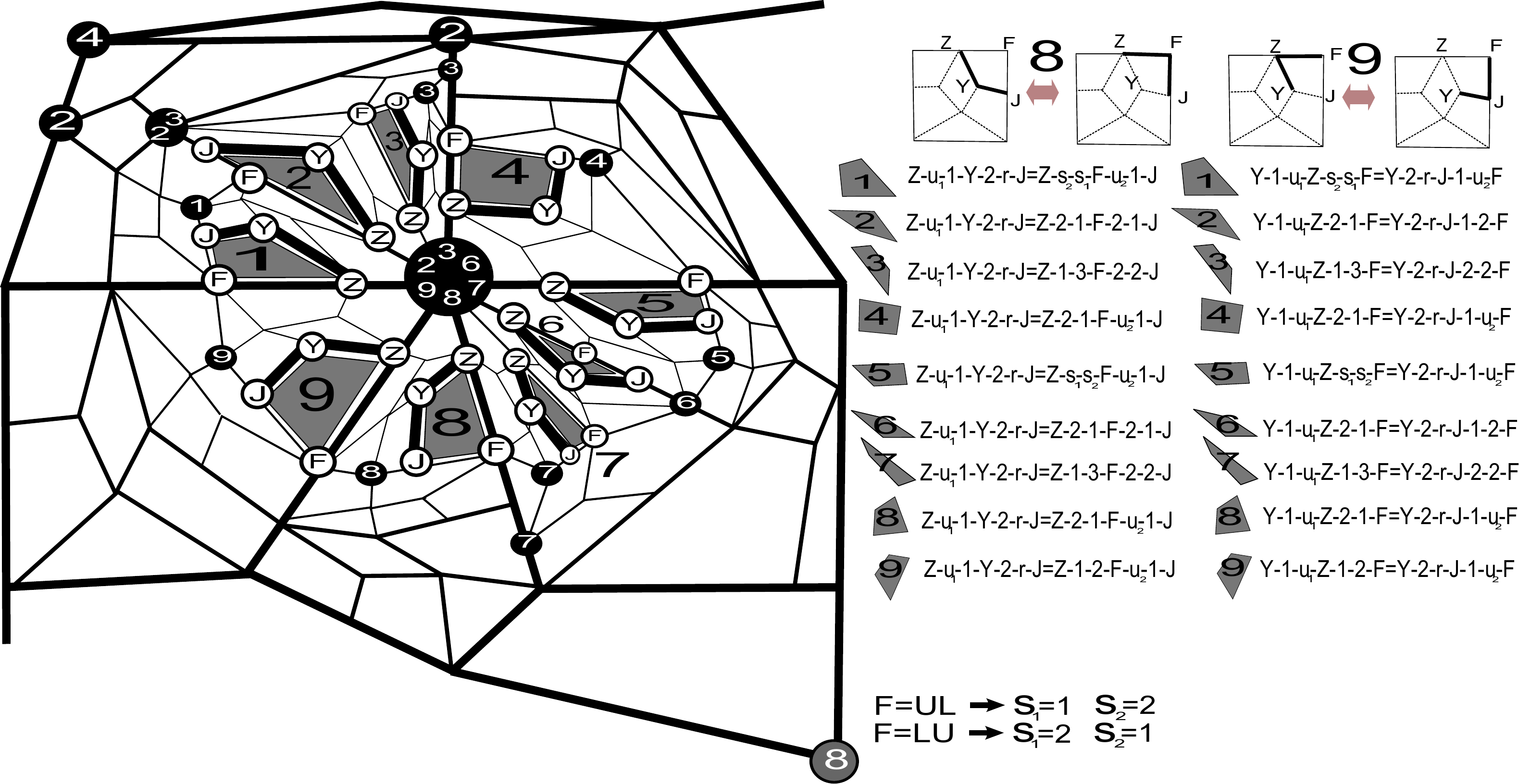}
\caption{Случаи расположения пути вокруг UL2-цепи и определяющие соотношения для локальных преобразований 8 и 9}
\label{UL2b}
\end{figure}

В правой части рисунка~\ref{UL2b} изображены локальные преобразования $8$ и $9$.

\medskip

\leftskip=-0.0cm

Зафиксируем вершину типа $\mathbb{UL}$ третьего уровня и ее начальников.
Черными кругами отметим вершины, являющиеся начальниками вершин  $Z$, $Y$, $F$, $J$.
Заметим, что зная окружение центрального $\mathbb{UL}$-узла, мы можем найти их типы, уровни и окружения. Так как начальники вершин  $Z$, $Y$, $F$, $J$ содержатся среди вершин, отмеченных черными кругами, то коды вершин $Z$, $Y$, $F$, $J$, во всех случаях мы можем назвать явно, с точностью до подклееных окружений. Это позволяет ввести определяющие отношения, записанные в правой части рисунка~\ref{UL2b}.

\medskip

{\bf Характеризация.} Составление таблиц полностью аналогично случаю $\mathbb{C}1$-цепи.

\medskip

{\bf Восстановление кода.}
Заметим, что у $Z$ и $F$ в каждом из случаев общие начальники, то есть код каждой из этих вершин может быть вычислен, исходя из кодов остальных трех вершин. $Y$ вычисляется также легко, так как во всех случаях первым начальником является $Z$, а тип второго в каждом случае легко виден из рисунка.

Начальники $J$ вычисляется следующим образом: в  случаях $1$, $4$, $5$, $8$, $9$ -- единственный начальник, узел $F$; в случаях $2$, $6$ -- начальники те же что и у $F$; в случаях  $3$ и $7$ -- первый и второй начальники как первый и третий у $F$.

Таким образом, зная код пути $ZYJ$, можно вычислить код пути $ZFJ$, и наоборот. А также зная код пути $FZY$, можно вычислить код пути $FJY$, и наоборот. Таким образом, мы можем осуществить локальное преобразование пути через операцию с его кодом.

\medskip

\subsection{Случай цепи $\mathbb{UL}2$; преобразования 7 и 10}

\begin{figure}[hbtp]
\centering
\includegraphics[width=1\textwidth]{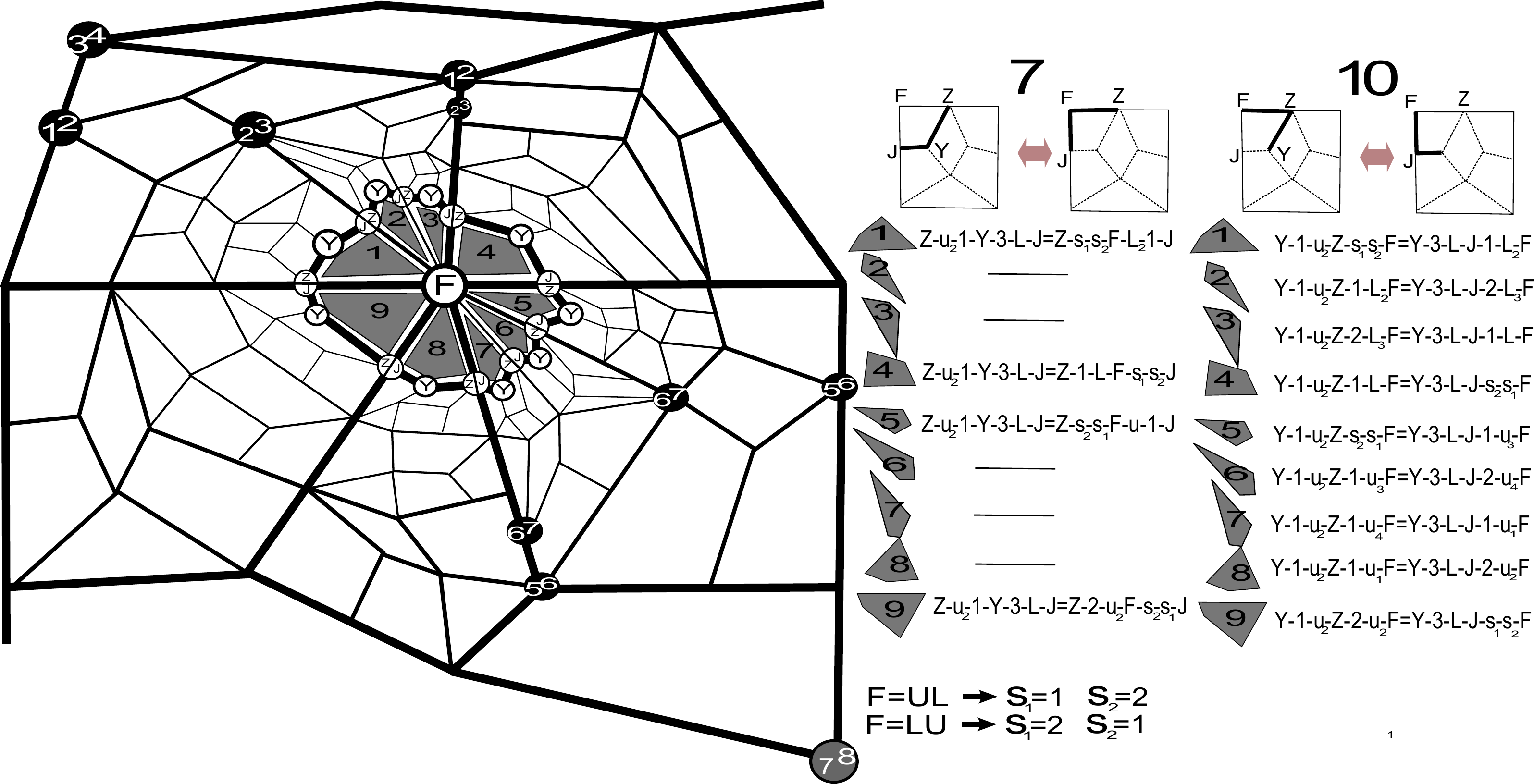}
\caption{Случаи расположения пути вокруг UL2-цепи и определяющие соотношения для локальных преобразований 7 и 10}
\label{UL2a}
\end{figure}

В правой части рисунка~\ref{UL2a} изображены локальные преобразования $7$ и $10$.

\medskip

Фиксируем некоторую вершину типа $\mathbb{UL}$ третьего уровня и ее начальников.
Черными и серыми кругами отметим вершины, являющиеся начальниками вершин  $Z$, $Y$, $F$, $J$.
Заметим, что зная окружение центральной вершины (совпадающей с $F$), мы можем вычислить коды вершин $Z$, $Y$, $F$, $J$, с точностью до подклееных окружений. Это позволяет ввести определяющие отношения, записанные в правой части рисунка~\ref{UL2a}.

Для  локального преобразования $7$ пути в случаях $2,3,6,7,8$ удовлетворяют условиям мертвого паттерна, и в этих случаях мы соотношения не вводим.

\medskip

{\bf Характеризация.} Составление таблиц полностью аналогично случаю $\mathbb{C}1$-цепи.

\medskip

{\bf Восстановление кода.}
В случаях $1,4,5,9$ вершина $F$ имеет общего начальника с одной из вершин $J$ или $Z$. То есть, зная код пути $ZYJ$, можно вычислить код пути $ZFJ$. В обратную сторону, а также для локального преобразования $10$: очевидно, что зная $F$, можно вычислить коды остальных вершин. Таким образом, мы можем осуществить локальное преобразование пути через операцию с его кодом.

\medskip

\subsection{Случай цепи $\mathbb{UL}3$; преобразования 7, 8, 9, 10}

Случай $\mathbb{UL}3$-цепи полностью аналогичен $\mathbb{UL}2$ случаю, соотношения выглядят идентично, только кодировки вершин $J$, $F$, $Z$, $Y$ отвечают $\mathbb{UL}3$-цепи. Все рассуждения о вычислении путей полностью аналогичны. Соотношений вводится столько же, сколько для $\mathbb{UL}2$ случая.

\medskip

\subsection{Случай цепи $\mathbb{UR}1$; преобразования 8 и 9}

\begin{figure}[hbtp]
\centering
\includegraphics[width=1\textwidth]{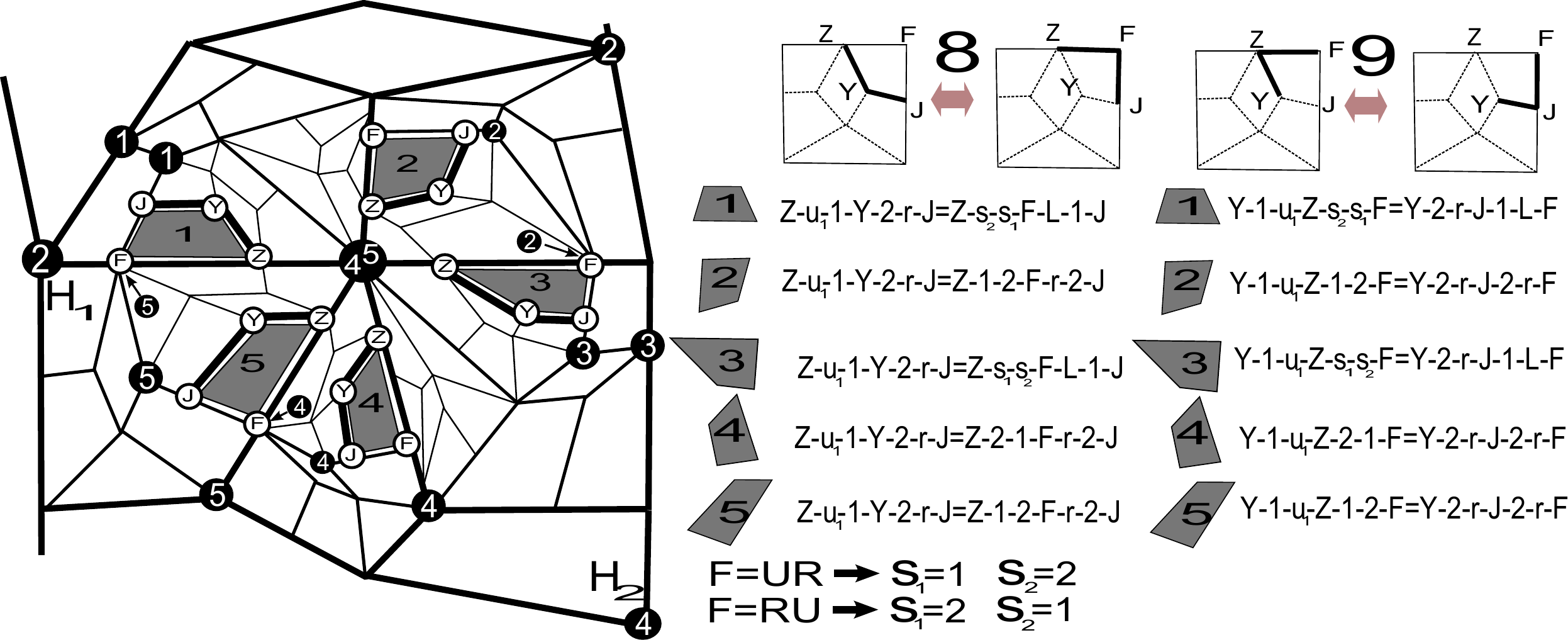}
\caption{Случаи расположения пути вокруг UR1-цепи и определяющие соотношения для локальных преобразований 8 и 9}
\label{UR1b}
\end{figure}

В правой части рисунка~\ref{UR1b} изображены локальные преобразования $8$ и $9$.

\medskip

Зафиксируем вершину типа $\mathbb{UR}$ третьего уровня и ее начальников.
Черными  кругами отметим вершины, являющиеся начальниками вершин  $Z$, $Y$, $F$, $J$. Три круга попадают в $F$-узлы других расположений, это отмечено стрелками.
Заметим, что зная окружение центрального $\mathbb{UR}$-узла, мы можем найти их типы, уровни и окружения (Кроме вершин $H_1$ и $H_2$, для которых можем найти тип). Так как начальники вершин  $Z$, $Y$, $F$, $J$ содержатся среди вершин, отмеченных черными кругами, то коды вершин $Z$, $Y$, $F$, $J$, во всех случаях мы можем назвать явно, с точностью до подклееных окружений. Это позволяет ввести определяющие отношения, записанные в правой части рисунка~\ref{UR1b}.

\medskip

{\bf Характеризация.} Составление таблиц полностью аналогично случаю $\mathbb{C}1$-цепи.

\medskip

{\bf Восстановление кода.}
Заметим, что у $Z$ и $F$ в каждом из случаев общие начальники, то есть код каждой из этих вершин может быть вычислен, исходя из кодов остальных трех вершин. $Y$ вычисляется также легко, так как во всех случаях первым начальником является $Z$, а тип второго в каждом случае виден из рисунка.

Начальники $J$ вычисляется следующим образом: случай $1$ -- это $\mathbf{Prev}(F)$; случай $2$ -- это $\mathbf{UpRightChain.FBoss}(Z)$ и $\mathbb{B}$-тип второго начальника; случай $3$ -- это $\mathbf{Prev}(F)$; случай $4$ -- это $0$-цепь с указателем $1$ вокруг $\mathbb{A}$, окружения которого как $U$-часть $\mathbf{Fboss}(Z)$, и тип второго начальника $\mathbb{B}$;
случай $5$ -- это $\mathbf{Plus.FBoss}(Z)$ и тип второго начальника $\mathbb{B}$.

Таким образом, зная код пути $ZYJ$, можно вычислить код пути $ZFJ$, и наоборот. А также зная код пути $FZY$, можно вычислить код пути $FJY$, и наоборот. Таким образом, мы можем осуществить локальное преобразование пути через операцию с его кодом.

\medskip

\subsection{Случай цепи $\mathbb{UR}1$; преобразования 7 и 10}

\begin{figure}[hbtp]
\centering
\includegraphics[width=1\textwidth]{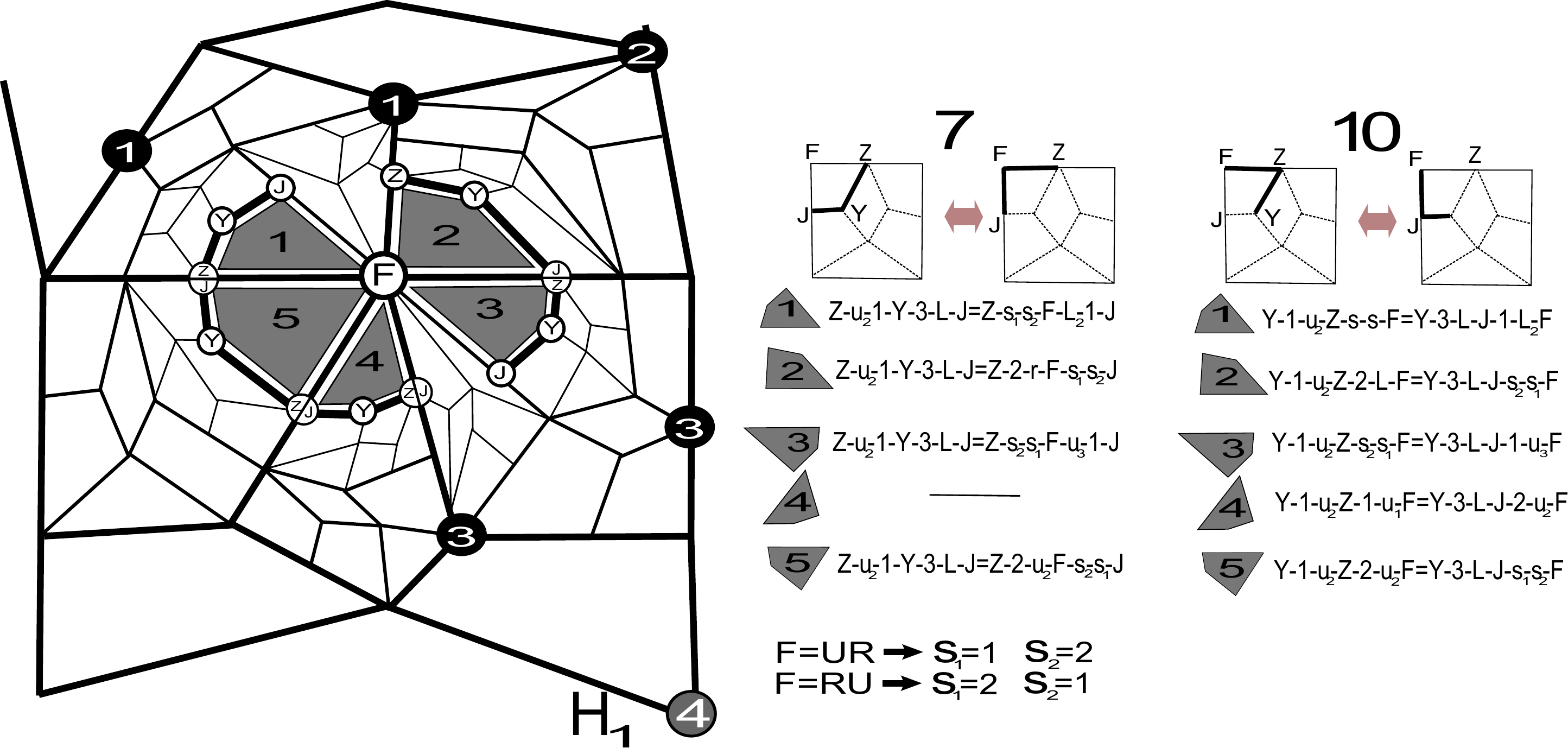}
\caption{Случаи расположения пути вокруг UR1-цепи и определяющие соотношения для локальных преобразований 7 и 10}
\label{UR1a}
\end{figure}

В правой части рисунка~\ref{UR1a} изображены локальные преобразования $7$ и $10$.

\medskip

Фиксируем некоторую вершину типа $\mathbb{UR}$ третьего уровня и ее начальников.
Черными кругами отметим вершины, являющиеся начальниками вершин  $Z$, $Y$, $F$, $J$.
Заметим, что зная окружение центральной вершины, мы можем вычислить коды вершин $Z$, $Y$, $F$, $J$, с точностью до подклееных окружений. Это позволяет ввести определяющие отношения, записанные в правой части рисунка~\ref{UR1a}.

\medskip

{\bf Характеризация.} Составление таблиц полностью аналогично случаю $\mathbb{C}1$-цепи.

\medskip

{\bf Восстановление кода.}
Вершина $F$ имеет общего начальника с одной из вершин $J$ или $Z$, в каждом из случаев, кроме случая $4$ (а в случае $4$ нам ее не надо вычислять). Таким образом, зная код пути $ZYJ$, можно вычислить код пути $ZFJ$. В обратную сторону, а также для  локального преобразования $10$: очевидно, что зная $F$, можно вычислить коды остальных вершин. То есть, мы можем осуществить локальное преобразование пути через операцию с его кодом.

\medskip

\subsection{Случай цепи $\mathbb{UR}2$; преобразования 8 и 9}

\begin{figure}[hbtp]
\centering
\includegraphics[width=1\textwidth]{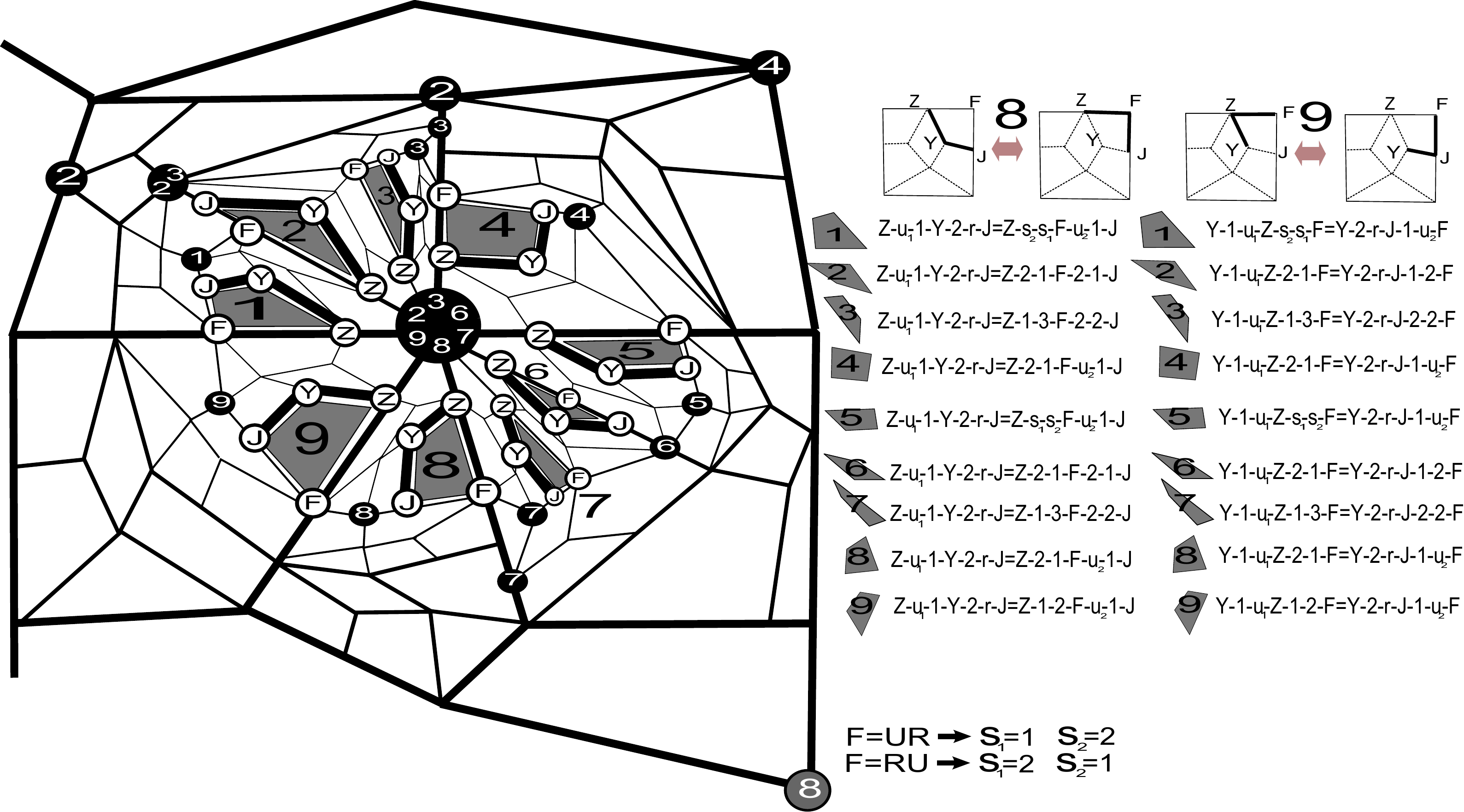}
\caption{Случаи расположения пути вокруг UR2-цепи и определяющие соотношения для локальных преобразований 8 и 9}
\label{UR2b}
\end{figure}

В правой части рисунка~\ref{UR2b} изображены локальные преобразования $8$ и $9$.

\medskip

Фиксируем вершину типа $\mathbb{UR}$ третьего уровня и ее начальников.
Черными кругами отметим вершины, являющиеся начальниками вершин  $Z$, $Y$, $F$, $J$.
Заметим, что зная окружение центрального $\mathbb{UR}$-узла, мы можем найти их типы, уровни и окружения (для вершины в правом нижнем углу -- только тип). Так как начальники вершин  $Z$, $Y$, $F$, $J$ содержатся среди вершин, отмеченных черными кругами, то коды вершин $Z$, $Y$, $F$, $J$, во всех случаях мы можем назвать явно, с точностью до подклееных окружений. Это позволяет ввести определяющие отношения, записанные в правой части рисунка~\ref{UR2b}.

\medskip

{\bf Характеризация.} Составление таблиц полностью аналогично случаю $\mathbb{C}1$-цепи.

\medskip

{\bf Восстановление кода.}
Заметим, что у $Z$ и $F$ в каждом из случаев общие начальники, то есть код каждой из этих вершин может быть вычислен, исходя из кодов остальных трех вершин. $Y$ вычисляется также легко, так как во всех случаях первым начальником является $Z$, а тип второго в каждом случае ясен.

Начальники $J$ вычисляется следующим образом: в  случаях $1$, $4$, $5$, $8$, $9$ -- единственный начальник, узел $F$; в  случаях $2$, $6$ -- начальники те же что и у $F$; в случаях $3$ и $7$ -- первый и второй начальники как первый и третий у $F$.

Таким образом, зная код пути $ZYJ$, можно вычислить код пути $ZFJ$, и наоборот. А также зная код пути $FZY$, можно вычислить код пути $FJY$, и наоборот. Таким образом, мы можем осуществить локальное преобразование пути через операцию с его кодом.

\medskip

\subsection{Случай цепи $\mathbb{UR}2$; преобразования 7 и 10}

\begin{figure}[hbtp]
\centering
\includegraphics[width=1\textwidth]{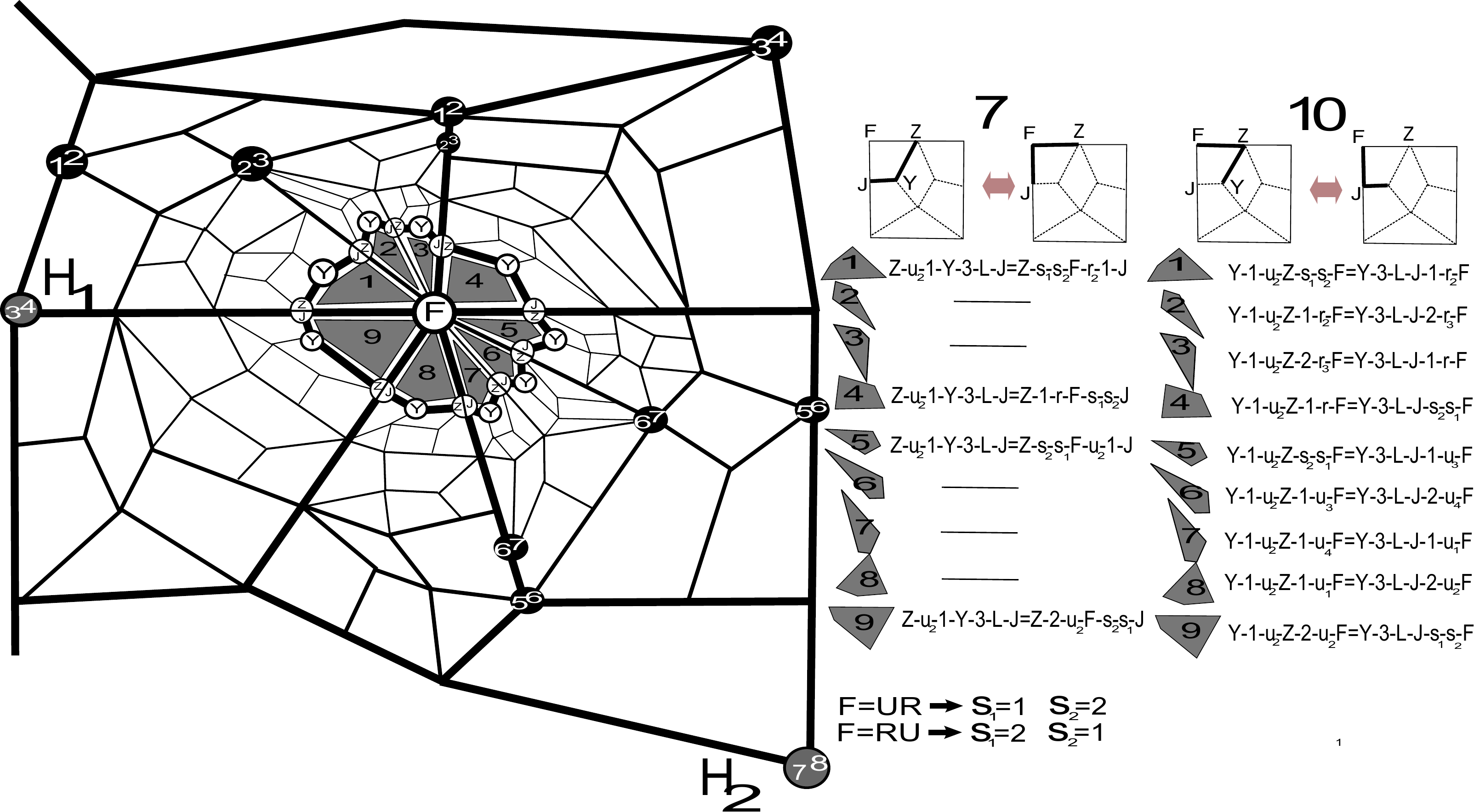}
\caption{Случаи расположения пути вокруг UR2-цепи и определяющие соотношения для локальных преобразований 7 и 10}
\label{UR2a}
\end{figure}

В правой части рисунка~\ref{UR2a} изображены локальные преобразования $7$ и $10$.

\medskip

Фиксируем некоторую вершину типа $\mathbb{UR}$ и ее начальников.
Черными кругами отметим вершины, являющиеся начальниками вершин  $Z$, $Y$, $F$, $J$.
Заметим, что зная окружение этой вершины, мы можем вычислить коды вершин $Z$, $Y$, $F$, $J$, с точностью до подклееных окружений (кроме вершин $H_1$ и $H_2$, для которых можно найти тип). Это позволяет ввести определяющие отношения, записанные в правой части рисунка~\ref{UR2a}.

Для локального преобразования $7$ в случаях $2$, $3$, $6$, $7$, $8$  пути удовлетворяют условиям мертвого паттерна, и в этих случаях мы соотношения не вводим.

\medskip

{\bf Характеризация.} Составление таблиц полностью аналогично случаю $\mathbb{C}1$-цепи.

\medskip

{\bf Восстановление кода.}
В случаях $1,4,5,9$ вершина $F$ имеет общего начальника с одной из вершин $J$ или $Z$. То есть, зная код пути $ZYJ$, можно вычислить код пути $ZFJ$. В обратную сторону, а также для  локального преобразования $10$: очевидно, что зная $F$, можно вычислить коды остальных вершин. Таким образом, мы можем осуществить локальное преобразование пути через операцию с его кодом.

\medskip

\subsection{Случай цепи $\mathbb{UR}3$; преобразования 7, 8, 9, 10}

Случай $\mathbb{UR}3$-цепи полностью аналогичен $\mathbb{UR}2$ случаю, соотношения выглядят идентично, только кодировки вершин $J$, $F$, $Z$, $Y$ отвечают $\mathbb{UR}3$-цепи. Все рассуждения о вычислении путей полностью аналогичны. Соотношений вводится столько же, сколько для $\mathbb{UR}2$ случая.

\medskip

\subsection{Случай цепи $\mathbb{DR}1$; преобразования 8 и 9}

\begin{figure}[hbtp]
\centering
\includegraphics[width=0.5\textwidth]{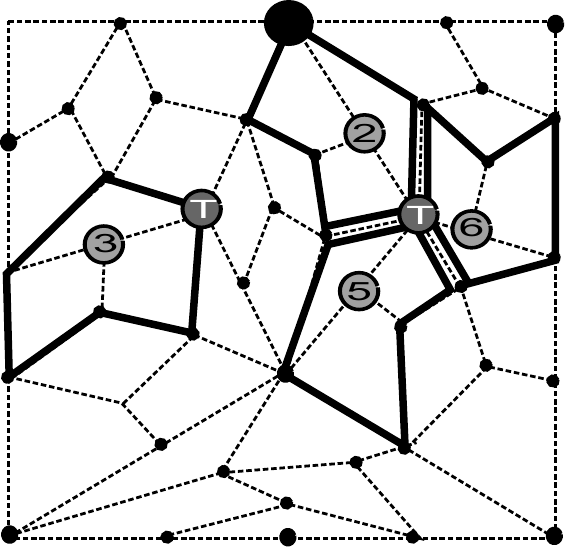}
\caption{Случаи расположения $\mathbb{DR}$-вершины}
\label{RDplacement}
\end{figure}

Аналогично предыдущим случаям, мы можем ввести определяющие соотношения. Сначала зафиксируем вершину типа $\mathbb{DR}$ третьего уровня и ее начальников. Всего возможно четыре случая расположения $\mathbb{DR}$-вершины, на четырех внутренних ребрах (типов $2$, $3$, $5$, $6$). Все эти случаи расположения показаны на рисунке~\ref{RDplacement}.

\medskip

\begin{figure}[hbtp]
\centering
\includegraphics[width=1\textwidth]{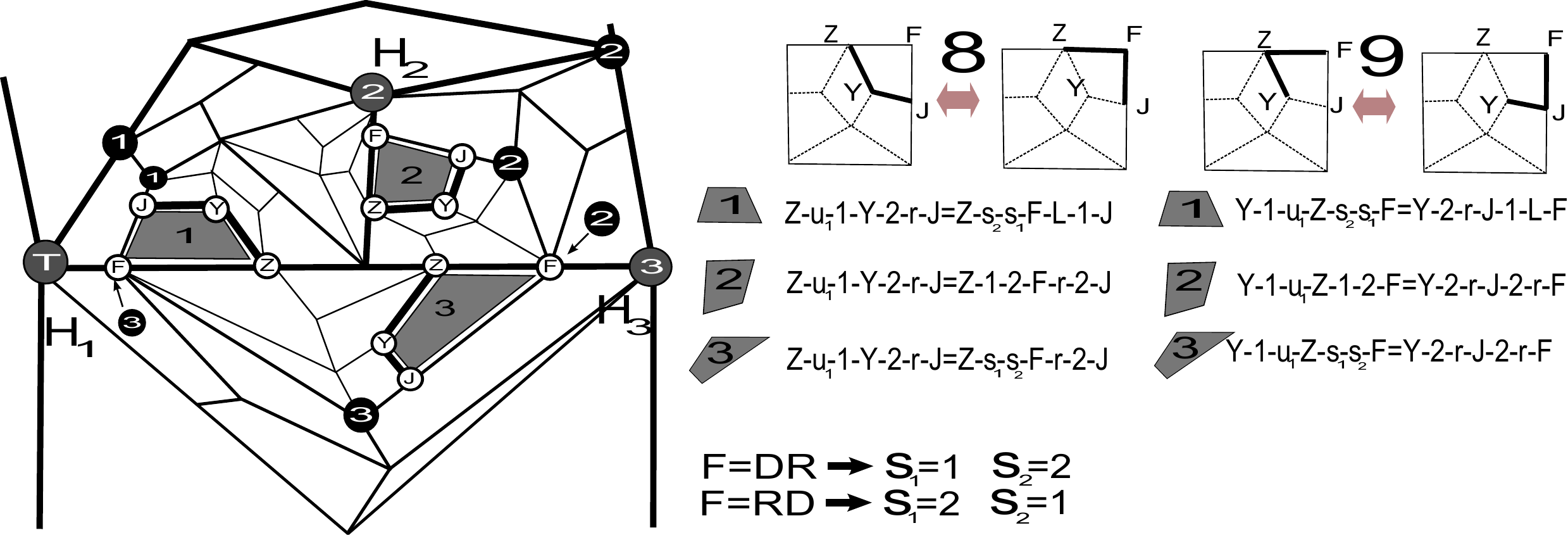}
\caption{Случаи расположения пути вокруг $\mathbb{DR}$1-цепи  и определяющие соотношения для локальных преобразований 8 и 9}
\label{DR1b}
\end{figure}

Первым начальником нашей $\mathbb{DR}$-вершины во всех случаях будет вершина в середине верхней стороны. Теперь непосредственно рассмотрим цепь с центром в этой $\mathbb{DR}$-вершине.
В правой части рисунка~\ref{DR1b} изображены локальные преобразования $8$ и $9$.

\medskip

Черными кругами отметим вершины, являющиеся начальниками вершин  $Z$, $Y$, $F$, $J$.  Два круга попадают в $F$-узлы других расположений, это отмечено стрелками.
Зная окружение центрального $\mathbb{DR}$-узла, мы можем найти их типы, уровни и окружения (для вершин $H_1$ $H_2$, $H_3$ можем найти только тип).

Тип узла, отмеченного как $T$, мы можем определить исходя из типа внутреннего ребра, на котором лежит центральный $\mathbb{DR}$-узел. Если это ребро $3$ -- то $\mathbb{A}$, в остальных случаях ($2$, $5$ или $6$) -- $\mathbb{B}$.

Так как начальники вершин  $Z$, $Y$, $F$, $J$ содержатся среди вершин, отмеченных черными кругами, то коды вершин $Z$, $Y$, $F$, $J$, во всех случаях мы можем назвать явно, с точностью до подклееных окружений. Это позволяет ввести определяющие отношения, записанные в правой части рисунка~\ref{DR1b}.

\medskip

{\bf Характеризация.} Составление таблиц полностью аналогично случаю $\mathbb{C}1$-цепи.

\medskip

{\bf Восстановление кода.}
Поскольку мы можем установить, с каким именно случаем расположения мы имеем дело, то окружение каждой вершины мы можем вычислить, зная код остальных трех.
Кроме того, заметим, что у $Z$ и $F$ в каждом из случаев общие начальники, то есть код каждой из этих вершин может быть вычислен, исходя из кодов остальных трех вершин. $Y$ вычисляется также легко, так как во всех случаях первым начальником является $Z$, а тип правого нижнего угла в каждом случае ясен.

Начальники $J$ вычисляется следующим образом:  случай $1$ -- это $\mathbf{Prev}(F)$;  случай $2$ -- первый начальник -- $\mathbf{UpRightChain.FBoss}(Z)$ и $\mathbb{B}$-тип второго;  случай $3$ -- это $1$-цепь вокруг узла $T$ c указателями в зависимости от типа ребра центрального $\mathbb{DR}$-узла: $3$ для $3$ типа, $1$ для $2$ типа, $2$ для $6$ типа, $3$ для $5$ типа.

Таким образом, зная код пути $ZYJ$, можно вычислить код пути $ZFJ$, и наоборот. А также зная код пути $FZY$, можно вычислить код пути $FJY$, и наоборот. Таким образом, мы можем осуществить локальное преобразование пути через операцию с его кодом.

\medskip

\subsection{Случай цепи $\mathbb{DR}1$; преобразования 7 и 10}

\begin{figure}[hbtp]
\centering
\includegraphics[width=1\textwidth]{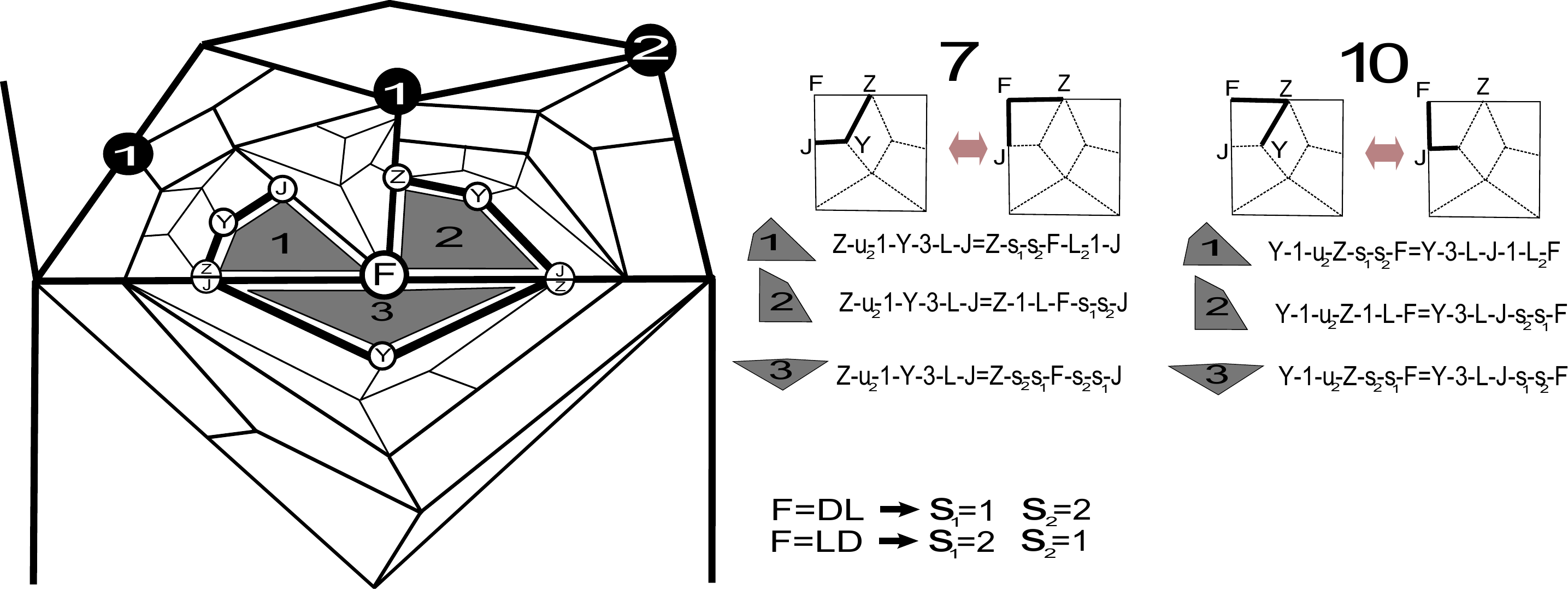}
\caption{Случаи расположения пути вокруг DR1-цепи и определяющие соотношения для локальных преобразований ~7 и 10}
\label{DR1a}
\end{figure}

В правой части рисунка~\ref{DR1a} изображены локальные преобразования $7$ и $10$.

\medskip

Фиксируем некоторую вершину типа $\mathbb{DR}$ третьего уровня и ее начальников.
Черными  кругами отметим вершины, являющиеся начальниками вершин  $Z$, $Y$, $F$, $J$.
Заметим, что зная окружение центральной $\mathbb{DR}$-вершины, мы можем в каждом случае вычислить коды вершин $Z$, $Y$, $F$, $J$ с точностью до подклееных окружений. Это позволяет ввести определяющие отношения, записанные в правой части рисунка~\ref{DR1a}.

\medskip

{\bf Характеризация.} Составление таблиц полностью аналогично случаю $\mathbb{C}1$-цепи.

\medskip

{\bf Восстановление кода.}
Вершина $F$ имеет общего начальника с одной из вершин $J$ или $Z$, в каждом из случаев.  То есть, зная код пути $ZYJ$, можно вычислить код пути $ZFJ$. В обратную сторону, а также для $10$ локального преобразования: очевидно, что зная $F$, можно вычислить коды остальных вершин. Таким образом, мы можем осуществить локальное преобразование пути через операцию с его кодом.

\medskip

\subsection{Случай цепи $\mathbb{DR}2$; преобразования 8 и 9}

\begin{figure}[hbtp]
\centering
\includegraphics[width=0.9\textwidth]{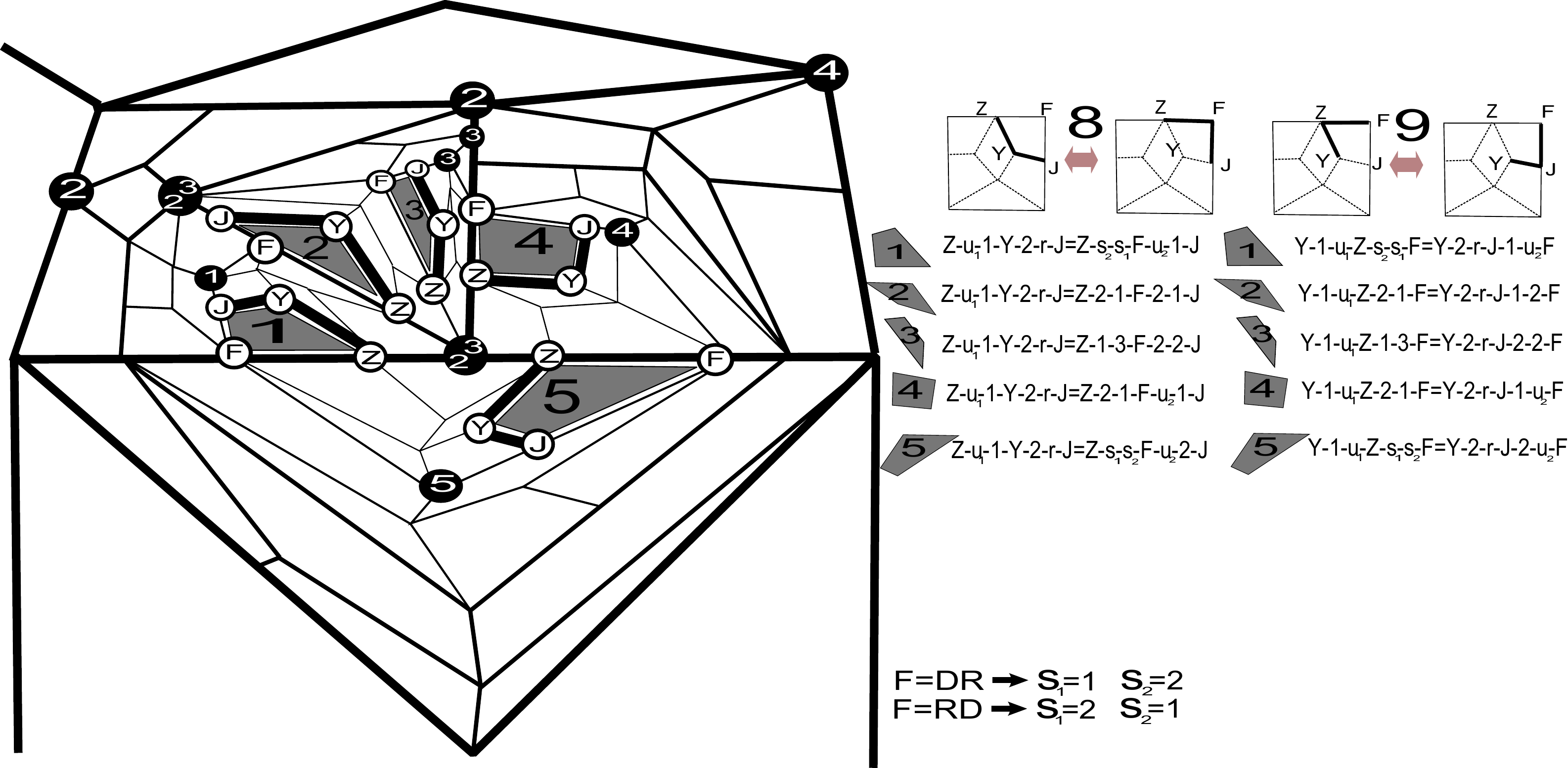}
\caption{Случаи расположения пути вокруг DR2-цепи  и определяющие соотношения для локальных преобразований 8 и 9}
\label{DR2b}
\end{figure}

В правой части рисунка~\ref{DR2b} изображены локальные преобразования $8$ и $9$.

\medskip

Зафиксируем вершину типа $\mathbb{DR}$ третьего уровня и ее начальников.
Черными кругами отметим вершины, являющиеся начальниками вершин  $Z$, $Y$, $F$, $J$.
Заметим, что зная окружение центрального $\mathbb{DR}$-узла, мы можем найти их типы, уровни и окружения. Так как начальники вершин  $Z$, $Y$, $F$, $J$ содержатся среди вершин, отмеченных черными кругами, то коды вершин $Z$, $Y$, $F$, $J$, во всех случаях мы можем назвать явно, с точностью до подклееных окружений. Это позволяет ввести определяющие отношения, записанные в правой части рисунка~\ref{DR2b}.

\medskip

{\bf Характеризация.} Составление таблиц полностью аналогично случаю $\mathbb{C}1$-цепи.

\medskip

{\bf Восстановление кода.}
Заметим, что у $Z$ и $F$ в каждом из случаев общие начальники, то есть код каждой из этих вершин может быть вычислен, исходя из кодов остальных трех вершин. $Y$ вычисляется также легко, так как во всех случаях первым начальником является $Z$, а тип второго в каждом случае ясен.

Начальники $J$ вычисляется следующим образом: в  случаях $1$, $4$ и $5$ -- единственный начальник, узел $F$; в случае $2$ -- начальники те же что и у $F$; в  случае $3$ -- первый и второй начальники как первый и третий у $F$.

Таким образом, зная код пути $ZYJ$, можно вычислить код пути $ZFJ$, и наоборот. А также зная код пути $FZY$, можно вычислить код пути $FJY$, и наоборот. Таким образом, мы можем осуществить локальное преобразование пути через операцию с его кодом.

\medskip

\subsection{Случай цепи $\mathbb{DR}2$; преобразования 7 и 10}

\begin{figure}[hbtp]
\centering
\includegraphics[width=1\textwidth]{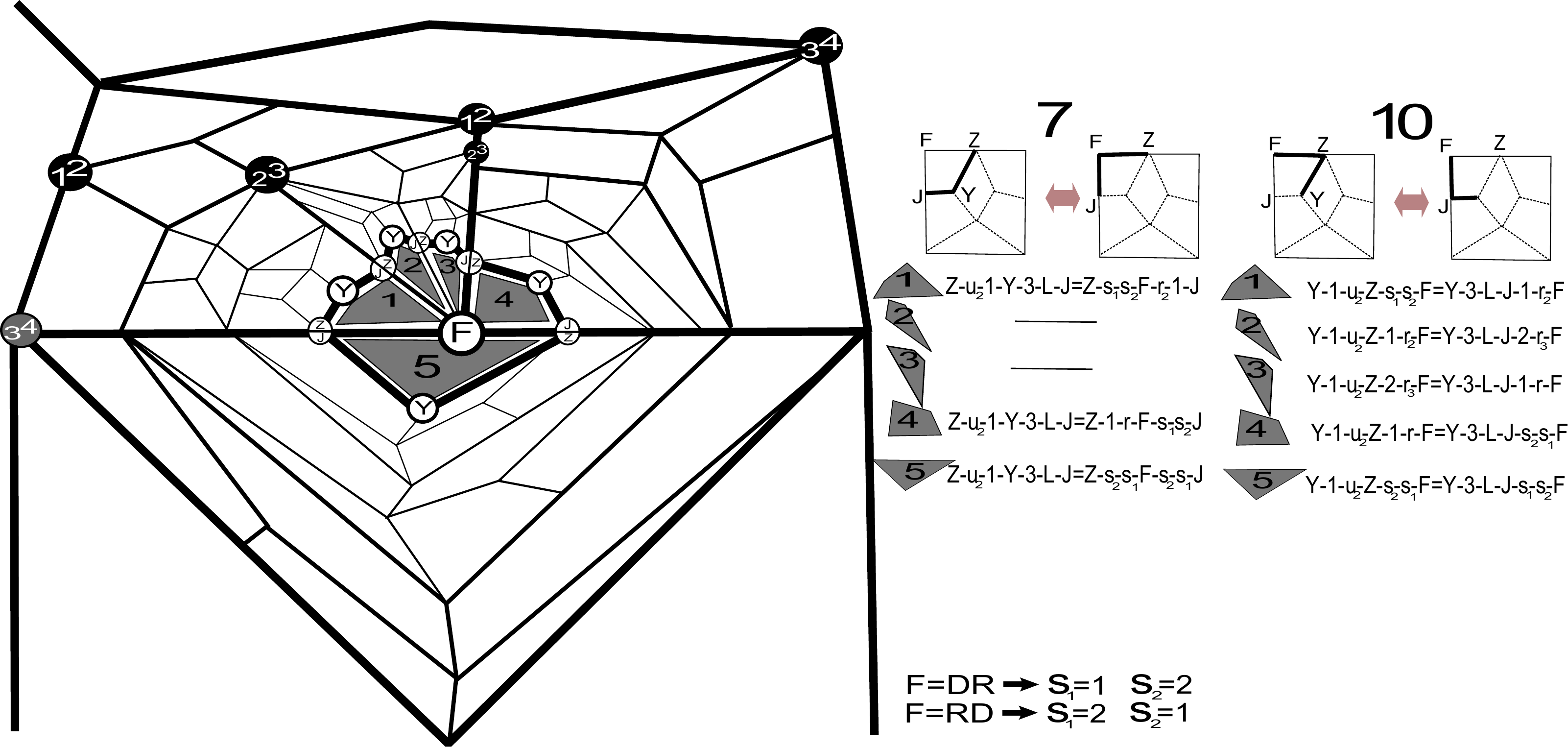}
\caption{Случаи расположения пути вокруг DR2-цепи и определяющие соотношения для локальных преобразований 7 и 10}
\label{DR2a}
\end{figure}

В правой части рисунка~\ref{DR2a} изображены локальные преобразования $7$ и $10$.

\medskip

Зафиксируем некоторую вершину типа $\mathbb{DR}$ третьего уровня и ее начальников.
Черными кругами отметим вершины, являющиеся начальниками вершин  $Z$, $Y$, $F$, $J$. Заметим, что зная окружение этой вершины, мы можем вычислить коды вершин $Z$, $Y$, $F$, $J$, с точностью до подклееных окружений. Это позволяет ввести определяющие отношения, записанные в правой части рисунка~\ref{DR2a}.

Для  локального преобразования $7$ в случаях $2$, $3$, $6$, $7$, $8$ пути удовлетворяют условиям мертвого паттерна, и в этих случаях мы соотношения не вводим.

\medskip

{\bf Характеризация.} Составление таблиц полностью аналогично случаю $\mathbb{C}1$-цепи.

\medskip

{\bf Восстановление кода.}
В случаях $1,4,5$ вершина $F$ имеет общего начальника с одной из вершин $J$ или $Z$. То есть, $ZYJ$, можно вычислить код пути $ZFJ$. В обратную сторону, а также для  локального преобразования $10$: очевидно, что зная $F$, можно вычислить коды остальных вершин. Таким образом, мы можем осуществить локальное преобразование пути через операцию с его кодом.

\medskip

\subsection{Случай цепи $\mathbb{DR}3$; преобразования 7, 8, 9, 10}

Случай $\mathbb{DR}3$-цепи полностью аналогичен $\mathbb{DR}2$ случаю, соотношения выглядят идентично, только кодировки вершин $J$, $F$, $Z$, $Y$ отвечают $\mathbb{DR}3$-цепи. Все рассуждения о вычислении путей полностью аналогичны. Соотношений вводится столько же, сколько для $\mathbb{DR}2$ случая.

\medskip

\subsection{Случай цепи $\mathbb{DL}1$; преобразования 8 и 9}

\begin{figure}[hbtp]
\centering
\includegraphics[width=0.9\textwidth]{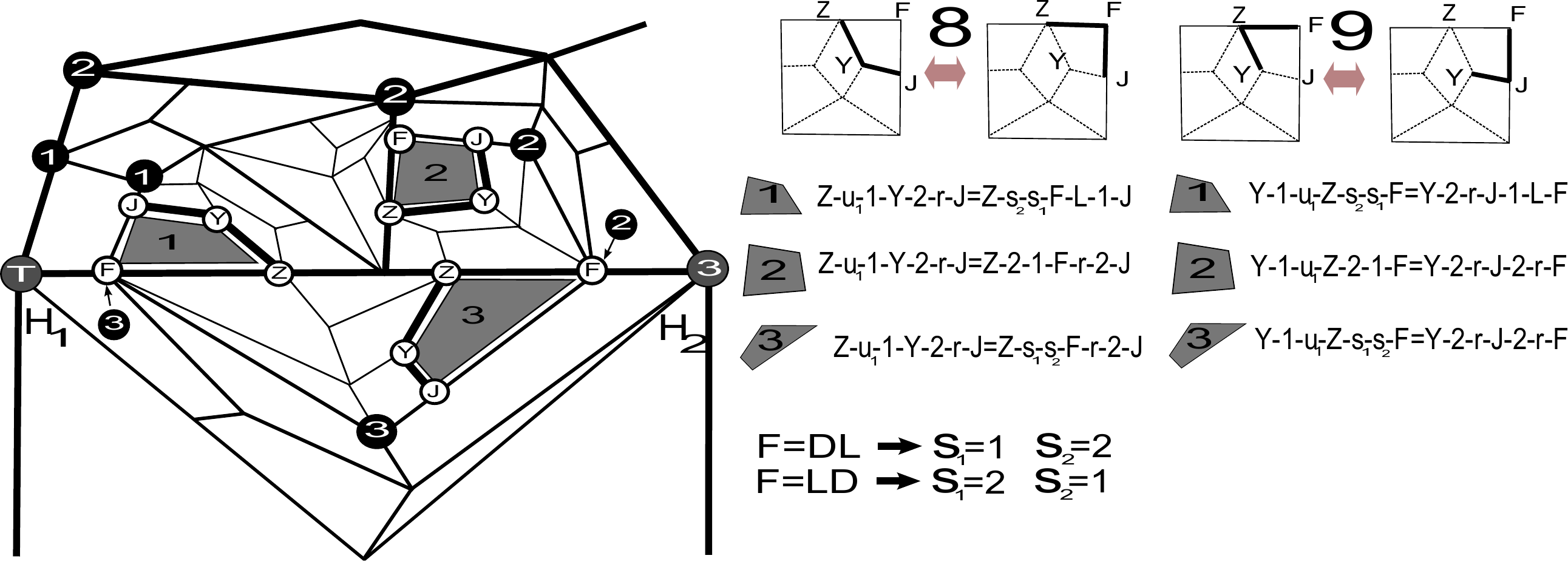}
\caption{Случаи расположения пути вокруг DL1-цепи  и определяющие соотношения для локальных преобразований 8 и 9}
\label{DL1b}
\end{figure}

В правой части рисунка~\ref{DL1b} изображены локальные преобразования $8$ и $9$.

\medskip

Зафиксируем вершину типа $\mathbb{DL}$ третьего уровня и ее начальников.
Черными и серыми кругами отметим вершины, являющиеся начальниками вершин  $Z$, $Y$, $F$, $J$. Два круга попадают в $F$-узлы других расположений, это отмечено стрелками.
Числа в них обозначают, начальниками вершин какого случая является данная вершина. Зная окружение центрального $\mathbb{DL}$-узла, мы можем найти их типы, уровни и окружения (для вершин $H_1$ и $H_2$ -- тип).

Так как начальники вершин  $Z$, $Y$, $F$, $J$ содержатся среди вершин, отмеченных черными кругами, то коды вершин $Z$, $Y$, $F$, $J$, во всех случаях мы можем назвать явно, с точностью до подклееных окружений. Это позволяет ввести определяющие отношения, записанные в правой части рисунка~\ref{DL1b}.

\medskip

{\bf Характеризация.} Составление таблиц полностью аналогично случаю $\mathbb{C}1$-цепи.

\medskip

{\bf Восстановление кода.}
Заметим, что у $Z$ и $F$ в каждом из случаев общие начальники, то есть код каждой из этих вершин может быть вычислен, исходя из кодов остальных трех вершин. $Y$ вычисляется также легко, так как во всех случаях первым начальником является $Z$, а тип второго в каждом случае ясен.

Начальники $J$ вычисляется следующим образом:  случай $1$ -- $\mathbf{Prev}(F)$;  случай $2$ --  первый -- $\mathbf{BottomLeftChain.FBoss}(Z)$ и $\mathbb{A}$-тип второго. Для  случая $3$ заметим, что тип узла, отмеченного символом $T$ может быть установлен, в зависимости от типа ребра, на котором лежат $Z$ и $F$. Так как это ребро, с центром в  $\mathbb{DL}$-вершине, то тип его может быть либо $4$ либо $7$, в соответствии со структурой разбиения макроплитки на подплитки. Если ребро имеет $4$ тип, то вершина $T$ имеет тип $\mathbb{A}$, а первый начальник $J$ -- это вершина из $\mathbb{A}1$ цепи с указателем $2$. Если же ребро имеет тип $7$, то $T$-вершина является вторым начальником $F$ и $Z$, а первый начальник $J$ -- это вершина из $2$-цепи вокруг него с указателем, соответствующим $7$ ребру.

Таким образом, зная код пути $ZYJ$, можно вычислить код пути $ZFJ$, и наоборот. А также зная код пути $FZY$, можно вычислить код пути $FJY$, и наоборот. Таким образом, мы можем осуществить локальное преобразование пути через операцию с его кодом.

\medskip

\subsection{Случай цепи $\mathbb{DL}1$; преобразования 7 и 10}

\begin{figure}[hbtp]
\centering
\includegraphics[width=1\textwidth]{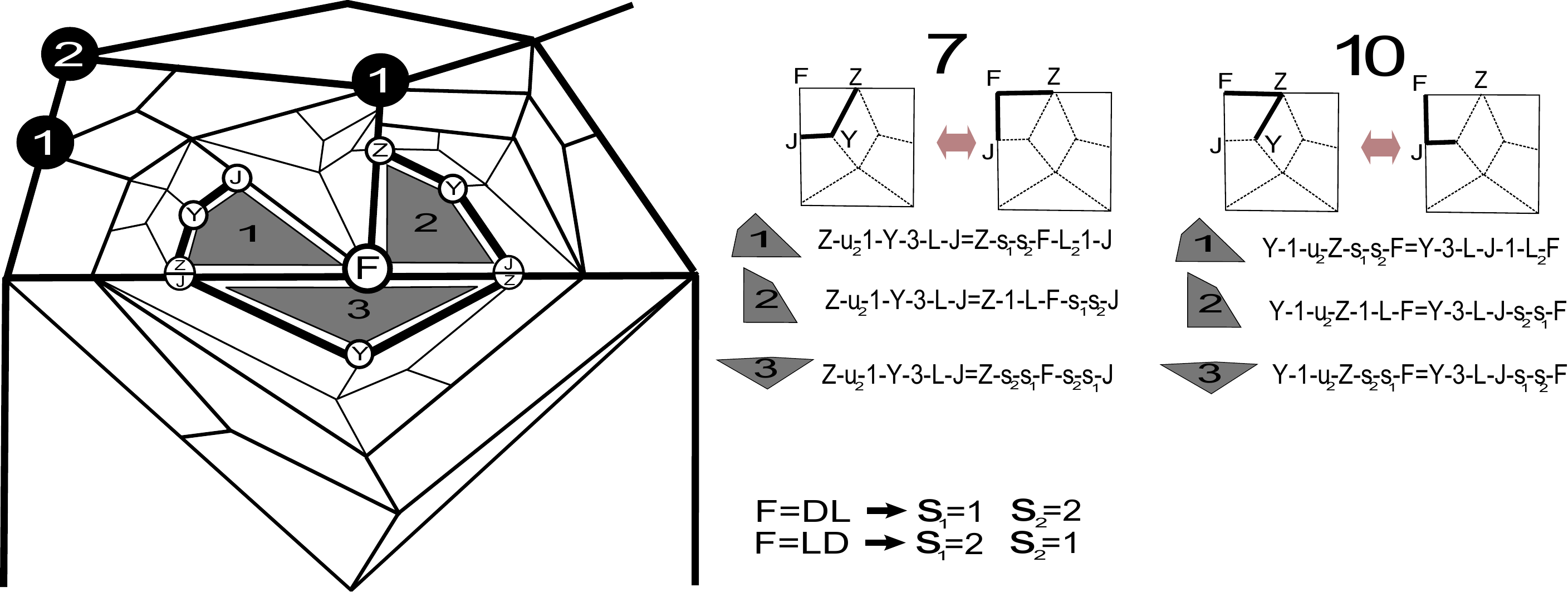}
\caption{Случаи расположения пути вокруг DL1-цепи и соответствующие им определяющие соотношения для локальных преобразований 7 и 10}
\label{DL1a}
\end{figure}

В правой части рисунка~\ref{DL1a} изображены локальные преобразования $7$ и $10$.

\medskip

Фиксируем некоторую вершину типа $\mathbb{DL}$ третьего уровня и ее начальника.
Черными кругами отметим вершины, являющиеся начальниками вершин  $Z$, $Y$, $F$, $J$.
Заметим, что зная окружение центральной вершины, мы можем вычислить коды вершин $Z$, $Y$, $F$, $J$, с точностью до подклееных окружений. Это позволяет ввести определяющие отношения, записанные в правой части рисунка~\ref{DL1a}.

\medskip

{\bf Характеризация.} Составление таблиц полностью аналогично случаю $\mathbb{C}1$-цепи.

\medskip

{\bf  Восстановление кода.}
Вершина $F$ имеет общего начальника с одной из вершин $J$ или $Z$, в каждом из случаев.  То есть, $ZYJ$, можно вычислить код пути $ZFJ$. В обратную сторону, а также для  локального преобразования $10$: очевидно, что зная $F$, можно вычислить коды остальных вершин. Таким образом, мы можем осуществить локальное преобразование пути через операцию с его кодом.

\medskip

\subsection{Случай цепи $\mathbb{DL}2$; преобразования 8 и 9}

\begin{figure}[hbtp]
\centering
\includegraphics[width=0.9\textwidth]{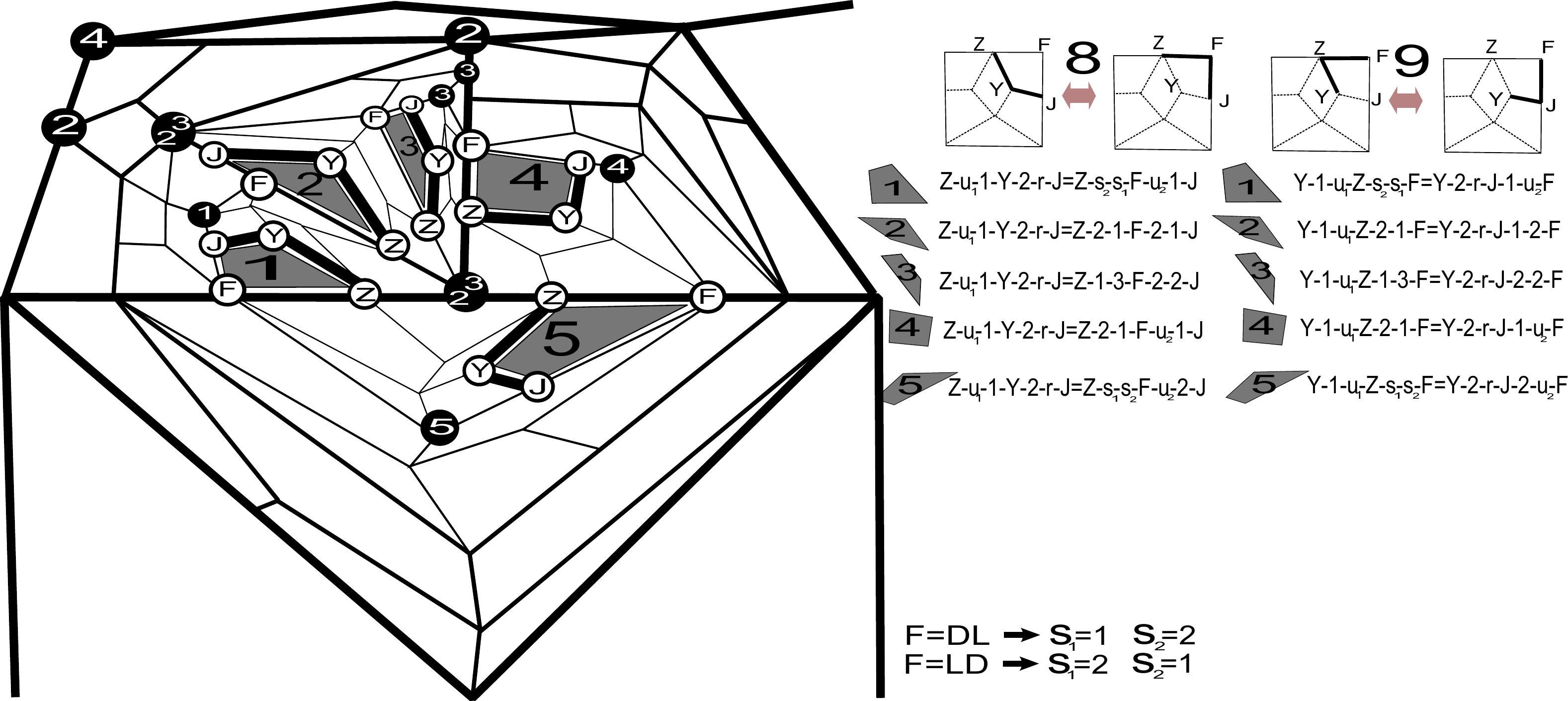}
\caption{Случаи расположения пути вокруг DL2-цепи и определяющие соотношения для локальных преобразований 8 и 9}
\label{DL2b}
\end{figure}

В правой части рисунка~\ref{DL2b} изображены локальные преобразования $8$ и $9$.

\medskip

Зафиксируем вершину типа $\mathbb{DR}$ третьего уровня и ее начальников.
Черными  кругами отметим вершины, являющиеся начальниками вершин  $Z$, $Y$, $F$, $J$.
Заметим, что зная окружение центрального $\mathbb{DL}$-узла, мы можем найти их типы, уровни и окружения. Так как начальники вершин  $Z$, $Y$, $F$, $J$ содержатся среди вершин, отмеченных черными кругами, то коды вершин $Z$, $Y$, $F$, $J$, во всех случаях мы можем назвать явно, с точностью до подклееных окружений. Это позволяет ввести определяющие отношения, записанные в правой части рисунка~\ref{DL2b}.

\medskip

{\bf Характеризация.} Составление таблиц полностью аналогично случаю $\mathbb{C}1$-цепи.

\medskip

{\bf Восстановление кода.}
Заметим, что у $Z$ и $F$ в каждом из случаев общие начальники, то есть код каждой из этих вершин может быть вычислен, исходя из кодов остальных трех вершин. $Y$ вычисляется также легко, так как во всех случаях первым начальником является $Z$, а тип второго в каждом случае ясен.

Начальники $J$ вычисляется следующим образом: в  случаях $1$, $4$, $5$ -- единственный начальник, узел $F$; в  случае $2$ -- начальники те же что и у $F$; в  случае $3$ -- первый и второй начальники как первый и третий у $F$.

Таким образом, зная код пути $ZYJ$, можно вычислить код пути $ZFJ$, и наоборот. А также зная код пути $FZY$, можно вычислить код пути $FJY$, и наоборот. Таким образом, мы можем осуществить локальное преобразование пути через операцию с его кодом.

\medskip

\subsection{Случай цепи $\mathbb{DL}2$; преобразования 7 и 10}

\begin{figure}[hbtp]
\centering
\includegraphics[width=1\textwidth]{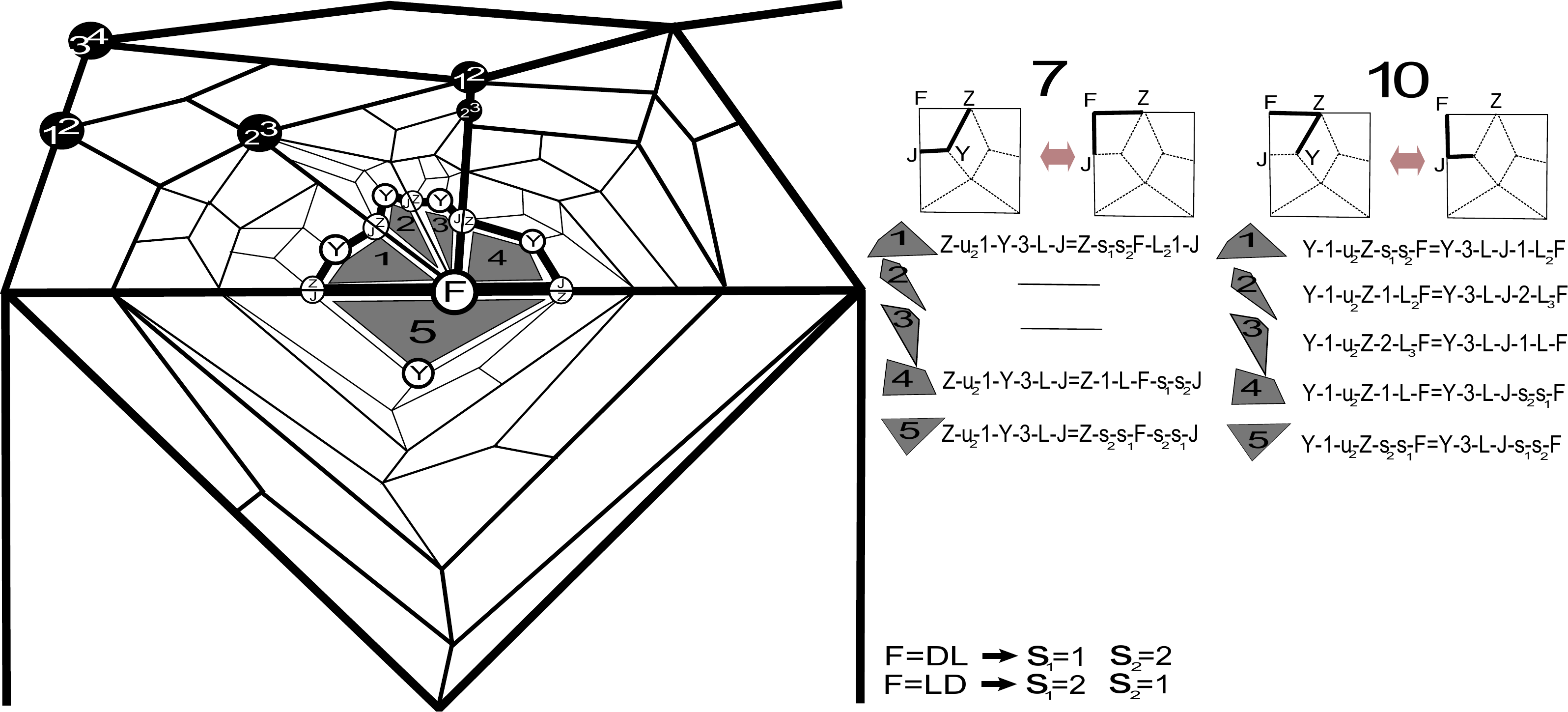}
\caption{Случаи расположения пути вокруг DL2-цепи и определяющие соотношения для локальных преобразований ~7 и 10}
\label{DL2a}
\end{figure}

В правой части рисунка~\ref{DL2a} изображены локальные преобразования $7$ и $10$.

\medskip

Фиксируем некоторую вершину типа $\mathbb{DL}$ третьего уровня и ее начальников.
Черными  кругами отметим вершины, являющиеся начальниками вершин  $Z$, $Y$, $F$, $J$.
Заметим, что зная окружение этой вершины, мы можем вычислить коды вершин $Z$, $Y$, $F$, $J$, с точностью до подклееных окружений. Это позволяет ввести определяющие отношения, записанные в правой части рисунка~\ref{DL2a}.

Для локального преобразования $7$ пути в случаях $2$ и $3$ удовлетворяют условиям мертвого паттерна, и в этих случаях мы соотношения не вводим.

\medskip

{\bf  Характеризация.} Составление таблиц полностью аналогично случаю $\mathbb{C}1$-цепи.

\medskip

{\bf Восстановление кода.}
В случаях $1,4,5$ вершина $F$ имеет общего начальника с одной из вершин $J$ или $Z$. То есть, $ZYJ$, можно вычислить код пути $ZFJ$. В обратную сторону, а также для  локального преобразования $10$: очевидно, что зная $F$, можно вычислить коды остальных вершин. Таким образом, мы можем осуществить локальное преобразование пути через операцию с его кодом.

\medskip

\subsection{Случай цепи $\mathbb{DL}3$; преобразования 7, 8, 9, 10}

Случай $\mathbb{DL}3$-цепи полностью аналогичен $\mathbb{DL}2$ случаю, соотношения выглядят идентично, только кодировки вершин $J$, $F$, $Z$, $Y$ отвечают $\mathbb{DL}3$-цепи. Все рассуждения о вычислении путей полностью аналогичны. Соотношений вводится столько же, сколько для $\mathbb{DL}2$ случая.

\medskip

Далее рассмотрим случаи цепей вблизи края макроплиток, то есть цепи с центрами в вершинах типов $\mathbb{D}$, $\mathbb{U}$, $\mathbb{R}$, $\mathbb{L}$. Рассмотрение этих случаев полностью аналогично случаям обычных боковых узлов, просто используются не все из локальных преобразований. Тем не менее, мы кратко приведем вводимые определяющие соотношения.

\medskip

\subsection{Случай цепи $\mathbb{D}1$; преобразования 8 и 9}

\begin{figure}[hbtp]
\centering
\includegraphics[width=0.9\textwidth]{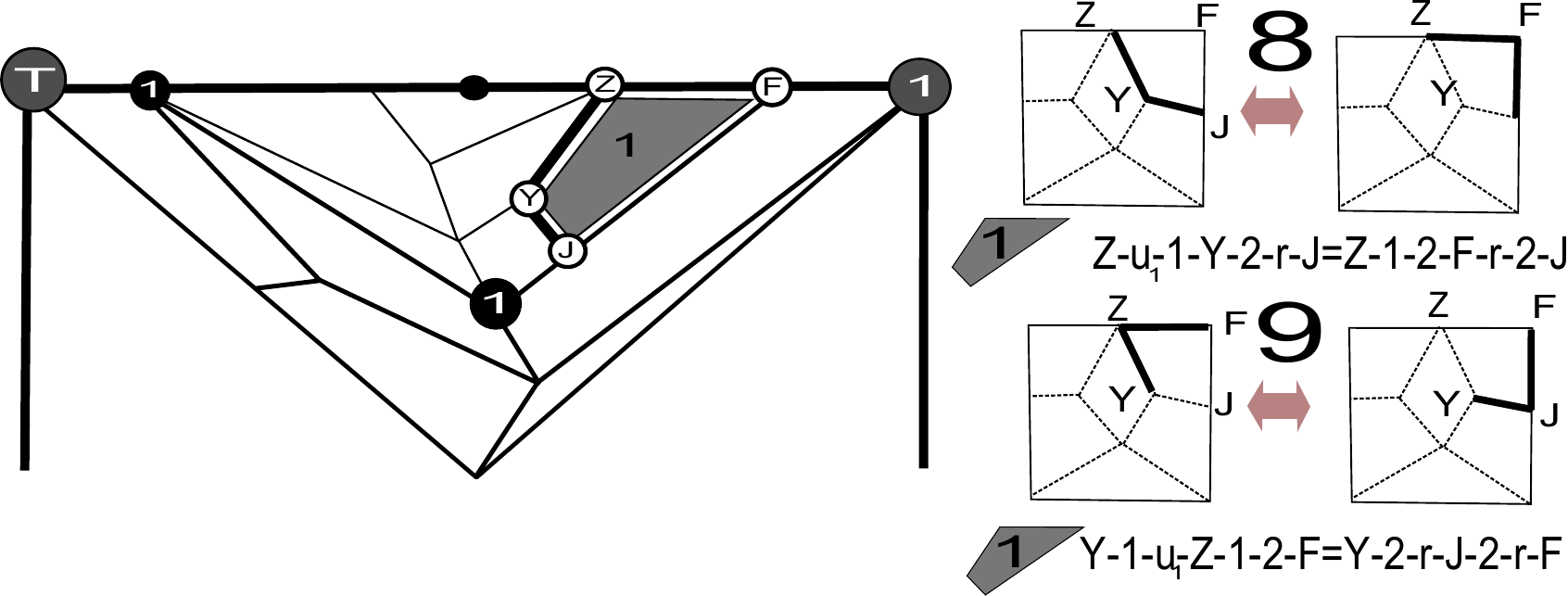}
\caption{Случаи расположения пути вокруг D1-цепи  и соответствующие им определяющие соотношения для локальных преобразований 8 и 9}
\label{D1b}
\end{figure}

В правой части рисунка~\ref{D1b} изображены локальные преобразования $8$ и $9$.

\medskip

Характеризация и восстановление кода полностью аналогично случаю $\mathbb{DL}1$ или  $\mathbb{DR}1$ цепи.

\subsection{Случай цепи $\mathbb{D}1$; преобразования 7 и 10}

\begin{figure}[hbtp]
\centering
\includegraphics[width=1\textwidth]{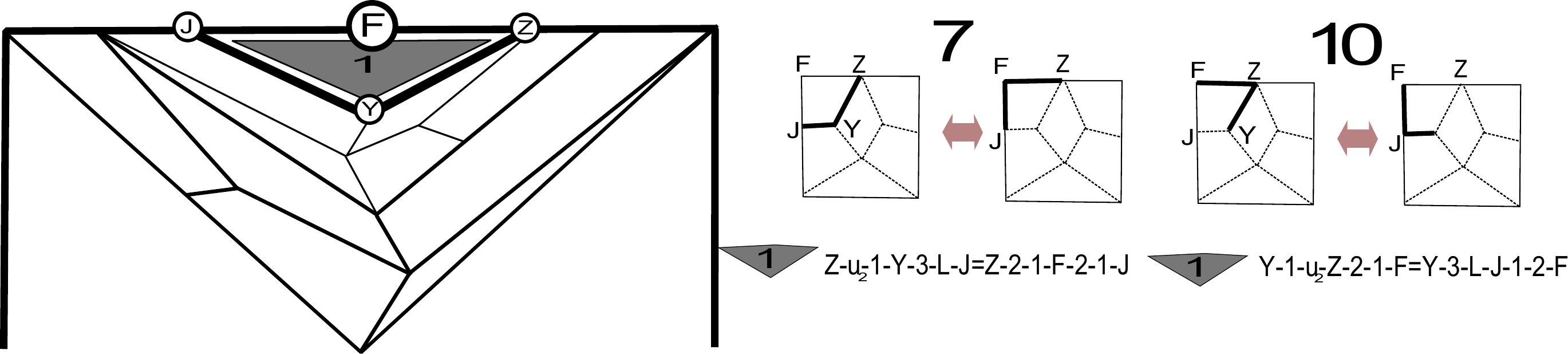}
\caption{Случаи расположения пути вокруг D1-цепи и соответствующие им определяющие соотношения для локальных преобразований 7 и 10}
\label{D1a}
\end{figure}

В правой части рисунка~\ref{D1a} изображены локальные преобразования $7$ и $10$.

\medskip

Характеризация и восстановление кода полностью аналогично случаю $\mathbb{DL}1$ или  $\mathbb{DR}1$ цепи.
\medskip

\subsection{Случай цепи $\mathbb{D}2$; преобразования 8 и 9}

\begin{figure}[hbtp]
\centering
\includegraphics[width=0.9\textwidth]{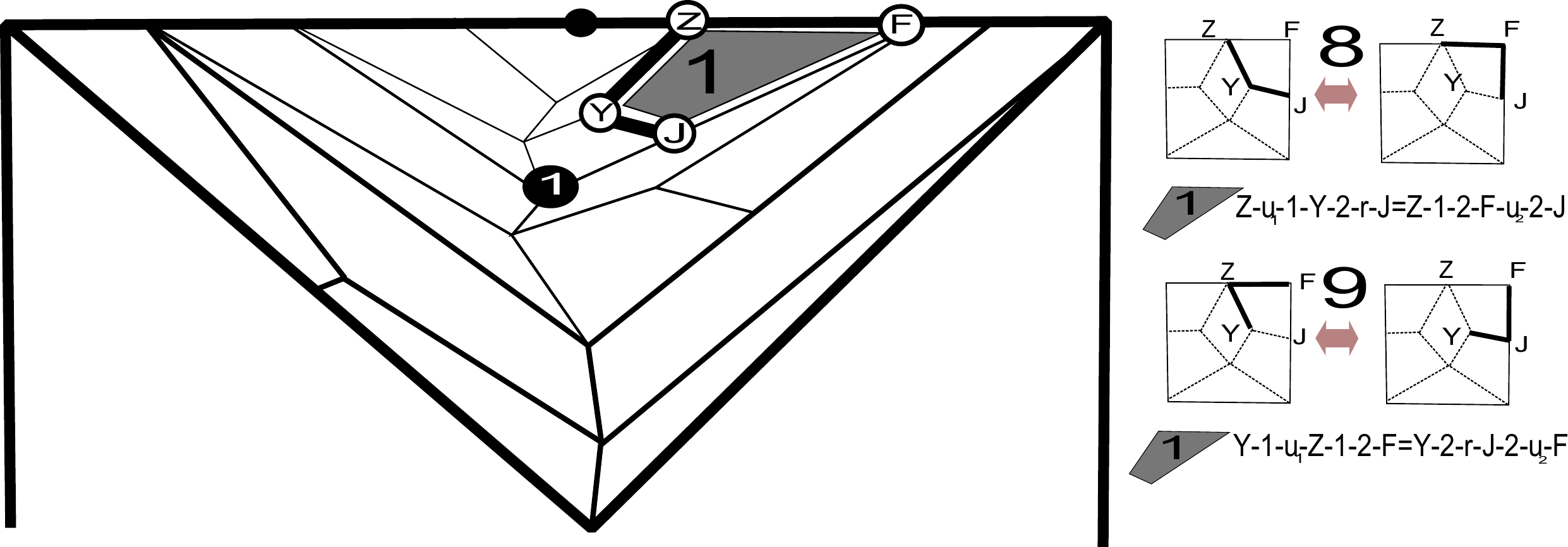}
\caption{Случаи расположения пути вокруг D2-цепи  и соответствующие им определяющие соотношения для локальных преобразований 8 и 9}
\label{D2b}
\end{figure}

В правой части рисунка~\ref{D2b} изображены локальные преобразования $8$ и $9$.

\medskip

Характеризация и восстановление кода полностью аналогично случаю $\mathbb{DL}2$ или  $\mathbb{DR}2$ цепи.

\subsection{Случай цепи $\mathbb{D}2$; преобразования 7 и 10}

\begin{figure}[hbtp]
\centering
\includegraphics[width=1\textwidth]{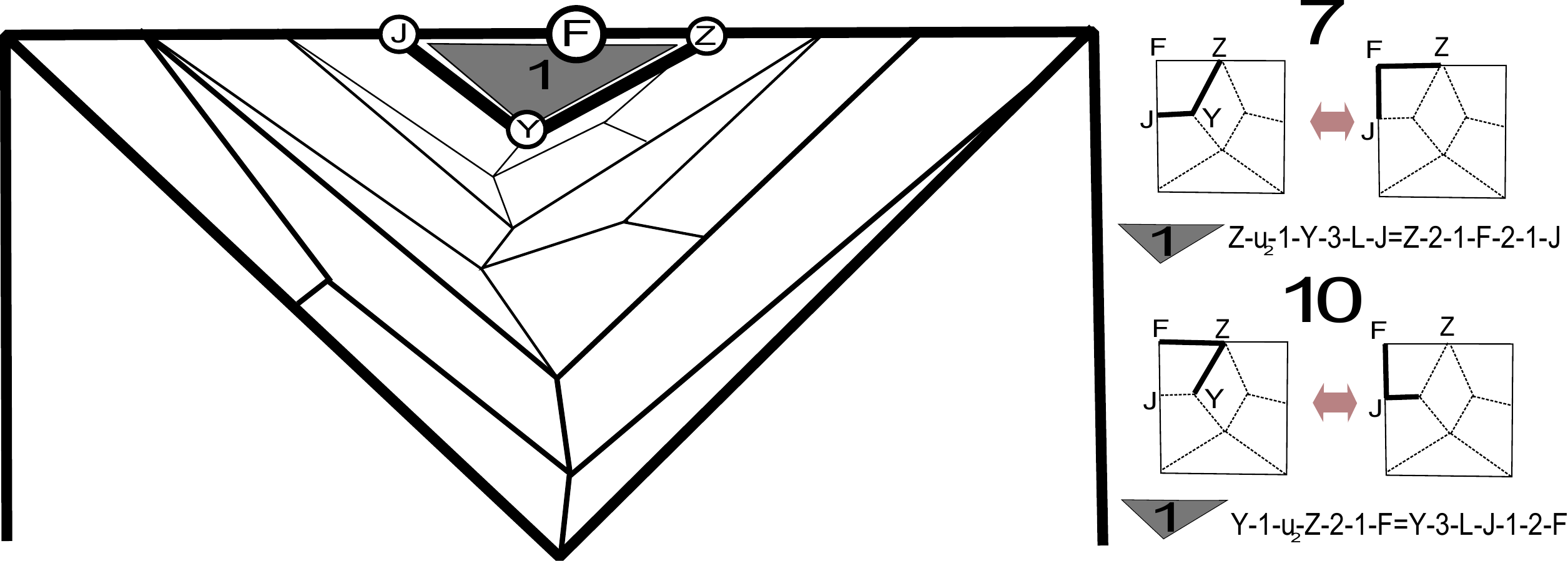}
\caption{Случаи расположения пути вокруг D2-цепи и соответствующие им определяющие соотношения для локальных преобразований 7 и 10}
\label{D2a}
\end{figure}

В правой части рисунка~\ref{D2a} изображены локальные преобразования $7$ и $10$.

\medskip

Характеризация и восстановление кода полностью аналогично случаю $\mathbb{DL}2$ или  $\mathbb{DR}2$ цепи.
\medskip

\subsection{Случай цепи $\mathbb{U}1$; преобразования 8 и 9}

\begin{figure}[hbtp]
\centering
\includegraphics[width=0.9\textwidth]{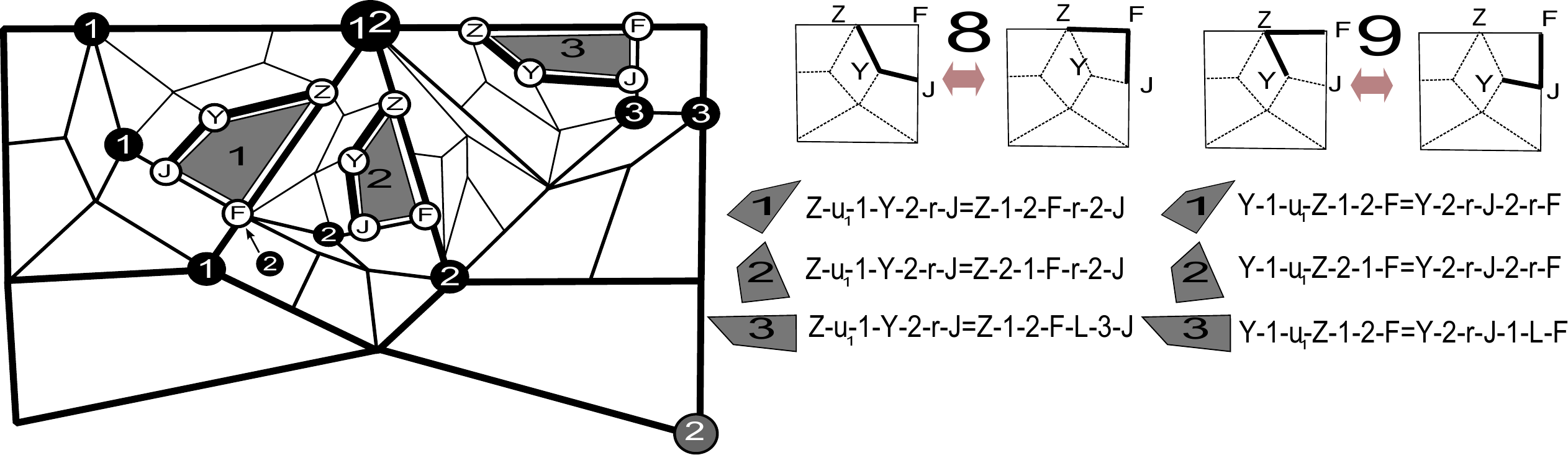}
\caption{Случаи расположения пути вокруг U1-цепи  и соответствующие им определяющие соотношения для 8 и 9 локальных преобразований}
\label{U1b}
\end{figure}

В правой части рисунка~\ref{U1b} изображены локальные преобразования $8$ и $9$.

\medskip

Характеризация и восстановление кода полностью аналогично случаю $\mathbb{UL}1$ или  $\mathbb{UR}1$ цепи.

\medskip

\subsection{Случай цепи $\mathbb{U}1$; преобразования 7 и 10}

\begin{figure}[hbtp]
\centering
\includegraphics[width=1\textwidth]{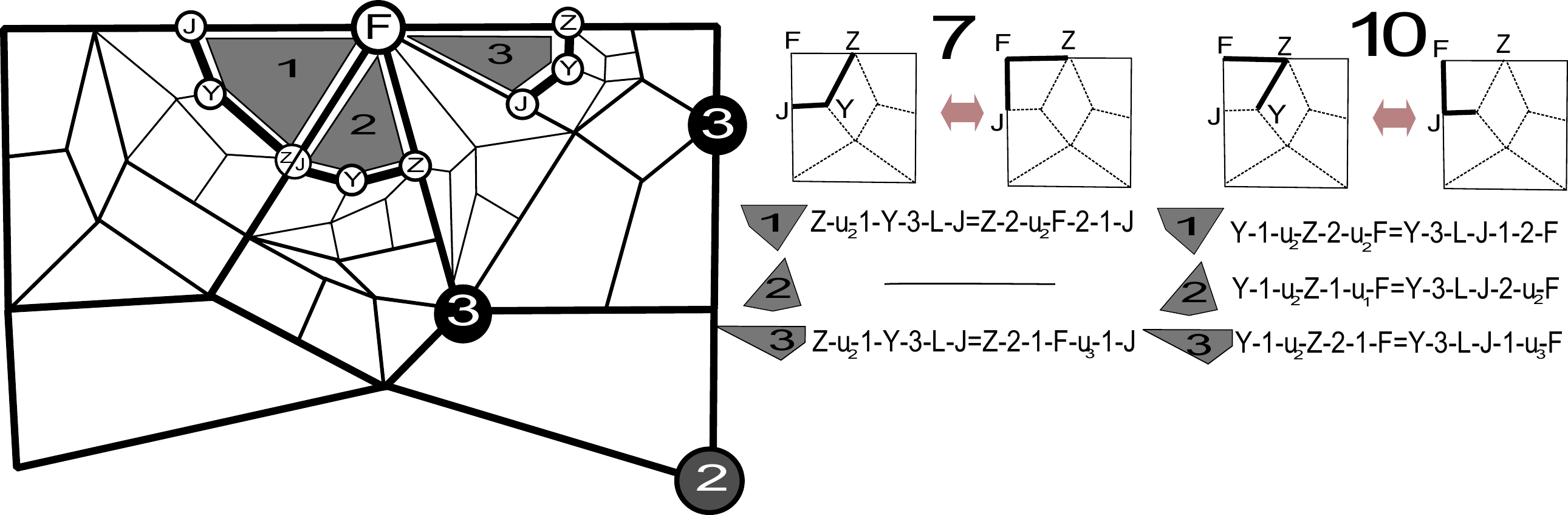}
\caption{Случаи расположения пути вокруг U1-цепи и соответствующие им определяющие соотношения для локальных преобразований 7 и 10}
\label{U1a}
\end{figure}

В правой части рисунка~\ref{U1a} изображены локальные преобразования $7$ и $10$.

\medskip

Характеризация и восстановление кода полностью аналогично случаю $\mathbb{UL}1$ или  $\mathbb{UR}1$ цепи.

\medskip

\subsection{Случай цепи $\mathbb{U}2$; преобразования 8 и 9}

\begin{figure}[hbtp]
\centering
\includegraphics[width=1\textwidth]{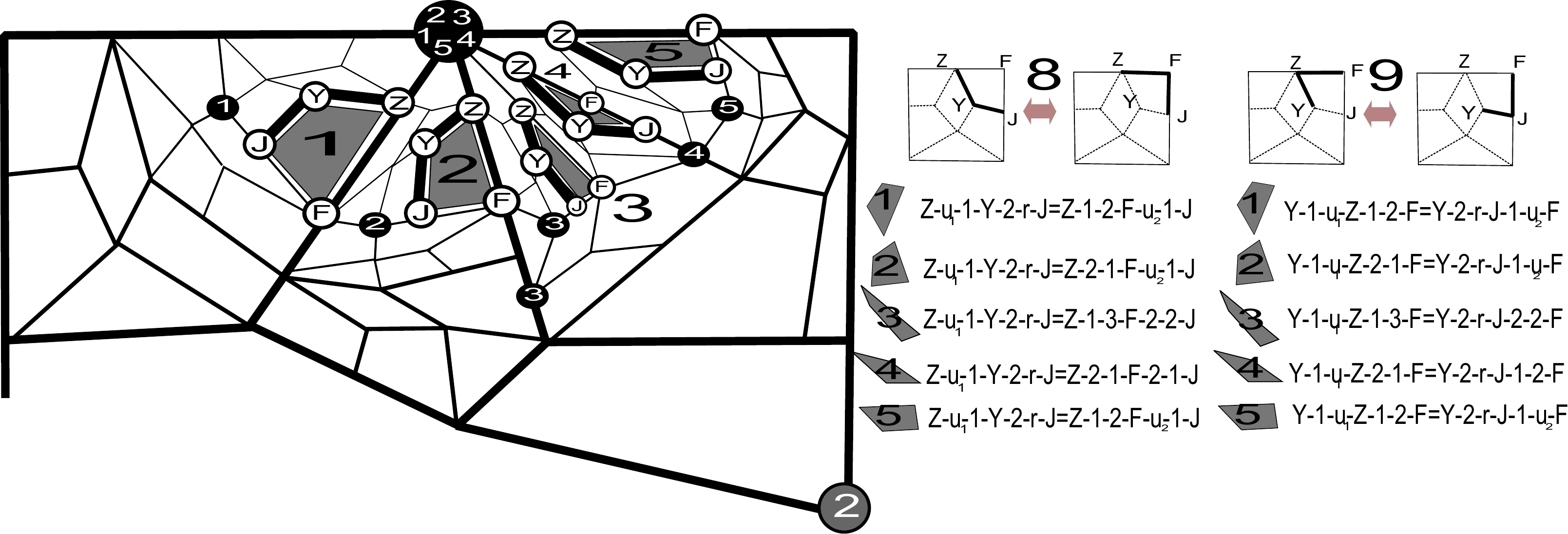}
\caption{Случаи расположения пути вокруг U2-цепи  и соответствующие им определяющие соотношения для локальных преобразований 8 и 9}
\label{U2b}
\end{figure}

В правой части рисунка~\ref{U2b} изображены локальные преобразования $8$ и $9$.

\medskip

Характеризация и восстановление кода полностью аналогично случаю $\mathbb{UL}2$ или  $\mathbb{UR}2$ цепи.

\medskip

\subsection{Случай цепи $\mathbb{U}2$; преобразования 7 и 10}

\begin{figure}[hbtp]
\centering
\includegraphics[width=0.9\textwidth]{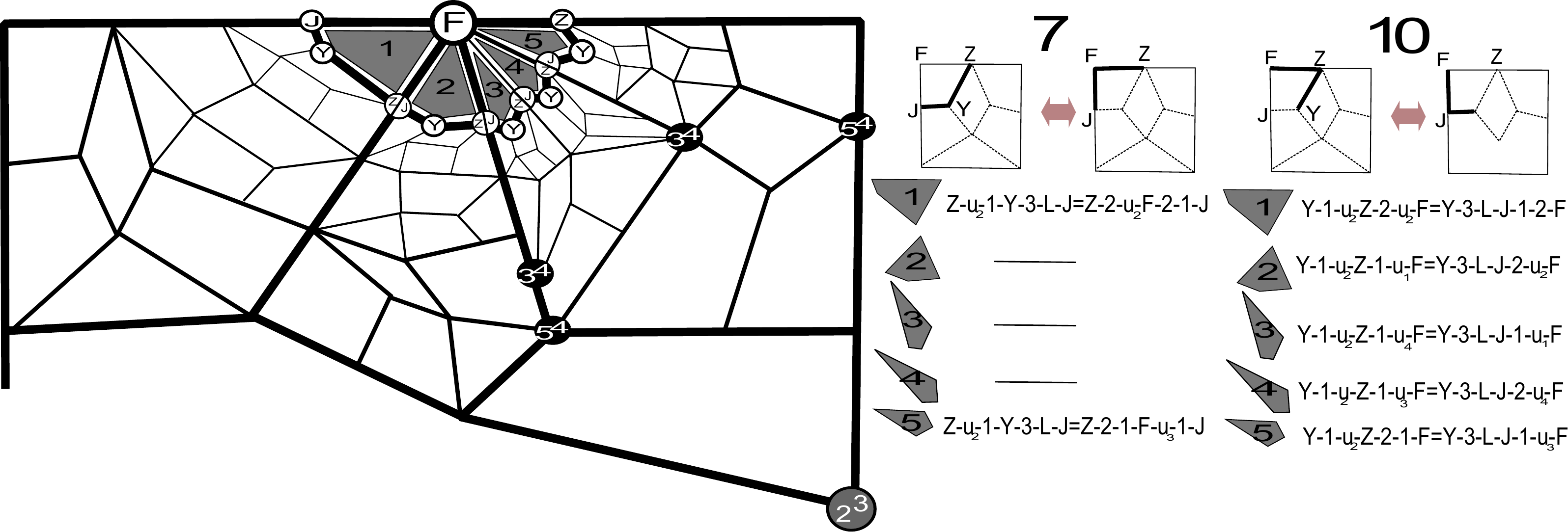}
\caption{Случаи расположения пути вокруг U2-цепи и соответствующие им определяющие соотношения для локальных преобразований 7 и 10}
\label{U2a}
\end{figure}

В правой части рисунка~\ref{U2a} изображены локальные преобразования $7$ и $10$.

\medskip

Характеризация и восстановление кода полностью аналогично случаю $\mathbb{UL}2$ или  $\mathbb{UR}2$ цепи.

\medskip

\subsection{Случай цепи $\mathbb{U}3$; преобразования 7, 8, 9, 10}

Случай $\mathbb{U}3$-цепи полностью аналогичен случаю $\mathbb{U}2$-цепи, соотношения выглядят идентично, только кодировки вершин $J$, $F$, $Z$, $Y$ отвечают $\mathbb{U}3$-цепи. Все рассуждения о вычислении путей полностью аналогичны. Соотношений вводится столько же, сколько для $\mathbb{U}2$ случая.

\medskip

\subsection{Случай цепи $\mathbb{L}1$; преобразования 8 и 9}

\begin{figure}[hbtp]
\centering
\includegraphics[width=0.9\textwidth]{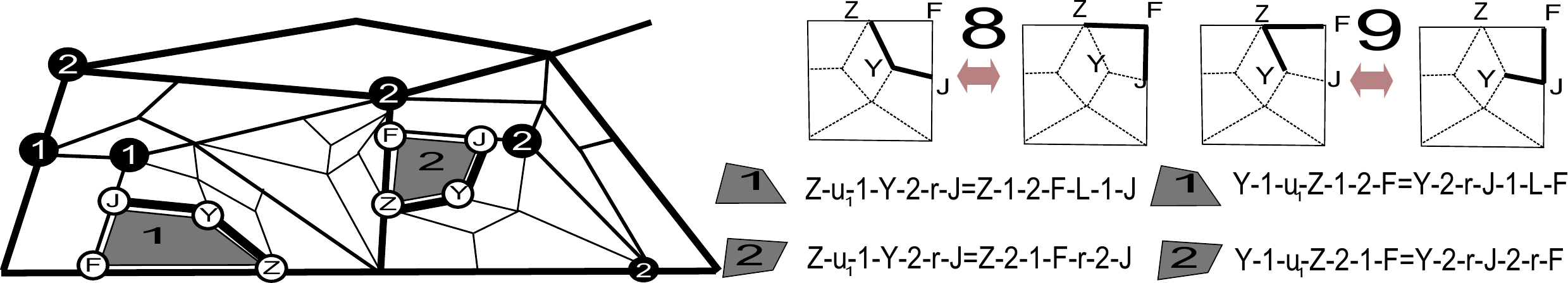}
\caption{Случаи расположения пути вокруг L1-цепи  и соответствующие им определяющие соотношения для локальных преобразований 8 и 9}
\label{L1b}
\end{figure}

В правой части рисунка~\ref{L1b} изображены локальные преобразования $8$ и $9$.

\medskip

Характеризация и восстановление кода полностью аналогично случаю $\mathbb{UL}1$ или  $\mathbb{DL}1$ цепи.

\medskip

\subsection{Случай цепи $\mathbb{L}1$; преобразования 7 и 10}

\begin{figure}[hbtp]
\centering
\includegraphics[width=1\textwidth]{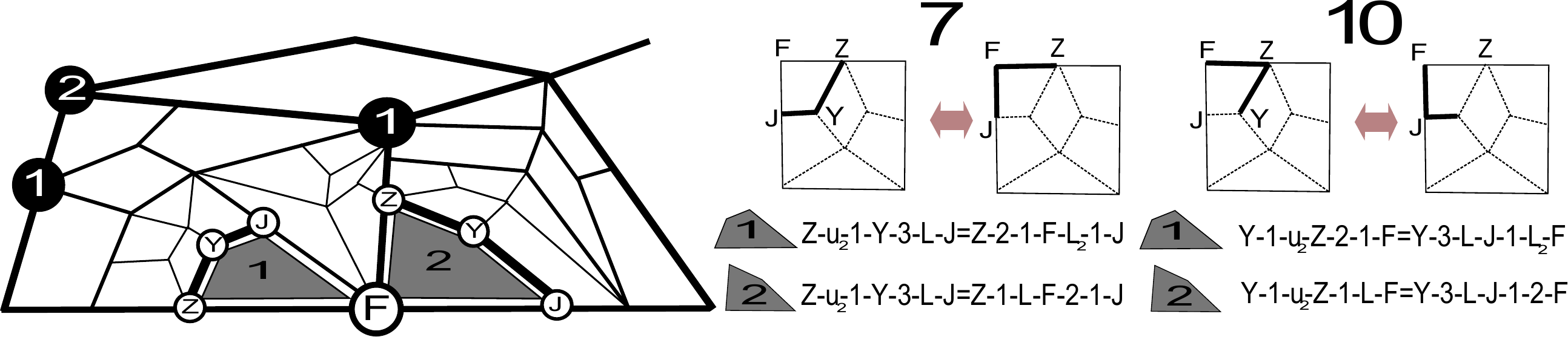}
\caption{Случаи расположения пути вокруг L1-цепи и соответствующие им определяющие соотношения для локальных преобразований 7 и 10}
\label{L1a}
\end{figure}

В правой части рисунка~\ref{L1a} изображены локальные преобразования $7$ и $10$.

\medskip

Характеризация и восстановление кода полностью аналогично случаю $\mathbb{UL}1$ или  $\mathbb{DL}1$ цепи.

\medskip

\subsection{Случай цепи $\mathbb{L}2$; преобразования 8 и 9}

\begin{figure}[hbtp]
\centering
\includegraphics[width=1\textwidth]{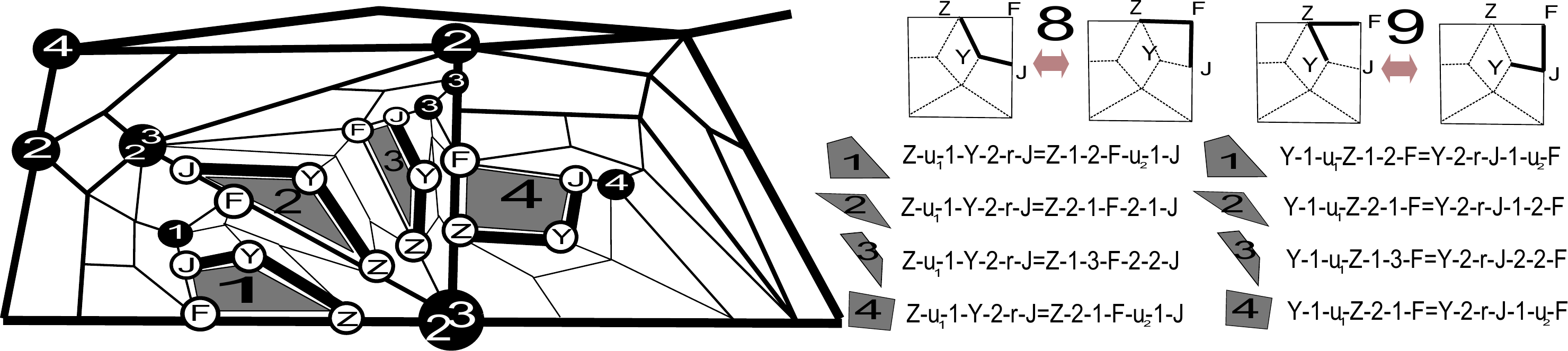}
\caption{Случаи расположения пути вокруг L2-цепи  и соответствующие им определяющие соотношения для локальных преобразований 8 и 9}
\label{L2b}
\end{figure}

В правой части рисунка~\ref{L2b} изображены локальные преобразования $8$ и $9$.

\medskip

Характеризация и восстановление кода полностью аналогично случаю $\mathbb{UL}2$ или  $\mathbb{DL}2$ цепи.

\medskip

\subsection{Случай цепи $\mathbb{L}2$, преобразования 7 и 10}

\begin{figure}[hbtp]
\centering
\includegraphics[width=1\textwidth]{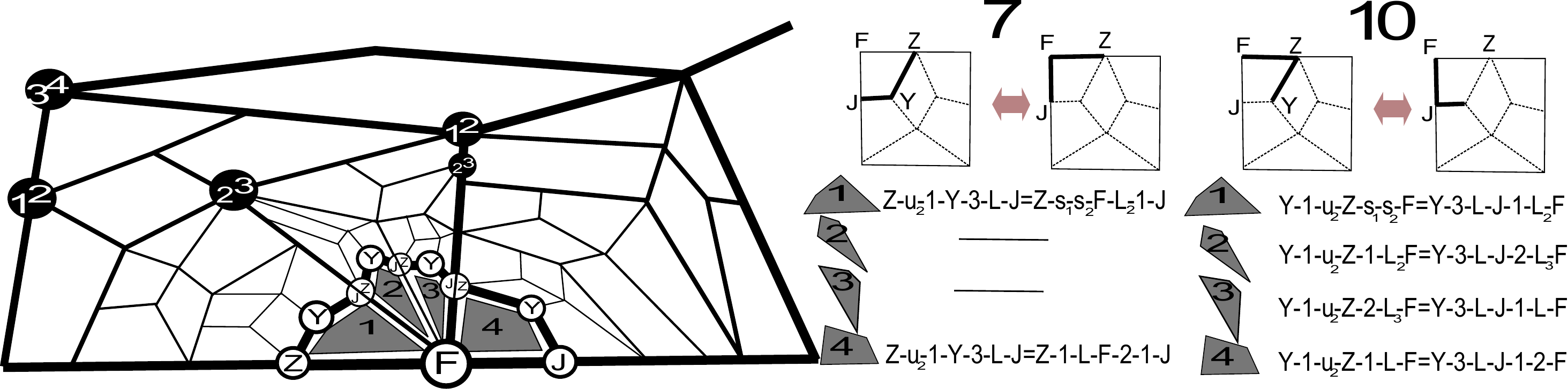}
\caption{Случаи расположения пути вокруг L2-цепи и определяющие соотношения для локальных преобразований 7 и 10}
\label{L2a}
\end{figure}

В правой части рисунка~\ref{L2a} изображены локальные преобразования $7$ и $10$.

\medskip

Характеризация и восстановление кода полностью аналогично случаю $\mathbb{UL}2$ или  $\mathbb{DL}2$ цепи.

\medskip

\subsection{Случай цепи $\mathbb{L}3$; преобразования 7, 8, 9, 10}

Случай $\mathbb{L}3$-цепи полностью аналогичен случаю $\mathbb{L}2$-цепи, соотношения выглядят идентично, только кодировки вершин $J$, $F$, $Z$, $Y$ отвечают $\mathbb{L}3$-цепи. Все рассуждения о вычислении путей полностью аналогичны. Соотношений вводится столько же, сколько для $\mathbb{L}2$ случая.

\medskip

\subsection{Случай цепи $\mathbb{R}1$; преобразования 8 и 9}

\begin{figure}[hbtp]
\centering
\includegraphics[width=1\textwidth]{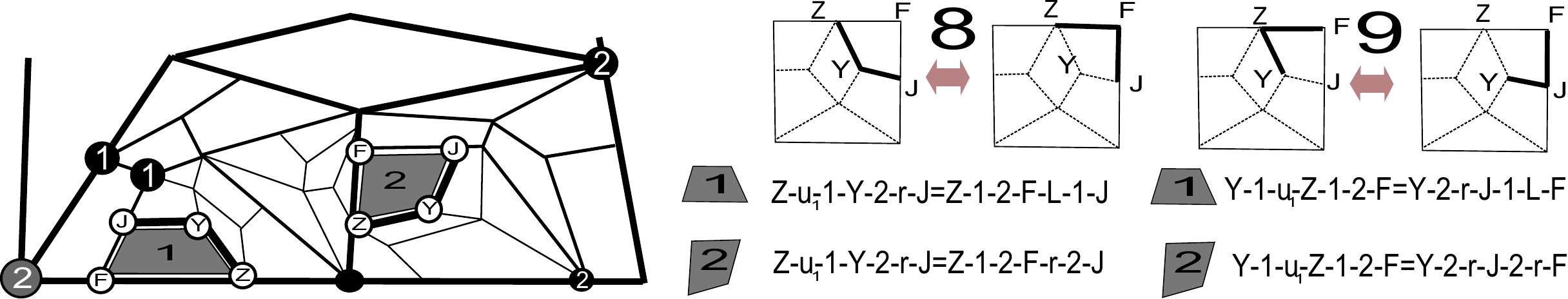}
\caption{Случаи расположения пути вокруг R1-цепи и определяющие соотношения для локальных преобразований 8 и 9}
\label{R1b}
\end{figure}

В правой части рисунка~\ref{R1b} изображены локальные преобразования $8$ и $9$.

\medskip

Характеризация и восстановление кода полностью аналогично случаю $\mathbb{DR}1$ или  $\mathbb{UR}1$ цепи.

\medskip

\subsection{Случай цепи $\mathbb{R}1$; преобразования 7 и 10}

\begin{figure}[hbtp]
\centering
\includegraphics[width=0.9\textwidth]{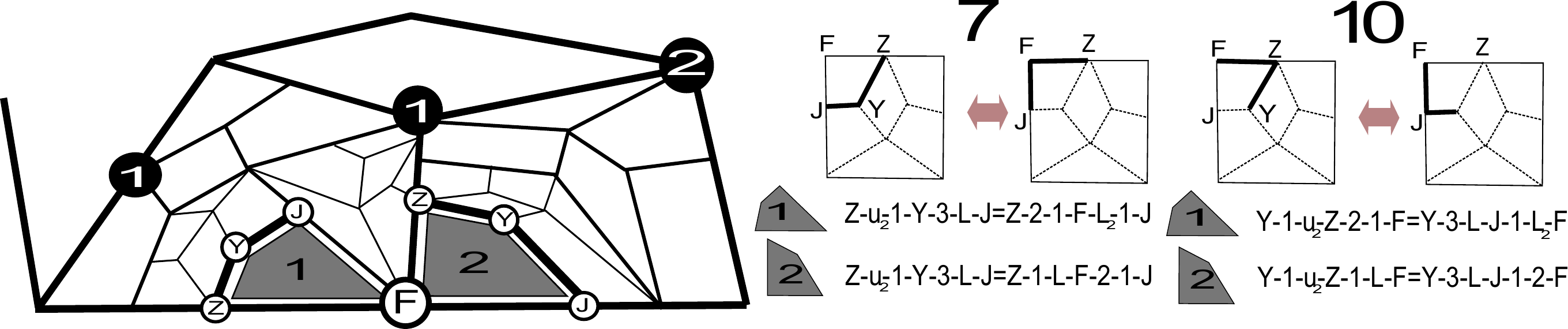}
\caption{Случаи расположения пути вокруг R1-цепи и определяющие соотношения для локальных преобразований 7 и 10}
\label{R1a}
\end{figure}

В правой части рисунка~\ref{R1a} изображены локальные преобразования $7$ и $10$.

\medskip

Характеризация и восстановление кода полностью аналогично случаю $\mathbb{DR}1$ или  $\mathbb{UR}1$ цепи.

\medskip

\subsection{Случай цепи $\mathbb{R}2$; преобразования 8 и 9}

\begin{figure}[hbtp]
\centering
\includegraphics[width=1\textwidth]{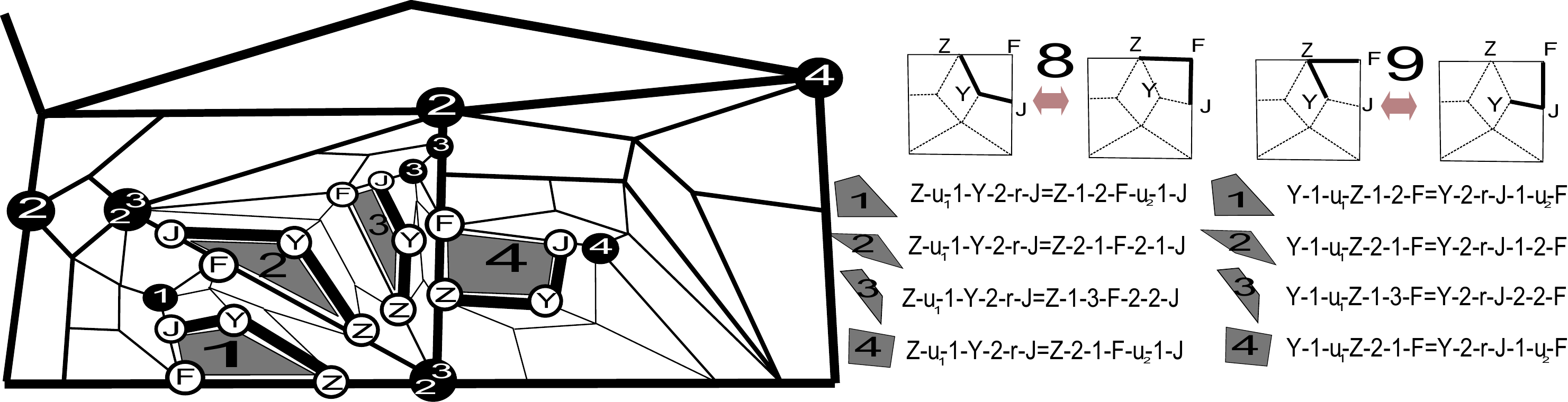}
\caption{Случаи расположения пути вокруг R2-цепи и определяющие соотношения для локальных преобразований 8 и 9}
\label{R2b}
\end{figure}

В правой части рисунка~\ref{R2b} изображены локальные преобразования $8$ и $9$.

\medskip

Характеризация и восстановление кода полностью аналогично случаю $\mathbb{DR}2$ или  $\mathbb{UR}2$ цепи.

\medskip

\subsection{Случай цепи $\mathbb{R}2$; преобразования 7 и 10}

\begin{figure}[hbtp]
\centering
\includegraphics[width=0.9\textwidth]{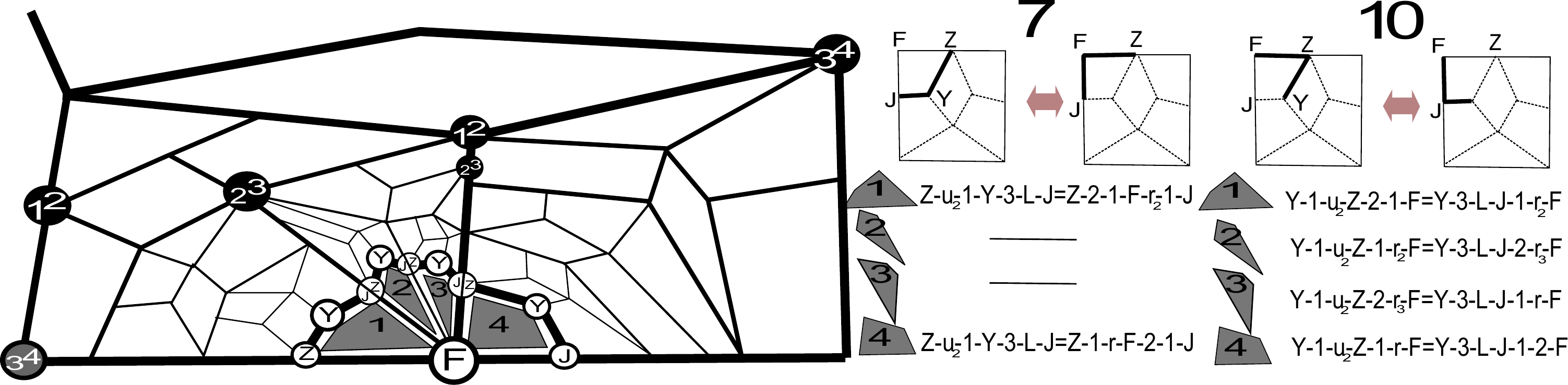}
\caption{Случаи расположения пути вокруг R2-цепи и определяющие соотношения для локальных преобразований 7 и 10}
\label{R2a}
\end{figure}

В правой части рисунка~\ref{R2a} изображены локальные преобразования $7$ и $10$.

\medskip

Характеризация и восстановление кода полностью аналогично случаю $\mathbb{DR}2$ или  $\mathbb{UR}2$ цепи.

\medskip

\subsection{Случай цепи $\mathbb{R}3$; преобразования 7, 8, 9, 10}

Случай $\mathbb{R}3$-цепи полностью аналогичен случаю $\mathbb{R}2$-цепи, соотношения выглядят идентично, только кодировки вершин $J$, $F$, $Z$, $Y$ отвечают $\mathbb{R}3$-цепи. Все рассуждения о вычислении путей полностью аналогичны. Соотношений вводится столько же, сколько для $\mathbb{R}2$ случая.

\medskip

\subsection{Случай цепи $\mathbb{CUL}0$; преобразования 7, 8, 9, 10}

\begin{figure}[hbtp]
\centering
\includegraphics[width=0.9\textwidth]{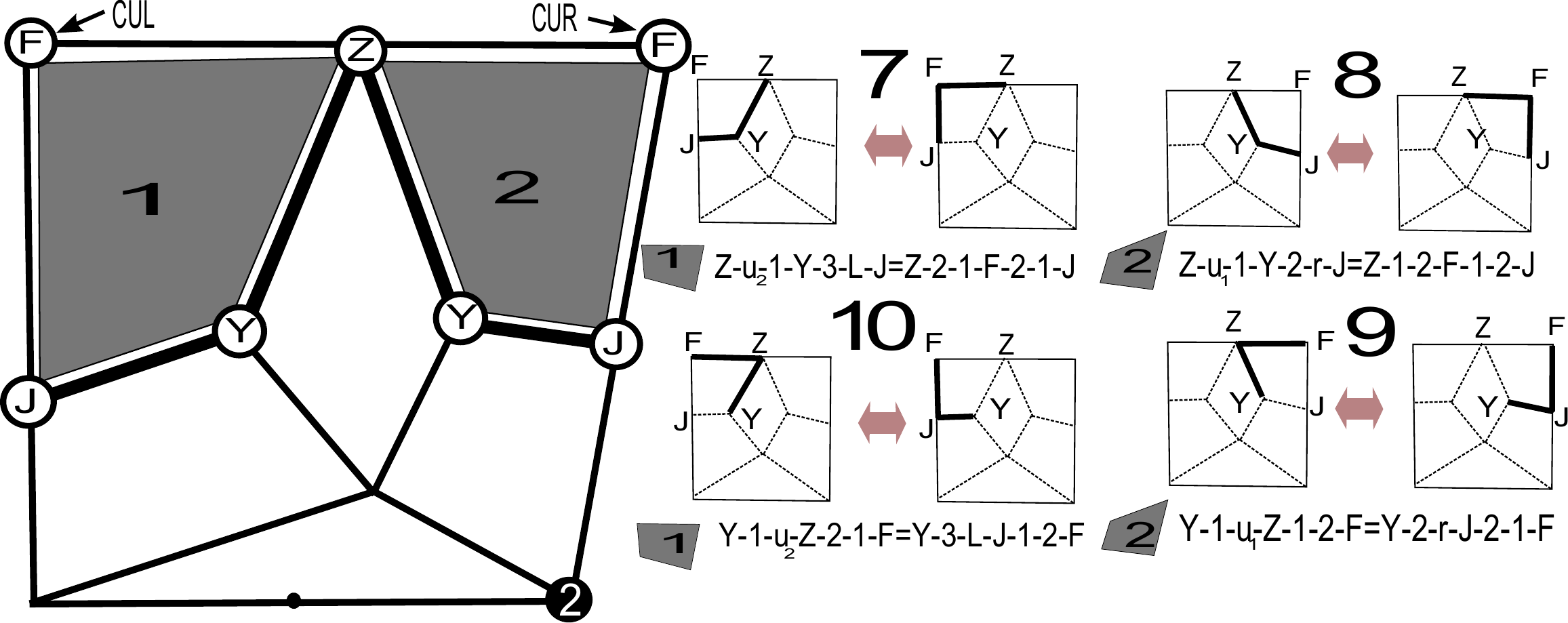}
\caption{Случаи расположения пути вокруг $\mathbb{CUL}0$-цепи и определяющие соотношения для локальных преобразований $7$, $8$, $9$, $10$}
\label{CUL0ab}
\end{figure}

В правой части рисунка~\ref{CUL0ab} изображены локальные преобразования $7$, $8$, $9$ и $10$.

\medskip

Характеризация очевидна.  Восстановление кода четвертой вершины по известным трем другим также ясно во всех случаях.

\medskip

\subsection{Случай цепи $\mathbb{CUL}1$; преобразования  7, 8, 9, 10}

\begin{figure}[hbtp]
\centering
\includegraphics[width=0.9\textwidth]{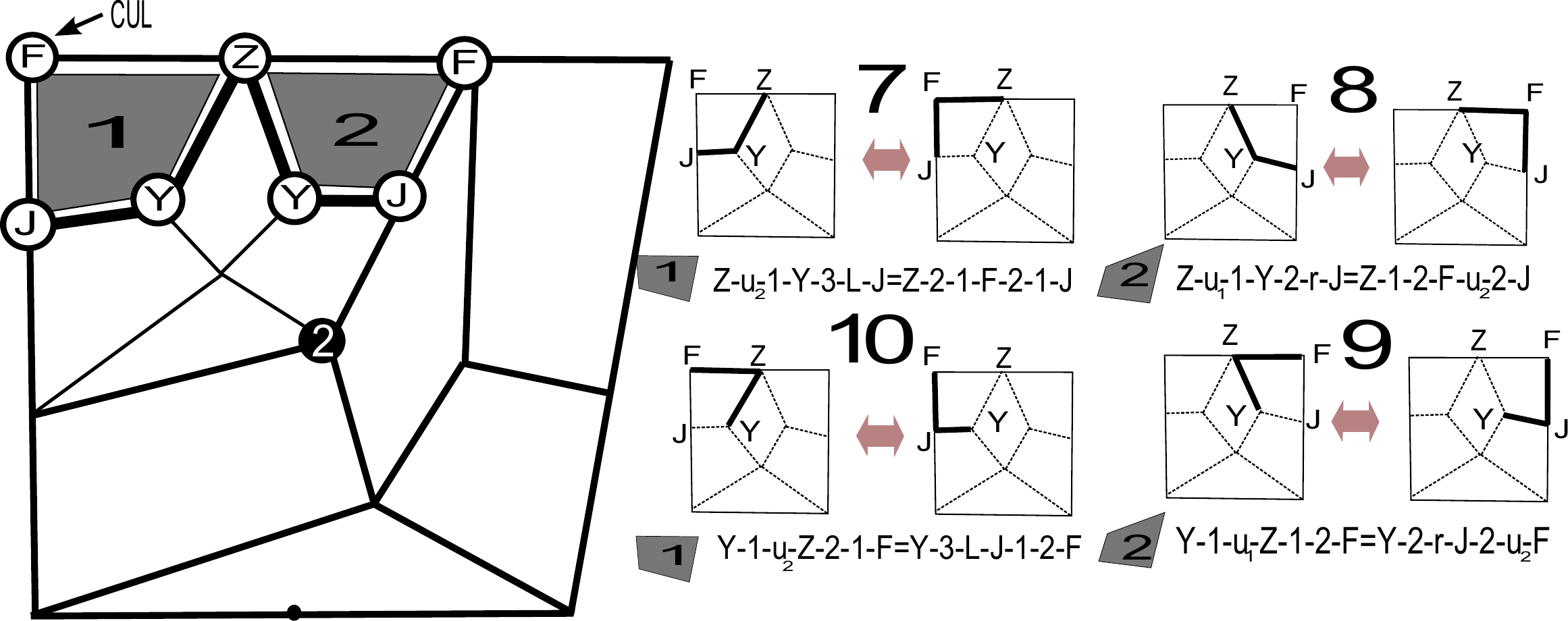}
\caption{Случаи расположения пути вокруг $\mathbb{CUL}1$-цепи и определяющие соотношения для локальных преобразований $7$, $8$, $9$, $10$}
\label{CUL1ab}
\end{figure}

В правой части рисунка~\ref{CUL1ab} изображены локальные преобразования $7$, $8$, $9$ и $10$.

\medskip

Характеризация очевидна.  Восстановление кода четвертой вершины по известным трем другим также ясно во всех случаях.

\medskip

\subsection{Случай цепи $\mathbb{CUR}1$; преобразования  7, 8, 9, 10}

\begin{figure}[hbtp]
\centering
\includegraphics[width=0.9\textwidth]{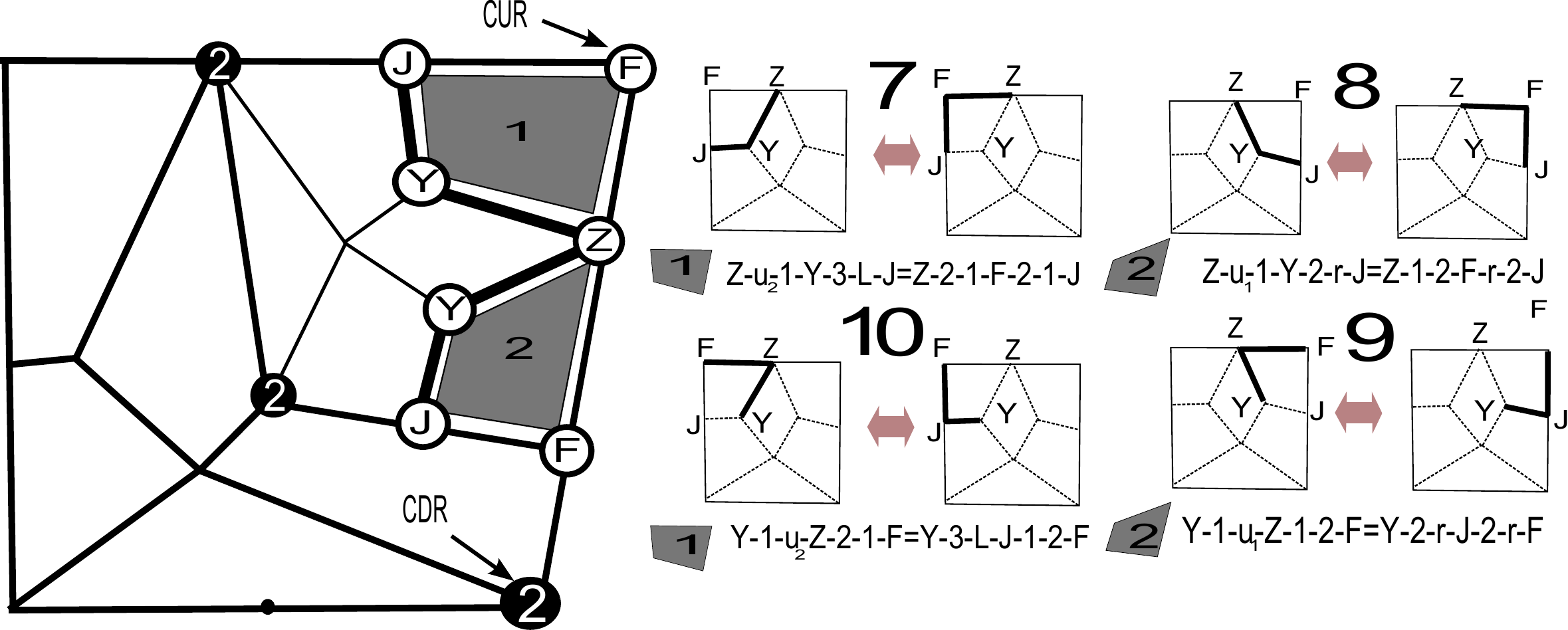}
\caption{Случаи расположения пути вокруг $\mathbb{CUR}1$-цепи и определяющие соотношения для локальных преобразований $7$, $8$, $9$, $10$}
\label{CUR1ab}
\end{figure}

В правой части рисунка~\ref{CUR1ab} изображены локальные преобразования $7$, $8$, $9$ и $10$.

\medskip

Характеризация очевидна.  Восстановление кода четвертой вершины по известным трем другим также ясно во всех случаях.

\medskip

\subsection{Случай цепи $\mathbb{CUR}2$; преобразования 7, 8, 9, 10}

\begin{figure}[hbtp]
\centering
\includegraphics[width=0.9\textwidth]{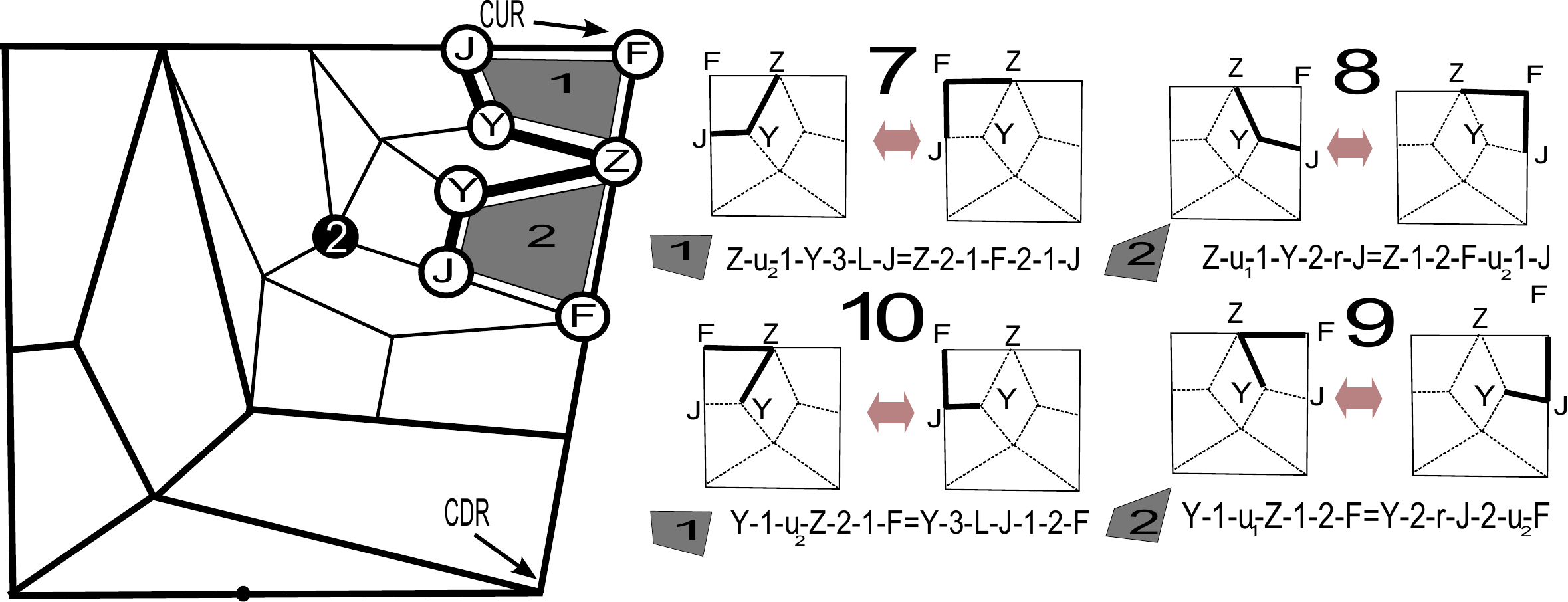}
\caption{Случаи расположения пути вокруг $\mathbb{CUR}2$-цепи  и соответствующие им определяющие соотношения для локальных преобразований $7$, $8$, $9$, $10$}
\label{CUR2ab}
\end{figure}

В правой части рисунка~\ref{CUR2ab} изображены локальные преобразования $8$ и $9$.

\medskip

Характеризация очевидна.  Восстановление кода четвертой вершины по известным трем другим также ясно во всех случаях.

\medskip

\subsection{Случай цепи $\mathbb{CDL}1$; преобразования  7, 8, 9, 10}

\begin{figure}[hbtp]
\centering
\includegraphics[width=0.9\textwidth]{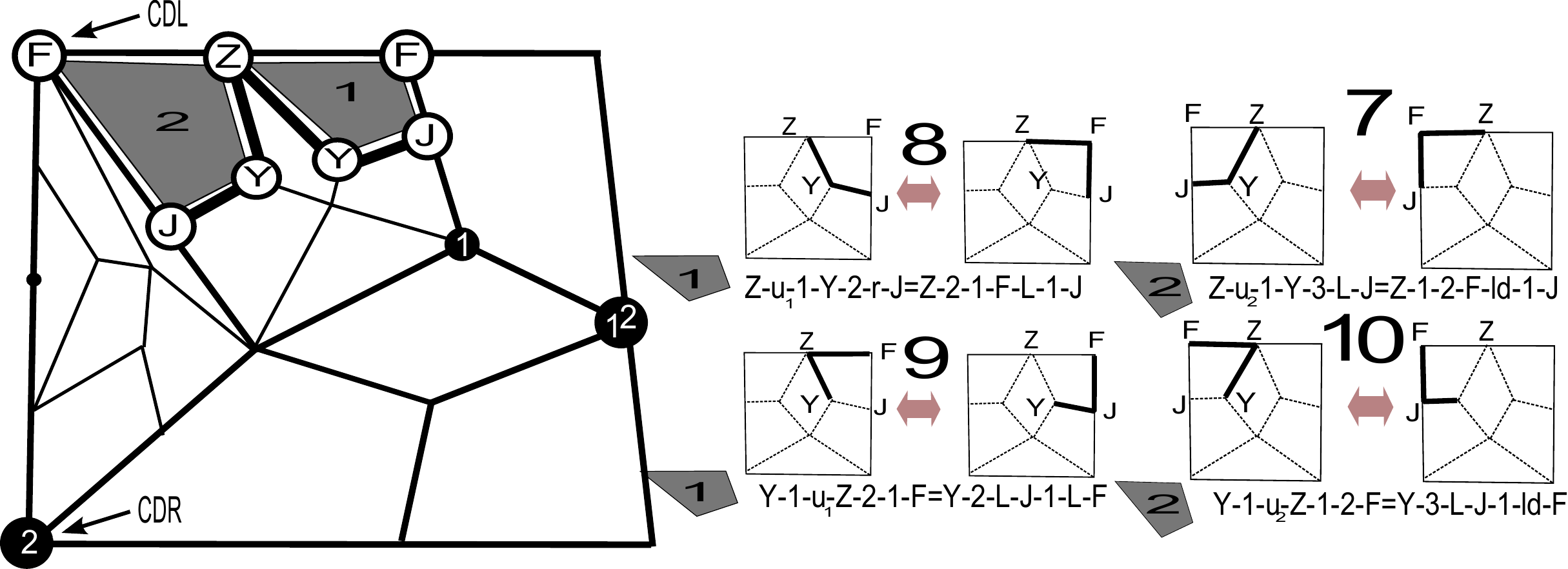}
\caption{Случаи расположения пути вокруг $\mathbb{CDL}1$-цепи и определяющие соотношения для локальных преобразований $7$, $8$, $9$, $10$}
\label{CDL1ab}
\end{figure}

В правой части рисунка~\ref{CDL1ab} изображены локальные преобразования $7$, $8$, $9$ и $10$.

\medskip

Характеризация очевидна. Восстановление кода ясно, так как у $F$ и $J$ -- общее множество начальников, у $Y$ первый начальник -- $Z$, и тип второго во втором случае  -- $\mathbb{A}$, а $Z$ -- краевая вершина без начальников.

\medskip

\subsection{Случай цепи $\mathbb{CDL}2$; преобразования 8 и 9}

\begin{figure}[hbtp]
\centering
\includegraphics[width=0.9\textwidth]{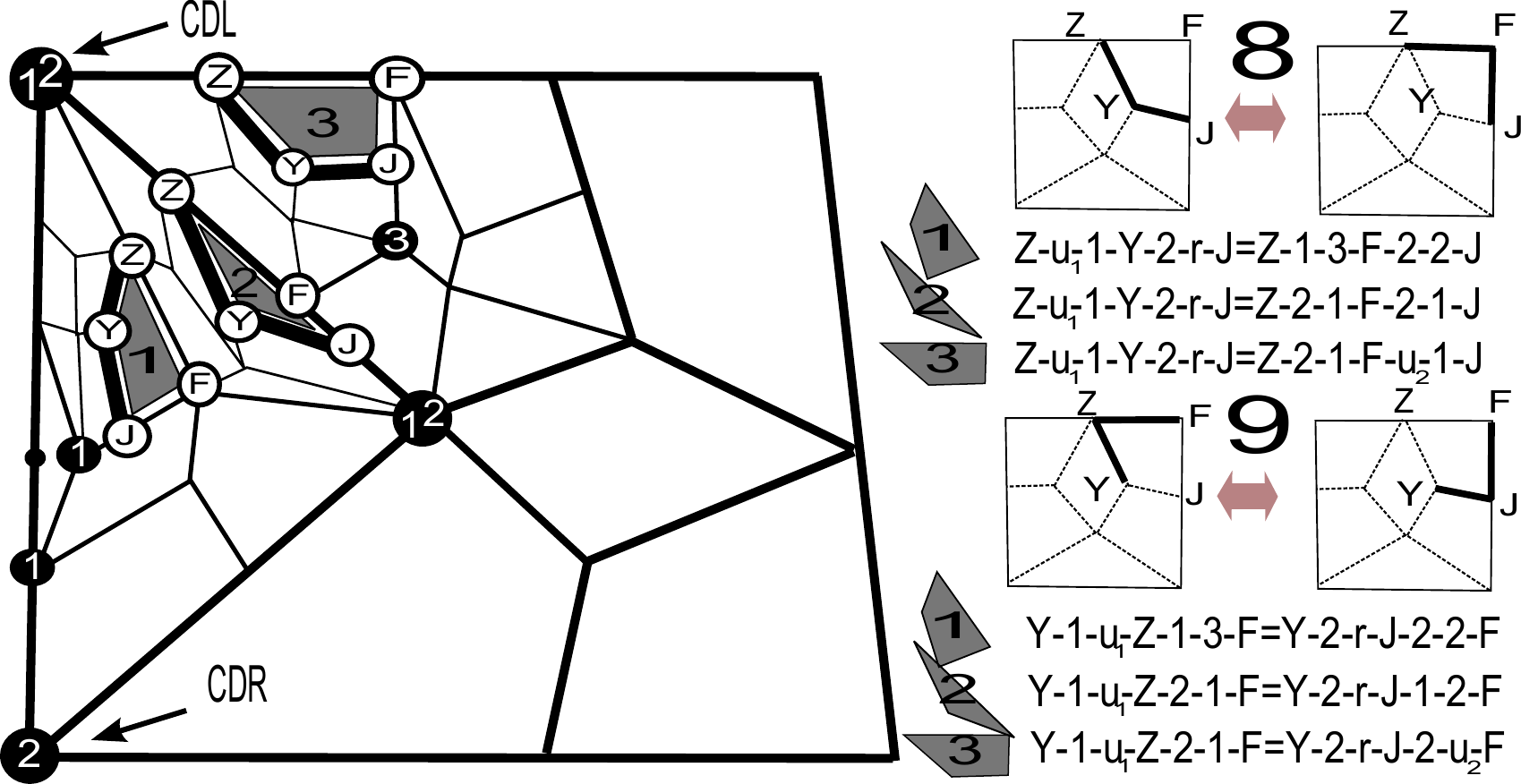}
\caption{Случаи расположения пути вокруг $\mathbb{CDL}2$-цепи и определяющие соотношения для  локальных преобразований 8 и 9}
\label{CDL2b}
\end{figure}

В правой части рисунка~\ref{CDL2b} изображены локальные преобразования $8$ и $9$.

\medskip

Характеризация очевидна. Восстановление кода также ясно, учитывая, что у $F$ и $Z$ всегда общее множество начальников, у $Y$ первый начальник  -- $Z$, а тип второго ясен из расположения. Начальники $J$ во всех случаях либо содержатся среди начальников $F$, либо сам $F$ является начальником.

\medskip

\subsection{Случай цепи $\mathbb{CDL}2$; преобразования 7 и 10}

\begin{figure}[hbtp]
\centering
\includegraphics[width=0.9\textwidth]{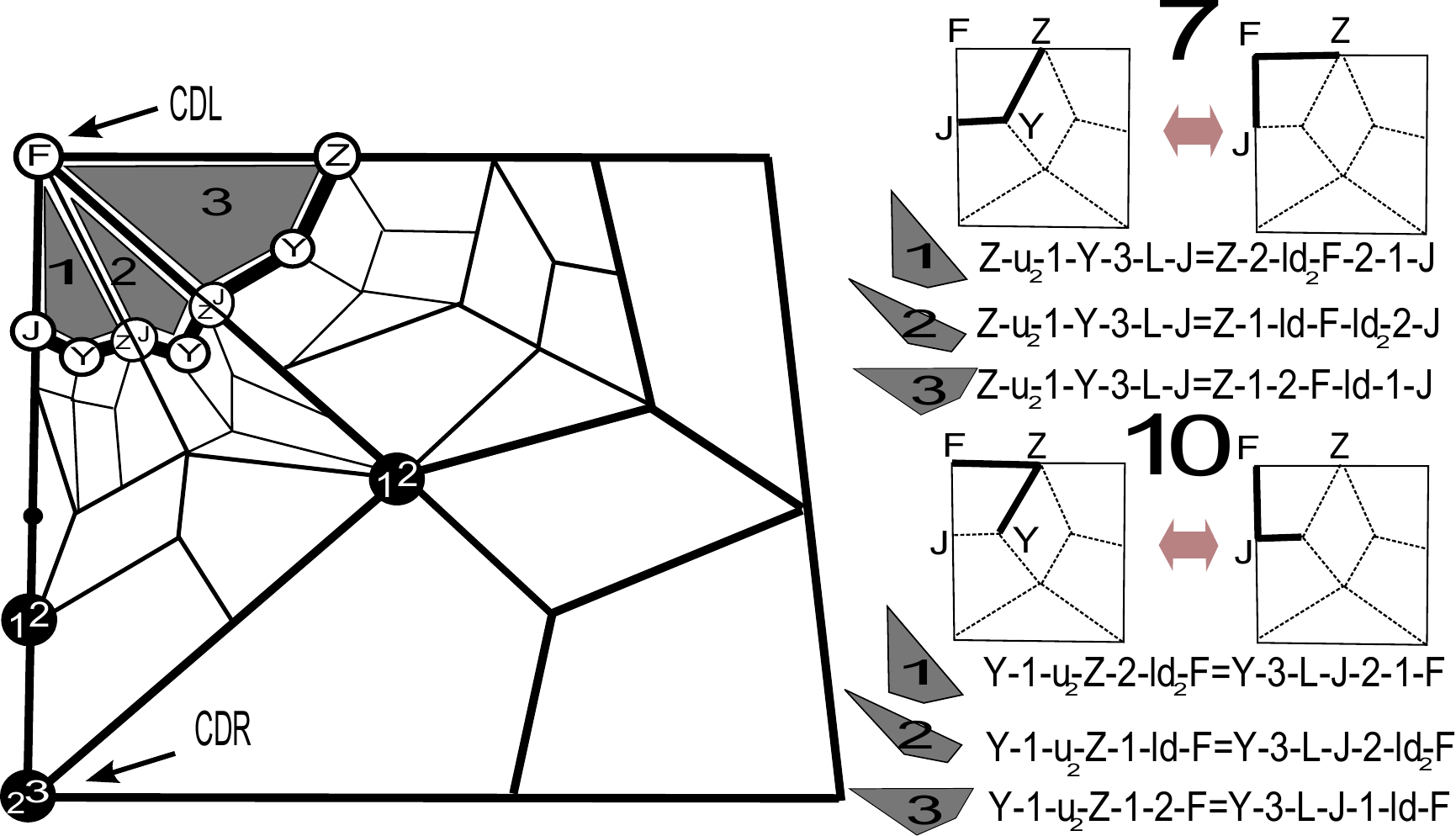}
\caption{Случаи расположения пути вокруг $\mathbb{CDL}2$-цепи и определяющие соотношения для локальных преобразований 7 и 10}
\label{CDL2a}
\end{figure}

В правой части рисунка~\ref{CDL2a} изображены локальные преобразования $7$ и $10$.

\medskip

Характеризация очевидна. Для доказательства восстановления кода достаточно заметить, что окружения четырех вершин в каждом случае очевидны, а окружение начальников (вершины, отмеченных черными кругами), можно явным образом написать.

\medskip

\subsection{Случай цепи $\mathbb{CDL}3$; преобразования 7 и 10}

Случай $\mathbb{CDL}3$-цепи полностью аналогичен случаю $\mathbb{CDL}2$-цепи, соотношения выглядят идентично, только кодировки вершин $J$, $F$, $Z$, $Y$ отвечают $\mathbb{CDL}3$-цепи. Все рассуждения о вычислении путей полностью аналогичны. Соотношений вводится столько же, сколько для $\mathbb{CDL}2$ случая.

\medskip

\subsection{Случай цепи $\mathbb{CDR}1$; преобразования 8 и 9}

\begin{figure}[hbtp]
\centering
\includegraphics[width=0.9\textwidth]{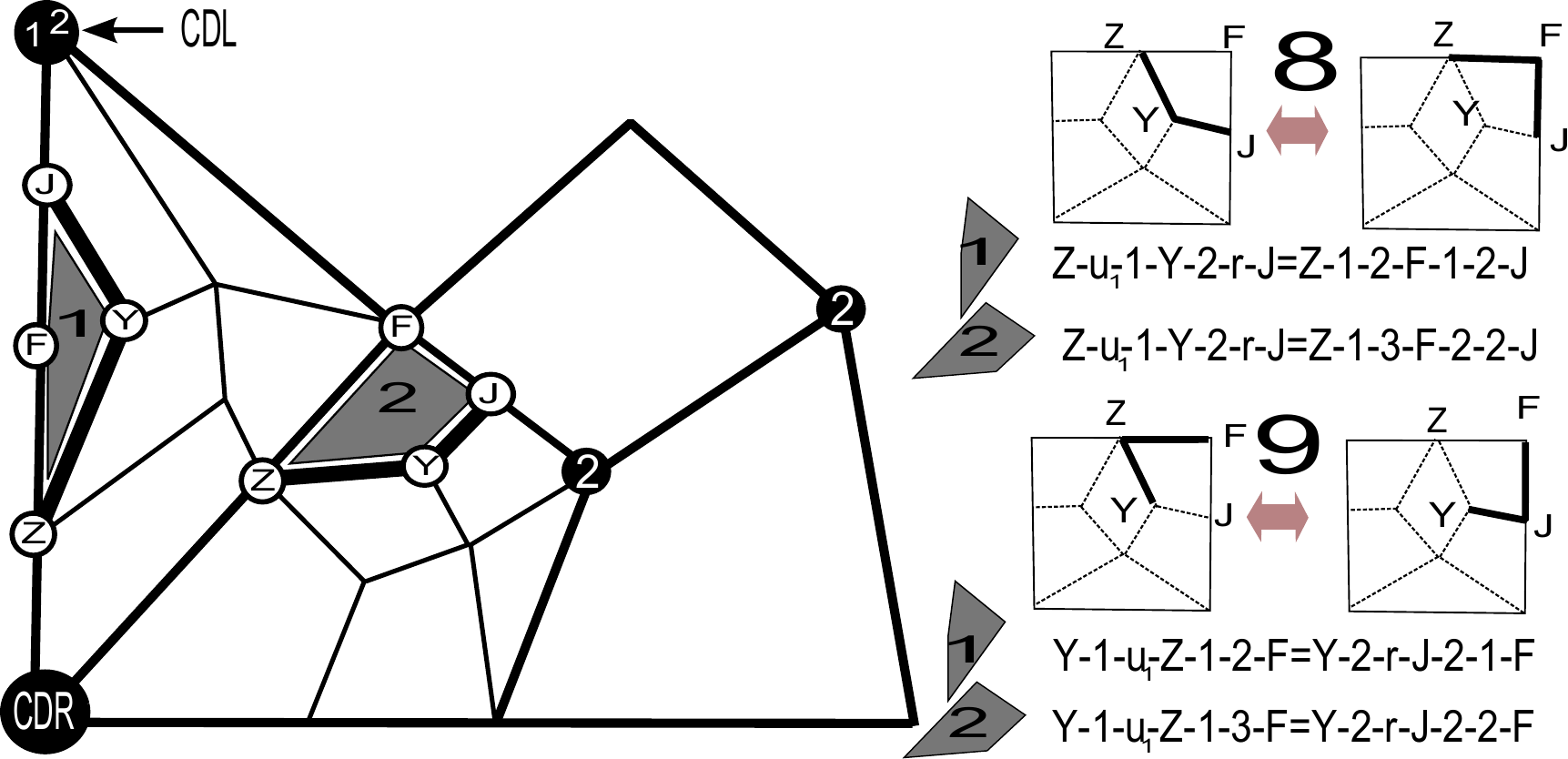}
\caption{Случаи расположения пути вокруг $\mathbb{CDR}1$-цепи и определяющие соотношения для локальных преобразований 8 и 9}
\label{CDR1b}
\end{figure}

В правой части рисунка~\ref{CDR1b} изображены локальные преобразования $8$ и $9$.

\medskip

Характеризация очевидна. Для восстановления кода достаточно заметить, что начальник каждой вершины либо является другой вершиной среди четырех, либо является начальником другой вершины из четырех.

\medskip

\subsection{Случай цепи $\mathbb{CDR}1$; преобразования 7 и 10}

\begin{figure}[hbtp]
\centering
\includegraphics[width=0.9\textwidth]{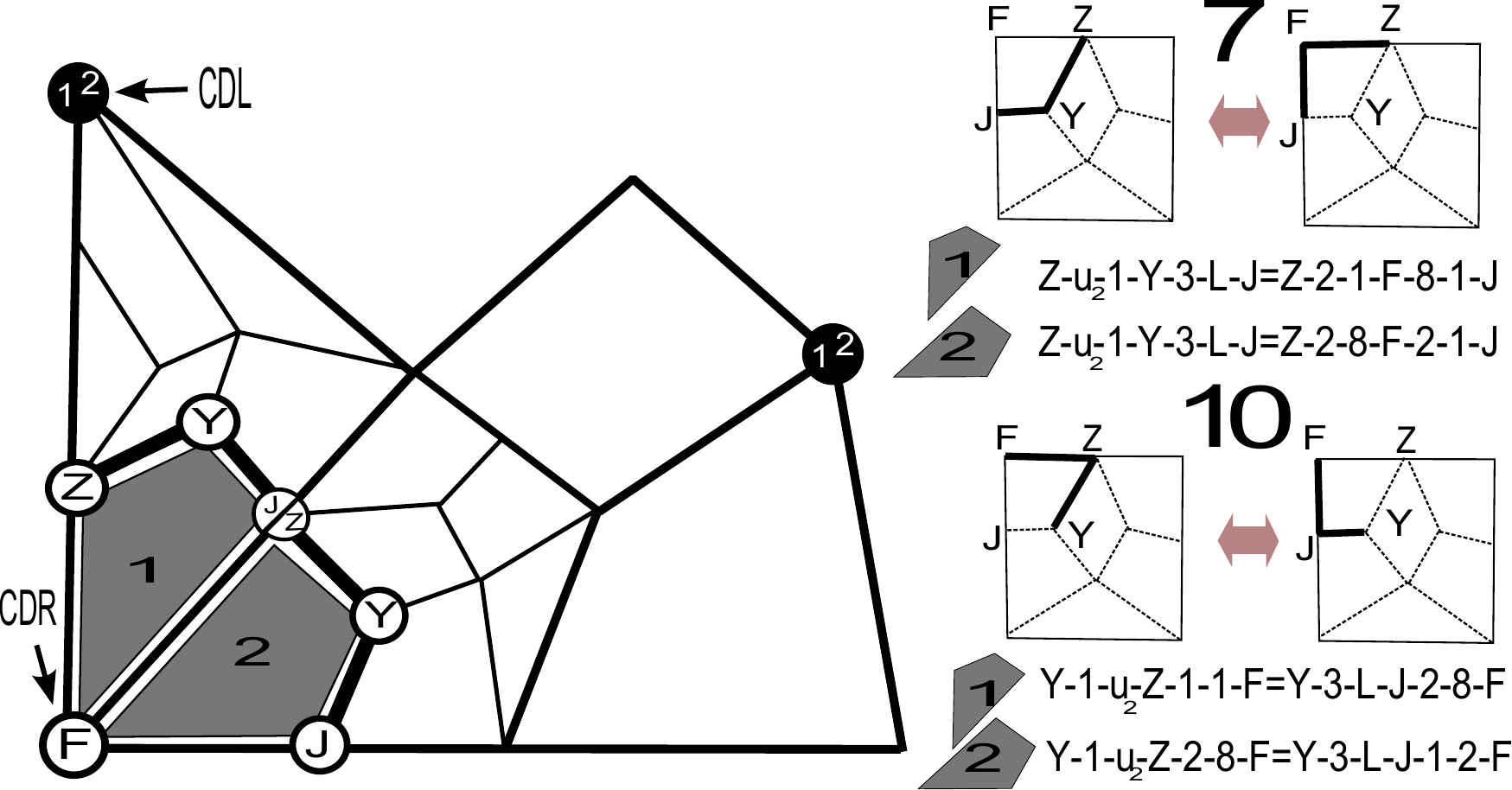}
\caption{Случаи расположения пути вокруг $\mathbb{CDR}1$-цепи и соответствующие им определяющие соотношения для локальных преобразований 7 и 10}
\label{CDR1a}
\end{figure}

В правой части рисунка~\ref{CDR1a} изображены локальные преобразования $7$ и $10$.

\medskip

Характеризация очевидна. Для восстановления кода достаточно заметить, что начальник каждой вершины либо является другой вершиной среди четырех, либо является начальником другой вершины из четырех.

\medskip

\subsection{Случай цепи $\mathbb{CDR}2$; преобразования 8 и 9}

\begin{figure}[hbtp]
\centering
\includegraphics[width=0.9\textwidth]{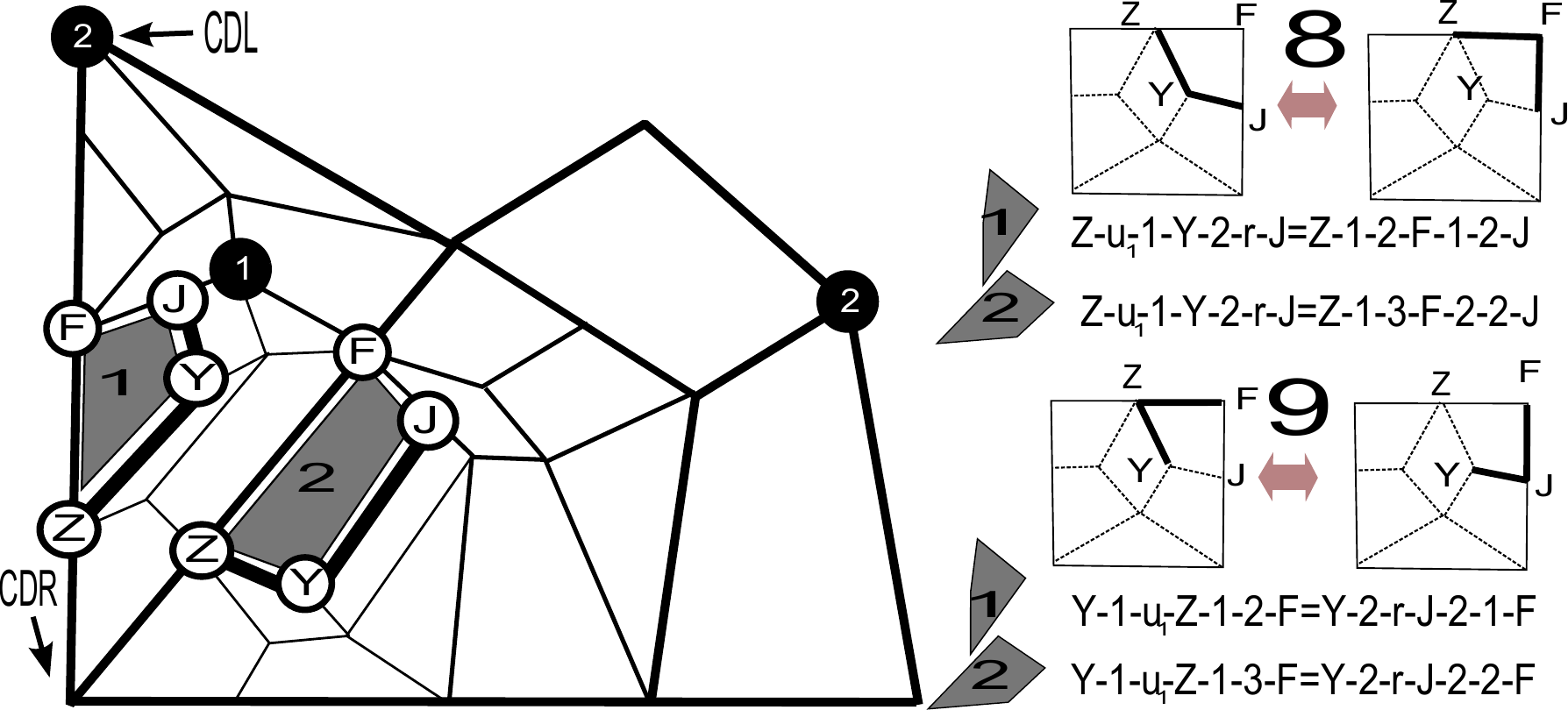}
\caption{Случаи расположения пути вокруг CDR2-цепи  и определяющие соотношения для локальных преобразований 8 и 9}
\label{CDR2b}
\end{figure}

В правой части рисунка~\ref{CDR2b} изображены локальные преобразования $8$ и $9$.

\medskip

Характеризация очевидна. Для восстановления кода достаточно заметить, что начальник каждой вершины либо является другой вершиной среди четырех, либо является начальником другой вершины из четырех.

\medskip

\subsection{Случай цепи $\mathbb{CDR}2$; преобразования 7 и 10}

\begin{figure}[hbtp]
\centering
\includegraphics[width=0.9\textwidth]{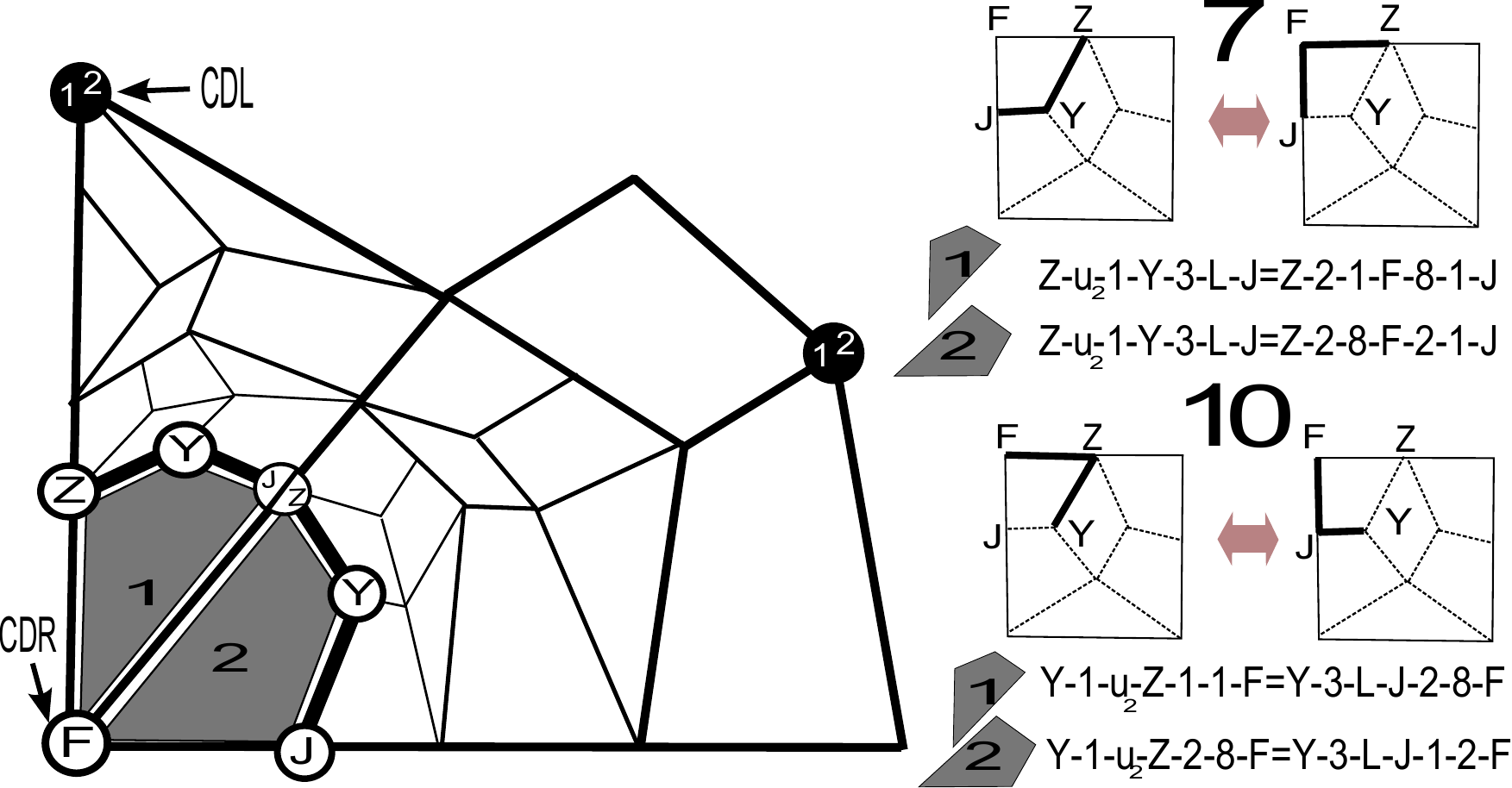}
\caption{Случаи расположения пути вокруг $\mathbb{CDR}2$-цепи и определяющие соотношения для локальных преобразований 7 и 10}
\label{CDR2a}
\end{figure}

В правой части рисунка~\ref{CDR2a} изображены локальные преобразования $7$ и $10$.

\medskip

Характеризация очевидна. Для восстановления кода достаточно заметить, что начальник каждой вершины либо является другой вершиной среди четырех, либо является начальником другой вершины из четырех.

\medskip

\subsection{Подсчет введенных соотношений}
Пусть $F$ -- число различных флагов макроплиток, $P$ -- число различных подклееных окружений, $I$ -- число различных информаций. Пусть также  $\mathbf{Num}(\mathbb{A})$ -- число базовых окружений вершин типа $\mathbb{A}$, аналогично для других типов.
Крайние вершины в пути могут иметь произвольное подклееное окружение, кроме того, соотношения вводятся для каждого заданного флага подклейки.

Для каждого из  локальных преобразований $(1)-(6)$ мы фиксировали один узел с окружением, флагом макроплитки и информацией. Также еще для одного узла могла быть зафиксирована информация. После чего, в каждом из шести случаев расположения мы вводили не более $2$ соотношений для каждого сочетания подклееных окружений крайних узлов пути и флага макроплитки. То есть общее число соотношений для локальных преобразований $(1)-(6)$ не превосходит $72 \cdot I\cdot F \cdot P\cdot \mathbf{Num}$, где $\mathbf{Num}$ -- число всевозможных сочетаний типа, уровня, окружения и информации плоской вершины.
То есть для локальных преобразований $(1)-(6)$ соотношений было введено не более

$$72\cdot 181 \cdot 10^6 \cdot 5 \cdot 10^{18} \cdot 39 \cdot 16 \cdot 10^{12} < 41 \cdot 10^{42}.$$

Теперь посчитаем число введенных соотношений для  локальных преобразований ~$7$, $8$, $9$, $10$.
Напомним, что мы рассматривали вершины из цепей, фиксировали для них информацию, флаг подклейки, а также два подклееных окружения. После этого вводилось не более двух соотношений для каждого локального преобразования.
Таким образом, число введенных соотношений для локальных преобразований ~$7$, $8$, $9$, $10$ не более $8 \cdot F \cdot P^2 \cdot N$, где $N$ -- общее число вершин всех типов, всех базовых окружений, всех возможных информаций.
Учитывая $F \le 16\cdot 10^{18}$ и $P=39$, получаем, что сочетаний базовых окружений и информаций не более $2749 \cdot 181 \cdot 10^6$. Значит число введенных соотношений для ~$7$, $8$, $9$, $10$ локальных преобразований не превосходит $13 \cdot 10^{33}$.

\section{Локальные преобразования при выходе в подклееную макроплитку} \label{pasting_section}

В предыдущей главе мы описали, как вводятся локальные преобразования (аналог определяющих соотношений) для плоских участков путей. В этом главе мы опишем, как это делается для путей, содержащих выход в подклееную макроплитку. То есть мы рассматриваем произвольный путь $X_1e_1e_2X_2e_3e_4X_3$, где $X_1$, $X_2$, $X_3$ -- буквы, отвечающие кодам вершин, а $e_1$, $e_2$, $e_3$, $e_4$ -- буквы, отвечающие ребрам входа и выхода. Хотя бы одно из ребер входа и выхода отвечает ребру, выходящему в подклейку.

\medskip

Для начала докажем следующее утверждение, позволяющее узнать ребра входа и выхода в пути из трех узлов, для случая когда средний узел является ядром некоторой подклейки.

\medskip

{\bf Ядро подклейки.} Пусть $Y$ -- ядро некоторой подклееной макроплитки $T$. Напомним, что в каждом узле $T$ хранится параметр {\it ``флаг подклейки''}, содержащий в себе полный код (тип, уровень, окружение и информацию) соответствующего ядра макроплитки $Y$, а также сочетания типов выходящих из $Y$ ребер, которые лежат на сторонах $T$. например, если $T$ подклеена по ребрам $1$ и $3$ к узлу типа $\mathbb{A}$ c окружением $(\mathbf{7A},\mathbf{7A},\mathbf{3B},\mathbf{4A})$, информацией $\mathbf{Info}$, то параметр {\it флаг подклейки} для узлов $T$ выглядит как $[\mathbb{A},\mathbf{Info},(\mathbf{7A},\mathbf{7A},\mathbf{3B},\mathbf{4A}),1,3]$.
Этот параметр (для узла $X$ внутри $T$) обозначается как $\mathbf{Core}(X)$. Тип выходящего из ядра ребра, соответствующего верхней стороне $T$ обозначается как $\mathbf{T.Core}(X)$. Соответственно, для ребра являющегося левой стороной $T$:  $\mathbf{L.Core}(X)$. Для окружения или типа $Y$ мы используем напрямую обозначение $\mathbf{Core}(X)$.

\medskip

\begin{lemma}[О ребрах около ядра макроплитки.]   \label{tilecore}

Пусть $XYZ$ -- некоторый плоский путь, причем $Y$ -- узел-ядро подклееной макроплитки, на сторонах которой лежат $X$ и $Z$. Тогда ребро выхода из $X$ (это всегда ребро $1$ или ребро $2$) можно установить, зная ребро входа в $Y$. Аналогично, ребро входа в $Z$ устанавливается по известному ребру выхода из $Y$.
\end{lemma}

Напомним, у каждого внутреннего ребра есть А-сторона и B-сторона. При этом первое главное ребро -- то, для которого A-сторона лежит справа.
В соответствии с этим правилом и с определением кодировок ребер входов и выходов можно выписать ребра выхода из $X$ и входа в $Z$ (приведены таблицы~\ref{tableAexits},~\ref{tableBexits},~\ref{tableCexits},~\ref{tableSideexits} для разных типов $Y$):

\begin{table}[hbtp]
\caption{Ребра выхода для $\mathbb{A}$-узла}
\centering
 \begin{tabular}{|c|c|c|c|c|}  \hline
ребро входа в $Y$    & ребро выхода из $X$  \cr
(или выхода из $Y$)  & (или входа в $Z$)  \cr  \hline
$1$ & $1$  \cr \hline
$2$ & $2$  \cr \hline
$3$ & $2$  \cr \hline
$\mathbf{lu}$  & $2$  \cr \hline
$\mathbf{ld}$ & $2$  \cr \hline
  \end{tabular}
\label{tableAexits}
\end{table}

\begin{table}[hbtp]
\caption{Ребра выхода для $\mathbb{B}$-узла}
\centering
 \begin{tabular}{|c|c|c|c|c|}  \hline
ребро входа в $Y$    & ребро выхода из $X$    \cr
(или выхода из $Y$)  & (или входа в $Z$)  \cr  \hline
$1$ & $2$  \cr \hline
$2$ & $1$  \cr \hline
$3$ & $1$  \cr \hline
$\mathbf{ru}$ & $2$  \cr \hline
$\mathbf{mid}$ & $2$  \cr \hline
$\mathbf{rd}$ & $2$  \cr \hline
  \end{tabular}
\label{tableBexits}
\end{table}

\begin{table}[hbtp]
\caption{Ребра выхода для $\mathbb{C}$-узла}
\centering
 \begin{tabular}{|c|c|c|c|c|}  \hline
ребро входа в $Y$    & ребро выхода из $X$    \cr
(или выхода из $Y$)  & (или входа в $Z$)  \cr  \hline
$1$ & $1$  \cr \hline
$2$ & $2$  \cr \hline
$3$ & $1$  \cr \hline
$4$ & $2$  \cr \hline
$\mathbf{ld}_1$ & $1$  \cr \hline
$\mathbf{ld}_2$ & $2$  \cr \hline
$\mathbf{mid}_1$ & $1$  \cr \hline
$\mathbf{mid}_2$ & $2$  \cr \hline
$\mathbf{d}_1$ & $1$  \cr \hline
$\mathbf{d}_2$ & $2$  \cr \hline
  \end{tabular}
\label{tableCexits}
\end{table}

\begin{table}[hbtp]
\caption{Ребра выхода для $\mathbb{UL}$, $\mathbb{UR}$, $\mathbb{DR}$, $\mathbb{DL}$ узлов}
\centering
 \begin{tabular}{|c|c|c|c|c|}  \hline
ребро входа в $Y$    & ребро выхода из $X$    \cr
(или выхода из $Y$)  & (или входа в $Z$)  \cr  \hline
1 & 2  \cr \hline
2 & 1  \cr \hline
$\mathbf{u}_1$ & $1$  \cr \hline
$\mathbf{u}_2$ & $2$  \cr \hline
$\mathbf{u}_3$ & $1$  \cr \hline
$\mathbf{u}_4$ & $2$  \cr \hline
$\mathbf{l}$ & $1$  \cr \hline
$\mathbf{l}_2$ & $1$  \cr \hline
$\mathbf{l}_3$ & $2$  \cr \hline
$\mathbf{r}$ & $2$  \cr \hline
$\mathbf{r}_2$ & $1$  \cr \hline
$\mathbf{r}_3$ & $2$  \cr \hline
  \end{tabular}
\label{tableSideexits}
\end{table}

Таким образом, зная хотя бы один узел подклееной макроплитки $T$ (из тех, которые лежат внутри нее, то есть не на верхней или левой стороне), мы можем узнать ребра входа и выхода для путей вдоль левой или верхней стороны $T$. В дальнейшем, если нам потребуется выписать типы ребер вдоль левой или верхней стороны подклееной макроплитки, мы будем использовать ссылку на предложение~\ref{tilecore}.

\medskip

\begin{lemma}[Об информации узла рядом с ядром макроплитки]   \label{tilecore2}

Пусть $XY$ -- некоторый плоский путь, причем $Y$ -- узел-ядро подклееной макроплитки, на одной из сторон которой лежит $X$. Будем считать, что нам известен код $Y$ (тип, уровень, окружение и информация), а также тип входящего в $Y$ ребра и тип, уровень и окружение $X$. Тогда информация у $X$ может быть восстановлена по этим данным.
\end{lemma}

Доказательство. Допустим сначала, что $X$ имеет один из внутренних типов ($\mathbb{A}$, $\mathbb{B}$, $\mathbb{C}$).

Пусть тип $X$ это $\mathbb{A}$. Рассмотрим ребро входа в $Y$. Имеется пять возможных типов ребер для узла типа  $\mathbb{A}$, это ребра $1$, $2$, $3$, $\mathbf{lu}$ , $\mathbf{ld}$. Если ребро входа в $Y$ это $1$, $2$ или $3$, информация у $X$ будет такой же, как у $Y$. Для ребра $\mathbf{lu}$ , первым начальником $X$ будет $\mathbf{LevelPlus.FBoss}(Y)$, вторым -- $\mathbf{Next.FBoss}(Y)$, третьим -- сам узел $Y$. В случае $\mathbf{ld}$, первым начальником $X$ будет $\mathbf{BottomLeftChain}(Y)$, вторым -- $\mathbb{C}$, с окружением $Y$, третьим -- сам узел $Y$.

 \medskip

Пусть тип $Y$ это $\mathbb{B}$. Имеется шесть возможных ребер из узла типа $\mathbb{B}$, это ребра $1$, $2$, $3$, $\mathbf{ru}$, $\mathbf{rd}$, $\mathbf{mid}$. Если ребро входа в $Y$ имеет тип $1$, $2$ или $3$, то информация у $X$ такая же как у $Y$.
Для ребра $\mathbf{rd}$, первым начальником будет узел в середине ребра 7 для макроплитки $T$, где $Y$ выступает в роли  $\mathbb{B}$-узла. То есть это $1$-цепь вокруг правого нижнего угла $T$, с указателем соответствующим входу по ребру 7. Тип правого нижнего угла $T$ нам известен, так как это второй начальник для $Y$ (это входит в информацию для $Y$). Второй начальник $X$ -- это узел в середине правой стороны $T$. То есть первый начальник $Y$ это $\mathbf{TopFromRight.SBoss}(X)$. Третьим начальником $X$ будет сам узел $Y$.
Для ребра $\mathbf{ru}$, первым начальником будет $1$-цепь вокруг правого верхнего угла $T$, с указателем соответствующим входу по правому ребру. Тип правого верхнего угла $T$ мы можем определить, применив функцию $\mathbf{TopRightType}(Y)$. Второй начальник $X$ совпадает с первым начальником $Y$.  Третьим начальником $X$ будет сам узел $Y$.
Для ребра $\mathbf{mid}$, первым начальником будет $1$-цепь вокруг $\mathbb{A}$-узла $T$ (то есть окружение как у $Y$), с указателем $1$. Второй начальник $X$ это $\mathbb{C}$-узел $T$ (окружение как у $Y$).  Третьим начальником $X$ будет сам узел~$Y$.

 \medskip

Пусть тип $Y$ это $\mathbb{C}$. Имеется десять возможных ребер из узла типа $\mathbb{C}$ ($4$ главных и $6$ неглавных). Можно заметить, что во всех случаях окружения начальников $X$ могут быть восстановлены по информации и окружению $Y$, этот процесс полностью аналогичен описанному выше.

\medskip

Пусть теперь $X$ имеет боковой тип. Если входящее в $Y$ ребро -- главное, то $X$ и $Y$ лежат на одном ребре в некоторой макроплитке, и поэтому у них одинаковая информация.
Пусть входящее в $Y$ ребро -- неглавное. Тогда тип ребра может быть одним из следующих: $\mathbf{u}_1$, $\mathbf{u}_1$, $\mathbf{u}_1$, $\mathbf{u}_1$, $\mathbf{r}$, $\mathbf{r}_2$, $\mathbf{r}_3$, $\mathbf{l}$, $\mathbf{l}_2$, $\mathbf{l}_3$.
Для $\mathbf{u}_1$, $\mathbf{u}_2$ первым начальником $X$ будет $Y$, при этом для $\mathbf{u}_1$ тип второго начальника для $X$ может быть восстановлен по информации и окружению $Y$.
Для $\mathbf{l}$, первым и единственным начальником $X$ будет $\mathbf{Prev}(Y)$, а для $\mathbf{r}$ первый начальник $X$ может быть восстановлен по процедуре  $\mathbf{TopFromRight}(Y)$, тип второго -- по процедуре  $\mathbf{BottomRightTypeFromRight}(Y)$.
Для остальных ребер все начальники $X$ восстанавливаются аналогично случаям внутренних вершин $Y$.

\medskip

\subsection{Обзор случаев}

Пусть есть путь, участвующий в локальном преобразовании, лежащий в макроплитке $T$, причем одно из его ребер -- входящее в подклееную макроплитку. Нужно показать, как происходит преобразование во всех случаях.

\medskip

Заметим, что один конец ребра в подклееную область всегда лежит на верхней или левой стороне макроплитки.

Итак, нужно рассмотреть случаи локальных преобразований, где макроплитка $T$ примыкает к границе левой или верхней стороны подклееной макроплитки. Таким образом, нужно рассмотреть те из случаев $B1-B20$, которые включают граничные стороны $\mathbf{top}$ или $\mathbf{left}$, таблица~\ref{P1P10}.

\medskip
\begin{table}[hbtp]
\caption{Десять окружений макроплиток, содержащих граничные стороны $\mathbf{top}$ или $\mathbf{left}$.}
\centering
 \begin{tabular}{|c|c|c|c|}   \hline
\x{старое \cr обозначение} & \x{окружение \cr макроплитки} & & \x{новое \cr обозначение} \cr \hline
B1.  & $(\mathbf{left},\mathbf{top},\mathbf{right},\mathbf{bottom})$  & & P1 \cr \hline
B2. & $(\mathbf{left},\mathbf{top},\mathbf{1A},\mathbf{3A})$   & & P2 \cr \hline
B3.  & $(\mathbf{7A},x,\mathbf{3B},\mathbf{4A})$,   &$x=\mathbf{left},\mathbf{top}$ & P3    \cr \hline
B4.  & $(\mathbf{7A},x,\mathbf{1A},\mathbf{3A})$,  & $x=\mathbf{left},\mathbf{top}$  & P4 \cr \hline
B5. & $(\mathbf{top},\mathbf{right},\mathbf{6A},\mathbf{2A})$  & & P5 \cr \hline
B6. & $(\mathbf{top},\mathbf{right},\mathbf{1A},\mathbf{3A})$  &  & P6  \cr \hline
B11. & $(x,\mathbf{3B},\mathbf{6A},\mathbf{2A})$,   & $x=\mathbf{left},\mathbf{top}$ & P7 \cr  \hline
B12. & $(x,\mathbf{3B},\mathbf{1A},\mathbf{3A})$,   & $x=\mathbf{left},\mathbf{top}$ & P8 \cr  \hline
B17. & $(x,\mathbf{1A},\mathbf{6A},\mathbf{2A})$,  & $x=\mathbf{left},\mathbf{top}$ & P9  \cr \hline
B18. & $(x,\mathbf{1A},\mathbf{1A},\mathbf{3A})$,  & $x=\mathbf{left},\mathbf{top}$ & P10 \cr \hline
  \end{tabular}
\label{P1P10}
\end{table}

\medskip

Пронумеруем их как $P1-P10$ и рассмотрим каждый отдельно.

\medskip

\begin{figure}[hbtp]
\centering
\includegraphics[width=0.5\textwidth]{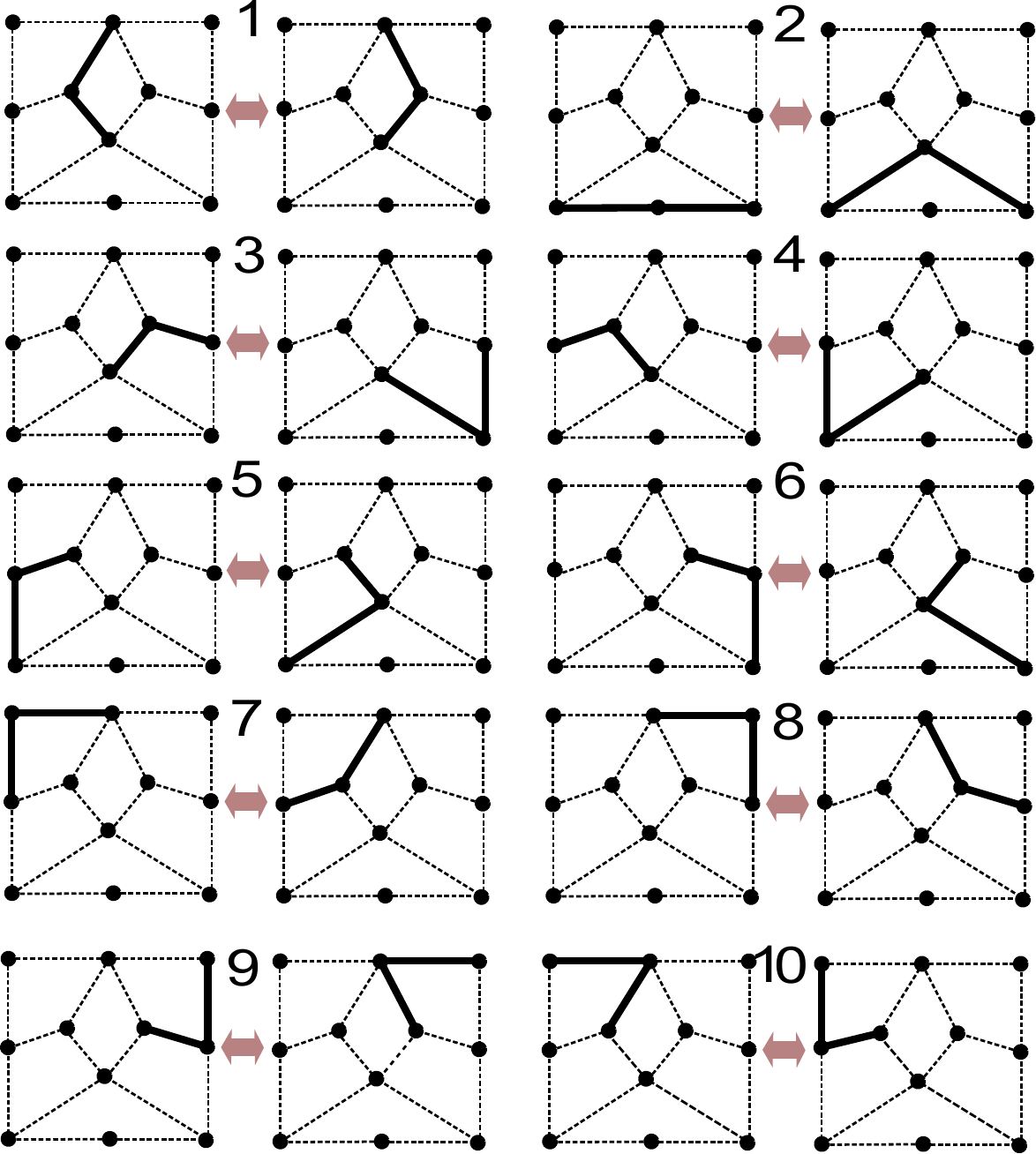}
\caption{Локальные преобразования}
\label{flippaths2}
\end{figure}

\medskip

Для каждого из этих случаев мы рассмотрим те локальные преобразования, где хотя бы один путь имеет узел на верхней или левой стороне. Это преобразования $1$, $2$, $4$, $5$, $7$, $8$, $9$, $10$ (рисунок~\ref{flippaths2}). В остальных случаях ($3$ и $6$) путь не может содержать ребро в подклееную область, так как пути не касаются левой или верхней стороны.

\medskip

{\bf Замечание.} Ниже мы не будем выписывать симметричное к введенному соотношение, отвечающее проходу пути в обратном порядке. Просто будем считать, что соотношений вводится в два раза больше.
Для экономии места в таблицах ниже, ``начальников'' будем называть ``боссами''.

\medskip

\subsection{P1: Макроплитка с окружением $(\mathbf{left},\mathbf{top},\mathbf{right},\mathbf{bottom})$}

В этом параграфе мы рассмотрим случай {\bf P1}, когда макроплитка $T$, внутри которой мы рассматриваем локальное преобразование пути, имеет окружение $(\mathbf{left},\mathbf{top},\mathbf{right},\mathbf{bottom})$. Это происходит, когда макроплитка сама по себе является только что подклееной макроплиткой.

Зафиксируем узел $Q$, являющийся ядром нашей подклееной макроплитки. Его подклееное окружение это $(\mathbf{left},\mathbf{top},\mathbf{right},\mathbf{bottom})$.  У него может быть произвольный тип, кроме углового, произвольное базовое окружение и произвольная информация. Если это боковой или краевой узел, то уровень у него должен быть третий (иначе бы к нему не была подклеена макроплитка).  Сторонами подклейки могут быть любые два ребра, выходящие из нашего узла.

Напомним, что в параметр {\it ядро подклейки} входит тип, уровень, расширенное окружение, информация и типы двух выходящих ребер, являющихся сторонами подклейки. То есть наше ядро подклейки задано.

\medskip

\begin{figure}[hbtp]
\centering
\includegraphics[width=0.9\textwidth]{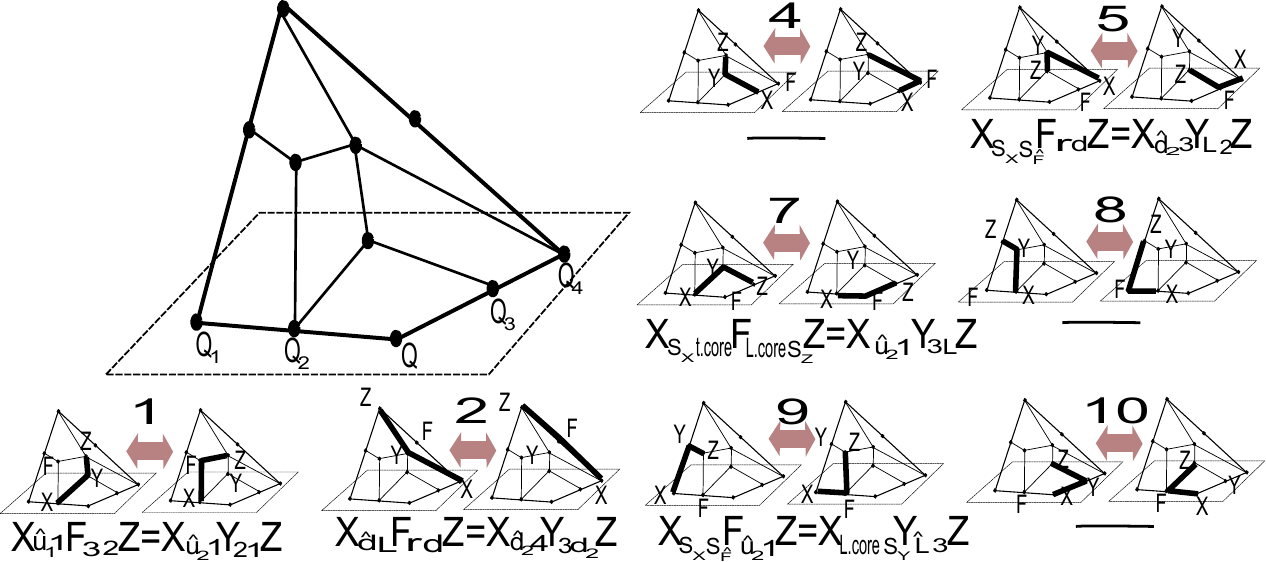}
\caption{Случай P1}
\label{P1}
\end{figure}

\medskip

Определив $Q$, можно также определить и узлы $Q_1$, $Q_2$, $Q_3$, $Q_4$, лежащие на сторонах подклейки. То есть, мы рассматриваем всевозможные сочетания кодов для пятерки вершин, существующих на комплексе. Для каждой такой пятерки, можно определить буквы коды вершин $X$, $Y$, $F$, $Z$, для каждого из указанных на рисунке~\ref{P1} восьми случаев. Также введем обозначения $e_{\mathbf{\mathbf{t.core}}}$ и $e_{\mathbf{\mathbf{l.core}}}$ для типов двух выходящих ребер, являющихся сторонами подклейки.

\medskip

{\bf Характеризация.} По выданному нам коду пути $P_1P_2P_3$, включающему переход по подклееному ребру, мы можем установить, встречается ли этот путь среди перечисленных на рисунке~\ref{P1}. Действительно, первый узел после перехода по подклееному ребру должен содержать параметр {\it ядро подклейки}, соответствующий рассматриваемому нами. Ясно также, что по конфигурации входящих-выходящих ребер и типов вершин можно выявить, какой именно случай из шестнадцати представленных на рисунке~\ref{P1} имеет место.

Итак, мы можем определить по коду пути к какому именно случаю расположения он относится. Теперь покажем, как по этому коду определить код парного к нему пути.

В случаях локальных преобразований $4$, $8$ и $10$ путь удовлетворяет условиям мертвого паттерна и не может быть частью достаточно длинного ненулевого пути, и для этих локальных преобразований мы соотношений не вводим.

\medskip

{\bf Вычисление кода $XFZ$.} Зная код пути $XYZ$ (то есть коды $X$, $Y$, $Z$, а также ребра входа в $Y$, $Z$ и ребра выхода из $X$ и $Y$ вдоль пути $XYZ$), можно выписать код $XFZ$.

Заметим, что коды $X$ и $Z$ мы и так знаем (из пути $XYZ$). Остается выписать ребра выхода из $X$, ребро входа в $F$, код $F$, ребро выхода из $F$ и ребро входа в $Z$ (таблица~\ref{P1a}).

\medskip

\begin{table}[hbtp]
\caption{Случай P1: вычисление кода $XFZ$. }
\centering
 \begin{tabular}{|c|c|c|c|c|c|c|}   \hline
  & ребро из $X$ & ребро в $F$ & \x{тип, уровень \cr и окружение $F$} & \x{информация \cr $F$} & \x{ ребро \cr из $F$ }& \x{ребро \cr в $Z$} \cr \hline
1. & $\widehat{\mathbf{u}_1}$ & $1$ & \x{ $\mathbb{B}$  (окруж. как у $Y$)} & 1 босс -- $X$ & $3$ & $2$  \cr \hline
2. & $\widehat{\mathbf{d}}$ & $\mathbf{l}$ & \x{$\mathbb{D}$, $1$ \cr (окруж. как у $Y$)} & -- & $\mathbf{r}$ & $\mathbf{d}$  \cr \hline
5. & Предл.\ref{tilecore} & Предл.\ref{tilecore} & $\mathbf{Next.FBoss}(Z)$ & как у $X$   & $\mathbf{l}$ & $3$  \cr \hline
7. & Предл.\ref{tilecore} & $\mathbf{T.Core}(Y)$ & $\mathbf{Core}(Y)$ & $\mathbf{Core}(Y)$ & $\mathbf{L.Core}(Y)$ & Предл.\ref{tilecore}   \cr \hline
9. & Предл.\ref{tilecore} & Предл.\ref{tilecore} & $\mathbf{FBoss}(Z)$& как у $X$  & $\widehat{\mathbf{u}_1}$ & $1$  \cr \hline
  \end{tabular}
\label{P1a}
\end{table}

\leftskip=0.0cm

\medskip

{\bf Примечание.} Через $\mathbf{T.Core}(Y)$ и $\mathbf{L.Core}(Y)$ мы обозначаем типы выходящих ребер, соответствующих сторонам подклееной макроплитки. Они содержатся в параметре флаг подклейки у вершины $Y$.

\medskip

{\bf Вычисление кода $XYZ$.} Зная код пути $XFZ$ (то есть коды $X$, $F$, $Z$, а также ребра входа в $F$, $Z$ и ребра выхода из $X$ и $F$ вдоль пути $XFZ$), можно выписать код $XYZ$.

Заметим, что коды $X$ и $Z$ мы и так знаем (из пути $XFZ$). Остается выписать ребра выхода из $X$, ребро входа в $Y$, код $Y$, ребро выхода из $Y$ и ребро входа в $Z$ (таблица~\ref{P1b}).

\medskip
\begin{table}[hbtp]
\caption{Случай P1: вычисление кода $XYZ$. }
\centering
 \begin{tabular}{|c|c|c|c|c|c|c|}   \hline
  & ребро из $X$ & ребро в $Y$ & \x{тип, уровень \cr и окружение $Y$} & информация $Y$ &  \x{ребро \cr из $Y$} & \x{ребро \cr в $Z$} \cr \hline
1. & $\widehat{\mathbf{u}_2}$ & $1$ & $\mathbb{A}$ (окружение как у $F$)& как у $F$ & $2$ & $1$  \cr \hline
2. & $\widehat{\mathbf{d}_2}$ & $4$ & $\mathbb{C}$ (окружение как у $F$)& $1$ босс -- как у $Z$  & $3$ & $\mathbf{d}_2$  \cr
&&&& $2$ босс -$X$, $3$ -$Z$  && \cr  \hline
5. & $\widehat{\mathbf{d}_2}$ & $3$ & $\mathbb{C}$ (окружение как у $Z$) & 1 босс как у $Z$   & $\mathbf{l}$ & $2$  \cr
&&&& $2$ босс -$X$, $3$ -- $\mathbb{CDR}$   &&\cr  \hline
7. & $\widehat{\mathbf{u}_2}$ & $1$ & $\mathbb{A}$, окружение -- & $1$ босс -- $X$ & $3$ & $\mathbf{l}$  \cr
&&& $(\mathbf{left},\mathbf{right},\mathbf{top},\mathbf{bottom})$  &&& \cr  \hline
9. &  $\widehat{\mathbf{r}}$ & $2$ & $\mathbb{R}$,$1$ окружение --  & --  & $\mathbf{r}$ & $2$  \cr
&&& $(\mathbf{left},\mathbf{right},\mathbf{top},\mathbf{bottom})$  &&& \cr  \hline
  \end{tabular}
\label{P1b}
\end{table}

\medskip

Когда все коды определены можно выписать определяющие соотношения, переводящие один путь в другой.

Обозначим результат применения процедуры, описанной в Предложении~\ref{tilecore} к вершине $Q$, связанной с этим ребром как $s_Q$. Как $s_{\widehat{Q}}$ обозначим противоположное главное ребро, то есть, если $s_Q=1$, то $s_{\widehat{Q}}=2$ и $s_Q=2$, то $s_{\widehat{Q}}=1$.  Обозначение $\mathbf{Past}(X)$ используем для обозначения подклееного окружения узла $X$.

\medskip

{\bf Соотношения.} Введем соотношения, реализующие описанные переходы.

\smallskip

$1)$  $X e_{\widehat{u_1}} e_1 F e_3 e_2 Z = X e_{\widehat{u_2}} e_1 Y e_2 e_1 Z   $

\medskip

$2)$  $X e_{\widehat{d}} e_l F e_r e_{d} Z = X e_{\widehat{d_2}} e_4 Y e_3 e_{d_2} Z   $

\medskip

$5)$  $X e_{s_X} e_{s_{\widehat{F}}} F e_l e_{3} Z = X e_{\widehat{d_2}} e_3 Y e_l e_{2} Z   $

\medskip

$7)$  $X e_{s_X} e_{t.core} F e_{l.core} e_{Z} Z = X e_{\widehat{u_2}} e_1 Y e_3 e_{l} Z   $

\medskip

$9)$  $X e_{s_X} e_{s_{\widehat{F}}} F e_{\widehat{u_1}} e_{1} Z = X e_{\widehat{l}} e_u Y e_r e_{2} Z   $

\medskip

\subsection{P2: Макроплитка с окружением $(\mathbf{left},\mathbf{top},\mathbf{1A},\mathbf{3A})$}

В этом параграфе мы рассмотрим случай {\bf P2}, когда макроплитка $T$, внутри которой мы рассматриваем локальное преобразование пути, имеет окружение $(\mathbf{left},\mathbf{top},\mathbf{1A},\mathbf{3A})$. Это происходит, когда макроплитка является прямым потомком подклееной макроплитки.

\begin{figure}[hbtp]
\centering
\includegraphics[width=1\textwidth]{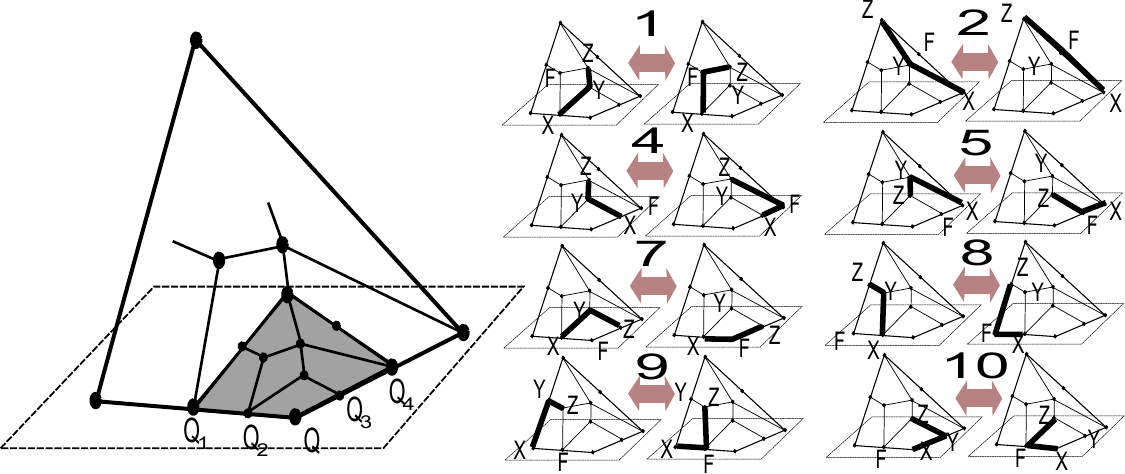}
\caption{Случай P2}
\label{P2}
\end{figure}

Аналогично случаю {\bf P1}, зафиксируем пять вершин $Q$, $Q_1$, $Q_2$, $Q_3$, $Q_4$, занимающих на комплексе положения как на рисунке~\ref{P2}. То есть $Q$ является ядром подклейки, а остальные вершины лежат на сторонах.

Для зафиксированного сочетания кодов этих пяти вершин, можно вычислить коды всех путей, указанных на рисунке~\ref{P2}. Кроме того, по выданному слову (коду пути) можно установить, действительно ли этот путь относится к случаю {\bf P2} и какую конфигурацию из изображенных на рисунке~\ref{P2} он имеет. Действительно по окружению узлов после перехода по подклееному ребру можно установить, что мы имеем дело именно со случаем $P2$, а последовательность входящих-выходящих ребер помогает установить нужную конфигурацию.

\medskip

\medskip

Теперь покажем, как по этому коду определить код парного к нему пути.

{\bf Вычисление кода $XFZ$.} Аналогично случаю {\bf P1}, зная код пути $XYZ$ (то есть коды $X$, $Y$, $Z$, а также ребра входа в $Y$, $Z$ и ребра выхода из $X$ и $Y$ вдоль пути $XYZ$), можно выписать код $XFZ$ (таблица~\ref{P2a}).

\medskip
\begin{table}[hbtp]
\caption{Случай P2: вычисление кода $XFZ$. }
\centering
 \begin{tabular}{|c|c|c|c|c|c|c|}   \hline
  & ребро из $X$ & ребро в $F$ & \x{тип, уровень \cr и окружение $F$} & \x{информация \cr $F$} &  \x{ребро \cr из $F$} & \x{ребро \cr в $Z$}
 \cr \hline
1. & $\widehat{\mathbf{u}_1}$ & $1$ & \x{$\mathbb{B}$  (окружение \cr как у $Y$)} & 1 босс -- $X$ & $3$ & $2$  \cr \hline
2. & $\widehat{\mathbf{l}}$ & $1$ & \x{ $\mathbb{DR}$,$1$, окружение \cr  $3$ ребро }& как у $Z$  & $2$ & $3$  \cr \hline
4. & Предл.\ref{tilecore} & Предл.\ref{tilecore} & $\mathbf{SBoss}(Z)$ & как у $X$  & $\mathbf{l}_2$ & $4$  \cr \hline
5. & Предл.\ref{tilecore} & Предл.\ref{tilecore} & $\mathbf{Next.FBoss}(Z)$ & как у $X$   & $\mathbf{l}$ & $3$  \cr \hline
7. & Предл.\ref{tilecore} & $\mathbf{T.Core}(Y)$ & $\mathbf{Core}(Y)$ & $\mathbf{Core}(Y)$  & $\mathbf{L.Core}(Y)$ & Предл.\ref{tilecore}  \cr \hline
8. & Предл.\ref{tilecore} & Предл.\ref{tilecore} & $\mathbf{FBoss}(Z)$& как у $X$  & $\mathbf{u}_2$ & $2$  \cr \hline
9. & Предл.\ref{tilecore} & Предл.\ref{tilecore} & $\mathbf{FBoss}(Z)$& как у $X$  & $\widehat{\mathbf{u}_1}$ & $1$  \cr \hline
10. & $\mathbf{Top.Core}(Z)$  & Предл.\ref{tilecore} & $\mathbf{FBoss}(Z)$ &  Предл.\ref{tilecore2}  & $\widehat{\mathbf{u}_2}$ & $1$  \cr \hline
  \end{tabular}
\label{P2a}
\end{table}

\medskip

{\bf Вычисление кода $XYZ$.} Зная код пути $XFZ$ (то есть коды $X$, $F$, $Z$, а также ребра входа в $F$, $Z$ и ребра выхода из $X$ и $F$ вдоль пути $XFZ$), можно выписать код $XYZ$ (таблица~\ref{P2b}).

\medskip

\begin{table}[hbtp]
  \caption{Случай P2: вычисление кода $XYZ$. }
\centering
 \begin{tabular}{|c|c|c|c|c|c|c|}   \hline
  & \x{ребро \cr из $X$} & \x{ребро \cr в $Y$} & \x{тип, уровень и \cr окружение $Y$} & информация $Y$ &  \x{ребро \cr из $Y$} & \x{ребро \cr в $Z$} \cr \hline
1. & $\widehat{\mathbf{u}_2}$ & $1$ & $\mathbb{A}$ (окружение как у $F$)& 1 босс -- $X$ & $2$ & $1$  \cr \hline
2. & $\widehat{\mathbf{d}_2}$ & $4$ & $\mathbb{C}$ (окружение как у $F$)& \x{  $1$ босс -- как у $Z$\cr $2$ босс -$X$, \cr $3$ босс -- $Z$  }  & $3$ & $\mathbf{d}_2$  \cr \hline
4. & $\widehat{\mathbf{l}}$ & $3$ & $\mathbb{A}$ (окружение как у $Z$)& 1 босс как у $Z$  & $2$ & $1$  \cr \hline
5. & $\widehat{\mathbf{d}_2}$ & $3$ & $\mathbb{C}$ (окружение как у $Z$) & \x{  1 босс как у $Z$ \cr  $2$ босс -$X$, \cr  $3$ босс -- $\mathbb{A}$, \cr окружение -- \cr $\mathbf{Past}(X)$ } & $\mathbf{l}$ & $2$  \cr \hline
7. & $\widehat{\mathbf{u}_2}$ & $1$ & \x{ $\mathbb{A}$, окружение -- \cr $(\mathbf{left},\mathbf{right},\mathbf{\mathbf{1A}},\mathbf{3A})$ }& $1$ босс -- $X$ & $3$ & $\mathbf{l}$  \cr \hline
8. & $\widehat{\mathbf{u}_1}$ & $1$ & \x{$\mathbb{B}$ окружение -- \cr $(\mathbf{left},\mathbf{right},\mathbf{1A},\mathbf{3A})$} & \x{ $1$ босс -- $X$ \cr тип $2$ -- $\mathbb{A}$ }  & $2$ & $\mathbf{r}$  \cr \hline
9. &  $\widehat{\mathbf{u}_2}$ & $2$ & \x{ $\mathbb{RU}$, окружение -- \cr $0$-цепь вокруг  $\mathbb{A}$ (ук $1$) \cr c окружением $\mathbf{Past}(X)$ \cr уровень 1 }& $1$ босс -- $X$  & $\mathbf{r}$ & $2$  \cr \hline
10. & $\mathbf{L.Core}(Z)$ & Предл.\ref{tilecore} & $\mathbf{Next.Past.FBoss}(Z)$& Предл.\ref{tilecore2}  & $\widehat{\mathbf{l}}$ & $3$  \cr \hline
  \end{tabular}
\label{P2b}
\end{table}

{\bf Примечание.}
Запись $\mathbf{Next.Past.FBoss}(Z)$ означает, что берется следующий узел в цепи для первого начальника $Z$, причем рассматривается его подклееное окружение (а не основное).

\medskip

{\bf Соотношения.} Вводятся аналогично случаю $P1$, коды вершин и ребер указаны в таблицах выше.

\subsection{P3: Макроплитка с окружением $(\mathbf{7A},x,\mathbf{3B},\mathbf{4A})$; $x=\mathbf{left},\mathbf{top}$}

В этом параграфе мы рассмотрим случай {\bf P3}, когда макроплитка $T$, внутри которой мы рассматриваем локальное преобразование пути, имеет окружение $(\mathbf{7A},x,\mathbf{1A},\mathbf{3A})$; $x=\mathbf{left},\mathbf{top}$. Это происходит, когда макроплитка является левой нижней подплиткой и верхней стороной выходит на левую или верхнюю сторону некоторой подклееной макроплитки.

\begin{figure}[hbtp]
\centering
\includegraphics[width=1\textwidth]{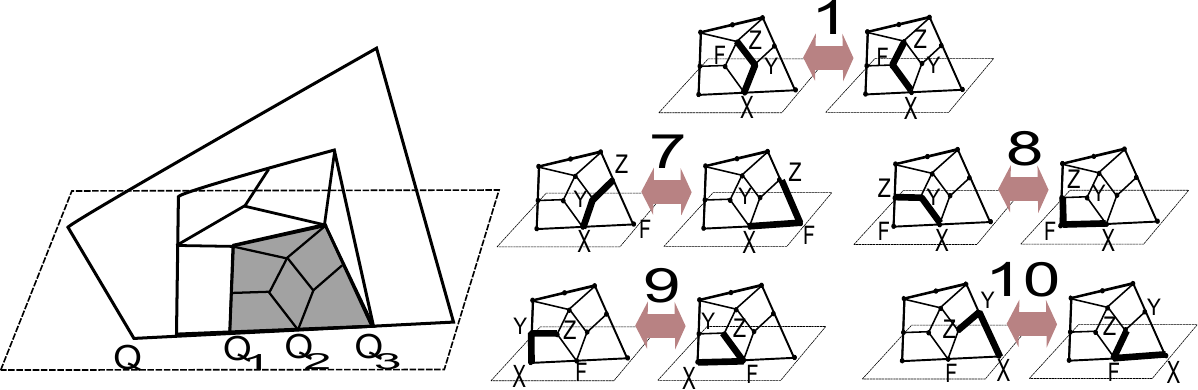}
\caption{Случай P3}
\label{P3}
\end{figure}

Зафиксируем пять вершин $Q$, $Q_1$, $Q_2$, $Q_3$, занимающих на комплексе положения как на рисунке~\ref{P3}. То есть $Q$ является ядром подклейки, а остальные вершины лежат на сторонах.

Для зафиксированного сочетания кодов этих пяти вершин, можно вычислить коды всех путей, указанных на рисунке~\ref{P3}. Кроме того, по выданному слову (коду пути) можно установить, действительно ли этот путь относится к случаю $P2$ и какую конфигурацию из изображенных на рисунке~\ref{P3} он имеет.

\medskip

\medskip

{\bf Вычисление кода $XFZ$.} Зная код пути $XYZ$ (то есть коды $X$, $Y$, $Z$, а также ребра входа в $Y$, $Z$ и ребра выхода из $X$ и $Y$ вдоль пути $XYZ$), можно выписать код $XFZ$. Обозначение $\mathbf{Past}(X)$ используем для обозначения подклееного окружения узла $X$ (таблица~\ref{P3a}).

\medskip
\begin{table}[hbtp]
  \caption{Случай P3: вычисление кода $XFZ$.}
\centering
 \begin{tabular}{|c|c|c|c|c|c|c|}   \hline
  & ребро из $X$ & ребро в $F$ & \x{тип, уровень и \cr окружение $F$} & информация $F$ &  ребро из $F$ & \x{ребро \cr в $Z$} \cr \hline
1. & $\widehat{\mathbf{u}_1}$ & $1$ & \x{$\mathbb{B}$ \cr (окруж как у $Y$) }& \x{ 1 босс -- $X$ \cr тип 2 -- $\mathbb{A}$ } & $3$ & $2$  \cr \hline
7. &Предл.\ref{tilecore} & Предл.\ref{tilecore} & $\mathbf{SBoss}(Z)$ & как у $X$  & \x{ $\widehat{\mathbf{u}_3}$ если \cr $\mathbf{Past}(F)=\mathbb{U}$ \cr  $\widehat{\mathbf{l}_2}$ если \cr $\mathbf{Past}(F)=\mathbb{L}$ } & $\mathbf{l}$   \cr \hline
8. & Предл.\ref{tilecore} & Предл.\ref{tilecore} & $\mathbf{Next.FBoss}(Z)$ & как у $X$  & $\mathbf{l}$ & $1$  \cr \hline
9. & Предл.\ref{tilecore} & Предл.\ref{tilecore} & $\mathbf{FBoss}(Z)$ & как у $X$  & $\widehat{\mathbf{u}_1}$ & $1$  \cr \hline
10. & Предл.\ref{tilecore} & Предл.\ref{tilecore} & $\mathbf{FBoss}(Z)$& как у $X$  & $\widehat{\mathbf{u}_2}$ & $1$  \cr \hline
  \end{tabular}
\label{P3a}
\end{table}

\medskip

{\bf Вычисление кода $XYZ$.} Зная код пути $XFZ$ (то есть коды $X$, $F$, $Z$, а также ребра входа в $F$, $Z$ и ребра выхода из $X$ и $F$ вдоль пути $XFZ$), можно выписать код $XYZ$ (таблица~\ref{P3b}).

\medskip
\begin{table}[hbtp]
 \caption{Случай P3: вычисление кода $XYZ$. }
\centering
 \begin{tabular}{|c|c|c|c|c|c|c|}   \hline
& \x{ребро \cr из $X$} & \x{ребро \cr в $Y$} & \x{тип, уровень \cr и окружение $Y$} & информация $Y$ & \x{ ребро \cr из $Y$} & \x{ребро \cr в $Z$ }\cr \hline
1. & $\widehat{\mathbf{u}_2}$ & $1$ & $\mathbb{A}$ (окруж. как у $F$)& 1 босс -- $X$ & $2$ & $1$  \cr \hline
7. & $\widehat{\mathbf{u}_2}$ & $1$ & \x{ $\mathbb{A}$, окружение -\cr  $(\mathbf{7A},x,\mathbf{3B},\mathbf{4A})$, \cr $x=\mathbf{left},\mathbf{top}$  }& $1$ босс -- $X$ & $3$ & $\mathbf{l}$  \cr  \hline
8. & $\widehat{\mathbf{u}_1}$ & $1$ &  \x{  $\mathbb{B}$, окружение -- \cr $(\mathbf{7A},x,\mathbf{3B},\mathbf{4A})$, \cr  $x=\mathbf{left},\mathbf{top}$  }& $1$ босс -- $X$  & $2$ & $\mathbf{r}$  \cr \hline
9. &  $\widehat{\mathbf{l}}$ & $2$ & \x{  $\mathbb{DR}$,1 окружение -- $3$ } & $1$ босс -- $Y$  & $\mathbf{r}$ & $2$  \cr \hline
10. & $\widehat{\mathbf{l}_2}$ или $\widehat{\mathbf{u}_1}$  & $1$ & $\mathbf{Next.Past.FBoss}(Z)$ & \x{ 1 босс -- \cr $\mathbf{TopFromCorner.Past}(X)$ \cr  $2$ босс -- $X$, $3$ босс -- \cr $\mathbf{RightCorner.Past}(X)$} & $\mathbf{l}$ &  $3$   \cr  \hline
  \end{tabular}
\label{P3b}
\end{table}

\medskip

{\bf Примечание}. Ребро из $X$ в случае $10$ имеет тип $\widehat{\mathbf{l}_2}$ если $\mathbf{Past}(X)=\mathbb{L}$ и $\widehat{\mathbf{u}_1}$ если $\mathbf{Past}(X)=\mathbb{U}$.

Запись $\mathbf{Next.Past.FBoss}(Z)$ означает, что берется следующий узел в цепи для первого начальника $Z$, причем рассматривается его подклееное окружение (а не основное).

\medskip

\subsection{P4: Макроплитка с окружением $(\mathbf{7A},x,\mathbf{1A},\mathbf{3A})$; $x=\mathbf{left},\mathbf{top}$}

В этом параграфе мы рассмотрим случай {\bf P4}, когда макроплитка $T$, внутри которой мы рассматриваем локальное преобразование пути, имеет окружение $(\mathbf{7A},x,\mathbf{1A},\mathbf{3A})$; $x=\mathbf{left},\mathbf{top}$. Это происходит, когда макроплитка является прямым потомком некоторой макроплитки, являющейся левой нижней подплиткой.  Верхней стороной $T$ выходит на левую или верхнюю сторону некоторой подклееной макроплитки.

Определение кодов вершин можно провести аналогично предыдущим случаям.

\medskip

\begin{figure}[hbtp]
\centering
\includegraphics[width=1\textwidth]{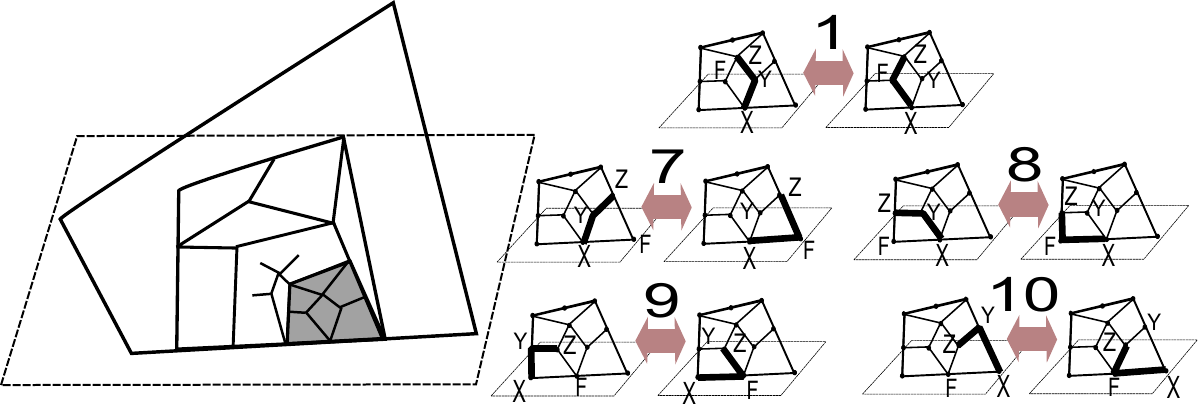}
\caption{Случай P4}
\label{P4}
\end{figure}

\medskip

{\bf Вычисление кода $XFZ$.} Зная код пути $XYZ$ (то есть коды $X$, $Y$, $Z$, а также ребра входа в $Y$, $Z$ и ребра выхода из $X$ и $Y$ вдоль пути $XYZ$), можно выписать код $XFZ$ (таблица~\ref{P4a}).

\medskip

\begin{table}[hbtp]
  \caption{Случай P4: вычисление кода $XFZ$. }
\centering
 \begin{tabular}{|c|c|c|c|c|c|c|}   \hline
  & ребро из $X$ & ребро в $F$ & \x{тип, уровень и \cr окружение $F$} & \x{информация \cr $F$} &  ребро из $F$ & \x{ребро \cr в $Z$}
 \cr \hline
1. & $\widehat{\mathbf{u}_1}$ & $1$ & \x{  $\mathbb{B}$ \cr  (окруж как у $Y$) }& \x{ 1 босс -- $X$ \cr тип 2 -- $\mathbb{A}$ }& $3$ & $2$  \cr \hline
7. & Предл.\ref{tilecore} & Предл.\ref{tilecore} & $\mathbf{SBoss}(Z)$ & как у $X$  & \x{ $\widehat{\mathbf{u}_3}$ если \cr $\mathbf{Past}(F)=\mathbb{U}$ \cr $\widehat{\mathbf{l}_2}$ если \cr $\mathbf{Past}(F)=\mathbb{L}$ } & $\mathbf{l}$   \cr \hline
8. & Предл.\ref{tilecore} & Предл.\ref{tilecore} & $\mathbf{FBoss}(Z)$ & как у $X$  & $\widehat{\mathbf{u}_2}$ & $2$  \cr \hline
9. & Предл.\ref{tilecore} & Предл.\ref{tilecore} & $\mathbf{FBoss}(Z)$ & как у $X$  & $\widehat{\mathbf{u}_1}$ & $1$  \cr \hline
10. & Предл.\ref{tilecore} & Предл.\ref{tilecore} & $\mathbf{FBoss}(Z)$& как у $X$  & $\widehat{\mathbf{u}_2}$ & $1$  \cr \hline
  \end{tabular}
\label{P4a}
\end{table}

\medskip

{\bf Вычисление кода $XYZ$.} Зная код пути $XFZ$ (то есть коды $X$, $F$, $Z$, а также ребра входа в $F$, $Z$ и ребра выхода из $X$ и $F$ вдоль пути $XFZ$), можно выписать код $XYZ$ (таблица~\ref{P4b}).

\medskip

\begin{table}[hbtp]
  \caption{Случай P4: вычисление кода $XYZ$. }
\centering
 \begin{tabular}{|c|c|c|c|c|c|c|}   \hline
& \x{ребро \cr из $X$} & \x{ребро \cr в $Y$} & \x{тип, уровень и \cr окружение $Y$} & информация $Y$ &  \x{ ребро \cr из $Y$ }& \x{ребро \cr в $Z$} \cr \hline
1. & $\widehat{\mathbf{u}_2}$ & $1$ & \x{$\mathbb{A}$, \cr (окружение как у $F$)} & 1 босс -- $X$ & $2$ & $1$  \cr \hline
7. & $\widehat{\mathbf{u}_2}$ & $1$ & \x{  $\mathbb{A}$, окружение -- \cr  $(\mathbf{7A},x,\mathbf{1A},\mathbf{3A})$, \cr $x=\mathbf{left},\mathbf{top}$ } & $1$ босс -- $X$ & $3$ & $\mathbf{l}$  \cr \hline
8. & $\widehat{\mathbf{u}_1}$ & $1$ & \x{ $\mathbb{B}$, окружение -- \cr  $(\mathbf{7A},x,\mathbf{1A},\mathbf{3A})$, \cr $x=\mathbf{left},\mathbf{top}$ }& $1$ босс -- $X$  & $2$ & $\mathbf{r}$  \cr  \hline
9. &  $\widehat{\mathbf{u}_2}$ & $2$ & \x{  $\mathbb{RU}$,1, окружение -- \cr $\mathbb{A}0$-цепь  (указ $1$) \cr c окруж $\mathbf{Past}(X)$ } & $1$ босс -- $X$  & $\mathbf{r}$ & $2$  \cr \hline
10. & $\widehat{\mathbf{l}_2}$ или $\widehat{\mathbf{u}_1}$  & $1$ & $\mathbf{Next.Past.FBoss}(Z)$ & \x{  1 босс -- \cr $\mathbf{TopFromCorner.Past}(X)$ \cr 2 босс -- $X$,  $3$ босс -- \cr $\mathbf{RightCorner.Past}(X)$ }& $\mathbf{l}$ & $3$  \cr  \hline
  \end{tabular}
\label{P4b}
\end{table}

\medskip

{\bf Примечание.} Ребро из $X$ в случае $10$ имеет тип $\widehat{\mathbf{l}_2}$ если $\mathbf{Past}(X)=\mathbb{L}$, и $\widehat{\mathbf{u}_1}$ если $\mathbf{Past}(X)=\mathbb{U}$.

Запись $\mathbf{Next.Past.FBoss}(Z)$ означает, что берется следующий узел в цепи для первого начальника $Z$, причем рассматривается его подклееное окружение (а не основное).

\medskip

\subsection{P5: Макроплитка с окружением $(\mathbf{top},\mathbf{right},\mathbf{6A},\mathbf{2A})$}

В этом параграфе мы рассмотрим случай {\bf P5}, когда макроплитка $T$, внутри которой мы рассматриваем локальное преобразование пути, имеет окружение $(\mathbf{top},\mathbf{right},\mathbf{6A},\mathbf{2A})$. Это происходит, когда макроплитка является правой верхней подплиткой некоторой подклееной макроплитки.

Определение кодов вершин можно провести аналогично предыдущим случаям.

\medskip

\begin{figure}[hbtp]
\centering
\includegraphics[width=1\textwidth]{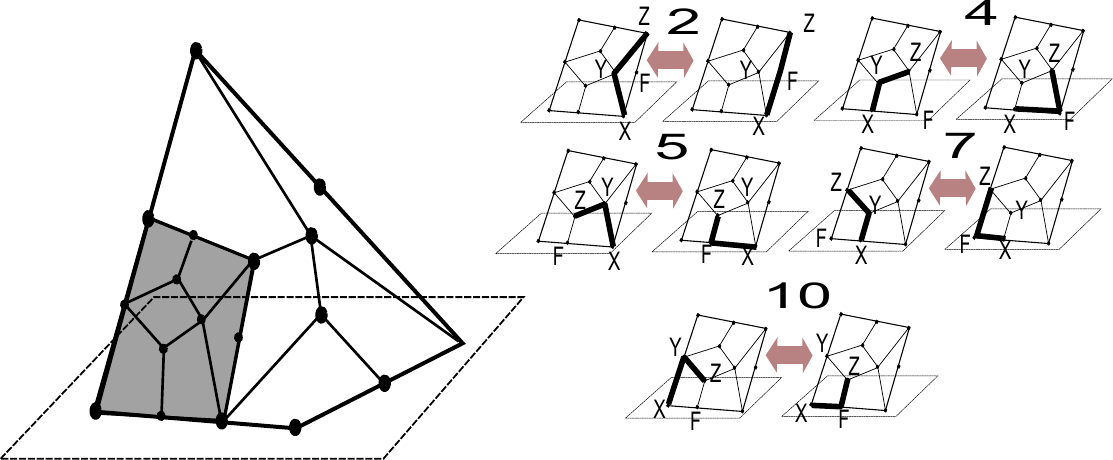}
\caption{Случай P5}
\label{P5}
\end{figure}

\medskip

{\bf Вычисление кода $XFZ$.} Зная код пути $XYZ$ (то есть коды $X$, $Y$, $Z$, а также ребра входа в $Y$, $Z$ и ребра выхода из $X$ и $Y$ вдоль пути $XYZ$), можно выписать код $XFZ$ (таблица~\ref{P5a}).

Для преобразования $7$ пути удовлетворяют условиям мертвого паттерна. Для них мы соотношения не вводим.

\medskip
\begin{table}[hbtp]
\caption{Случай P5: вычисление кода $XFZ$. }
\centering
 \begin{tabular}{|c|c|c|c|c|c|c|}   \hline
  & ребро из $X$ & ребро в $F$ & \x{тип, уровень и \cr окружение $F$} & информация $F$ &  ребро из $F$ & ребро в $Z$
 \cr  \hline
2. & $\widehat{\mathbf{u}_1}$ & $1$ & \x{ $\mathbb{DR}$,1, \cr окружение -- $3$ } & \x{ 1 босс -- $X$ \cr тип 2 -- $\mathbb{CDR}$  } & $1$ & $2$  \cr \hline
4. & Предл.\ref{tilecore} &  Предл.\ref{tilecore} & $\mathbf{SBoss}(Z)$ & как у $X$  & $\widehat{\mathbf{u}_3}$  & $4$   \cr \hline
5. &  Предл.\ref{tilecore} &  Предл.\ref{tilecore} & $\mathbf{Next.FBoss}(Z)$ & как у $X$  & $\widehat{\mathbf{l}}$ & $3$  \cr \hline
10. &  Предл.\ref{tilecore} &  Предл.\ref{tilecore} & $\mathbf{Next.FBoss}(Z)$ & как у $X$  & $\widehat{\mathbf{l}}$ & $3$  \cr \hline
  \end{tabular}
\label{P5a}
\end{table}

\medskip

{\bf Вычисление кода $XYZ$.} Зная код пути $XFZ$ (то есть коды $X$, $F$, $Z$, а также ребра входа в $F$, $Z$ и ребра выхода из $X$ и $F$ вдоль пути $XFZ$), можно выписать код $XYZ$ (таблица~\ref{P5b}).

\medskip

\begin{table}[hbtp]
\caption{Случай P5: вычисление кода $XYZ$. }
\centering
 \begin{tabular}{|c|c|c|c|c|c|c|}   \hline
& \x{ребро \cr из $X$} & \x{ребро \cr в $Y$} & \x{тип, уровень и \cr окружение $Y$} & информация $Y$ &  \x{ребро \cr из $Y$} & \x{ребро \cr в $Z$} \cr \hline
2. & $\widehat{\mathbf{u}_3}$ & $4$ & \x{ $\mathbb{C}$, окружение как у $T$} & \x{1 босс -- $\mathbb{U}$, \cr окружение как у $T$, \cr $(\mathbf{top},\mathbf{right},\mathbf{6A},\mathbf{2A})$, \cr  2 босс -- $X$, \cr 3 босс -- $\mathbb{B}$, окружение \cr $(\mathbf{left},\mathbf{top},\mathbf{right},\mathbf{bottom})$}  & $2$ & $1$  \cr \hline
4. & $\widehat{\mathbf{l}}$ & $3$ & $\mathbb{A}$, окружение как у $T$ & $1$ босс как у $Z$ & $2$ & $1$  \cr \hline
5. & $\widehat{\mathbf{u}_3}$ & $4$ & $\mathbb{C}$, окружение как у $T$ & \x{ $1$ босс -- как у $Z$ \cr 2 босс -- $X$,\cr 3 босс --  $\mathbb{B}$, окружение \cr $(\mathbf{left},\mathbf{top},\mathbf{right},\mathbf{bottom})$ }  & $1$ & $2$  \cr  \hline
10. & $\widehat{\mathbf{r}}$ & $2$ & $\mathbf{FBoss}(Z)$ & --   & $\mathbf{u}_2$ & $1$  \cr \hline
  \end{tabular}
\label{P5b}
\end{table}

\medskip

\subsection{P6: Макроплитка с окружением $(\mathbf{top},\mathbf{right},\mathbf{1A},\mathbf{3A})$}

В этом параграфе мы рассмотрим случай {\bf P6}, когда макроплитка $T$, внутри которой мы рассматриваем локальное преобразование пути, имеет окружение $(\mathbf{top},\mathbf{right},\mathbf{1A},\mathbf{3A})$. Это происходит, когда макроплитка является прямым потомком некоторой макроплитки, являющейся правой верхней подплиткой некоторой подклееной макроплитки.

Определение кодов вершин можно провести аналогично предыдущим случаям.

\medskip

\begin{figure}[hbtp]
\centering
\includegraphics[width=1\textwidth]{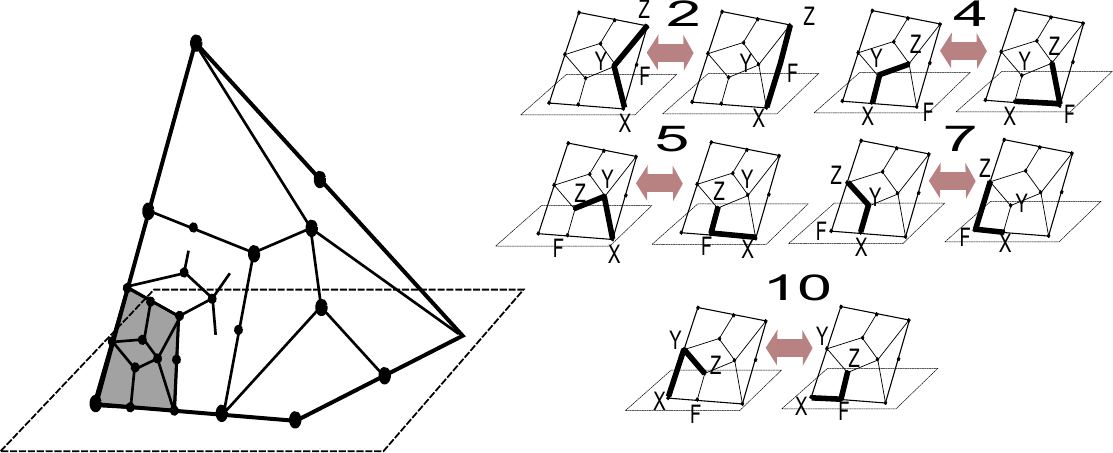}
\caption{Случай P6}
\label{P6}
\end{figure}

\medskip

{\bf Вычисление кода $XFZ$.} Зная код пути $XYZ$ (то есть коды $X$, $Y$, $Z$, а также ребра входа в $Y$, $Z$ и ребра выхода из $X$ и $Y$ вдоль пути $XYZ$), можно выписать код $XFZ$ (таблица~\ref{P6a}).

Для преобразования $7$ пути удовлетворяют условиям мертвого паттерна. Для них мы соотношения не вводим.

\medskip

\begin{table}[hbtp]
  \caption{Случай P6: вычисление кода $XFZ$. }
\centering
\begin{tabular}{|c|c|c|c|c|c|c|}   \hline
  & ребро из $X$ & ребро в $F$ & \x{тип, уровень и \cr окружение $F$} & информация $F$ &  ребро из $F$ & ребро в $Z$ \cr \hline
2. & $\widehat{\mathbf{l}}$ & $1$ & \x{ $\mathbb{DR}$,1, \cr окружение -- $3$ } & 1 босс как у $Z$ & $2$ & $3$  \cr \hline
4. &  Предл.\ref{tilecore} &  Предл.\ref{tilecore} & $\mathbf{SBoss}(Z)$ & как у $X$  & $\widehat{\mathbf{l}_3}$  & $4$   \cr \hline
5. &  Предл.\ref{tilecore} &  Предл.\ref{tilecore} & $\mathbf{Next.FBoss}(Z)$ & как у $X$  & $\widehat{\mathbf{l}}$ & $3$  \cr \hline
10. &  Предл.\ref{tilecore} &  Предл.\ref{tilecore} & $\mathbf{Next.FBoss}(Z)$ & как у $X$  & $\widehat{\mathbf{l}}$ & $3$  \cr \hline
  \end{tabular}
\label{P6a}
\end{table}

\medskip

{\bf Вычисление кода $XYZ$.} Зная код пути $XFZ$ (то есть коды $X$, $F$, $Z$, а также ребра входа в $F$, $Z$ и ребра выхода из $X$ и $F$ вдоль пути $XFZ$), можно выписать код $XYZ$ (таблица~\ref{P6b}).

\medskip

\begin{table}[hbtp]
\caption{Случай P6: вычисление кода $XYZ$. }
\centering
 \begin{tabular}{|c|c|c|c|c|c|c|}   \hline
& \x{ребро \cr из $X$} & \x{ребро \cr в $Y$} & \x{тип, уровень и \cr окружение $Y$} & информация $Y$ &  \x{ребро \cr из $Y$} & \x{ребро \cr в $Z$} \cr \hline
2. & $\widehat{\mathbf{l}_2}$ & $4$ & \x{$\mathbb{C}$, \cr окружение -- как у $T$} & \x{ 1 босс -- $\mathbb{U}$, \cr окружение как у $T$\cr 2 босс -- $X$, \cr 3 босс -- $\mathbb{A}$, окружение \cr как у $\mathbf{FBoss}(F)$  } & $2$ & $1$  \cr \hline
4. & $\widehat{\mathbf{l}}$ & $3$ & $\mathbb{A}$, окружение как у $T$ & $1$ босс как у $Z$ & $2$ & $1$  \cr \hline
5. & $\widehat{\mathbf{u}_3}$ & $4$ & $\mathbb{C}$, окружение как у $T$ & \x{ $1$ босс -- как у $Z$\cr 2 босс -- $X$, \cr 3 босс --  $\mathbb{A}$,окружение \cr $\mathbf{RightCorner.Past}(X)$  }  & $1$ & $2$  \cr  \hline
10. & $\widehat{\mathbf{r}}$  & $2$ & $\mathbf{FBoss}(Z)$ & 1 босс -- $X$   & $\mathbf{u}_2$ & $1$  \cr \hline
  \end{tabular}
\label{P6b}
\end{table}

\medskip

\subsection{P7: Макроплитка с окружением $(x,\mathbf{3B},\mathbf{6A},\mathbf{2A})$; $x=\mathbf{left},\mathbf{top}$}

В этом параграфе мы рассмотрим случай {\bf P7}, когда макроплитка $T$, внутри которой мы рассматриваем локальное преобразование пути, имеет окружение $(x,\mathbf{3B},\mathbf{6A},\mathbf{2A})$; $x=\mathbf{left},\mathbf{top}$. Это происходит, когда макроплитка является правой верхней подплиткой некоторой макроплитки, являющейся левой нижней подплиткой подклееной макроплитки.

Определение кодов вершин можно провести аналогично предыдущим случаям.

\medskip

\begin{figure}[hbtp]
\centering
\includegraphics[width=1\textwidth]{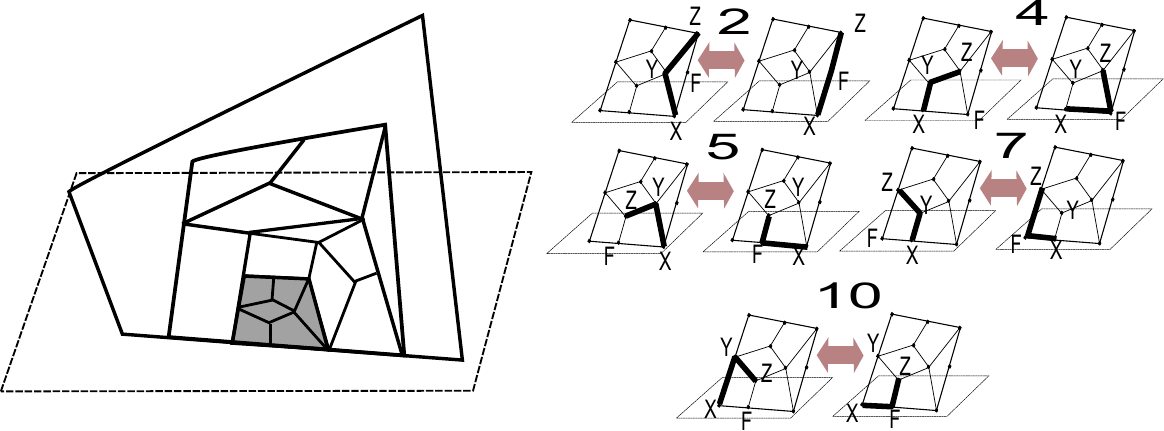}
\caption{Случай P7}
\label{P7}
\end{figure}

\medskip

{\bf Вычисление кода $XFZ$.} Зная код пути $XYZ$ (то есть коды $X$, $Y$, $Z$, а также ребра входа в $Y$, $Z$ и ребра выхода из $X$ и $Y$ вдоль пути $XYZ$), можно выписать код $XFZ$ (таблица~\ref{P7a}).

\medskip

\begin{table}[hbtp]
  \caption{Случай P7: вычисление кода $XFZ$. }
\centering
 \begin{tabular}{|c|c|c|c|c|c|c|}   \hline
  & ребро из $X$ & ребро в $F$ & \x{тип, уровень и \cr окружение $F$} & информация $F$ &  ребро из $F$ & ребро в $Z$ \cr \hline
2. & $\widehat{\mathbf{u}_1}$ & $1$ &  \x{ $\mathbb{DR}$,1, \cr окружение -- $2$  } & \x{ 1 босс -- $X$ \cr  тип 2 -- $\mathbb{A}$} & $2$ & $1$  \cr \hline
4. & Предл.\ref{tilecore} & Предл.\ref{tilecore}& $\mathbf{SBoss}(Z)$ & как у $X$  & $\widehat{\mathbf{u}_3}$  & $4$   \cr \hline
5. & Предл.\ref{tilecore} & Предл.\ref{tilecore} & $\mathbf{Next.FBoss}(Z)$ & как у $X$  & $\widehat{\mathbf{l}}$ & $3$  \cr \hline
7. & Предл.\ref{tilecore} & Предл.\ref{tilecore} & $\mathbf{Next.FBoss}(Z)$ & как у $X$  & $\widehat{\mathbf{l}}$ & $1$  \cr \hline
10. & Предл.\ref{tilecore} & Предл.\ref{tilecore} & $\mathbf{Next.FBoss}(Z)$ & как у $X$  & $\widehat{\mathbf{l}}$ & $3$  \cr \hline
  \end{tabular}
\label{P7a}
\end{table}

\medskip

{\bf Вычисление кода $XYZ$.} Зная код пути $XFZ$ (то есть коды $X$, $F$, $Z$, а также ребра входа в $F$, $Z$ и ребра выхода из $X$ и $F$ вдоль пути $XFZ$), можно выписать код $XYZ$ (таблица~\ref{P7b}).

\medskip

\begin{table}[hbtp]
\caption{Случай P7: вычисление кода $XYZ$. }
\centering
\begin{tabular}{|c|c|c|c|c|c|c|}   \hline
& ребро  & ребро & \x{тип, уровень и \cr окружение $Y$} & информация $Y$ &  ребро & ребро  \cr
&из $X$ &в $Y$&&& из $Y$ & в $Z$ \cr \hline
2. & $\widehat{\mathbf{u}_3}$ & $4$ & $\mathbb{C}$, окружение -- как у $T$ & \x{ 1 босс -- \cr $\mathbf{TopFromCorner.Past}(X)$,\cr 2 босс -- $X$, \cr 3 босс -- $Z$  } & $3$ & $\mathbf{ru}$  \cr \hline
4. & $\widehat{\mathbf{l}}$ & $3$ & $\mathbb{A}$ окружение как у $T$ & $1$ босс как у $Z$ & $2$ & $1$  \cr \hline
5. & $\widehat{\mathbf{u}_3}$ & $4$ & $\mathbb{C}$ окружение как у $T$ &  \x{$1$ босс -- как у $Z$\cr 2 босс -- $X$, \cr 3 босс --  $\mathbb{B}$,окружение \cr  $(\mathbf{7A},x,\mathbf{3B},\mathbf{4A})$  } & $2$ & $1$  \cr \hline
7. &  $\widehat{\mathbf{l}}$ & $3$ & $\mathbb{A}$ окружение как у $T$  & $1$ босс -- $Z$  & $1$ & $\mathbf{u}_2$  \cr\hline
10. & $\widehat{\mathbf{l}}$  & $1$ & $\mathbf{FBoss}(Z)$ & 1 босс -- $\mathbf{Prev}(X)$  & $\mathbf{u}_2$ & $1$  \cr \hline
  \end{tabular}
\label{P7b}
\end{table}

\medskip

\subsection{P8: Макроплитка с окружением $(x,\mathbf{3B},\mathbf{1A},\mathbf{3A})$; $x=\mathbf{left},\mathbf{top}$}

В этом параграфе мы рассмотрим случай {\bf P8}, когда макроплитка $T$, внутри которой мы рассматриваем локальное преобразование пути, имеет окружение $(x,\mathbf{3B},\mathbf{1A},\mathbf{3A})$; $x=\mathbf{left},\mathbf{top}$. Это происходит, когда макроплитка является прямым потомком правой верхней подплитки некоторой макроплитки, являющейся левой нижней подплиткой подклееной макроплитки.

Определение кодов вершин можно провести аналогично предыдущим случаям.

\medskip

\begin{figure}[hbtp]
\centering
\includegraphics[width=1\textwidth]{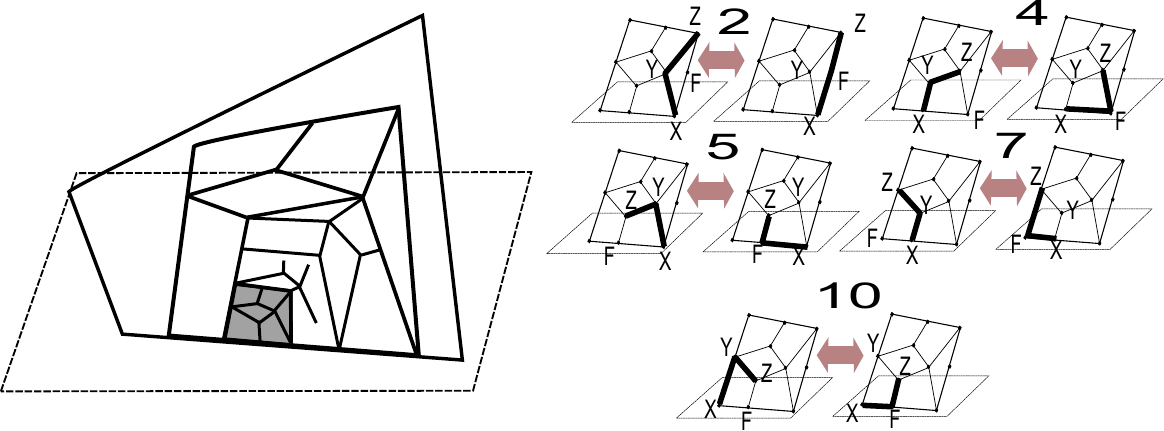}
\caption{Случай P8}
\label{P8}
\end{figure}

\medskip

{\bf Вычисление кода $XFZ$.} Зная код пути $XYZ$ (то есть коды $X$, $Y$, $Z$, а также ребра входа в $Y$, $Z$ и ребра выхода из $X$ и $Y$ вдоль пути $XYZ$), можно выписать код $XFZ$ (таблица~\ref{P8a}).

\medskip

\begin{table}[hbtp]
\caption{Случай P8: вычисление кода $XFZ$. }
\centering
 \begin{tabular}{|c|c|c|c|c|c|c|}   \hline
  & ребро из $X$ & ребро в $F$ & \x{тип, уровень и \cr окружение $F$} & информация $F$ &  ребро из $F$ & ребро в $Z$ \cr \hline
2. & $\widehat{\mathbf{l}}$ & $1$ & \x{  $\mathbb{DR}$,1, \cr окружение -- $3$  }& 1 босс -- как у $Z$ & $2$ & $3$  \cr \hline
4. & Предл.\ref{tilecore} & Предл.\ref{tilecore} & $\mathbf{SBoss}(Z)$ & как у $X$  & $\widehat{\mathbf{l}_2}$  & $4$   \cr \hline
5. & Предл.\ref{tilecore} & Предл.\ref{tilecore} & $\mathbf{Next.FBoss}(Z)$ & как у $X$  & $\widehat{\mathbf{l}}$ & $3$  \cr \hline
7. & Предл.\ref{tilecore} & Предл.\ref{tilecore} & $\mathbf{Next.FBoss}(Z)$ & как у $X$  & $\widehat{\mathbf{l}}$ & $1$  \cr \hline
10. & Предл.\ref{tilecore} & Предл.\ref{tilecore} & $\mathbf{Next.FBoss}(Z)$ & как у $X$  & $\widehat{\mathbf{l}}$ & $3$  \cr \hline
  \end{tabular}
\label{P8a}
\end{table}

\medskip

{\bf Вычисление кода $XYZ$.} Зная код пути $XFZ$ (то есть коды $X$, $F$, $Z$, а также ребра входа в $F$, $Z$ и ребра выхода из $X$ и $F$ вдоль пути $XFZ$), можно выписать код $XYZ$ (таблица~\ref{P8b}).

\medskip

\begin{table}[hbtp]
 \caption{Случай P8: вычисление кода $XYZ$. }
\centering
\begin{tabular}{|c|c|c|c|c|c|c|}   \hline
& \x{ребро \cr из $X$} & \x{ребро \cr в $Y$} & \x{тип, уровень и \cr окружение $Y$} & информация $Y$ &  \x{ребро \cr из $Y$} & \x{ребро \cr в $Z$} \cr \hline
2. & $\widehat{\mathbf{l}_2}$ & $4$ & $\mathbb{C}$, окружение -- как у $T$ & \x{1 босс -- \cr  $\mathbf{TopFromCorner.Past}(X)$, \cr 2 босс -- $X$ \cr 3 босс -- $Z$   } & $3$ & $\mathbf{ru}$  \cr \hline
4. & $\widehat{\mathbf{l}}$ & $3$ & $\mathbb{A}$, окружение как у $T$ & $1$ босс как у $Z$ & $2$ & $1$  \cr \hline
5. & $\widehat{\mathbf{l}_2}$ & $4$ & $\mathbb{C}$, окружение как у $T$ & \x{ $1$ босс -- как у $Z$\cr 2 босс -- $X$, \cr 3 босс --  $\mathbb{A}$, окружение \cr $\mathbf{RightCorner.Past}(X)$ } & $1$ & $2$  \cr \hline
7. &  $\widehat{\mathbf{l}}$ & $3$ & $\mathbb{A}$, окружение как у $T$  & $1$ босс -- $Z$  & $1$ & $\mathbf{u}_2$  \cr\hline
10. & $\widehat{\mathbf{l}}$  & $2$ & $\mathbf{FBoss}(Z)$ & 1 босс -- $\mathbf{Prev}(X)$  & $\mathbf{u}_2$ & $1$  \cr \hline
  \end{tabular}
\label{P8b}
\end{table}

\medskip

\subsection{P9: Макроплитка с окружением $(x,\mathbf{1A},\mathbf{6A},\mathbf{2A})$; $x=\mathbf{left},\mathbf{top}$}

В этом параграфе мы рассмотрим случай {\bf P9}, когда макроплитка $T$, внутри которой мы рассматриваем локальное преобразование пути, имеет окружение $(x,\mathbf{3B},\mathbf{6A},\mathbf{2A})$; $x=\mathbf{left},\mathbf{top}$. Это происходит, когда макроплитка является правой верхней подплиткой некоторой макроплитки, являющейся левой верхней подплиткой подклееной макроплитки.

Определение кодов вершин можно провести аналогично предыдущим случаям.

\medskip

\begin{figure}[hbtp]
\centering
\includegraphics[width=1\textwidth]{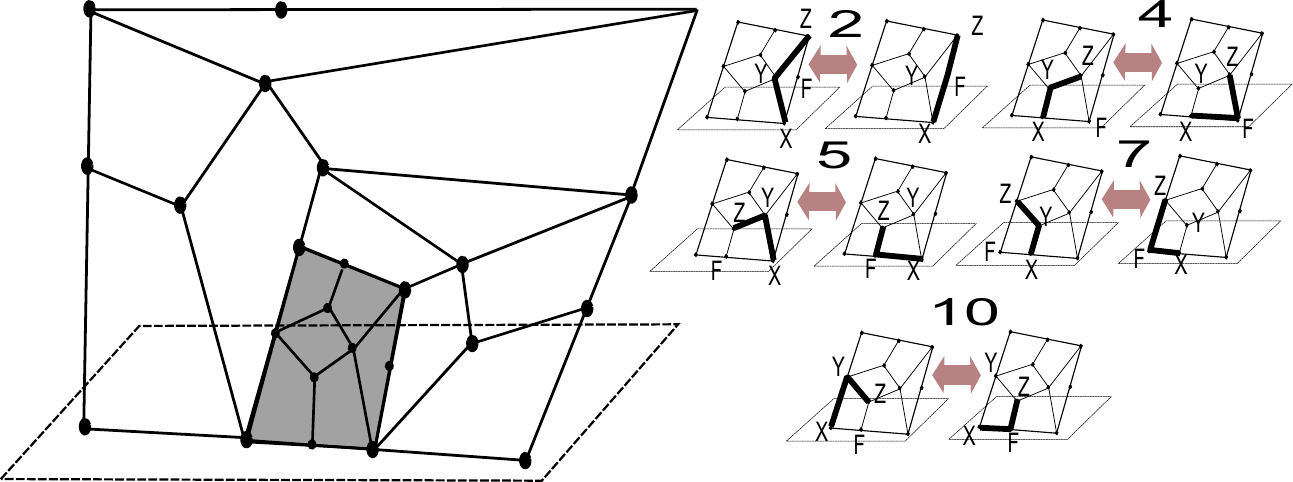}
\caption{Случай P9}
\label{P9}
\end{figure}

\medskip

{\bf Вычисление кода $XFZ$.} Зная код пути $XYZ$ (то есть коды $X$, $Y$, $Z$, а также ребра входа в $Y$, $Z$ и ребра выхода из $X$ и $Y$ вдоль пути $XYZ$), можно выписать код $XFZ$ (таблица~\ref{P9a}).

\medskip
\begin{table}[hbtp]
  \caption{Случай P9: вычисление кода $XFZ$. }
\centering
 \begin{tabular}{|c|c|c|c|c|c|c|}   \hline
  & ребро из $X$ & ребро в $F$ & \x{тип, уровень и \cr окружение $F$} & информация $F$ &  ребро из $F$ & ребро в $Z$ \cr \hline
2. & $\widehat{\mathbf{u}_1}$ & $1$ & \x{  $\mathbb{DR}$,1, \cr окружение -- $2$ }& \x{ 1 босс -- $X$ \cr  тип 2 -- $\mathbb{A}$} & $2$ & $1$  \cr \hline
4. & Предл.\ref{tilecore} & Предл.\ref{tilecore} & $\mathbf{SBoss}(Z)$ & как у $X$  & $\widehat{\mathbf{u}_3}$  & $4$   \cr \hline
5. & Предл.\ref{tilecore} & Предл.\ref{tilecore} & $\mathbf{Next.FBoss}(Z)$ & как у $X$  & $\widehat{\mathbf{l}}$ & $3$  \cr \hline
7. & Предл.\ref{tilecore} & Предл.\ref{tilecore} & $\mathbf{Next.FBoss}(Z)$ & как у $X$  & $\widehat{\mathbf{u}_2}$ & $1$  \cr \hline
10. & Предл.\ref{tilecore} & Предл.\ref{tilecore} & $\mathbf{Next.FBoss}(Z)$ & как у $X$  & $\widehat{\mathbf{l}}$ & $3$  \cr \hline
  \end{tabular}
\label{P9a}
\end{table}

\medskip

{\bf Вычисление кода $XYZ$.} Зная код пути $XFZ$ (то есть коды $X$, $F$, $Z$, а также ребра входа в $F$, $Z$ и ребра выхода из $X$ и $F$ вдоль пути $XFZ$), можно выписать код $XYZ$ (таблица~\ref{P9b}).

\medskip

\begin{table}[hbtp]
\caption{Случай P9: вычисление кода $XYZ$.}
\centering
\begin{tabular}{|c|c|c|c|c|c|c|}   \hline
& ребро  & ребро  & \x{тип, уровень и \cr окружение $Y$} & информация $Y$ &  ребро  & ребро  \cr
&из $X$ &в $Y$&&& из $Y$ & в $Z$ \cr \hline
2. & $\widehat{\mathbf{u}_3}$ & $4$ & $\mathbb{C}$, окружение -- как у $T$ & \x{ 1 босс -- \cr  $\mathbf{TopFromCorner.Past}$(X), \cr 2 босс -- $X$, \cr 3 босс -- $Z$ } & $3$ & $\mathbf{ru}$  \cr \hline
4. & $\widehat{\mathbf{l}}$ & $3$ & $\mathbb{A}$, окружение как у $T$ & $1$ босс как у $Z$ & $2$ & $1$  \cr \hline
5. & $\widehat{\mathbf{u}_3}$ & $4$ & $\mathbb{C}$, окружение как у $T$ & \x{  $1$ босс -- как у $Z$\cr 2 босс -- $X$, \cr 3 босс --  $\mathbb{B}$, окружение \cr $(\mathbf{left},\mathbf{top},\mathbf{1A},\mathbf{3A})$ } & $1$ & $2$  \cr \hline
7. &  $\widehat{\mathbf{l}}$ & $3$ & $\mathbb{A}$, окружение как у $T$  & $1$ босс -- $Z$  & $1$ & $\mathbf{u}_2$  \cr\hline
10. & $\widehat{\mathbf{l}}$  & $1$ & $\mathbf{FBoss}(Z)$ & 1 босс -- $X$  & $\mathbf{u}_2$ & $1$  \cr \hline
  \end{tabular}
\label{P9b}
\end{table}

\medskip

\subsection{P10: Макроплитка с окружением $(x,\mathbf{1A},\mathbf{1A},\mathbf{3A})$; $x=\mathbf{left},\mathbf{top}$}

В этом параграфе мы рассмотрим случай {\bf P10}, когда макроплитка $T$, внутри которой мы рассматриваем локальное преобразование пути, имеет окружение $(x,\mathbf{3B},\mathbf{1A},\mathbf{3A})$; $x=\mathbf{left},\mathbf{top}$. Это происходит, когда макроплитка является прямым потомком правой верхней подплитки некоторой макроплитки, являющейся левой верхней подплиткой подклееной макроплитки.

Определение кодов вершин можно провести аналогично предыдущим случаям.

\medskip

\begin{figure}[hbtp]
\centering
\includegraphics[width=1\textwidth]{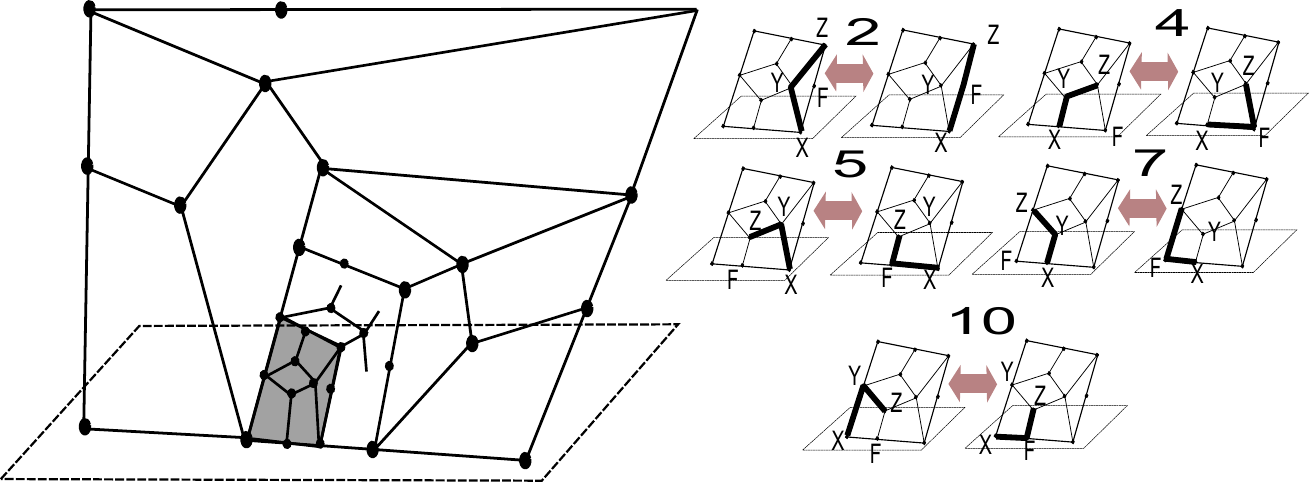}
\caption{Случай P10}
\label{P10}
\end{figure}

\medskip

{\bf Вычисление кода $XFZ$.} Зная код пути $XYZ$ (то есть коды $X$, $Y$, $Z$, а также ребра входа в $Y$, $Z$ и ребра выхода из $X$ и $Y$ вдоль пути $XYZ$), можно выписать код $XFZ$ (таблица~\ref{P10a}).

\medskip

\begin{table}[hbtp]
\caption{Случай P10: вычисление кода $XFZ$. }
\centering
 \begin{tabular}{|c|c|c|c|c|c|c|}   \hline
  & ребро из $X$ & ребро в $F$ & \x{тип, уровень и \cr окружение $F$} & информация $F$ &  ребро из $F$ & ребро в $Z$ \cr \hline
2. & $\widehat{\mathbf{l}}$ & $1$ & \x{  $\mathbb{DR}$ окружение -- $3$ }& 1 босс как у $Z$ & $2$ & $3$  \cr \hline
4. & Предл.\ref{tilecore} & Предл.\ref{tilecore} & $\mathbf{SBoss}(Z)$ & как у $X$  & $\widehat{\mathbf{l}_2}$  & $4$   \cr \hline
5. & Предл.\ref{tilecore} & Предл.\ref{tilecore} & $\mathbf{Next.FBoss}(Z)$ & как у $X$  & $\widehat{\mathbf{l}}$ & $3$  \cr \hline
7. & Предл.\ref{tilecore} & Предл.\ref{tilecore} & $\mathbf{FBoss}(Z)$ & как у $X$  & $\widehat{\mathbf{u}_2}$ & $2$  \cr \hline
10. & Предл.\ref{tilecore} & Предл.\ref{tilecore} & $\mathbf{Next.FBoss}(Z)$ & как у $X$  & $\widehat{\mathbf{l}}$ & $3$  \cr \hline
  \end{tabular}
\label{P10a}
\end{table}

\medskip

{\bf Вычисление кода $XYZ$.} Зная код пути $XFZ$ (то есть коды $X$, $F$, $Z$, а также ребра входа в $F$, $Z$ и ребра выхода из $X$ и $F$ вдоль пути $XFZ$), можно выписать код $XYZ$ (таблица~\ref{P10b}).

\medskip
\begin{table}[hbtp]
  \caption{Случай P10: вычисление кода $XYZ$. }
\centering
\begin{tabular}{|c|c|c|c|c|c|c|}   \hline
& \x{ребро \cr из $X$} & \x{ребро \cr в $Y$} & \x{тип, уровень и \cr окружение $Y$} & информация $Y$ &  \x{ребро \cr из $Y$ }& \x{ребро \cr в $Z$} \cr \hline
2. & $\widehat{\mathbf{l}_2}$ & $4$ & $\mathbb{C}$, окружение -- как у $T$ & \x{ 1 босс -- \cr $\mathbf{LevelPlus.FBoss}(Z)$, \cr 2 босс -- $X$, \cr 3 босс -- $Z$} & $3$ & $\mathbf{ru}$  \cr \hline
4. & $\widehat{\mathbf{l}}$ & $3$ & $\mathbb{A}$, окружение как у $T$ & $1$ босс как у $Z$ & $2$ & $1$  \cr \hline
5. & $\widehat{\mathbf{l}_2}$ & $4$ & $\mathbb{C}$, окружение как у $T$ & \x{ $1$ босс -- как у $Z$\cr 2 босс -- $X$, \cr 3 босс --  $\mathbb{A}$, окружение \cr $\mathbf{RightCorner.Past}(X)$ } & $1$ & $2$  \cr  \hline
7. &  $\widehat{\mathbf{l}}$ & $3$ & $\mathbb{A}$, окружение как у $T$  & $1$ босс -- $Z$  & $1$ & $\mathbf{u}_2$  \cr\hline
10. & $\widehat{\mathbf{u}_2}$  & $2$ & $\mathbf{FBoss}(Z)$ & 1 босс -- $X$  & $\mathbf{u}_2$ & $1$  \cr \hline
  \end{tabular}
\label{P10b}
\end{table}

\medskip

Таким образом, во всех случаях можно провести локальное преобразование, как и для плоских путей. Это дает нам возможности преобразовывать пути, меняя локальные участки на им эквивалентные и приводить слово к канонической форме.

\medskip

\subsection{Оценка числа введенных соотношений}

В каждом случае мы выбирали ядро подклейки, после чего выбирали не более четырех вершин ($Q_1$, $Q_2$, $Q_3$, $Q_4$) на сторонах подклееной макроплитки, для которых фиксировали их коды. После чего записывали не более $16$ соотношений, для каждого выбранного случая. Учитывая конструкцию подклейки, можно заметить, что информация вершин $Q_i$ может быть вычислена по информации ядра подклейки и известным ребрам-сторонам. То есть соотношений будет не более чем $16 \cdot F \cdot H^4$, где $F$ -- количество флагов подклейки, $H$ -- число расширенных окружений. То есть не более чем $16 \cdot 82485^4 \cdot 5 \cdot 10^{18}<4\cdot 10^{39}$ соотношений.

\medskip

%\section{Библиография}

%\addcontentsline{toc}{chapter}{\Numline {}Библиография}

%\markboth{}{Библиография}

%\begin{enumerate}

% Доклады, сделанные по материалам работы, приведены в конце списка литературы.

\end{document}